\documentclass[12pt]{article}
\usepackage[a4paper,left=2cm,right=1.2cm,top=3cm,bottom=4cm]
{geometry}

\usepackage{amsmath,amstext, amsthm, empheq, extpfeil,  
amsbsy,amssymb,marvosym,fancyhdr,graphicx,amscd,amsfonts,latexsym,delarray,stackrel,lineno,color,cite,appendix,xcolor,url,hyperref}
\usepackage{wasysym}

\usepackage{srcltx}
\usepackage{mathrsfs}
\usepackage{float} 
\usepackage[all]{xy}
\usepackage{t1enc}
\usepackage{mathrsfs}
\usepackage{pifont}

\usepackage{mathpple}
\usepackage[T1]{fontenc}

\newcommand*\patchAmsMathEnvironmentForLineno[1]{%
  \expandafter\let\csname old#1\expandafter\endcsname\csname #1\endcsname
  \expandafter\let\csname oldend#1\expandafter\endcsname\csname end#1\endcsname
  \renewenvironment{#1}%
     {\linenomath\csname old#1\endcsname}%
     {\csname oldend#1\endcsname\endlinenomath}}%
\newcommand*\patchBothAmsMathEnvironmentsForLineno[1]{%
  \patchAmsMathEnvironmentForLineno{#1}%
  \patchAmsMathEnvironmentForLineno{#1*}}%
\AtBeginDocument{%
\patchBothAmsMathEnvironmentsForLineno{equation}%
\patchBothAmsMathEnvironmentsForLineno{align}%
\patchBothAmsMathEnvironmentsForLineno{flalign}%
\patchBothAmsMathEnvironmentsForLineno{alignat}%
\patchBothAmsMathEnvironmentsForLineno{gather}%
\patchBothAmsMathEnvironmentsForLineno{multline}%
}

\definecolor{Green}{rgb}{0,1,0}
\definecolor{Blue}{RGB}{0,0,191}
\definecolor{mathmodecolor}{RGB}{0,102,0}
\definecolor{keywordcolor}{RGB}{0,51,151}
\definecolor{sourcebackgroundcolor}{RGB}{255,247,223}
\definecolor{unixagred}{RGB}{255,0,0}
\definecolor{lightgray}{RGB}{191,191,191}
\definecolor{green}{RGB}{1,191,191}

\newtheorem{thm}{Theorem}[section]
\newtheorem{prop}[thm]{Proposition}
\newtheorem{cor}[thm]{Corollary}
\newtheorem{lem}[thm]{Lemma}

\newtheorem{defn}[thm]{Definition}
\newtheorem{rem}[thm]{Remark}
\newtheorem{example}[thm]{Example}

\def\Aut{{\rm Aut}}
\def\Coker{{\rm Coker}}

\def\End{{\rm End}}

\def\Hom{{\rm Hom}}
\def\Homi{{\underline{\rm Hom}}}
\def\Endi{{\underline{\rm End}}}
\def\ext{{\rm Ext}}
\def\id{{\rm id}}

\def\Ker{{\rm Ker}}
\def\Kpr{{\rm Ker}_p\,}
\def\Cpr{{\rm Coker}_p\,}
\def\sker{{\rm ker}}
\def\equ{{\rm Equ}}
\def\coker{{\rm Coker}}
\def\scoker{{\rm cok}}
\def\coequ{{\rm Coequ}}

\def\Range{{\rm Range}}
\def\rad{{\rm Rad}}

\def\B{{\mathbb B}}

\def\N{{\mathbb N}}

\def\R{{\mathbb R}}
\def\Z{{\mathbb Z}}

\def\aarith{{\mathfrak A}}

\def\hweak{H^{\rm weak}}
\def\imm{\overline{\rm Im}}
\def\im{{\rm Im}}
\def\coimm{{\rm Coim}}
\def\nsb{{\rm Nsb}}

\def\mC{{\mathscr C}}

\def\exxx{{doubly exact}~}
\def\strg{strict }
\def\Strg{Strict }
\def\strgly{strictly }
\def\Strgly{Strictly }
\def\cA{{\mathcal A}}
\def\cB{{\mathcal B}}
\def\cC{{\mathcal C}}

\def\cE{{\mathcal E}}

\def\cH{{\mathcal H}}
\def\cI{{\mathcal I}}
\def\cF{{\mathcal F}}
\def\cJ{{\mathcal J}}
\def\cK{{\mathcal K}}
\def\cL{{\mathcal L}}
\def\cN{{\mathcal N}}

\def\cO{{\mathcal O}}

\def\cR{{\mathcal R}}
\def\cS{{\mathcal S}}

\def\cU{{\mathcal U}}

\def\aarith{{\mathscr A}}
\def\scal{{(\rnt,\cO)}}
\def\scal1{{\hat \aarith}}
\def\dar[#1]{\ar@<2pt>[#1]\ar@<-2pt>[#1]}
\def\qqq{\,,\,~\forall}

\def\bs{{\B[s]}}

\newcommand{\ie}{{\it i.e.\/}\ }
\newcommand{\eg}{{\it e.g.\/}\ }
\newcommand{\cf}{{\it cf.}}
\newcommand{\opcit}{{\it op.cit.\/}\ }

\def\id{{\mbox{Id}}}

\def\ker{{\mbox{Ker}}}

\def\Hom {{\mbox{Hom}}}

\def\End{{\mbox{End}}}

\def\chain{{\rm Ch_\bullet}}
\def\chainp{{\rm Ch_+}}
\def\sequ{{\rm Sh}}

\def\yb{\mathfrak{s}}

\def\rma{\R_{\rm max}}

\def\bmod{{\rm \B mod}}
\def\bm2{{\rm \B mod^2}}
\def\b2{{\rm \B mod^{\mathfrak s}}}

\def\Se{\frak{ Sets}}

\def\coim{{\rm Coim}}

\newcommand{\nil}[1]{}

\def\rnt{{[0,\infty)\rtimes{\N^{\times}}}}

\title
{Homological algebra in characteristic one}

\author{Alain Connes, Caterina Consani\thanks{Partially supported by the Simons Foundation collaboration grant n.~353677. This author would also like to thank the Coll\`ege de France for some financial support.}}

\date{}

\begin{document}

\maketitle

\begin{abstract}This article  develops several main results for a general theory of homological algebra in categories such as the category of  sheaves of idempotent modules over a topos. In the analogy with the development of homological algebra   for abelian categories the present paper should be viewed as the analogue  of the development of homological algebra for abelian groups. Our selected prototype, the category $\bmod$ of modules over the Boolean semifield $\B:=\{0,1\}$  is the replacement for the category of abelian groups. We show that the semi-additive category $\bmod$ fulfills analogues of the axioms AB1 and AB2 for abelian categories. By introducing a precise comonad on $\bmod$ we obtain the conceptually related  Kleisli and Eilenberg-Moore categories. The latter category $\b2$ is simply $\bmod$ in the topos of sets endowed with an involution and as such it shares with $\bmod$ most of its abstract categorical properties. The three main results of the paper are the following. First, when endowed with the natural ideal of null morphisms, the category $\b2$ is a semi-exact, homological category in the sense of M. Grandis. Second, there is a far reaching analogy between $\b2$ and the category of operators in Hilbert space, and in particular results relating null kernel and injectivity for morphisms. The third fundamental result is  that, even for finite objects of $\b2$ the resulting homological algebra is non-trivial and  gives rise to a computable Ext functor. We determine explicitly this functor in the case provided by the diagonal morphism of the Boolean semiring into its square.
	
\end{abstract}

\vspace{0.5in}

\tableofcontents

\section{Introduction}

This article lays the ground for the development of a cohomological theory of sheaves of $\rma$-modules over a topos. The need for such a theory is obvious in tropical geometry and arises naturally in the development of our work on the Riemann-Roch theorem on the scaling site and its square (\!\!\cite{CCss, CCss1}).  The expected new cohomology theory  ought to be described in terms of $\rma$-modules. However, sheaves of $\rma$-modules do not form an abelian category and thus abstract homological algebra cannot be applied in its classical form. To by-pass this difficulty, we decided to follow the strategic path developed by Grothendieck in his Tohoku paper \cite{Grtohoku}, where abelian sheaf cohomology is realized as an important example of homological algebra on abelian categories. Therefore this article advocates the existence of a meaningful theory of homological algebra in categories, such as the category of  sheaves of $\rma$-modules over a topos, which share the same abstract properties of our selected prototype namely the category $\bmod$ of modules over the Boolean semifield $\B:=\{0,1\}$ (\cf~\S\ref{sect bmodules}). In the analogy with the development of homological algebra in abelian categories, $\bmod$ is the replacement for the category of abelian groups and the results of the present paper should be viewed as the analogue of the development of homological algebra for abelian groups as in \cite{CE}. The key property that allows to transfer our main results on $\bmod$ to the category of sheaves of $\B$-modules is provided by the well known fact that a constructive proof in a category still holds true in any topos. \newline
The category $\bmod$ is usually referred in books as the category of join semi-lattices with a least element, or idempotent semi-modules and has a long history: we refer to \cite{Cohen} for an overview. In particular there is a well known equivalence between $\bmod$ and the category of algebraic lattices. This equivalence together with the crucial notion of Galois connection are recalled in \S\ref{galois connect}. This review is based on the duality results of \cite{Cohen} whose proofs are given, for the sake of completeness, in the simplest case of $\bmod$. Incidentally, the equivalence of $\bmod$ with the category of algebraic lattices provides also the reader with a wealth of concrete and interesting examples even by restricting to the case of finite objects. In fact  $\bmod$ plays a universal role taking into account the transfer functor \cite{gra1} which associates to an object of a category the lattice of its subobjects.  \newline While $\bmod$ is not an additive category it is  semi-additive \ie  it has a zero object (both initial and final), all finite products and coproducts exist,  and for any pair of objects $M,N$ in $\bmod$ the canonical morphism $M \amalg N\to M\times N$
from the coproduct  to the product 
 is an isomorphism.
One sets $M\oplus N=M \amalg N\simeq M\times N$, and  lets $p_j:M \oplus M\to M$ be the canonical projections and $s_j:M\to M\oplus M$  the canonical inclusions. In a semi-additive category $\mC$ the set of the endomorphisms  of an object in $\mC$ is canonically endowed with the structure of semiring. Moreover,  the analogues of the Axioms AB1 and AB2 for abelian categories (\!\cite{Grtohoku}, \S 1.4) take the following form\vspace{0.2cm}

	AB1': $\mC$ admits equalizers $\equ$ and coequalizers $\coequ$. \newline
	Assuming AB1' we let, for $f\in \Hom_\mC(L,M)$, $
f^{(2)}:=(f\circ p_1,f\circ p_2)
$ and $\equ f^{(2)}$ be the equalizer of $(f\circ p_1,f\circ p_2)$. It is  the categorical kernel of $f$. Similarly, with $
f_{(2)}:=(s_1\circ f,s_2\circ f)
$ the coequalizer  $\coequ f_{(2)}$ is the categorical cokernel of $f$. The analogue of the axiom AB2 is\newline
	AB2': For $f\in \Hom_\mC(L,M)$, the natural morphism $\coim f:=\coequ(\equ f^{(2)})\stackrel{\underline f}{\to} \im(f):=\equ(\coequ f_{(2)})$ is an isomorphism.

 Then, for any morphism $f:L\to M$ in $\mC$ one derives the sequence
\begin{equation}\label{defncokergenintro}
	\equ f^{(2)} \stackbin[\iota_2]{\iota_1}{\rightrightarrows} L\stackrel{f}{\to} M \stackbin[\gamma_2]{\gamma_1}{\rightrightarrows}\coequ f_{(2)}
\end{equation}

In \S\ref{sect abstract bmod} we prove that the category $\bmod$ fulfills the above Axioms $AB1'$ and $AB2'$. In this part  a key role is played by  a  comonad $\perp$ on the category $\bmod$ which encodes the endofunctor $M\longrightarrow M^2$ involved in the above construction of $f^{(2)}$ and $f_{(2)}$. The  comonad $\perp$ is defined by
\begin{enumerate}
	\item The endofunctor $\perp:\bmod\longrightarrow \bmod$,~~ $\perp M = M^2$, $\perp f=(f,f)$.
	\item The counit $\epsilon: \perp \to 1_{\bmod}$,~~$\epsilon_M = p_1$,~ $p_1:M^2 \to M$, $p_1(a,b)=a$.
	\item The coproduct $\delta:\perp \to \perp\circ\perp$,~~ $\delta_M= (M^2 \to (M^2)^2)$, ~$(x,y)\mapsto (x,y,y,x)$. 
\end{enumerate}	
   In Proposition \ref{Kleisli} we determine the  Kleisli and Eilenberg-Moore categories associated to the comonad $\perp$. The Kleisli category $\bm2$ has the same objects as $\bmod$ whereas the morphisms are  pairs of morphisms in $\bmod$ composing  following the rule
\begin{equation}\label{paircompfirstintro}
	(f,g)\circ (f',g'):=(f\circ f'+g\circ g', f\circ g'+g\circ f').
\end{equation}
The category $\bm2$ is the very natural enlargement of $\bmod$ where pairs of morphisms  are taken into account so that the notions of kernels and cokernels become similar to their incarnations in abelian categories.   In \S\ref{weak and strong} we study, as a technical tool revealing the properties of the category $\bm2$, the following \strg form of exactness 
\begin{defn}\label{exsequ intro} We say that the sequence  $L \stackbin[\alpha_2]{\alpha_1}{\rightrightarrows}M\stackbin[\beta_2]{\beta_1}{\rightrightarrows}N$ in $\bm2$ is {\em \strgly exact} at $M$ if  $ B(\alpha_1,\alpha_2)+\Delta=Z(\beta_1,\beta_2)$, where $\Delta\subset M\times M$ is the diagonal and one sets
\begin{equation}	
\label{exactseqstr}
B(\alpha_1,\alpha_2):=\{(\alpha_1(x)+\alpha_2(y),\alpha_2(x)+\alpha_1(y))\mid x,y\in L\}
\end{equation}
\begin{equation}	
\label{exactseqstr2}
Z(\beta_1,\beta_2):=\{(u,v)\in M\mid \beta_1(u)+\beta_2(v)=\beta_2(u)+\beta_1(v)\}.
\end{equation}

	\end{defn}
	In \S \ref{sectveryshort}, Proposition \ref{iso9}, we show that a sequence $0 \stackbin[0]{0}{\rightrightarrows}L\stackbin[g]{f}{\rightrightarrows}M\stackbin[0]{0}{\rightrightarrows}0$ in $\bm2$ is  \strgly exact if and only if $\phi=(f,g)$ is an isomorphism in $\bm2$ (\ie is invertible). Moreover, we find that the  group $\Aut_\bm2(M)$ of automorphisms of an object in $\bm2$ is  the semi-direct product of the $2$-group  of decompositions $M=M_1\times M_2$ of $M$ as a product, by the group of automorphisms of the $\B$-module $M$. In \S\S \ref{sectmeaningst}, \ref{sectepistrong},  we compare \strg exactness with the categorical notions of   monomorphisms and epimorphisms in $\bm2$. The conclusion is that while being an epimorphism is equivalent to  the \strg exactness of $L\stackbin[g]{f}{\rightrightarrows}M\stackbin[0]{0}{\rightrightarrows}0$,  being a monomorphism is a more restrictive notion than claiming the \strg exactness of  $0\stackbin[0]{0}{\rightrightarrows}L\stackbin[g]{f}{\rightrightarrows}M$ (\cf~Example \ref{example not mono} in \S \ref{sectmeaningst}).
	
Example \ref{notmonoepi} displays a critical defect of $\bm2$ by showing that some specific morphism in the category does not admit the epi-mono factorization ``monomorphism $\circ$ epimorphism". Additionally, while a submodule $N\subset M$ of a $\B$-module $M$ defines a natural covariant functor 
$$
M/N:\bm2\longrightarrow {\rm Bmod},\quad M/N(X):=\{(f,g)\in \Hom_{\bm2}(M,X)\mid f(x)=g(x)\qqq x\in N\}
$$ this functor is in general not representable in $\bm2$. The general theory of \cite{barr}  suggests that these defects  are eliminated by  enlarging  $\bm2$ to the Eilenberg-Moore category of  the comonad $\perp$. This extension turns out to be exactly the category $\b2$ of $\bs$-modules which was first introduced in the thesis of S.~Gaubert \cite{Gaubert} (see also \cite{Akian}). We refer to  \cite{Rowen} for an update account of this approach. 
	 \begin{defn} We let $\b2$ be the category of $\B$-modules endowed with an involution $\sigma$. The morphisms in $\b2$ are the morphisms of $\B$-modules commuting with $\sigma$, \ie equivariant for the action of $\Z/2\Z$. 
\end{defn} 
It should be noted from the start that the category $\b2$ is simply the category $\bmod$ in the topos of sets endowed with an involution and as such it shares with $\bmod$ most of its abstract categorical properties. Moreover 
$\b2$ is the same as the category of $\bs$-modules where $\bs$ is the semiring which is the symmetrization of $\B$ introduced in \cite{Gaubert}, p 71 (see also \cite{Akian}, definition 2.18). One has $\bs=\B\oplus s \B$, $s^2=1$ and the element $p=1+s$ fulfills $p^2=p$ and plays the role of ``null" element since $N=\{0,p\}$ is an ideal in $\bs$.
 
A main finding of this paper  is that the Eilenberg-Moore category $\b2$  falls  within the framework of the homological categories developed by M. Grandis in \cite{gra1}, where the theory of lattices of subobjects in an ambient category is tightly linked to the existence and coherence of fundamental constructions in  homological algebra. The astonishing outcome is that even in the simplest case of the category $\b2$, the resulting homological algebra is highly non-trivial and gives rise to a computable Ext functor.\newline 
In \S\ref{sectnullmorphisms} we prove  the fundamental result  stating   that the pair $(\b2,\cN)$, where $\cN$ is the ideal (closed as in \S 1.3.2 of \cite{gra1}) of the null morphisms  $f$ in $\b2$ (\ie  $\sigma\circ f=f$)  is a semiexact homological category.
\begin{thm}\label{homologicalintro}
The pair given by the category $\b2$ and the null morphisms: $\cN\subset \Hom_\b2(L,M)$
  forms a homological category in the sense of \cite{gra,gra1}. 
\end{thm}
 It turns out that $\cN=\cN(\cO)$ \ie $\cN$ is generated by the class $\cO$ of  null objects in $\b2$ \ie the objects $N$ such that $\id_N\in \cN$ (equivalently the $N$'s for which the involution $\sigma_N$ is the identity).  This class of objects is stable under retracts.
   The pair $(\b2,\cN)$ is a semiexact category and as such is provided with kernels and cokernels with respect to the ideal $\cN$.
     \begin{defn}\label{ker&coker intro} For $h\in \Hom_\b2(L,M)$  one sets
 $\Ker(h):=h^{-1} (M^\sigma) $
endowed with the induced involution.  One also defines $\coker(h):=M/\!\!\sim$ ~where
 \begin{equation}\label{cokernelb2intro}
 b\sim b'\iff f(b)=f(b')~\&~\forall X, \forall f\in \Hom_\b2(M,X):~ {\rm Range}(h)\subset \Ker(f).
\end{equation}
\end{defn}
The homological structure of $\b2$ allows for the definition of subquotients and their induced morphisms, for the construction of the connecting morphism and for the introduction of the homology sequence associated to a short exact sequence of chain complexes over $\b2$.\newline By applying the functor $y_\B(-)=\Hom_\bm2(\B,-)$,
the category $\bm2$ becomes  a full subcategory of $\b2$ and the new objects in this latter category provide both a factorization of any morphism as ``monomorphism $\circ$ epimorphism" and a representation of the functors $M/N$. 
The notion of cokernel in  $\b2$ gives an interpretation of the coequalizer $M \stackbin[\gamma_2]{\gamma_1}{\rightrightarrows}\coequ f_{(2)}
$ of the sequence \eqref{defncokergenintro} as a cokernel.
The notion of \strg exactness continues to make sense in $\b2$. A sequence   $L\stackrel{f}{\to} M\stackrel{g}{\to}N$ in $\b2$ is {\em \strgly exact} at $M$ when 
$
	{\rm Range}(f)+M^\sigma=\ker(g) 
$. In $\b2$,  not all subobjects $E\subset F$ of an object $F$, with $F^\sigma\subset E$, qualify as kernels of morphisms. This fact suggests to introduce the following new  notion of exactness which turns out to be the same as the notion of exactness given in \cite{gra1} for semiexact categories.
\begin{defn}\label{exactnessb2bisintro}  $(i)$~For $f\in \Hom_\b2(L,M)$, we let the {\em normal image}  $\imm(f)\subset M$ be the kernel of its cokernel.\newline
  $(ii)$~A sequence  of $\b2$: $L\stackrel{f}{\to} M\stackrel{g}{\to}N$ is {\em  exact} at $M$ when $\imm(f)=\ker(g)$.\newline
  $(iii)$~A subobject $E\subset F$ of an object in $\b2$ is {\em normal} if $E$ is the kernel of some morphism in $\b2$.
 \end{defn}
  Note that \strg exactness  implies exactness in $\b2$, and the definition of the normal image of a morphism as the kernel of its cokernel corresponds to the standard definition of the image of a morphism in an abelian category.  By applying these concepts we find
      \begin{prop}\label{tentative1intro} Let $f\in\Hom_\b2(L,M)$. 
 The sequence 
$$
 	0 \to \Ker(f)\to L\stackrel{f}{\to}M\to \coker(f)\to 0.
 $$
 is exact in $\b2$.
\end{prop}
The slight difference between the notions of \strg exactness and exactness in $\b2$ is due to the fact that not all subobjects $E\subset F$, with $F^\sigma\subset E$, are normal. For this reason it is crucial to determine the ``closure" of a subobject in $\b2$ \ie the intersection of all kernels containing this subobject or equivalently find the normal  image $\imm(E)$ of the inclusion $E\subset F$.  This part is developed in Proposition \ref{existphi1}. We report here the output of this computation. Let $p(x):=x+\sigma(x)$ be the projection on null elements.
 \begin{prop}\label{existphi1intro} Let $E\subset F$ be a subobject in $\b2$, with $F^\sigma\subset E$. Then for $\xi \in F$ one has  $\xi \in \imm(E)$ (the normal image of the inclusion) if and only if there exists a finite sequence $a_0,a'_0,a_1,a'_1,\ldots, a_n,a'_n$ of elements of $E$ such that 
 $
 \xi=a_0+\xi,\ p(a_0)=p(a'_0),\ a'_0+\xi=a_1+\xi,\ p(a_1)=p(a'_1), \ldots, p(a_n)=\xi+\sigma(\xi)
 $.
  \end{prop} 
  A corresponding explicit description of the cokernel of a morphism in $\b2$ is given in Proposition \ref{compute cokernel} and states as follows
  \begin{prop}\label{compute cokernel intro} Let $E\subset F$ be a subobject in $\b2$ containing $F^\sigma$. Then the cokernel of the inclusion $E\subset F$ is the quotient of $F^\sigma\cup E^c$ ($E^c:=F\setminus E$), by the smallest equivalence relation such that 
   \begin{equation}\label{simple add kernel intro}
   \xi \in E^c,\ u,v\in E,\ p(u)=p(v)~\Longrightarrow~ \xi +u\sim \xi +v.
  \end{equation} 
   The cokernel map is the quotient map on $E^c$ and the projection $p$ on $E$. 
   \end{prop}


In \S\ref{section analogy operators} we develop a second fundamental topic of this paper namely the analogy between the category $\b2$ and the category of linear operators in Hilbert space. This parallel relies on a duality property  holding in $\b2$ which we shall now review. Let $\yb: \bmod \to \b2$ be the functor acting on objects as $M \mapsto (M^2,\sigma_M)$, where $\sigma_M(x,y) = (y,x)$ and on morphisms $f:M\to N$ as $\yb(f)(x,y)=(f(x),f(y))$. One defines the orthogonal object of a subobject $E$ of an object $F$ in $\b2$ by implementing the following natural pairing connecting $F$ and $F^*:=\Hom_\b2(F,\yb\B)$: 
\[
\langle x,y\rangle_\sigma:=y(x)\qqq x\in F,y\in F^*, \ \ E^\perp:=\{y\mid  \langle x,y\rangle_\sigma \in (\yb\B)^\sigma,~\forall x\in E\}.
\] 
Then the following fact holds 
\begin{prop}\label{bidualcoker intro} $(i)$~Let $E\subset F$ be a subobject of the object $F$ of $\b2$. Then the least normal subobject $\imm(E)$ containing $E$ satisfies: $\imm(E)=(E^\perp)^\perp$.\newline
$(ii)$~Let $\phi\in \Hom_\b2(E,F)$. Then the dual of the cokernel of $\phi$ is canonically isomorphic to the kernel of $\phi^*\in \Hom_\b2(F^*,E^*)$ \ie  $\coker(\phi)^*=\Ker(\phi^*)$.
 \end{prop}
 In the paper we construct a large class of objects $X$ of $\b2$ such that for any morphism $\phi\in \Hom_\b2(X,Y)$ the nullity of $\ker(\phi)$ is equivalent to the injectivity of $\phi$. More precisely we obtain the following
\begin{thm}\label{null kernel thm intro} Let $M$ be a finite object of $\bmod$ whose dual $M^*$ is generated by its minimal non-zero elements. Then, for any object $N$ of $\b2$ a morphism $\phi\in \Hom_\b2(\yb M,N)$ is a monomorphism if and only if its kernel is null.	
\end{thm}

In \S\ref{homological grandis}, we finally apply to the semiexact, homological category $\b2$ the general theory developed in \cite{gra1}. The key notion is the analogue of the ``kernel-cokernel pairs'' in additive categories \cite{Buhler}, which is introduced in \S1.3.5 of \cite{gra1}. In our context such short \exxx sequences are given by a pair of maps $$L\stackrel{m}{\rightarrowtail}M\stackrel{p}{\xtwoheadrightarrow{}}N:~~ m= \sker(p)\ \& \  p= \scoker(m)$$ 
  where  
  $L=\Ker(p)$ is a normal subobject of $M$ and $N=\coker(m)$ is a normal quotient of $M$. Notice that the definition of a short \exxx sequence  as above is more restrictive than requiring the exactness of  $0\to L\stackrel{m}{\to}M\stackrel{p}{\to}N\to 0$ and the two notions become equivalent if one requires the maps $m,p$ to be exact (Proposition \ref{shortex}).  
  
  Connected sequences of functors and satellites exist for semiexact categories and their limit-construction is given as pointwise Kan extensions. For the definition of a cohomology theory of $\rma$-modules, we are mainly interested in the construction of the right satellite of the left exact functor defined by the internal Hom functor: $\Homi(M,-)$.  In an abelian category this functor is known to be left exact  and one obtains
the Ext functors as right derived functors. In our context, the right satellite is constructed as a left Kan extension (\cf  \cite{MacLane} p. 240). The great advantage to work within a semiexact, homological category $\mC$ (such as $\b2$) is that one can construct a homology sequence for any short \exxx sequence of chain complexes over $\mC$. This homology sequence is always of order two, meaning that the composition of two consecutive morphisms is null, and a very simple condition on the middle complex gives a partial exactness property for the homology sequence, sufficient for studying the universality of chain homology. The main obstacle for a straightforward application of the results of \cite{gra1} is that the categories we work with do not satisfy the modularity requirement which allows one  to bypass the role of the difference to relate for instance the injectivity of a morphism with the nullity of its kernel. This drawback inputs a further obstruction  to a direct proof of the main condition $(a)$ of  Theorem 4.2.2 of \cite{gra1}, which entails, if satisfied, to compute satellite functors using semi-resolutions. The key result of \S\ref{homological grandis}  is a sleek solution of this problem stating that for the right satellite of  $\Homi(M,-)$ one has the following result
   \begin{thm}\label{reduction toendo intro} Let $\iota:I'\stackrel{i'}{\rightarrowtail{}} I \stackrel{i''}{\xtwoheadrightarrow{}} I''$ be a short \exxx sequence in $\b2$ with the middle term $I$ being  both injective and projective. 
Then the functor $F:=\Homi_\b2(I'',-)$ satisfies the condition $a)$ of \cite{gra1} with respect to an arbitrary morphism of short \exxx sequences $c\to \iota$, provided this holds for all  endomorphisms of $I'$ in $\b2$.
 \end{thm}

 The final \S\ref{sectcokerdiag}  is a crucial last preliminary before attacking the construction of  sheaf cohomology for sheaves of $\rma$-modules on a topos. We apply Theorem \ref{reduction toendo intro} to determine explicitly the Ext functor associated to the  key example provided by the cokernel of the diagonal.  After reviewing the role of the cokernel of the diagonal in ordinary sheaf cohomology with coefficients in abelian groups, we analyze the diagonal in the context of $\B$-modules and investigate the simplest case provided by the diagonal for the module $\B$ itself. The cokernel of the diagonal is given in this simplest example by a computable object $Q$ of $\b2$, which is described in details in \S \ref{sectrepfunctor}. 
It fits in the following short \exxx sequence
\begin{equation}\label{short diag exact intro}
	\alpha:K=\ker(\scoker(\yb\Delta))\stackrel{\subset}{\rightarrowtail} \B^2\times \B^2\stackrel{\scoker(\yb\Delta)}{\xtwoheadrightarrow{}}Q
\end{equation}
and the main point is to detect the obstruction to lift back from $Q$ to $\B^2\times \B^2$. Thus we  consider the representable functor $F:=\Homi_\b2(Q,-)$, viewed as a covariant endofunctor in $\b2$  (using the natural internal Hom) and we compute its right derived functor.
 Lemma \ref{condition to lift} shows that the short \exxx sequence \eqref{short diag exact intro},  gives rise, after applying the functor  $F$, to a non-null cokernel for the map  $F(\scoker(\yb\Delta))$. This yields naturally an element of the satellite functor $SF(\ker(\scoker(\yb\Delta)))$. In Lemma \ref{coker compute} we compute $\coker(F(\scoker(\yb\Delta)))$: its graph is displayed in Figure \ref{graphcokernel}. In \S \ref{subsubcorrespondence} we investigate the correspondence between endomorphisms of $K:=\ker(\scoker(\yb\Delta))$ and of $Q$ as provided by the two maps from the endomorphisms $\End(\alpha)$ of the short \exxx sequence $\alpha$ of \eqref{short diag exact intro} \begin{equation}\label{quotient map intro}
	\End_\b2(K)\stackrel{{ res}}{\leftarrow} \End(\alpha)\stackrel{{ quot}}{\to} \End_\b2(Q)
\end{equation}
 We   show that this correspondence  is multivalued. The remaining issue then is to prove that after applying the functor $F$ it becomes single valued so that one obtains functoriality. This result is achieved in \S \ref{subsecfunctact}  (Theorems \ref{sequence alpha}, \ref{reduction toendo cor2}) and provides one with the following third fundamental result of this paper  
 \begin{thm}\label{sequence alpha intro} The short \exxx sequence $\alpha:\Ker(\scoker(\yb\Delta))\stackrel{\subset}{\rightarrowtail} \B^2\times \B^2\stackrel{\scoker(\yb\Delta)}{\xtwoheadrightarrow{}}Q$ satisfies condition $a)$ of \cite{gra1} with respect to the functor $F=\Homi_\b2(Q,-)$.  	
 The right satellite functor $SF$  is non-null and $SF(K)=\coker(F(\scoker(\yb\Delta)))$.\end{thm}
This result shows that the Ext functor is non-trivial already in the case of finite objects of $\b2$. For this reason the present paper provides a strong motivation for a systematic development of homological algebra for categories of sheaves of modules over semirings of characteristic one over a topos. 
 
\vspace{0.1in}
{\bf{Acknowledgement.}}
 We are grateful to Marco Grandis and Stephane Gaubert for several helpful comments in the elaboration of this paper.

\section{$\B$-modules and algebraic lattices}\label{sect bmodules}

A $\B$-module $M$ is a commutative monoid written additively such that $x+x=x$ $\forall x\in M$. The idempotent addition is often denoted by $x\vee y$ instead of $x+y$. The action of $\B$ on $M$ is set to be: $0\cdot x=0$  and $1\cdot x=x$, $\forall x\in M$. In \S\ref{sect bmodules1} we recall why $\B$-modules correspond to join-semilattices with a least element.  In \S\ref{sectbmodual} we review basic material on duality \cite{Cohen}, and \S\ref{galois connect} surveys the classical equivalence with algebraic lattices. In \S\ref{sect dec morphisms} we implement Galois connections \cite{Cohen} to  give a canonical decomposition of morphisms in $\bmod$. This construction will play an important role in the analogy with operators developed in \S\ref{section analogy operators}.  Finally, in \S\ref{section radical} we discuss for future application in Theorem \ref{kernel and injective}, the condition on a $\B$-module $E$  that ensures that its dual $E^*$ is generated by its minimal non-zero elements.

\subsection{$\B$-modules as partially ordered sets}\label{sect bmodules1}
 
The traditional terminology for $\B$-modules is of join-semilattices with a least element (written $0$), where ``upper-semilattice" is also used in place of join-semilattice. The join of two elements is denoted $x\vee y$ and the morphisms $f:S\to T$ in the category satisfy $f(x\vee y)=f(x)\vee f(y)$ and $f(0)=0$. The following statement shows that this category is the same as that of $\B$-modules.
\begin{prop}\label{Bmod1}$(i)$~Let $M$ be a $\B$-module, then the condition
\begin{equation}\label{Bmod2}
	x\leq y\iff x\vee y =y
\end{equation}
defines a partial order on $M$ such that any two elements of $M$ have a join, \ie there is a smallest element $z=x\vee y$ among the elements larger than $x$ and $y$, and moreover $0$ is the smallest element.\vspace{.005in}

Conversely:\newline
	$(ii)$~Let $(E,\leq)$ be a partially ordered set with a smallest element in which any two elements $x,y\in E$ have a join $x\vee y$. Then $E$ endowed with the addition $(x,y)\mapsto x\vee y$ is a $\B$-module.
	\newline
	$(iii)$~The two functors defined in $(i)$ and $(ii)$ are inverse of each other.
\end{prop}
 \proof $(i)$~Evidently $x\leq x$ since $x\vee x=x$. Assume $x\leq y$ and $y\leq z$ then it follows from the associativity of the operation that
 $
 x\vee z= x\vee (y\vee z)=(x\vee y)\vee z=y\vee z=z.
 $
 Thus the relation $\leq$ is transitive. Moreover if $x\leq y$ and $y \leq x$ one has $x=x \vee y =y$ thus $\leq$ is a partial order relation. One also has $x\leq x\vee y$ since $x\vee (x\vee y)=x\vee y$. Let $z\in M$ be such that $x\leq z$ and $y\leq z$. Then 
 $
 (x\vee y)\vee z=x\vee (y\vee z)=x\vee z=z,
 $
 thus $z\geq x\vee y$ and hence $x\vee y$ is the smallest among the elements larger than $x$ and $y$. Finally $0$ is the smallest element since $0\vee x=x$ $\forall x\in M$.\newline
 $(ii)$~We show that in $(E,\le)$ the operation $(x,y)\mapsto x\vee y$ is associative. Indeed, let $x,y,z\in E$ and let $t=(x\vee y)\vee z$. Let $u$ with $u\geq x$, $u\geq y$, $u\geq z$ then $u\geq x\vee y$ and thus $u\geq t$. Thus $t$ is the smallest element in the set of $u$'s satisfying the above conditions. The same holds for $x\vee (y\vee z)$ and thus one gets the equality. The operation $(x,y)\mapsto x\vee y$ is commutative by construction and idempotent, thus one gets a $\B$-module with the smallest element as the $0$. \newline
	$(iii)$~Let $M$ be a $\B$-module and $\leq$ the associated partial order on $M$. We have shown that the join of any pair in $M$ is $x\vee y$ and thus the $\B$-module associated by $(ii)$ to the partial order $\leq$ is $(M,\vee)$. Let $E$ be a partially ordered set as in $(ii)$ and $\vee$ the join. Then one has for the partial order of $E$: 
	$
	x\leq y\iff x\vee y=y
	$.  
This shows that the partial order for the associated $\B$-module $(E,\vee)$ is the same as the partial order of $E$.\endproof
In terms of partial orders the morphisms in the category of $\B$-modules are increasing maps of semi-lattices preserving the sup and the minimal element.
 
 \subsection{Duality in $\bmod$}\label{sectbmodual}

This section recalls in the simple case of the semiring $\B$ the duality results of \cite{Cohen}. We give complete proofs for convenience.
The  functor $\Hom_\B$ with one of the two entries fixed,
determines an endofunctor in $\bmod$. The simplest case to consider is the duality 
determined by the contravariant endofunctor
\begin{equation}\label{the dual}
	M\longrightarrow M^*:=\Hom_\B(M,\B).
\end{equation}
The terminology ``ideal" is often used in the context of partially ordered sets to denote a hereditary submodule. We prefer to adopt the latter to avoid  confusion  with the theory of semirings.
\begin{defn}\label{heredit} $(i)$~A subset $H\subset E$ of a partially ordered set $E$ is hereditary if \[
y\in H\, \& \, x\leq y~\Longrightarrow~ x\in H.
\]
$(ii)$~A subset $H\subset M$ of a $\B$-module $M$ is hereditary if it is so for the partial order \eqref{Bmod2} on $M$.
	\end{defn}
	Next proposition is easily deduced from \cite{Cohen} (\cf~\S 4.2 and Corollaries 37 and 40). In particular $(iii)$ implies that $\B$ is injective in the category of $\B$-modules
\begin{prop}\label{inj1}$(i)$~The map $M^*\to \{H\subset M\}$, $\phi\mapsto\phi^{-1}(0)$ is a bijection between $M^*$ and the set of non-empty hereditary submodules $H\subset M$.\newline
$(ii)$~A subset $H\subset M$ is an hereditary submodule if and only if it is a filtering union of intervals $I_x:=\{y\in M\mid y\leq x\}$.\newline
$(iii)$~Let $N\subset M$ be a submodule, then the restriction map $M^*\to N^*$ is surjective.	
\newline
$(iv)$~The duality between $M$ and $M^*$  is separating.\newline
$(v)$~Let $E,F$ be arbitrary objects  in $\bmod$ and $f\in \Hom_\B(E,F)$, then 
   \begin{equation}\label{orthogonality1}
  	f(\xi)=f(\eta)\iff 
  	 \langle\xi,u\rangle=\langle\eta,u\rangle\qqq u\in {\rm Range}(f^*).
  \end{equation}
\end{prop}
\proof $(i)$~For $\phi\in M^*$ the subset $\phi^{-1}(0)$ is a hereditary submodule of $M$. Conversely, if $H\subset M$ is a non-empty hereditary submodule of $M$ the equality $\phi(x)=0 \iff x\in H$ defines a map $\phi:M\to \B$ satisfying $\phi(x\vee y)=\phi(x)\vee\phi(y)$. In fact, either $x,y\in H$ in which case $x\vee y\in H$ and the two sides are $0$, or say $x\notin H$ and then $x\vee y \notin H$ since $H$ is hereditary and then both sides are equal to $1$.\newline
$(ii)$~Let $H\subset M$. If $H$ is a hereditary submodule one has $H=\cup_{x\in H} I_x$ and the inclusion $I_x\cup I_y\subset I_{x\vee y}$ shows that the union is filtering. Conversely, any interval $I_x$ is a hereditary submodule, any union of hereditary subsets is hereditary and any filtering union of submodules is a submodule. \newline
$(iii)$~Let $N\subset M$ be a submodule. Let $\phi\in N^*$. Let $$
H\subset M,\ \ H=\cup_{x\in N\mid \phi(x)=0} \, I_x=\{y\in M\mid \exists x\in N,\, y\leq x\,\&\, \phi(x)=0\}.
$$
Then $H\subset M$ is a filtering union of intervals and hence a hereditary submodule, thus 
 the equality $\psi(x)=0 \iff x\in H$ defines a map $\psi:M\to \B$ and $\psi\in M^*$. Moreover one has $N\cap H=\phi^{-1}(0)$ and thus the restriction of $\psi$ to $N$ is equal to $\phi$.\newline
$(iv)$~Let $x,y\in M$, such that $\phi(x)=\phi(y)$ $\forall\phi\in M^*$, then $x\in I_y$ and $y\in I_x$ so $x=y$. \newline
$(v)$~The duality of $F$ with $F^*$ is separating by  $(iv)$, and this shows the implication $\Leftarrow$. The implication $\Rightarrow$ is straightforward.\endproof

Note that by construction the map $\psi: M\to \B$ defined in the proof of $(iii)$ is the largest extension of $\phi: N\to \B$ to $M$ since for any extension $\psi'$ one has $\psi'(x)=0$ for any $x\in H$ for $H$ as in $(iii)$. Applying $(iii)$ to $N=0$ one thus obtains the largest element of $M^*$, as the map $\alpha:M\to \B$ such that $\alpha(x)=1$, $\forall x\neq 0$.\vspace{.05in}

Next proposition shows that $\B$ is a cogenerator in the category of $\B$-modules.

\begin{prop}\label{inj3}$(i)$~The product $\B^X$ (resp. the coproduct $\B^{(X)}$) of any number of copies of $\B$ is an injective (resp. projective) object in the category of $\B$-modules.\newline
$(ii)$~Let $M$ be a $\B$-module then $M$ is isomorphic to a submodule of a product $\B^X$ and there exists a surjective map of $\B$-modules of the form $\B^{(X)}\to M$.	
\end{prop}
\proof $(i)$~A morphism $\phi:N\to \B^X$ can be equivalently described by a family of morphisms $\phi_x\in \Hom_\B(N,\B)=N^*$. Since when $N\subset M$ each $\phi_x$ admits an extension to $M$, one obtains that the naturally associated restriction map $\Hom_\B(M,\B^X)\to\Hom_\B(N,\B^X)$ is surjective. Similarly, a morphism $\phi:\B^{(X)}\to M$ is equivalent to a family of morphisms $\phi_x\in \Hom_\B(\B,M)=M$ and one derives in this way that $\B^{(X)}$ is projective.\newline
$(ii)$~Let $\iota: M\to (M^*)^*$ be defined by $\iota(x)(\phi):=\phi(x)$. One has
$$
\iota(x\vee y)(\phi)=\phi(x\vee y)= \phi(x)\vee \phi(y)=\iota(x)(\phi)\vee \iota(y)(\phi).
$$
Then $\iota$ defines a morphism $M\to\B^{M^*}$. This morphism is injective since for $x\in M$ the element $\phi_x\in M^*$ defined by $\phi_x^{-1}(0)=I_x=\{y\mid y\leq x\}$ satisfies the rule $\phi_x(y)=\phi_x(x)\iff y\leq x$. Thus if $x\neq y$ one has either that $\phi_x(y)\neq\phi_x(x)$ or that $\phi_y(y)\neq\phi_y(x)$. The existence of a surjective map $\B^{(X)}\to M$ is clear by taking $X=M$ and using the addition. \endproof 

\begin{rem}\label{worry} {\rm It might seem contradictory that the lattice $\cL$ of subspaces $E\subset V$ of a finite dimensional vector space $V$ over a field $K$ can be embedded as a submodule of a product $\B^X$. To see why this construction is meaningful one first considers the orthogonal $E^\perp\subset V^*$ and this step replaces the lattice operation $E_1\vee E_2$ by the dual operation \ie the intersection:
$
(E_1\vee E_2)^\perp= E_1^\perp \cap E_2^\perp.
$
In this way $\cL$ embeds as a sub-lattice of the lattice of subsets of the set $V^*$ endowed with the operation of intersection. Finally, by using the map from a subset to its complement, this lattice is the same as the product $\B^X$, for $X=V^*$. 	
}\end{rem}
We end this subsection with the simple form of the Separation Theorem of \cite{Cohen} in the case of $\B$-modules. It gives the following analogue of the Hahn-Banach Theorem. Since the proof of \cite{Cohen}  takes a simple form in our context we include it here below.
\begin{lem} \label{imlem} 
Let $N$ be a $\B$-module and $E\subset N$ be a submodule: we denote $\iota:E\to N$ the inclusion. Then  
\begin{equation}\label{separ}
x\in  E\iff \phi(x)=\psi(x) \qqq  \phi,\psi \in \Hom_\B(N,\B)~{\rm s.t.}~ \phi\circ \iota=\psi \circ \iota.
\end{equation}
\end{lem}
\proof Let $\xi\in N$, $\xi\notin E$. It is enough to show that there are two morphisms $\phi,\psi\in \Hom_\B(N,\B)$ which agree on $E$ but take different values at $\xi$.  Let $\phi=\phi_\xi$ so that $\phi^{-1}(0)=\{\eta\mid \eta\leq \xi\}$. We let 
$$
F=\{\alpha \in N\mid \exists \eta\in E,\ \alpha\leq \eta\leq \xi\}.
$$
One has $\alpha,\alpha'\in F\Rightarrow \alpha\vee \alpha'\in F$ since for $\eta,\eta'\in E$, $\alpha\leq \eta\leq \xi$, $\alpha'\leq \eta'\leq \xi$, $\eta\vee \eta'$ is an element of $E$, and one has  $\alpha\vee \alpha'\leq \eta\vee \eta'\leq \xi$. Moreover by construction one has $\alpha \in F\ \& \ \beta\leq \alpha\Rightarrow \beta\in F$. It follows that the equality $\psi(\alpha)=0\iff \alpha \in F$ defines an element  $\psi\in\Hom_\B(N,\B)$.  Moreover $\xi \notin F$ since $\xi \notin E$, thus $\psi(\xi)=1$.  One has $\phi(\xi)=0\neq \psi(\xi)=1$. Since $F\cap E=[0,\xi]\cap E$ one has $\phi(\eta)=\psi(\eta)$, $\forall \eta \in E$. \endproof
Lemma \ref{imlem} allows one to determine the epimorphisms in $\bmod$ and one gets the following corollary.
\begin{prop}\label{monoepi} Let $f\in \Hom_\B(M,N)$ be a morphism in $\bmod$.\newline
$(i)$~$f$ is a monomorphism if and only if the underlying map is injective.\newline
$(ii)$~$f$ is an epimorphism if and only if the underlying map is surjective.	
\end{prop}
\proof $(i)$~This follows from the identification $M=\Hom_\B(\B,M)$ which associates  to any $x\in M$ the unique morphism $\B\to M$  sending $1$ to $x$.\newline
$(ii)$~This follows from Lemma \ref{imlem}.\endproof

\subsection{Algebraic lattices}\label{galois connect}

In this subsection we describe the well-known equivalence between the category $\bmod$ and the category of  algebraic  lattices. The notion of  algebraic  lattice is recalled in the following
\begin{defn} \label{defn complete lattice} $(i)$~A complete lattice $\cL$ is a partially ordered set in which all subsets have both a supremum (join) and an infimum (meet).	\newline
$(ii)$~An element $x\in \cL$ of a complete lattice is {\em compact} (one writes $x\in K(\cL)$) if the following rule applies 
$$
x\leq \vee_A y~\Longrightarrow~ \exists F\subset A,~F\ {\rm finite,~ s.t.}~ x\leq \vee_F y.
$$
$(iii)$~A  lattice $\cL$ is {\em algebraic} iff it is  complete and every element $x$ is the supremum of compact elements below $x$. \end{defn}

Examples of algebraic  lattices are:
\begin{itemize}
\item Subgroups$(G)$ for a group $G$, where the compact elements are the finitely generated subgroups
\item Subspaces$(V)$ for a vector space $V$, where the compact elements are the finite dimensional subspaces
\item Id$(R)$ \ie ideals for a ring $R$, where the compact elements are the finitely generated ideals
\end{itemize}

A morphism  of algebraic (complete) lattices is a compactness-preserving complete join homomorphism \ie  
\begin{equation}\label{compactness-preserving}
	f:S\to T~{\rm s.t.}~ f(\vee_A y)=\vee_A f(y)\qqq A\subset S,\ \ f(K(S))\subset K(T).
\end{equation}
To understand this notion we show
\begin{lem}\label{algebraic lattice inj} Let $f:S\to T$ be a morphism of algebraic lattices. Then $f$ is injective if and only if its restriction $f:K(S)\to K(T)$ is injective.	
\end{lem}
\proof We assume that the restriction to $f:K(S)\to K(T)$ is injective. Let then $s,s'\in S$ such that $f(s)=f(s')$. Let $A=\{x\in K(S)\mid x\leq s\}$, $B=\{y\in K(S)\mid y\leq s'\}$. One has $s=\vee_A x$, $s'=\vee_B y$, and $A,B\subset K(S)$. One has $\vee_A f(x)=\vee_B f(y)$ and the compactness of $f(y)$ for $y\in B$ shows that there exists $x\in A$ with $f(x)\geq f(y)$. This means $f(x)+f(y)=f(x)$ and hence, as $x+y,x\in K(S)$, that $x+y=x$, \ie $y\leq x$. This proves that $A$ and $B$ are cofinal and thus $s=s'$. \endproof

The notion of compactness of elements is best understood in terms of the small category  $\cC(\cL)$ whose objects are the elements of $\cL$ and there is a single morphism from $x$ to $y$ iff $x\leq y$, while otherwise $\Hom(x,y)=\emptyset$. In these terms, an element $x\in \cL$ is compact iff it is {\em finitely presentable} in the sense that the functor $\Hom(x,-):\cC(\cL)\longrightarrow \Se$ commutes with filtered colimits. In this language and for any ordered set $\cL$ one has the equivalence

\vspace{.1in}

\centerline{$\cL$ is a complete algebraic lattice $\iff$ $\cC(\cL)$ is locally finitely presentable}
\vspace{.1in}

By definition a locally finitely presentable category has all small colimits and any object is a filtered colimit of the canonical diagram of finitely presentable objects mapping into it.
There is a classical equivalence  of categories between join-semilattices with a least element, \ie the category $\bmod$ and the category $\cA$ of complete algebraic lattices with compactness-preserving complete join-homomorphisms. This equivalence is a special case of the Gabriel-Ulmer duality \cite{GU}. 

To a $\B$-module  $M$ one associates the algebraic complete lattice $\id(M)$ of hereditary submodules of Definition \ref{heredit}. The join operation in the lattice is the hereditary submodule generated by the union. The join of the intervals $I_x:=\{z\in M\mid z\leq x\}$ and $I_y$ is the interval $I_{x\vee y}$.   Proposition   \ref{inj1} then shows that one gets an algebraic lattice and that the intervals  are the most general compact elements of $\id(M)$. To a morphism $f:M\to N$ in $\bmod$ one associates the morphism $\id(f):\id(M)\to \id(N)$ obtained by mapping  an hereditary submodule $I\subset M$ to the hereditary submodule generated by $f(I)\subset N$. This defines a compactness-preserving complete join-homomorphism but since in general $f(I\cap J)\neq f(I)\cap f(J)$, such homomorphism does not respect the intersections. By applying Proposition   \ref{inj1} one has a canonical identification of $\id(M)$ with $M^*:=\Hom_\B(M,\B)$ but this identification is not compatible with the join since 
$$
(\phi\vee \psi)^{-1}(0)=\phi^{-1}(0)\cap \psi^{-1}(0)\qqq \phi, \   \psi \in M^*.
$$
Moreover the functor $\id:\bmod\longrightarrow\cA$ is covariant while the endofunctor $M\mapsto M^*$ on $\bmod$ is contravariant. The precise link between these two functors is derived from the following notion of Galois connection: we refer to \cite{DP} Chapter VII  \cite{gra1} 1.1.3, and \cite{Gierz} to read more details.

\begin{defn}\label{galois connection}
 A covariant Galois connection $f \dashv g$ between two  partially ordered sets $E,F$ consists of two monotone  functions: $f : E \to F$ and $g: F \to E$, such that 
\begin{equation}\label{equ galois connection}
f(a) \leq b \iff  a \leq  g(b).	
\end{equation}
\end{defn}
Note that in fact one does not have to assume that the maps $f$ and $g$ are monotone since this property follows from \eqref{equ galois connection}. Indeed, we show that $x\leq y\Rightarrow f(x)\leq f(y)$. It follows from \eqref{equ galois connection} that  $f(x)\leq f(y)\iff x\leq g(f(y))$ and moreover using $f(y)\leq f(y)$ one derives $y\leq  g(f(y))$. Thus if $x\leq y$ one gets $x\leq g(f(y))$ and hence $f(x)\leq f(y)$.\newline 
There is a natural composition law for Galois connections:  given $f\dashv g$ with $f : E \to F$ and $h\dashv k$ with  $h: F\to G$, the pair $(h\circ f, g\circ k)$ is a Galois connection $h\circ f\dashv  g\circ k$, with $h\circ f:E\to G$. Indeed, one has for $a\in E$, $v\in G$ 
  $$
  (h\circ f)(a)\leq v\iff f(a)\leq k(v)\iff a\leq g(k(v))=(g\circ k)(v).
  $$
  Moreover a Galois connection $f\dashv g$ with  $f:E\to F$ automatically fulfills the equalities: $f=fgf$ and $g=gfg$ since $f(a)\leq f(a)\Rightarrow a\leq (g\circ f)(a)\Rightarrow f(a)\leq (f\circ g\circ f)(a)$, while $g(b)\leq g(b)\Rightarrow (f\circ g)(b)\leq b$ and taking 
  $b=f(a)$ one gets  $(f\circ g\circ f)(a)\leq a$ and thus $(f\circ g\circ f)(a)= a$.
  
For a Galois connection $f\dashv g$, $f$ is called the left adjoint of $g$ and $g$ is called the right adjoint of $f$. This terminology fits with  that of adjoint functors by considering the category $\cC(E)$.  The fact that $f$ is monotone means that it defines a covariant functor $\cC(f):\cC(E)\to \cC(F)$. The relation $f \dashv g$ means exactly that $\cC(f) \dashv \cC(g)$ as a pair of adjoint functors. 
Any homomorphism between complete lattices which preserves all joins  is the left adjoint of some Galois connection. This fact holds in particular for $\id(f):\id(M)\to \id(N)$ with the above notations (\ie with $f\in \Hom_\B(M,N)$) and one derives the following conclusion
\begin{prop}\label{galois connection1} Let $M,N$ be objects of $\bmod$ and $f\in \Hom_\B(M,N)$. Then, under the canonical identifications of partially ordered sets $\id(M)\simeq (M^*)^{\rm op}$ and $\id(N)\simeq (N^*)^{\rm op}$, the following pair defines a monotone Galois connection
$$
\id(f)\dashv f^*.
$$	
\end{prop}
\proof The map $f^*:N^*\to M^*$ is defined as composition with $f$, \ie $\psi\in N^*\mapsto \psi\circ f\in M^*$. One has for $\phi\in M^*$, $\psi\in N^*$
$$
\id(f)(\phi^{-1}(0))\subset \psi^{-1}(0)\iff f(\phi^{-1}(0))\subset \psi^{-1}(0).
$$
In turns this is equivalent to $\phi(x)=0 \Rightarrow \psi(f(x))=0$ \ie to $\phi \geq f^*\psi$. Since the identification $\id(M)\simeq (M^*)^{\rm op}$ is order reversing one gets the required equivalence of Definition \ref{galois connection}.\endproof 
The next statement is a corollary of  the theorem  asserting the equivalence of  $\bmod$ with the category $\cA$
\begin{cor}\label{galois connection2} Let $M$ be a $\B$-module, and $M^*=\Hom_\B(M,\B)$.
Then $(M^*)^{\rm op}$ is a complete algebraic lattice and the map which to $x\in M$ associates $\phi_x\in M^*$, with $\phi_x^{-1}(0)=[0,x]$ is an order preserving bijection of $M$ with the compact elements of $(M^*)^{\rm op}$.
\end{cor}
\proof Since the partially ordered set $(M^*)^{\rm op}$ is isomorphic to $\id(M)$ it is a complete algebraic lattice. The compact elements of $\id(M)$ are the intervals $[0,x]$ and the corresponding elements of $M^*=\Hom_\B(M,\B)$ are the $\phi_x$. \endproof 
Next, we compare the complete algebraic lattice $\id(M)$ with the bidual $(M^*)^*$. In general the latter is strictly larger than $\id(M)$ and there is a simple notion of complete join homomorphism between complete lattices, \ie  those morphisms which preserve arbitrary $\vee$. It is natural to compare them with normal positive maps in the von Neumann algebra context. Thus we let 
$$
N^*_{\rm norm}:=\{\phi\in \Hom_\B(N,\B)\mid \phi(\vee x_\alpha)=\vee(\phi(x_\alpha)\}.
$$
 \begin{prop}\label{bidual} $(i)$~The canonical map $\epsilon: M\to (M^*)^*$ given by evaluation $\epsilon(x)(\psi)=\psi(x)$ extends to an isomorphism $\tilde \epsilon: \id(M)\simeq (M^*)^*_{\rm norm}$. 	\newline
 $(ii)$~The $\B$-module $M$ is the submodule of $(M^*)^*_{\rm norm}$ of the compact elements of this complete algebraic lattice.
 \end{prop}
\proof $(i)$~The $\B$-module $N=M^*$ is a complete lattice. An element $\phi \in N^*$ is characterized by the hereditary submodule $J=\phi^{-1}(\{0\})$ and $\phi\in N^*_{\rm norm}$ iff $J$ is closed under arbitrary $\vee$. This fact holds iff $J=[0,n]$ for some element $n\in N$. One thus obtains a canonical order preserving bijection $(M^*)^*_{\rm norm}\simeq (M^*)^{\rm op}=\id(M)$ which associates to $\phi\in (M^*)^*_{\rm norm}$ the unique $\psi\in M^*$ such that $
\phi(\alpha)=0\iff \alpha\leq \psi
$. 
Let $x\in M$ and $\phi=\epsilon_x$. One has $\phi\in (M^*)^*_{\rm norm}$ and the associated $\psi\in M^*$ is such that $\alpha(x)=0\iff \alpha\leq \psi$ and thus $\psi=\phi_x$ with $\phi_x^{-1}(0)=[0,x]$ as in Corollary \ref{galois connection2}. This shows that the canonical order preserving bijection $(M^*)^*_{\rm norm}\simeq (M^*)^{\rm op}=\id(M)$ is the extension by completion of the  map $\epsilon: M\to (M^*)^*$ given by evaluation.\newline
$(ii)$~ follows from $(i)$ and the identification of $M$ as the submodule of compact elements in $\id(M)$. \endproof 

\begin{rem}\label{lattice implicit}{\rm The additive structure of a semifield $F$ of characteristic one is that of a lattice as a partially ordered set. Moreover, the morphisms of such semifields preserve not only the join $\vee$ but also the meet $\wedge$. Indeed, while for the natural partial order on $F$ given by $x\leq y\iff x+y=y$ the sum gives the $\vee$, the map $x\mapsto x^{-1}$ is an order reversing isomorphism of $F^*=F\setminus \{0\}$ onto itself since one has
$$
x\leq y\iff x/y+1=1\iff 1/y+1/x=1/x\iff 1/y\leq 1/x.
$$ 
In fact, one obtains in this way an equivalence of categories between the category of $\ell$-groups (\ie lattice ordered abelian groups) and the category of idempotent semifields. 	
}\end{rem}

\subsection{Decomposition of morphisms in $\bmod$}\label{sect dec morphisms}

The construction of the Galois connection reviewed  in \S\ref{galois connect} suggests that one may associate to a subset $E\subset F$ of an object in $\bmod$ a projection $q_E:F\to E$ defined by the formula of \cite{Cohen} Remark 14
  \begin{equation}\label{try proj}
  	 q_E(\xi):=\wedge \{a\in E\mid a\geq \xi\}.
  \end{equation}
  Here $E$ is assumed to be a lattice for the  partial order induced by $F$, and the $\wedge$ is taken with respect to this induced order.  For simplicity we restrict to the case of finite objects and obtain the following decomposition of morphisms which plays an important role in \S\ref{section analogy operators}.
  \begin{prop}\label{support} Let $f\in \Hom_\B(F,G)$ be a morphism of finite objects in $\bmod$. Let $S=\{z\in F\mid f(y)\leq f(z)\Rightarrow y\leq z\}$. Then 
  \begin{enumerate}
  \item $S$ is a lattice for the partial order induced by $F$.
  \item Formula \eqref{try proj} defines a surjective morphism $q_S\in \Hom_\B(F,S)$, for $S$ endowed with the operation $\vee$ induced by its order, and the lowest element $0_S$.
  \item The restriction of $f$ to $S$ defines an injective morphism $f\vert_S\in \Hom_\B(S,G)$.
  \item One has $f=f\vert_S\circ q_S$.
  \item The inclusion $\iota:S\to F$ preserves the $\wedge$:~ $\iota(x\wedge_S y)=\iota(x)\wedge \iota(y)$, $\forall x,y \in S$.  	
  \end{enumerate}  	
  \end{prop}
    \proof Since $F$ and $G$ are finite they are lattices and thus there exists a unique monotone map $g:G\to F$ such that $f\dashv g$ is a Galois connection. By construction one has, for $x\in F$ and $y\in G$,
    \begin{equation}\label{galconnect}
  	 f(x)\leq y\iff x\leq g(y).
  \end{equation} 
It follows from the general theory on Galois connections that $g(x\wedge y)=g(x)\wedge g(y)$ for all $x,y\in G$ and that  $c(x):=g(f(x))$ is a closure operator, \ie $c$ fulfills (using the equality $f=fgf$) the following facts
\begin{enumerate}
  \item $x\leq y\Rightarrow c(x)\leq c(y)$.
  \item $c(c(x))=c(x)\qqq x\in F$.
  \item $x\leq c(x) \qqq x\in F$.  	
  \end{enumerate}
  Since \eqref{galconnect} with $y=f(z)$ means that $f(x)\leq f(z)\iff x\leq c(z)$, one has  $S=\{z\in F\mid f(y)\leq f(z)\iff y\leq z\}=\{z\in F \mid c(z)=z \}$. Moreover, using the above relations one derives  
 \begin{equation}\label{closure op}
  c(x)\leq y\iff x\leq y\qqq x\in F, \ y\in S.
  \end{equation}
  This shows that $c\dashv \iota$ where $\iota:S\to F$ is the inclusion of $S$ endowed with the induced partial order in $F$. Thus $c$ is the left adjoint of $\iota$,  it only depends upon $S$ and is defined by 
  $$
  \iota^{-1}([x,\infty))=[c(x),\infty)\subset S\qqq x\in F.
  $$
  One has $[x\vee y,\infty)=[x,\infty)\cap [y,\infty)$ and thus in the partially ordered set $S$ one gets $$
  [c(x\vee y),\infty)=[c(x),\infty)\cap [c(y),\infty).
  $$
   This formula shows that any pair $a,b\in S$ has a least upper bound $a\vee_S b=c(a\vee b)\in S$  and that the map $c:F\to S$ fulfills $c(x\vee y)=c(x)\vee_S c(y)$, $\forall x,y\in F$. Thus since $S$ is finite it is a lattice and moreover $O_S:=c(0)$ is the smallest element of $S$ since $\iota^{-1}([0,\infty))=[c(0),\infty)$. Hence $S$ endowed with the operation $\vee_S$ and the zero element $0_S\in S$ is an object of $\bmod$ and $c:F\to S$ is a morphism in $\bmod$. Since $c\circ \iota=\id$ the morphism $c$ is surjective. By construction and using \eqref{try proj} one has $c=q_S$. Let $h=f\vert_S$ be the restriction of $f$ to $S$. We prove that $h\in \Hom_\B(S,G)$. First, one has $h(0_S)=f(c(0))=f(0)=0$ since $c=gf$ and $fgf=f$. Next, let $a,b\in S$ then 
   $$
   h(a \vee_S b)=h(c(a\vee b)=fgf(a\vee b)=f(a\vee b)=f(a) \vee f(b)=h(a)\vee h(b).
   $$
   Moreover the map $h:S\to G$ is injective since for $x\in S$ one has $x=c(x)=g(f(x))$. \newline
   We show 4). One has $q_S=c$ and $c=gf$ so 4) follows from $f=fgf$. We have thus proven that $S\subset F$ fulfills the first 4 conditions. The last one can be derived from the adjunction $c\dashv \iota$.\endproof
      The subset $S\subset F$ is called the {\em support} of $f: F \to G$ and denoted ${\rm Support}(f)$. To stress the relation between the operator $q_S: F \to S$ of Proposition \ref{support} and an orthogonal projection we state the following
   
   \begin{lem}\label{ortho proj0} With the notations of  Proposition \ref{support}, let $\hat S\subset F^*=\Hom_\B(F,\B)$ be the submodule of $F^*$ given by the $\phi_s$, $s\in S$. Then one has $\hat S={\rm Range}(f^*)$ and
   \begin{equation}\label{orthogonality0}
  	q_S(\xi)=q_S(\eta)\iff \langle\xi,\phi_s\rangle=\langle\eta,\phi_s\rangle\qqq s\in S\iff 
  	 \langle\xi,u\rangle=\langle\eta,u\rangle\qqq u\in \hat S
  \end{equation}  
     \end{lem}
   \proof We show first that $\hat S={\rm Range}(f^*)$. Since $c=gf$ and $g=gfg$,  $S$ is the range of $g$ and moreover the maps $\phi_s$, $s\in S$, fulfill $\phi_{g(y)}=f^*(\phi_y)$ since
   $$
   \phi_{g(y)}(x)=0\iff f(x)\leq y \iff f^*(\phi_y)(x)=0.
   $$
   We now prove \eqref{orthogonality0}. By construction $\phi_s(\zeta)=0\iff \zeta\leq s$, thus the equality on the right hand side of \eqref{orthogonality0} means: $\xi\leq s\iff \eta\leq s$ for any $s\in S$. One has $q_S=c$ and thus by \eqref{closure op}: $\xi\leq s\iff \eta\leq s$ $\forall s\in S$ if and only if $q_S(\xi)=q_S(\eta)$.
     \endproof 

\begin{rem}\label{projonfix} Let $F$ be an object of $\bmod$, $\sigma\in \Aut_\bmod(F)$ an involution (\ie $\sigma^2=\id$) and $E=F^\sigma\subset F$ the fixed subset by $\sigma$.  Then formula \eqref{try proj} applied to the subset $E\subset F$ defines the projection $p: F\to E,~p(x):=x+\sigma(x)$. The equality $q_E=p$ shows that $p$ only depends upon the submodule $F^\sigma$.	
\end{rem}

\subsection{The radical of an object of $\bmod$}\label{section radical}

In this subsection we discuss for future use (\cf~Theorem \ref{kernel and injective}) the condition on   an object $E$ of $\bmod$ equivalent to state  that the dual $E^*$ is generated by its minimal non-zero elements. To make clear the analogy with the notion of the radical  we shall use in the following the terminology ``ideal" in place of ``hereditary submodule" as in Definition \ref{heredit}. 
 For an object $E$ of $\bmod$ the correspondence between ideals $J\subset E$ and elements $\phi\in E^*$ given by: $J=\phi^{-1}(\{0\})$ \& $\phi(x)=0 \Leftrightarrow x\in J$, for $J\subset E$ hereditary submodule (\cf~Proposition~\ref{inj1}) fulfills the rule: $J\subset J'\Leftrightarrow \phi\geq \phi'$. Thus minimal non-zero elements of $E^*$ correspond to maximal ideals of $E$, where $E$ itself is not counted as an ideal by convention. We can thus reformulate, when $E$ is finite, the condition: "$E^*$ is generated by its minimal non-zero  elements", in terms of ideals of $E$ as follows\vspace{.05in}
 
\centerline{\em Every ideal $J\subset E$ is the intersection of maximal ideals containing $J$.}\vspace{.05in} 

Indeed, the above statement means exactly that any element of $E^*$ is a supremum of minimal non-zero elements. This condition makes sense in general, without any finiteness hypothesis and in fact it also suggests to define, for any object $E$ of $\bmod$, the following congruence relation on $E$ 
\begin{defn}\label{defn radical} The radical $\rad(E)$ of an object $E$ of $\bmod$ is defined as the following congruence 
$$
x\sim_{\rad(E)} y\iff \left(\forall \ \text{maximal}\  \text{ideal} \ J\subset E,  \ x\in J\iff y\in J\right).
$$	
\end{defn}
An equivalent formulation of this definition can be given in terms of the quotients $E/J$, for maximal ideals $J\subset E$. Here, the quotient $E/J$ is defined as the $\B$-module of equivalence classes for the relation: $u\sim v\iff \exists i,j\in J\mid u+i=v+j$. One checks that this is an equivalence relation compatible with $+$ and that $u\sim v$ holds iff $f(u)=f(v)$, for any morphism $f$ with $J\subset f^{-1}(\{0\})$. The following lemmas describe   properties of $\rad(E)$. 
\begin{lem}\label{lem radical} The radical $\rad(E)$ is the same as the relation 
	stating that $x,y\in E$  have the same image in $E/J$, for all maximal ideals $J\subset E$. 
\end{lem}
\proof To prove the statement one shows that $E/J=\B$ for any maximal ideal $J\subset E$.  Let $J\subset E$ be a maximal ideal, then any proper ideal in $F=E/J$ is reduced to $\{0\}$ since its inverse image  by the quotient map is an ideal of $E$ containing $J$. It follows that $F=\B$ since all intervals $[0,\xi]\neq F$ are reduced to $\{0\}$. Moreover the class of an element $u\in E$ is $0\in E/J= \B$ iff $u+i=j\in J$ for some $i\in J$ and this implies $u\in J$. Thus the image of $u$ in $E/J$ is entirely determined as $0$ if $u\in J$ and $1$ if $u\notin J$. This provides the equivalence with Definition \ref{defn radical}.\endproof
 \begin{lem}\label{max ideals} Let $E$ be an object of $\bmod$. An ideal $J\subset E$ is maximal if and only if the quotient $E/J\simeq \B$.   	
   \end{lem}
   \proof The proof of the previous lemma shows that if $J$ is maximal then the quotient $E/J\simeq \B$. Conversely, assume that the quotient  $E/J$  is $\B$. Let $J'\supset J$ be an ideal containing $J$ with $J'\neq E$.  Let $u\sim_J v$ be the equivalence relation defined as: $u\sim_J v\iff \exists i,j\in J\mid u+i=v+j$. One has $u\sim_J v \Rightarrow u\sim_{J'} v$ and thus one gets an induced surjective map  $s:\B=E/J\to E/J'$. The class of any $x\notin J'$ in $E/J'$ is non zero since $u\sim_{J'} 0\Rightarrow u\in J'$. Thus $E/J'$ contains two distinct elements,  $s$ is also injective and $J=J'$. \endproof  
\begin{lem}\label{lem radical1} 
Let  $E$ be an object of $\bmod$. Then the congruence $\rad(E)$ is trivial (\ie  $x\sim_{\rad(E)} y\Rightarrow x=y$) if and only if every principal ideal $J\subset E$ is the intersection of maximal ideals containing $J$.
\end{lem}
\proof  Assume that every principal ideal $J\subset E$ is the intersection of maximal ideals containing $J$. Apply this for $x\in E$ to the ideal $J_x=[0,x]\subset E$. Since this ideal uniquely determines $x$ (as its largest element) and since for any maximal ideal $K\subset E$ one has $J_x\subset K\iff x\in K$, it follows that the congruence $\rad(E)$ is trivial \ie that $x\sim y\Rightarrow x=y$. Conversely, assume that the congruence $\rad(E)$ is trivial. Let $J_x=[0,x]\subset E$ be a principal ideal. Let $J\subset E$ be the intersection of all maximal ideals containing $J_x$. One has $J_x\subset J$. Let $y\in J$: we show that $x+y\sim_{\rad(E)} x$.  For any maximal ideal $K$ one has: $
x\in K\Rightarrow x+y\in K \Rightarrow x\in K
$.
Thus since the congruence $\rad(E)$ is trivial one gets $x+y=x$, \ie $y\leq x$ so that $y\in J_x$. \endproof
Let  $E$ be an object of $\bmod$.  We let $\kappa:E^*\to \id(E)$ be the bijection which associates to $\phi\in E^*$  the ideal $\phi^{-1}(\{0\})=J$.
\begin{lem}\label{linear forms prep} Let $E$ be an object of $\bmod$ and ${\rm Max}(E)$ the set of maximal ideals of $E$. \newline
$(i)$~The map $\kappa$ induces a bijection of the set of minimal non-zero elements $M\subset E^*$ with ${\rm Max}(E)\subset \id(E)$. \newline
  $(ii)$~For $\phi\in E^*$, the following conditions are equivalent
   \begin{enumerate}       
   \item $\phi$ belongs to the complete submodule of $E^*$ generated by the minimal elements.
    \item $J=\kappa(\phi)$ is an intersection of maximal ideals.  	
   \end{enumerate} 	
   \end{lem}
   \proof The equivalence follows from the fact that the bijection $\kappa:E^*\to \id(E)$   transforms the operation $\vee$ on $\phi$ into the intersection of the ideals $J$.\endproof 
 We now replace the condition that the dual $E^*$ is generated by its minimal non-zero elements by the weaker condition given by  the triviality of the congruence $\rad(E)$. In view of Lemma 
\ref{lem radical1} the triviality of $\rad(E)$ is analogous to the statement of the vanishing of the radical in ordinary algebra.
\begin{prop}\label{quotient by radical} Let $E$ be a finite object of $\bmod$. Then
\newline
$(i)$~The dual of the quotient $E/\rad(E)$ is the submodule $S\subset E^*$ generated by the minimal non-zero elements.\newline
$(ii)$~The quotient $E/\rad(E)$ has a trivial radical congruence.\newline
$(iii)$~For $\phi\in E^*$, $J=\kappa(\phi)$ is an intersection of maximal ideals if and only if $x\sim_{\rad(E)} y\Rightarrow \phi(x)=\phi(y)$.
\end{prop}
\proof  $(i)$~Let $S$ be the submodule of $E^*$ generated by the minimal non-zero elements $M\subset E^*$ and $j:S\to E^*$ the inclusion. The morphism $j^*:E=(E^*)^*\to S^*$ is surjective. Moreover $j^*(x)=j^*(y)\iff x\sim_{\rad(E)} y$   since $M$ generates $S$. Thus,  one has $E/\rad(E)=S^*$ and $(E/\rad(E))^*=(S^*)^*=S$. \newline
$(ii)$~Since $(E/\rad(E))^*=S$ is generated by its minimal non-zero elements,    $E/\rad(E)$ has a trivial radical congruence.\newline
$(iii)$~This follows from $(i)$ and Lemma \ref{linear forms prep}.  \endproof 

   \begin{prop}\label{good morphisms} Let $E,F$ be finite objects of $\bmod$ and $f\in \Hom_\B(E,F)$. Then one has $x\sim_{\rad(E)} y\Rightarrow f(x)\sim_{\rad(F)}f(y)$ if and only if for every maximal ideal $J$ of $F$, the ideal $f^{-1}(J)$ is an intersection of maximal ideals. 
   \end{prop}
   \proof Let us first assume that $x\sim_{\rad(E)} y\Rightarrow f(x)\sim_{\rad(F)}f(y)$. Let $J$ be a maximal ideal in $F$ and $\phi\in F^*$, $\phi^{-1}(\{0\})=J$. Then by Proposition \ref{quotient by radical} one has  $x\sim_{\rad(F)} y\Rightarrow \phi(x)=\phi(y)$, and thus $u\sim_{\rad(E)} v\Rightarrow \phi(f(u))=\phi(f(v))$. Thus the ideal $f^{-1}(J)$ is an intersection of maximal ideals. Conversely, if this fact holds for any maximal ideal $J\subset F$, then if $x\sim_{\rad(E)} y$ one has $x\in f^{-1}(J)\iff y\in f^{-1}(J)$ or equivalently $f(x)\in J\iff f(y)\in J$, \ie $f(x)\sim_{\rad(F)}f(y)$.\endproof 
   Next  we provide the simplest example of a morphism  which does not fulfill the condition of Proposition \ref{good morphisms}. Example \ref{rad not functorial} shows also that in general a morphism $f$ in $\bmod$  does not induce a morphism on the quotients by the radical congruence.
   \begin{example}\label{rad not functorial}  Let $N=\{0,m,n\}$ with $0<m<n$, and  $M=\{0,m,x,n\}$. The idempotent addition is defined as follows: $n\vee m=n$, $x\vee m=n$, $x\vee n=n$. Let $f\in \Hom_\B(N,M)$ be the natural inclusion. One sees that in $N$, $\{0,m\}$ is the only maximal ideal and thus
   $m\sim_{\rad(N)} 0$. On the other hand, in $M$ the maximal ideal $\{0,x\}$ contains $0$ but not $m$ so that $m\nsim_{\rad(M)} 0$.   	
   \end{example}

\section{The category $\bmod$ of $\B$-modules and the comonad $\perp$}\label{sect abstract bmod}

  The category $\bmod$  is a symmetric, closed monoidal category. The object $\{0\}$ is both initial and final. Thus for any pair of objects $M, N$ in $\bmod$ there is  a natural morphism
$
\gamma_{M,N}:M\amalg N\to M\times N
$
from the coproduct of  $M$ and $N$ to their product. Indeed, by  construction of $M\amalg N$, a morphism $f:M\amalg N\to P$ is a pair of morphisms $M\to P$, $N\to P$, then for $P=M\times N$ one takes the morphisms $(\id, 0):M\to P$ and $(0,\id):N\to P$. We recall the proof of the  following well known Lemma (\!\cite{MacLane})
\begin{lem}\label{alab} In the category of $\B$-modules the morphisms $\gamma_{M,N}$ are isomorphisms.	
\end{lem}
\proof Let $M$ and $N$ be two $\B$-modules.  By definition  $M \amalg N$ is the initial object for pairs of morphisms  $\alpha:M\to X$, $\beta:N\to X$, where $X$ is any $\B$-module. Let $P=M\times N$ endowed with the operation 
$$
(x,y)+(x',y')=(x+x',y+y')\qqq x,x'\in M, \forall y,y'\in N.
$$
$P$ is a $\B$-module with zero element  $(0,0)$. Moreover one has canonical morphisms $s:M\to P$, $s(x)=(x, 0)$ and $t:N\to P$, $t(y)=(0,y)$.
 Given a pair of morphisms  $\alpha:M\to X$ and $\beta:N\to X$, there exists a unique morphism   $\rho: M\times N\to X$, $\rho(x,y)= \alpha(x)+\beta(y)$ such that $\rho\circ s=\alpha$ and $\rho\circ t=\beta$. This proves that $(M\times N, s,t)$ is also an initial object for pairs of morphisms $\alpha,\beta$ as above. 
 Given a pair of morphisms $\phi:X\to M$ and $\psi:X\to N$, the map $a\mapsto (\phi(a),\psi(a))$ defines a morphism $\sigma:X\to P$ uniquely characterized by $\phi=p_1\circ \sigma$ and $\psi=p_2\circ \sigma$, where the $p_j$ are the projections. This shows that $(M\times N,p_1,p_2)$ is the product in the category of $\B$-modules.  
 The  morphism $\gamma_{M,N}:M\amalg N\to M\times N$ is defined by implementing the pair of morphisms 
 $\alpha=(\id, 0):M\to P$ and $\beta=(0,\id):N\to P$. The corresponding morphism $\rho$ is the identity map and 
 this shows that the coproduct $M\amalg N$ is isomorphic to the product by means of the morphism $\gamma_{M,N}$.\endproof
 
In view of Lemma \ref{alab} we shall denote by $M\oplus N$ the $\B$-module $M\amalg N\simeq M\times N$. Given two morphisms $\alpha,\beta:L\rightrightarrows M$ of $\B$-modules, their equalizer is given by the submodule $\equ(\alpha,\beta)=\{x\in L\mid \alpha(x)=\beta(x)\}$. We now describe the coequalizer.

 \subsection{Coequalizer of two morphisms of $\B$-modules}
We give the explicit construction of the coequalizer in the category of $\B$-modules.  
 By definition, a congruence on a $\bmod$ $M$ is a submodule of $M\times M$ characterized as the graph $\cC$ of an equivalence relation. Thus $\cC\subset M \times M$ fulfills the following conditions
\begin{enumerate}
\item $(f_i,g_i)\in \cC, i=1,2\implies (f_1+f_2,g_1+g_2)\in \cC$
\item $(f,f)\in \cC$, $\forall f\in M$
\item $(f,g)\in \cC\implies (g,f)\in\cC$
\item $(f,g)\in \cC, \ (g,h)\in \cC\implies (f,h)\in\cC$.
\end{enumerate}
These conditions ensure that  the quotient $M/\cC$ is a $\B$-module. 
We  use this notion to construct the coequalizer of two morphisms $\alpha,\beta:L\rightrightarrows M$ of $\B$-modules. We let $\cC$ be the intersection of all submodules of $M\times M$ which fulfill the above four conditions and also contain $\{(\alpha(x),\beta(x))\mid x\in L\}$. Let $E=M/\cC$  and $\rho:M\to E$ be the quotient map. By construction one has $\rho\circ \alpha=\rho\circ \beta$. 
\begin{lem} \label{doublelem0} 
$(i)$~The pair $(E,\rho)$ is the coequalizer of the morphisms $\alpha,\beta:L\rightrightarrows M$:~~
$
L \stackbin[\beta]{\alpha}{\rightrightarrows} M \stackrel{\rho}{\to} E
$.\newline
$(ii)$~The coequalizer of the morphisms $\alpha,\beta:L\rightrightarrows M$ is the quotient of $M$ by the equivalence relation
\begin{equation}\label{coequequiv}
	x\sim x'\iff h(x)=h(x') \qqq X,\, h:M\to X\mid h\circ \alpha=h\circ \beta
\end{equation}

\end{lem}
\proof $(i)$~We have seen that  $\rho\circ \alpha=\rho\circ \beta$. To test the universality, we let $\phi:M\to X$ be a morphism of $\B$-modules such that $\phi\circ \alpha=\phi\circ \beta$. Let then
$
\cC_\phi:=\{(x,y)\in M\times M\mid \phi(x)=\phi(y)\}
$. 
One gets $\cC\subset \cC_\phi$ since by construction $\cC_\phi$ fulfills the above four conditions and also contains $\{(\alpha(x),\beta(x))\mid x\in S\}$. Thus one sees that $\phi(x)$ only depends on the image $\rho(x)\in E$ and moreover the map factors through $\rho$ as required by the universal property of the coequalizer. \newline
$(ii)$~This follows from $(i)$ since any $h:M\to X\mid h\circ \alpha=h\circ \beta$ factors through $\rho$.
\endproof
\begin{rem}\label{coequ}{\rm In some contexts the following equivalence relation on $M$ is introduced in the presence of two  morphisms $\alpha,\beta:L\rightrightarrows M$ of $\B$-modules
$$
x,y\in M\quad x\sim y\iff \exists u, v\in L~~{\rm s.t.}~     \ x+\alpha(u)+\beta(v)=y+\alpha(v)+\beta(u).
$$
This is an additive congruence which coequalizes $\alpha$ and $\beta$, however next example shows that in general this equivalence relation is not the coequalizer of $\alpha,\beta$. Let $\alpha=\beta=\id:M\to M$. Then the coequalizer as in Lemma \ref{doublelem0} is the identity map $\id:M\to M$ whereas the above equivalence relation reads as: $x\sim y\iff \exists a\in M~~{\rm s.t.}~ x+a=y+a$. In characteristic $1$ (\ie for idempotent structures) this gives $x\sim 0$ $\forall x\in M$. Thus the above congruence does not provide  in   general the coequalizer of two maps.
}\end{rem}

\subsection{Kernel  and coimage of morphisms of $\B$-modules}
The category $\bmod$  is a particular example for the general category RSmd of semimodules over a unital semiring $R$ considered in \cite{gra1} \S 1.6.2. In $\bmod$, the naive notion of kernel given, for  a morphism $f\in \Hom_\B(M,N)$, by $f^{-1}(\{0\})$ (\ie the equalizer of $f$ and the $0$-morphism) is not adequate. In fact, using this notion one can find plenty of examples of morphisms with trivial kernels but failing to be monomorphisms. In fact one has in general
\begin{lem}\label{maxform} Let $M$ be an object of $\bmod$ and $f\in \Hom_\B(M,\B)$ be defined by $f(x)=1$, $\forall x\neq 0$. Then $f^{-1}(\{0\})=\{0\}$ and $f$ is not a monomorphism unless $M=\B$.	
\end{lem}
\proof It is enough to show that $f$ is additive and this follows from $x+y=0\Rightarrow x=y=0$. \endproof 

The problem of the inadequacy of the naive notion of kernel is of the same nature as the inadequacy  
of the notion of ideal for semirings which needs to be replaced by the notion of a congruence.

For a morphism $f\in \Hom_\B(M,N)$, our goal is to construct an exact sequence which replaces, in this context, the sequence 
\begin{equation}\label{seq}
0\to \Ker f\to M\stackrel{f}{\to} N\to \coker f\to 0
\end{equation}
holding in an abelian category. To this end, we introduce the extended notion of morphism obtained by considering pairs $(f,g)$, with $f,g\in \Hom_\B(M,N)$ which compose 
as follows:
\begin{equation}\label{paircomp}
	(f,g)\circ (f',g'):=(f\circ f'+g\circ g', f\circ g'+g\circ f').
\end{equation}
This law makes sense since addition of morphisms makes sense in $\bmod$. This set-up determines the new category $\bm2$ that will be considered in details in \S\ref{weak and strong}. Any ordinary morphism $f\in \Hom_\B(M,N)$ is therefore seen as the pair $(f,0)$ of $\bm2$. Moreover, for pairs of morphisms one introduces the following definition
\begin{defn}\label{paircomp1} Let $f,g\in \Hom_\B(M,N)$ then 
\begin{enumerate}
\item $\equ(f,g)$ is the equalizer of $(f,g)$: ~~$\equ(f,g)\to M\stackbin[g]{f}{\rightrightarrows} N$.
\item $\coequ(f,g)$ is the 	coequalizer of $(f,g)$:~~$M\stackbin[g]{f}{\rightrightarrows} N \to \coequ(f,g)$.
\end{enumerate}	
\end{defn}
Diagonal pairs such as $(f,f)$ play a special role since $\forall f\in \Hom_\B(M,N)$ one has: 
$
\equ(f,f)=(M,\id_M),$ and $\coequ(f,f)=(N,\id_N)
$.
Thus the result is independent of the choice of $f$ and the outcome is similar to what happens for the morphism $0\in\Hom_\B(M,N)$.
It would however be too naive to simplify by such pairs and the outcome would correspond to the over-simplification described in Remark \ref{coequ}.
 To properly define the replacement for $\Ker f$ of \eqref{seq} we introduce the pair 
$$
f^{(2)}:=(f\circ p_1,f\circ p_2), \  \  f\circ p_j\in \Hom_\B(M^2,N),~j=1,2
$$
where the $p_j:M^2\to M$ are the two canonical projections. By definition $\equ f^{(2)}$ is the equalizer of $(f\circ p_1,f\circ p_2)$: 
$
\equ f^{(2)} \to M^2 \stackbin[f\circ p_2]{f\circ p_1}{\rightrightarrows} N
$, hence it is, by construction, a subobject of $M^2$ (By Proposition \ref{monoepi} any monomorphism in $\bmod$ is an injective map). This determines two maps $\equ f^{(2)}\stackbin[\iota_2]{\iota_1}{\rightrightarrows} M$ by composition with the $p_j$'s. Thus one obtains: \begin{equation}\label{defnker}
	 \equ f^{(2)} \stackbin[\iota_2]{\iota_1}{\rightrightarrows} M\stackrel{f}{\to} N
\end{equation}
where by construction $f\circ \iota_1=f\circ \iota_2$. 
Next, recall that in an abelian category one defines the coimage  $\coim f$ of a morphism $f$ as the cokernel of its kernel. In our context we use the kernel pair.
\begin{defn}\label{paircomp3} Let $f\in \Hom_\B(M,N)$ then 
\begin{enumerate}
\item $\Kpr f=\equ f^{(2)}\subset M^2$  (as in \eqref{defnker}).
\item $\coim f=\coequ(\iota_1,\iota_2)=M/\!\!\sim$ (as in \eqref{coequequiv}).
\end{enumerate}	
\end{defn}
It follows from Lemma \ref{doublelem0} that $\coim f$ is the quotient of $M$ by the equivalence relation fulfilling the $4$ conditions of that lemma  and containing the pairs $(\iota_1(x,y),\iota_2(x,y))$. These pairs are characterized by the equality $f(x)=f(y)$ and thus one derives
\begin{lem} \label{coimlem} 
Let $f\in \Hom_\B(M,N)$, then $\coim f$ is the quotient of $M$ by the congruence 
\[
x\sim y\iff f(x)=f(y).
\]
\end{lem}

As a corollary one obtains that $\coim f$ is isomorphic to the naive notion of image of $f$, \ie as the submodule of $N$ given by 
\begin{equation}\label{naiveim}
	{\rm Range}\ (f):=\{f(x)\mid x\in M\}\subset N.
	\end{equation}
	
\subsection{Cokernel  and image of morphisms of $\B$-modules}
Next, we investigate the cokernel and the image of a morphism in $\bmod$. Likewise for the kernel, we associate to $f\in \Hom_\B(M,N)$ the pair 
$$
f_{(2)}:=(s_1\circ f,s_2\circ f), \  \  s_j\circ f\in \Hom_\B(M,N\oplus N), j=1,2
$$
where the $s_j:N\to N\oplus N$ are the two canonical inclusions. By definition, $\coequ f_{(2)}$, which plays the role of the cokernel, is the coequalizer of $(s_1\circ f,s_2\circ f)$ and is hence a quotient of $N\oplus N$. This provides two maps $\gamma_j:N\to\coequ f_{(2)}$ by composition with $s_j$. Thus one obtains: 
\begin{equation}\label{defncoker}
	M\stackrel{f}{\to} N \stackbin[\gamma_2]{\gamma_1}{\rightrightarrows}\coequ f_{(2)}
\end{equation}
where by construction $ \gamma_1\circ f=\gamma_2\circ f$ and the coequalizer is the universal solution of this equation. This equation only involves the naive notion of image of $f$ \ie the submodule of $N$ given by \eqref{naiveim}
since the equality $ \alpha_1\circ f=\alpha_2\circ f$ is equivalent to $ \alpha_1\circ \iota=\alpha_2\circ \iota$, where $\iota:{\rm Range}\ (f)\to N$ is the inclusion. \newline
As in universal algebra we define the cokernel pair $\Cpr f$ of the morphism $f\in \Hom_\B(M,N)$ as the $\B$-module $\coequ f_{(2)}$ together with the morphism $(\gamma_1,\gamma_2):N\oplus N \to \coequ f_{(2)}$, and the categorical notion of image is defined as follows
\begin{defn}\label{propimagedefn} Let $f\in \Hom_\B(M,N)$ then 
\begin{enumerate}
\item $\Cpr f = (\coequ f_{(2)})=(N\oplus N)/\!\!\sim$ (as in \eqref{defncoker}).
\item $\im (f)=\equ(\gamma_1,\gamma_2)\subset N$.
\end{enumerate}	
\end{defn}
 
Thus we derive  the sequence  (analogue of \eqref{seq})
\[
\Kpr f\, \stackbin[\iota_2]{\iota_1}{\rightrightarrows} M\stackrel{f}{\to} N \stackbin[\gamma_2]{\gamma_1}{\rightrightarrows}\,\Cpr f.
\] 
Next, we shall compare the naive image of a morphism as defined in \eqref{naiveim} with the categorical notion given in Definition \ref{propimagedefn}. Both notions depend only upon the submodule ${\rm Range}\ (f)\subset N$ and by construction one has ${\rm Range}\ (f)\subset\im (f)$. 

\begin{prop} \label{imlem1} 
Let $N$ be a $\B$-module, $E\subset N$  a submodule and let $\iota:E\to N$ the inclusion. Then \newline
$(i)$~$\im(\iota)=\tilde E$, where $\tilde E\subset N$ is defined using arbitrary $\B$-modules as follows
\begin{equation}\label{separ1}
x\in\tilde E\iff \phi(x)=\psi(x) \qqq X,\ \phi,\psi \in \Hom_\B(N,X)~{\rm s.t.}~ \phi\circ \iota=\psi \circ \iota.
\end{equation}
$(ii)$~One has $\tilde E=E$.
\end{prop}
\proof $(i)$~Note that formula \eqref{separ1} defines the kernel (\ie the equalizer) of the cokernel pair of the inclusion $\iota:E\to N$. Indeed, one lets $s_j:N\to N\oplus N$ be the canonical inclusions $(\id,0)$ and $(0,\id)$ and $\rho:N\oplus N\to N$ the coequalizer of $(s_1\circ \iota,s_2\circ \iota)$. Then one has in general 
\begin{equation}\label{separbis}
x\in\tilde E\iff \rho\circ s_1(x)=\rho\circ s_2(x). 
\end{equation}
$(ii)$~follows from Lemma \ref{imlem}. 
\endproof

As a bi-product, one derives for $\bmod$  the key property AB2 holding for abelian categories
\begin{prop}\label{corinj}
For a morphism $f\in \Hom_\B(M,N)$, the natural map from $\coim f$ to $\im(f)$ is an isomorphism.	
\end{prop}
\proof  By Lemma \ref{coimlem} the coimage of $f$ is the quotient of $M$ by the congruence  $f(x)=f(y)$. Thus the coimage of $f$ is isomorphic to the submodule $E={\rm Range}\ (f)\subset N$. The image of $f$ as in Definition \ref{propimagedefn} is $\tilde E$ which is equal to $E$ by Proposition  \ref{imlem1}.  \endproof

\subsection{The comonad $\perp$ and its Eilenberg-Moore and Kleisli categories}\label{grandis remark} 

In the above discussion of kernels and cokernels of morphisms of $\bmod$ we made repeated use of the operation which replaces an object $M$ of $\bmod$ by its square $M^2$. This operation is an endofunctor and leads one to consider the following comonad (or cotriple).
\begin{prop}\label{comonad} The following rules define a comonad:
\begin{enumerate}
	\item The endofunctor $\perp:\bmod\longrightarrow \bmod$,~ $\perp\! M=M^2$, $(\perp\!f):=(f,f)$.
	\item The counit $\epsilon: \perp \to 1_{\bmod}$,~~$\epsilon_M = p_1$,~ $p_1:M^2 \to M$, $p_1(x,y)=x$.
	\item The coproduct $\delta:\perp \to \perp\circ\perp$,~~ $\delta_M= (M^2 \to (M^2)^2)$, ~$(x,y)\mapsto (x,y,y,x)$. \end{enumerate}	
\end{prop}
\proof  We first check the co-associativity, \ie the commutativity of the following diagram:
 $$
\xymatrix@C=45pt@R=55pt{
\perp\! M \ \ar[d]_{\delta_M}\ar[rr]^{\delta_M} && \perp\!(\perp\! M)\ar[d]^{\delta_{\perp\!M}}
\\
\perp\!(\perp\! M) \ \ar[rr]^{\perp\!(\delta_M)}&&\perp \perp\perp\! M}
$$
One has $\delta_{\perp\!M}\delta_M((x,y))=\delta_{\perp\!M}(((x,y),(y,x)))=(((x,y),(y,x)),((y,x),(x,y)))$ and \newline
$\perp\!(\delta_M)\delta_M((x,y))=(\delta_M((x,y)),\delta_M((y,x)))=((x,y,y,x),(y,x,x,y))$ so one gets the required equality. We now check the defining property of the counit, \ie the commutativity
$$
\xymatrix@C=45pt@R=55pt{
\perp\! M \ \ar[d]_{\delta_M}\ar[rr]^{\delta_M}\ar[drr]^\id && \perp\!(\perp\! M)\ar[d]^{\epsilon_{\perp\!M}}
\\
\perp\!(\perp\! M) \ \ar[rr]^{\perp\!(\epsilon_M)}&&\perp\! M}
$$
One has $\epsilon_{\perp\!M}\delta_M((x,y))=\epsilon_{\perp\!M}(((x,y),(y,x)))=(x,y)$ and similarly $\perp\!(\epsilon_M)\delta_M((x,y))=(\epsilon_M((x,y)),\epsilon_M((y,x)))=(x,y)$ so one gets the required equality.
\endproof 
A comonad gives rise to two categories, its Kleisli category and its Eilenberg-Moore category (\!\cite{barr} \S 3.2) which we now determine for the comonad $\perp$.
\begin{prop}\label{Kleisli} 
$(i)$~The Kleisli category of the comonad $\perp$ is the category $\bm2$ whose objects are  $\B$-modules and morphisms are pairs of morphisms in $\bmod$ with composition assigned by the formula
\begin{equation}\label{paircompbm2}
	(f,g)\circ (f',g'):=(f\circ f'+g\circ g', f\circ g'+g\circ f').
\end{equation}
$(ii)$~The Eilenberg-Moore category of the comonad $\perp$ is the category $\b2$ of $\B$-modules endowed with an involution $\sigma$.
\end{prop} 
\proof $(i)$~By construction of the Kleisli category $\cK_\perp$ its objects are the objects of $\bmod$. The morphisms $M\to N$ in $\cK_\perp$  are given by $\Hom_\B(\perp\! M,N)$, \ie, since $M^2=M\oplus M$, by pairs $(f,g)$ of morphisms $f,g\in \Hom_\B(M,N)$. The composition $(f',g')\circ (f,g)$ is given by: 
$$
\perp\! M\stackrel{\delta_M}{\to}\perp\perp\! M\stackrel{\perp((f,g))}{\to}\perp\! N\stackrel{(f',g')}{\to}P.
$$
One has $\delta_M((x,y))=(x,y,y,x)$ and  $\perp\!(f,g)(\delta_M((x,y)))=(f(x)+g(y),f(y)+g(x))\in \perp\! N$.
By applying $(f',g')$ one obtains 
$$
(f',g')\circ \perp\!(f,g)(\delta_M((x,y)))=f'(f(x)+g(y))+g'(f(y)+g(x))
$$ 
which coincides with $(f'\circ f+g'\circ g)(x)+(f'\circ g+g'\circ f)(y)$.\newline
$(ii)$~By construction the Eilenberg-Moore category of the comonad $\perp$ is the category of coalgebras for this comonad.
 A coalgebra in this context is given by an object $M$ of $\bmod$ and a morphism $\alpha:M\to \perp\! M$ that makes the following diagrams commutative
\begin{equation}\label{counit diagram}
\xymatrix@C=15pt@R=25pt{
 & \ M\ar[ld]_{\alpha}\ar[rd]^{\id_M} 
\\
\perp\! M \ \ar[rr]^{\epsilon_M}&& M}
\end{equation} 
	
	\begin{equation}\label{coassociative}
\xymatrix@C=25pt@R=35pt{
 M \ \ar[d]_{\alpha}\ar[rr]^{\alpha} && \perp\! M \ar[d]^{\perp(\alpha)}
\\
\perp\! M \ \ar[rr]^{\delta_M}&& \perp\perp\! M}
\end{equation}
 The commutativity of the diagram \eqref{counit diagram} means that $\alpha(x)=(x,\sigma(x))$ for some morphism $\sigma\in \End_\B(M)$. In the diagram \eqref{coassociative} one gets $\delta_M(\alpha(x))=(x,\sigma(x),\sigma(x),x)$ while $\perp\!\alpha(\alpha(x))=(x,\sigma(x),\sigma(x),\sigma^2(x))$. Thus the commutativity of the diagram \eqref{coassociative} means that $\sigma^2=\id$. \endproof 
Proposition \ref{Kleisli} gives the conceptual meaning of the categories $\bm2$ and $\b2$ which are studied in \S\S\ref{weak and strong},\ref{secthomologbmods}. A number of their properties  are corollaries   of  general properties of Kleisli and Eilenberg-Moore categories and for instance Lemma \ref{yb} $(ii)$ which identifies $\bm2$ as a full subcategory of $\b2$ is a special case of Proposition 13.2.11 of \cite{lambdacalc}. 

There is a natural notion of a projective object in a category $\cC$ endowed with a comonad $\perp$: see \eg \cite{Weibel} Definition 8.6.5.
\begin{defn}\label{projective def} An object $P$ of $\cC$ is $\perp$-projective if the counit map $\epsilon_P:\perp\! P\to P$ has a section, \ie there is a map $s:P\to \perp\! P$ such that $\epsilon_P\circ s=\id_P$.	
\end{defn}
In our setup we have
\begin{lem}\label{projective sect} Any object $P$ of $\bmod$ is $\perp$-projective for the monad $\perp=I\circ \yb$.	
\end{lem}
\proof Let $\iota_P:P\to \perp\! P$ be given by $\iota_P(x)=(x,0)$. Then one has $\epsilon_P\circ \iota_P=\id$. In fact the section $\iota_P$ can be characterized as the smallest one in the sense that for any other section $s:P\to \perp\! P$  one has $s+\iota_P=s$. \endproof

\section{The Kleisli  category $\bm2$}\label{weak and strong}

In this section we study the Kleisli category $\bm2$ of Proposition \ref{Kleisli} whose introduction can be justified independently as follows. 
The lack of the additive inverse for morphisms of the category $\bmod$ leads one to consider formal pairs $(f,g)$ of morphisms in $\bmod$ as a substitute for $f-g$. More precisely one introduces (see Proposition \ref{Kleisli} $(i)$)

\begin{defn} We denote by $\bm2$ the category whose objects are  $\B$-modules and morphisms are pairs of morphisms in $\bmod$ with composition assigned by the formula
\begin{equation}\label{paircompbm2biis}
	(f,g)\circ (f',g'):=(f\circ f'+g\circ g', f\circ g'+g\circ f').
\end{equation}
\end{defn}
By construction $\bm2$ is enriched over the category $\bmod$. In this section (\cf \S \ref{sectstrongexact}) we study a provisional notion of \strg exactness of sequences in $\bm2$. This definition 
 will be refined later on, in \S\ref{sectnullmorphisms}. In \S\S\ref{sectmeaningst}, \ref{sectveryshort}, \ref{sectepistrong}, we investigate the link of \strg exactness respectively with monomorphisms,  isomorphisms  and  epimorphisms. Finally in \S \ref{sectquotbm2} we consider the issue of defining quotients in $\bm2$.

\subsection{\Strgly exact sequences of $\B$-modules}\label{sectstrongexact}

The category $\bmod$ embeds as a subcategory of $\bm2$ by applying the functor $\kappa:\bmod\longrightarrow \bm2$  which is defined as  the identity on objects while  on morphisms one sets 
\begin{equation}\label{functoriota}
	\kappa(f):=(f,0)\in \Hom_\bm2(M,N)\qqq f\in \Hom_\B(M,N).
\end{equation}
 Definition \ref{paircomp3} gives the kernel in the case of pairs of the form $(f,0)$ thus the next step is to extend the construction of the kernel to an arbitrary pair $(f,g)$ of morphisms $M\to N$ in $\bmod$. The main idea is to re-interpret a congruence relation involving differences of the form $f(x)-g(x)=f(y)-g(y)$ as follows
$$
f(x)-g(x)=f(y)-g(y)\iff f(x)+g(y)=f(y)+g(x).
$$
This means that in the product $M\times M$ one looks for  all pairs $(x,y)$ such that $f(x)+g(y)=f(y)+g(x)$.  
\begin{prop}\label{ker1} Let $f,g\in \Hom_\B(M,N)$. \newline
$(i)$~$Z(f,g):=\{(x,y)\in M\times M\mid f(x)+g(y)=f(y)+g(x)\}$ is a submodule of $M\times M$ that we call the algebraic kernel of the pair $(f,g)$. \newline 
$(ii)$~Let $\iota_j:Z(f,g)\to M$ be the restrictions of the canonical projections $p_j:M^2\to M$, $j=1,2$. The sequence
\begin{equation}\label{defnkerbett}
	Z(f,g) \stackbin[\iota_2]{\iota_1}{\rightrightarrows} M  \stackbin[g]{f}{\rightrightarrows} N
\end{equation}	
satisfies the relation:~ $f\circ \iota_1+g\circ \iota_2=f\circ \iota_2+g\circ \iota_1$.\newline
$(iii)$~Let $\alpha,\beta\in \Hom_\B(L,M)$ be such that $f\circ \alpha+g\circ \beta=f\circ \beta+g\circ \alpha$. Then the range of the map $(\alpha,\beta):L\to M^2$ is included in $Z(f,g)\subset M^2$. \newline
$(iv)$~$Z(f,0)=\Kpr f$ (\cf~Definition \ref{paircomp3}). \end{prop}
\proof $(i)$~The pairs $(x,y)\in M\times M$ such that $f(x)+g(y)=f(y)+g(x)$ form a submodule  since 
$$
(x,y),(x',y')\in Z~\Rightarrow~ f(x+x')+g(y+y')=f(y+y')+g(x+x').
$$
$(ii)$~By definition $Z(f,g)$ is  a subobject of $M^2$, thus one defines two maps $Z(f,g)\stackbin[\iota_2]{\iota_1}{\rightrightarrows} M$ by composing with the $p_j$'s. From this one derives the sequence \eqref{defnkerbett}. By definition of $Z(f,g)$ one also has $f\circ \iota_1+g\circ \iota_2=f\circ \iota_2+g\circ \iota_1$.\newline
$(iii)$~For any $x\in L$ one has $(\alpha(x),\beta(x))\in Z(f,g)$.
\newline
$(iv)$~Both sides of the equality are defined as the submodule of $M^2$ given by the equation $f(x)=f(y)$. \endproof 
The statement $(iii)$ of the above proposition means that when the composition $(f,g)\circ (\alpha,\beta)$  of two successive pairs is equivalent to $0$, \ie given by a diagonal pair, one derives a factorization through the kernel $Z(f,g)$. Thus, by requiring that this factorization is onto $Z(f,g)$, one gets a first hint for the notion of \strg exactness. One has
\begin{prop}\label{coker1} Let $f,g\in \Hom_\B(L,M)$. \newline
$(i)$~$B(f,g)=\{(f(x)+g(y),f(y)+g(x))\mid x,y\in L\}$ is a submodule of $M\times M$. \newline 
$(ii)$~Let $\alpha,\beta\in \Hom_\B(M,N)$, then one has
\begin{equation}\label{exactcomp}
	B(f,g)\subset Z(\alpha,\beta)\iff \alpha\circ f+ \beta\circ g= \beta\circ f+ \alpha\circ g
\end{equation}	
$(iii)$~Let $\phi=(f,g)\in \Hom_\bm2(L,M)$, $\psi=(h,k)\in \Hom_\bm2(M,N)$,
then one has 
\begin{equation}\label{iso8}
	Z(f,g)=Z(\phi)\subset Z(\psi\circ \phi), \ \ B(h,k)=B(\psi)\supset B(\psi\circ \phi).
\end{equation}	

\end{prop}
\proof $(i)$~is straightforward. \newline
$(ii)$~$B(f,g)\subset Z(\alpha,\beta)$ if and only if for any $x,y\in L$ one has
$$
\alpha(f(x)+g(y))+\beta(f(y)+g(x))=\beta(f(x)+g(y))+\alpha(f(y)+g(x))
$$
or equivalently 
$$
(\alpha\circ f+ \beta\circ g)(x)+(\beta\circ f+ \alpha\circ g)(y)=
(\alpha\circ f+ \beta\circ g)(y)+(\beta\circ f+ \alpha\circ g)(x)
$$
This equality clearly holds if $\alpha\circ f+ \beta\circ g= \beta\circ f+ \alpha\circ g$. Conversely,  taking $y=0$, one gets $\alpha\circ f+ \beta\circ g= \beta\circ f+ \alpha\circ g$.\newline
$(iii)$~By Proposition \ref{ker1}, $Z(\psi\circ \phi)$ is given as the set of pairs $(x,y)\in L\times L$ such that 
$$
(h\circ f +k\circ g)(x)+ (h\circ g +k\circ f)(y)=(h\circ f +k\circ g)(y)+(h\circ g +k\circ f)(x).
$$
 Then the first inclusion in \eqref{iso8} follows by rewriting the above equality  as 
 $$
 h(f(x)+g(y))+k(g(x)+f(y))=h(f(y)+g(x))+k(g(y)+f(x)).
 $$
 By $(i)$, $B(\psi\circ \phi)$ is given as the set of pairs 
 $$
 \left( (h\circ f +k\circ g)(x)+ (h\circ g +k\circ f)(y),(h\circ f +k\circ g)(y)+(h\circ g +k\circ f)(x)\right).
 $$ 
 The second inclusion in \eqref{iso8} then follows by rewriting such pairs as 
 $$
 \left( h(X)+k(Y),h(Y)+k(X)\right), \  \  X=f(x)+g(y), \ Y=f(y)+g(x).
 $$
 \endproof 
 It is clear that in general one  has 
$\Delta\subset Z(\alpha,\beta)$, where $\Delta\subset M\times M$ is the diagonal, \ie $\Delta=\{(x,x)\mid x\in M\}$. Keeping in mind this fact, we introduce the notion of \strgly exact sequence as follows

	\begin{defn}\label{exsequ} The sequence  $L \stackbin[\alpha_2]{\alpha_1}{\rightrightarrows}M\stackbin[\beta_2]{\beta_1}{\rightrightarrows}N$ in $\bm2$ is \strgly exact at $M$ if $ B(\alpha_1,\alpha_2)+\Delta=Z(\beta_1,\beta_2)$.
	\end{defn}
	
	In the following we shall test this new notion in several cases. As a first case, we assume that $N=0$ so that both arrows $\beta_j=0$. Then, $Z(\beta_1,\beta_2)=M\times M$ and one needs to find out the relation between the coequalizer of the $\alpha_j$ and the \strg exactness of the sequence at $M$. This is provided by the following
	\begin{prop}\label{surjtest} Consider the sequence $L \stackbin[\alpha_2]{\alpha_1}{\rightrightarrows}M\stackbin[0]{0}{\rightrightarrows}0$ in $\bm2$.  \newline
$(i)$~The following three conditions are equivalent\newline a) The sequence is  \strgly exact at $M$,\newline b) $
\{\alpha_1(x)+\alpha_2(y)\mid x,y\in L, \ \alpha_2(x)+\alpha_1(y)=0\}=M
$ \newline c) $B(\alpha_1,\alpha_2)=M\times M$.
\newline 
$(ii)$~If $\alpha_2=0$  then the sequence is  \strgly exact at $M$ if and only if $\alpha_1$ is surjective. \newline
$(iii)$~If the $\alpha_j$'s have a non-trivial coequalizer then the sequence is not  \strgly exact at $M$. 
\end{prop}
\proof $(i)$~One has $Z(0,0)=M\times M$. By symmetry, one has $ B(\alpha_1,\alpha_2)+\Delta=M\times M$ if and only if $M\times \{0\}\subset B(\alpha_1,\alpha_2)+\Delta$. But an equality $a+b=0$ for $a,b\in M$ implies $a=0$ and $b=0$ since $a=a+a+b=a+b=0$. Thus for $x,y\in L$ and $z\in M$ such that 
$
(\alpha_1(x)+\alpha_2(y)+z,\alpha_2(x)+\alpha_1(y)+z)=(t,0)
$, 
one has $z=0$ and $\alpha_2(x)+\alpha_1(y)=0$.\newline
$(ii)$~If $\alpha_2=0$, $(i)$ states that the \strg exactness at $M$ holds if and only if
$
\{\alpha_1(x)\mid x,y\in L, \ \alpha_1(y)=0\}=M
$. This requirement is nothing but the surjectivity of $\alpha_1$.\newline
$(iii)$~Consider a non-trivial morphism $\phi:M\to N$ such that $\phi\circ \alpha_1=\phi\circ \alpha_2$. Let $t\in M$ with $\phi(t)\neq 0$. Then the pair $(t,0)$ cannot belong to $B(\alpha_1,\alpha_2)+\Delta$ since 
$
\phi(\alpha_1(x)+\alpha_2(y)+z)=\phi(\alpha_2(x)+\alpha_1(y)+z)~ \forall x,y\in L,\ z\in M
$.
\endproof
Next we test the notion of \strg exactness of $
0 \stackbin[0]{0}{\rightrightarrows}M\stackbin[\beta_2]{\beta_1}{\rightrightarrows}N
$.

\begin{prop}\label{injtest} Consider the sequence $0 \stackbin[0]{0}{\rightrightarrows}M\stackbin[\beta_2]{\beta_1}{\rightrightarrows}N$ in $\bm2$. \newline
$(i)$~The sequence is  \strgly exact at $M$ if and only if $$\beta_1(x)+\beta_2(y)=\beta_2(x)+\beta_1(y)~\Leftrightarrow~ x=y$$  
$(ii)$~If $\beta_2=0$ then the sequence is  \strgly exact at $M$ if and only if $\beta_1$ is injective. 
\end{prop}
\proof $(i)$~One has $B(0,0)=0$ and thus the \strg exactness holds if and only if $Z(\beta_1,\beta_2)=\Delta$.
\newline
$(ii)$~One has $\beta_1(x)=\beta_1(y)~\Leftrightarrow~ x=y$ if and only if $\beta_1$ is injective. \endproof 
\begin{cor}\label{isotest}
	The sequence $0 \stackbin[0]{0}{\rightrightarrows}M\stackbin[0]{\beta}{\rightrightarrows}N\stackbin[0]{0}{\rightrightarrows}0$ in $\bm2$ is  \strgly exact if and only if $\beta$ is an isomorphism in $\bmod$.
\end{cor}
\proof This follows from Propositions \ref{surjtest} and \ref{injtest}.\endproof

\subsection{Monomorphisms and \strg exact sequences}\label{sectmeaningst}

In this subsection we investigate the meaning of a  \strgly exact sequence as in Proposition \ref{injtest}.
\begin{prop}\label{injN} Let  $0 \stackbin[0]{0}{\rightrightarrows}M\stackbin[g]{f}{\rightrightarrows}N$ be a sequence in $\bm2$  \strgly exact at $M$, then:\newline
 $(i)$~The pair $(f,g)$ embeds $M$ as a submodule of $N^2$.\newline
 $(ii)$~The map $f+g$ embeds $M$ as a submodule of $N$.	
\end{prop}
\proof By Proposition \ref{injtest} one has 
\begin{equation}\label{exactleft}
	f(x)+g(y)=g(x)+f(y)\iff x=y
\end{equation}
$(i)$~Assume that for elements $x,y\in M$ one has $(f(x),g(x))=(f(y),g(y))$. Then $f(x)+g(y)=g(x)+f(y)$ thus by hypothesis $x=y$.\newline
$(ii)$~Assume that for elements $x,y\in M$ one has $f(x)+g(x)=f(y)+g(y)$. Then one has $f(x)+g(x+y)=f(x+y)+g(y)$ since adding $g(y)$ or $f(x)$ to $f(x)+g(x)=f(y)+g(y)$ does not change the result. Then by \eqref{exactleft} one gets $x=x+y$ and similarly $y=x+y$ so that $x=y$.\endproof
Thus, Proposition \ref{injN} $(i)$ shows that to understand the \strg exactness of  $0 \stackbin[0]{0}{\rightrightarrows}M\stackbin[g]{f}{\rightrightarrows}N$ at $M$ we can reduce to the case where $M\subset N\times N$ is a submodule of $N\times N$, while $f=p_1$ and $g=p_2$ are the two projections restricted to this submodule. \vspace{.05in}
The condition of \strg exactness on the submodule $M\subset N\times N$ reads as the implication for $x,y\in M\subset N\times N$
\begin{equation}\label{diagcond}
 x+\sigma(y)=\sigma(x)+y\Rightarrow y=x	
\end{equation}
where $\sigma:M\times M\to M\times M$ is the involution $\sigma(a,b):=(b,a)$.
One might guess at first  that \eqref{diagcond} implies the injectivity of the map  $\iota:N\times N\to M\times M$ given by
\begin{equation}\label{137}
	\iota(x,y):=x+\sigma(y)
\end{equation}
 But this fails as shown by the following
\begin{example}\label{example not mono}
	Let $X=\{1,2\}$ and  $N=2^X$. Consider the submodule $M\subset N\times N$ defined as
\begin{equation}\label{138}
M:=\{(\emptyset,\emptyset), (\emptyset,\{1\}),(\{1\},\{1,2\})\}.
\end{equation}
Then the sequence $0 \stackbin[0]{0}{\rightrightarrows}M\stackbin[p_2]{p_1}{\rightrightarrows}N$ is \strgly exact at $M$, but the map  $\iota:M\times M\to N\times N$ given by \eqref{137} is not injective.
\end{example}
 \proof 
The structure of $\B$-module of $M$ is the same as that of the totally ordered set $0<1<2$. For $x,y\in M$ the sum $x+\sigma(y)$ cannot be symmetric if $x$ or $y$ is $(\emptyset,\emptyset)$ and $x\neq y$. Moreover the sum of the two non-zero elements is $(\emptyset,\{1\})\vee (\{1\},\{1,2\})=(\{1\},\{1,2\})$ which is not symmetric. Thus $M$ fulfills the condition of \strg exactness at $M$ of the sequence $0 \stackbin[0]{0}{\rightrightarrows}M\stackbin[p_2]{p_1}{\rightrightarrows}N$. We now consider the map $\iota:M\times M\to N\times N$ as in \eqref{137}. We claim that its range has $7$ elements, precisely:
$$
\iota(M\times M)=\{(\emptyset,\emptyset), (\emptyset,\{1\}),(\{1\},\{1,2\}),(\{1\},\emptyset),(\{1,2\},\{1\}),(\{1\},\{1\}),(\{1,2\},\{1,2\})\}
$$
and the element $(\{1\},\{1,2\})$ is obtained twice since 
$$
(\{1\},\{1,2\})\vee \sigma((\emptyset,\{1\}))=(\{1\},\{1,2\})=(\{1\},\{1,2\})\vee \sigma((\emptyset,\emptyset)).
$$
Note that by applying $\sigma$ on both sides one gets an equality of the form $\iota(x,y)=\iota(x',y)$ with $x\neq x'$.\endproof 

\begin{figure}[H]
\begin{center}
\includegraphics[scale=0.4]{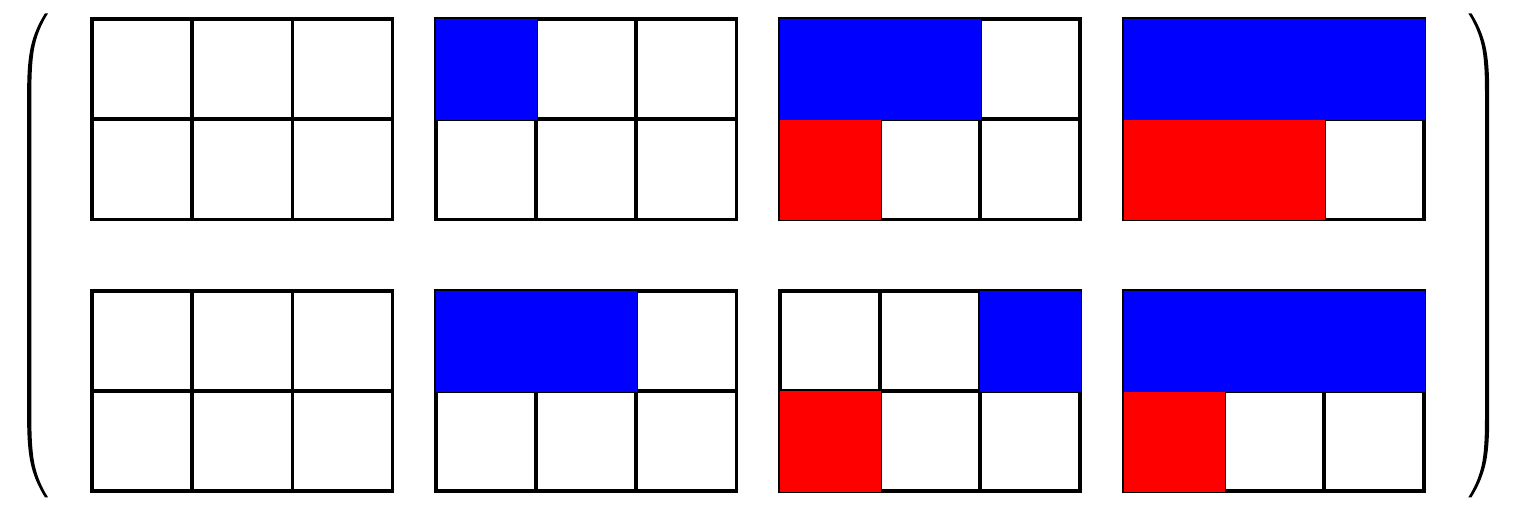}
\end{center}
\caption{Reduction to two cases $\vert M\vert=4$. \label{m3n4red} }
\end{figure}

Example \ref{example not mono} was the simplest case to consider. To understand the various choices of $M\subset N\times N$ which fulfill the \strg exactness of Lemma \ref{injN}, we consider the next case where  $N=2^X$, $\vert X\vert =3$. In this case we shall display all the maximal choices of $M$ up to permutations and symmetry.  The condition of \strg exactness \eqref{diagcond}
 is preserved by any automorphism $\alpha$ of $N\times N$ which commutes with the 
symmetry $\sigma$ and we simplify further using such automorphisms. Since $N=2^X$ one has $N\times N=2^{X \cup Y}$ where we represent $X$ and $Y$ as the first (in blue) and second line (in red) of the rectangle, and the elements of $N\times N$ as subsets of the rectangle. When the cardinality of $M$ is  $\vert M\vert =4$  the number of maximal cases up to symmetries is two. They are displayed in Figure \ref{m3n4red} where each line gives a choice of $M\subset N\times N$. 
Notice that these cases are not isomorphic since in the first line the module $M$ is totally ordered. Similarly for $\vert M\vert =5$ the number of maximal cases up to symmetries is two. 
\begin{figure}[H]
\begin{center}
\includegraphics[scale=0.4]{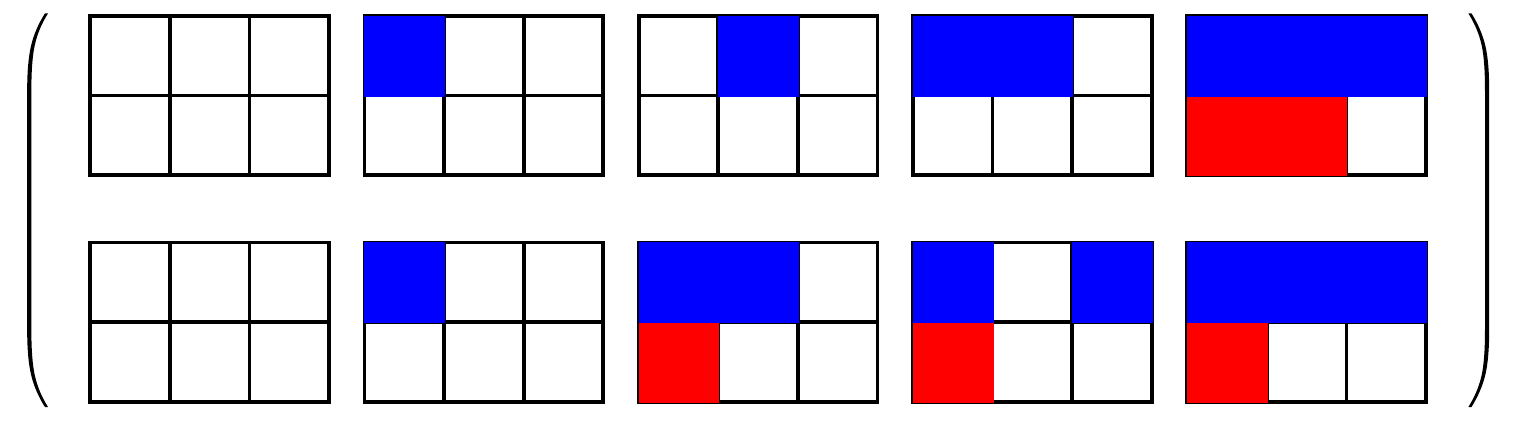}
\end{center}
\caption{Reduction to two cases $\vert M\vert=5$. \label{m3n5red} }
\end{figure}
For $\vert M\vert =6$ the number of maximal cases up to symmetries is one.
\begin{figure}[H]
\begin{center}
\includegraphics[scale=0.6]{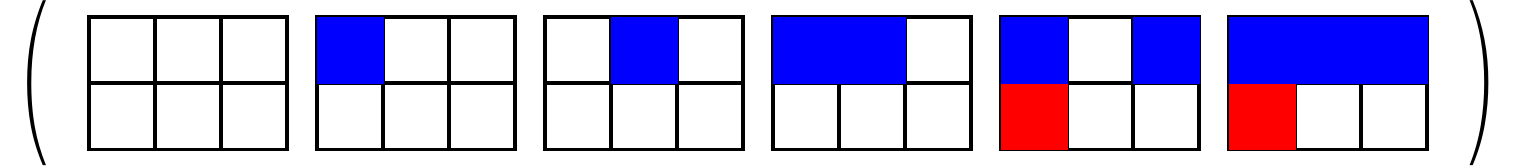}
\end{center}
\caption{Reduction to one case $\vert M\vert=6$. \label{m3n6red} }
\end{figure}
as well as the case $\vert M\vert =8$
\begin{figure}[H]
\begin{center}
\includegraphics[scale=0.7]{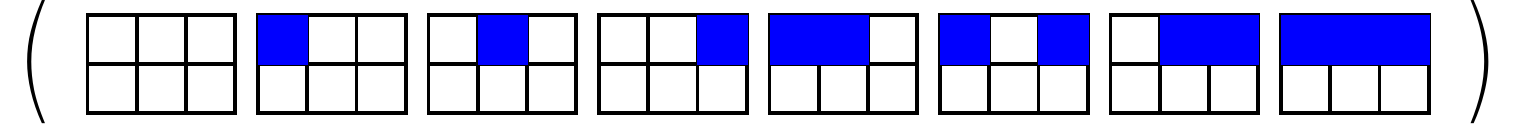}
\end{center}
\caption{Reduction to one case $\vert M\vert=8$. \label{m3n8red} }
\end{figure}
 Example \ref{example not mono} shows that \strg exactness of a sequence of the form $0 \stackbin[0]{0}{\rightrightarrows}M\stackbin[g]{f}{\rightrightarrows}N$ does not imply that the morphism $\phi=(f,g)$ is  a monomorphism. More precisely one has

      \begin{prop}\label{monomorphism}$(i)$~The  morphism $\phi:M\stackbin[g]{f}{\rightrightarrows}N$, is a monomorphism in the category $\bm2$
if and only if the map $M^2\to N^2,~(x,y)\mapsto (f(x)+g(y), g(x)+f(y))$ is injective.\newline 
$(ii)$~If $\phi=(f,g)$ is  a monomorphism, the sequence $0 \stackbin[0]{0}{\rightrightarrows}M\stackbin[g]{f}{\rightrightarrows}N$ is  \strgly exact at $M$.
\end{prop}
   \proof $(i)$~If the map $(x,y)\mapsto (f(x)+g(y),g(x)+f(y))$ fails to be injective, consider $(x,y)\in M^2$, $(x',y')\in M^2$ such that $(x,y)\neq (x',y')$ and
   $
   f(x)+g(y)=f(x')+g(y'),~g(x)+f(y)=g(x')+f(y')
   $.
  Let $\psi\in \Hom_\bm2(\B,M)$ be given by  $\B \stackbin[y]{x}{\rightrightarrows}M$, where $x$ stands for the unique morphism in $\bmod$ such that $1\mapsto x$. One uses similar notations  for $y$ and $\psi'$. The above equality then means that $\phi\circ \psi=\phi \circ \psi'$, where $\phi=(f,g)\in \Hom_\bm2(M,N)$. Thus $\phi$ fails to be a monomorphism in $\bm2$. Conversely, if the map $(x,y)\mapsto (f(x)+g(y), g(x)+f(y))$ is injective then the composition law in $\bm2$ shows that $\phi=(f,g)$ is a monomorphism. Indeed, the equality 
  $
  (f\circ \alpha+g\circ \beta, f\circ \beta+ g\circ \alpha)=
  (f\circ \alpha'+g\circ \beta', f\circ \beta'+ g\circ \alpha')
  $
  implies, by restriction to diagonal pairs $(a,a)$ that $\alpha(a)=\alpha'(a)$ and $\beta(a)=\beta'(a)$, $\forall a$.\newline 
$(ii)$~When $\phi=(f,g)$ is  a monomorphism, the above sequence is  \strgly exact
since the equality $f(x)+g(y)= g(x)+f(y)$ implies that $(x,y)$ and $(y,x)$ have the same image by the map $(x,y)\mapsto (f(x)+g(y),g(x)+f(y)):M^2\to N^2$.
   \endproof

  \subsection{Isomorphisms and \strg exact sequences}\label{sectveryshort}
 
 To understand how restrictive the notion of  \strgly exact sequence is we investigate the exact sequences in $\bm2$ of the form $0 \stackbin[0]{0}{\rightrightarrows}M\stackbin[g]{f}{\rightrightarrows}N\stackbin[0]{0}{\rightrightarrows}0$. 
 \begin{prop}\label{iso4}  Let the sequence $0 \stackbin[0]{0}{\rightrightarrows}M\stackbin[g]{f}{\rightrightarrows}N\stackbin[0]{0}{\rightrightarrows}0$ be  \strgly exact in $\bm2$. \newline
 Then there exists a unique decomposition $N=N_1\times N_2$ and a unique isomorphism of $\B$-modules $\alpha:M\to N$ such that 
 \begin{equation}
 	\label{iso5}
 f=(p_{N_1}\circ \alpha,0), \ \ g=(0,p_{N_2}\circ \alpha), \ \ p_{N_j}:N_1\times N_2\to N_j
 \end{equation}	
 \end{prop}
\proof By the \strg exactness of the sequence at $M$ and Lemma \ref{injN}, we can reduce to the case where $M\subset N\times N$, $f=p_1$ and $g=p_2$. By applying  Proposition \ref{surjtest} the \strg exactness of the sequence at $N$ means that 
\begin{equation}
 	\label{iso31}
\{p_1(x)+p_2(y)\mid x,y\in M, \ p_2(x)+p_1(y)=0\}=N
\end{equation}
 One then obtains the following two hereditary submodules $J_k\subset M$
$$
J_1:=\{x\in M\mid p_2(x)=0\}= M\cap (N\times 0), \ \ 
J_2:=\{x\in M\mid p_1(x)=0\}= M\cap (0\times N).
$$
The addition gives an injection $(N\times 0)\times (0\times N)\to N\times N$ and hence by restriction an injection $j:J_1\times J_2\to M$. Let $J=j(J_1\times J_2)=J_1+J_2\subset M$. The map $p:J\to N$, $p(x):=p_1(x) +p_2(x)$ is surjective  by \eqref{iso31}. By Lemma \ref{injN} $(ii)$ the \strg exactness at $M$ shows that the map $\tilde p:M\to N$, $\tilde p(x):=p_1(x) +p_2(x)$ is injective. Since its restriction $p:J\to N$ to $J\subset M$ is surjective this shows that $M=J$ and that $p$ is bijective and is an isomorphism.  To the decomposition of $J$ as $J_1\times J_2$ corresponds a decomposition $N=N_1\times N_2$ where $N_k=p(J_k)$.  This shows the existence of the decomposition \eqref{iso5}. We finally prove its uniqueness. Assuming \eqref{iso5}, one has $N_1=f(M)$, $N_2=g(M)$. This determines the decomposition $N=N_1\times N_2$ uniquely. Moreover both $p_{N_1}\circ \alpha$ and $p_{N_2}\circ \alpha$ are uniquely determined by $f$ and $g$ and thus $\alpha$ is unique.
\endproof
Proposition \ref{iso4} suggests that a morphism $M\stackbin[g]{f}{\rightrightarrows}N$ in the category $\bm2$ is an isomorphism iff the sequence $0 \stackbin[0]{0}{\rightrightarrows}M\stackbin[g]{f}{\rightrightarrows}N\stackbin[0]{0}{\rightrightarrows}0$ is  \strgly exact. 
 Let $M$ be a $\B$-module, and consider two decompositions $M=M_1\times M_2$, $M=M'_1\times M'_2$  of $M$ as a product. For $x\in M$ we denote by $x=x_1+x_2$ (resp. $x=x_{1'}+x_{2'}$) its unique decomposition with $x_j\in M_j$ (resp. $x_{j'}\in M'_j$). Since $M_j$ is a hereditary submodule one has, with $x_j=(x_j)_{1'}+(x_j)_{2'}$ that 
 $
 (x_j)_{k'}\in M_j\cap M'_k,$ and that any element $x\in M$ is uniquely decomposed as a sum $x=\sum x_{jk'}$ where $x_{jk'}\in  M_j\cap M'_k$. It follows that the projections $p_j$ and $p'_k$ associated to the two decompositions commute pairwise and that the composition in the category $\bm2$ of the morphisms $(p_1,p_2)$ and $(p'_1,p'_2)$ is given by the pair corresponding to the decomposition 
 $$
 M=\left((M_1\cap M'_1)\oplus (M_2\cap M'_2)\right)\times 
 \left((M_1\cap M'_2)\oplus (M_2\cap M'_1)\right).
 $$
 It follows that such pairs $(p_1,p_2)$ form a subgroup $\Aut^{(1)}_\bm2(M)\subset \Aut_\bm2(M)$ of the group of automorphisms $\Aut_\bm2(M)$. By construction this subgroup is abelian and every element is of order two.
 
  \begin{prop}\label{iso9}  $(i)$~The sequence $0 \stackbin[0]{0}{\rightrightarrows}M\stackbin[g]{f}{\rightrightarrows}N\stackbin[0]{0}{\rightrightarrows}0$ is  \strgly exact if and only if the pair $\phi=(f,g)$ is an isomorphism in $\bm2$.
 
 $(ii)$~The group of automorphisms $\Aut_\bm2(M)$ is  the semi-direct product of the $2$-group $\Aut^{(1)}_\bm2(M)$ of decompositions $M=M_1\times M_2$ by the group of automorphisms of the $\B$-module $M$.
 \end{prop}
 \proof $(i)$~Assume first that \eqref{iso5} holds, with $\alpha:M\to N$ an isomorphism in the category of $\B$-modules. The composition law $(u,v)\circ(\alpha,0)=(u\circ \alpha,v\circ \alpha)$ 
shows that to prove that the morphism $\phi=(f,g)$ is an isomorphism in $\bm2$ one can assume that $N=M$ and $\alpha=\id$. Thus we can assume that 
$
f=(p_{M_1},0),~ g=(0,p_{M_2}),~ p_{M_j}:M_1\times M_2\to M_j
$. 
One then has, since both $f$ and $g$ are idempotent and with vanishing product:~ 
$
(f,g)\circ (f,g)=(f\circ f+g\circ g,f\circ g +g\circ f)=(\id,0).
$
This shows that \strg exactness implies isomorphism. \newline
Conversely, if $\phi=(f,g)$ is an isomorphism in $\bm2$ it admits a left inverse $\psi$ such that $\psi\circ \phi=(\id,0)$ and Proposition \ref{coker1} $(iii)$ shows that $Z(f,g)\subset Z(\psi\circ \phi)=Z((\id,0))=\Delta$. Similarly, $\phi$ admits a right inverse $\psi$ such that $\phi\circ \psi=(\id,0)$ and Proposition \ref{coker1} $(iii)$ shows that $B(f,g)=B(\phi)\supset B(\phi\circ \psi)=B((\id,0))=N\times N$.\newline
$(ii)$~The map $\alpha\mapsto (\alpha,0)=\rho(\alpha)$ defines an injective group homomorphism $\rho:\Aut_\bmod(M)\to \Aut_\bm2(M)$. Moreover by $(i)$ together with  Proposition \ref{iso4}, we derive that any element of $\Aut_\bm2(M)$ is uniquely a product $\beta\circ \rho(\alpha)$, with $\beta\in  \Aut^{(1)}_\bm2(M)$ and $\alpha \in \Aut_\bmod(M)$. Thus one obtains the equality:~ 
$
\Aut_\bm2(M)=\Aut^{(1)}_\bm2(M)\rtimes \Aut_\bmod(M)
$, 
using the natural action of $\Aut_\bmod(M)$ by conjugation on $\Aut^{(1)}_\bm2(M)$. \endproof

\subsection{Epimorphisms and  \strg exact sequences}\label{sectepistrong}

 There is a direct relation between \strg exactness and the categorical  notion of epimorphism in $\bm2$, it is given by the following 
   \begin{prop}\label{epimorphism} Let $\phi=(f,g)\in \Hom_\bm2(M,N)$. The following conditions are equivalent\newline
   1. $B(f,g)=N\times N$.
   
  2. The sequence
    $M\stackbin[g]{f}{\rightrightarrows}N\stackbin[0]{0}{\rightrightarrows}0$
    is  \strgly exact at $N$.
    
    3.	The morphism $\phi$ is an epimorphism in the category $\bm2$.
     \end{prop}
    \proof By Proposition \ref{surjtest} the \strg exactness of the sequence at $N$ is equivalent to $B(f,g)=N\times N$ and the addition of $\Delta$ is irrelevant in that case. By Lemma \ref{imlem} the equality  $B(f,g)=N\times N$ holds iff for any $\B$-module $X$ and any morphisms $\psi,\psi'\in \Hom_\B(N\times N, X)$ one has 
    $$
    \psi=\psi'\iff \psi\vert_{B(f,g)}=\psi'\vert_{B(f,g)}.
    $$
   Let $\psi_j$ (resp. $\psi'_j$) be the components of $\psi$ so that $\psi((a,b))=\psi_1(a)+\psi_2(b)$ (resp. $\psi'((a,b))=\psi'_1(a)+\psi'_2(b)$). One then derives 
   $$
   \psi((f(x)+g(y),f(y)+g(x))=\psi_1(f(x)+g(y))+\psi_2(g(x)+f(y))\qqq x,y\in N.
   $$
   Thus (taking $y=0$ or $x=0$) the condition $\psi\vert_{B(f,g)}=\psi'\vert_{B(f,g)}$ is equivalent to 
   $$
   \psi_1\circ f+\psi_2\circ g =\psi'_1\circ f+\psi'_2\circ g,\ \ \ \psi_1\circ g+\psi_2\circ f =\psi'_1\circ g+\psi'_2\circ f
   $$
   which means exactly, using composition in $\bm2$, that $\psi\circ \phi=\psi'\circ \phi$, where $\psi$ is viewed as the element $(\psi_1,\psi_2)\in \Hom_\bm2(N,X)$ and similarly for $\psi'$. Thus \strg exactness holds iff the morphism $\phi=(f,g)$ is an epimorphism in the category $\bm2$.\endproof

 By using the inclusion $f(N)\subset M$ and Proposition \ref{monoepi} every morphism  $f$ in $\bmod$ admits a factorization ``monomorphism $\circ$ epimorphism".  This fact no longer holds in the category $\bm2$ since the size of the image of a monomorphism $M\stackbin[g]{f}{\rightrightarrows}N$, when viewed as the  submodule $E=\{(f(x),g(x))\mid x\in M\}\subset N\times N$, is limited by the condition of injectivity of the map $(x,y)\mapsto (f(x)+g(y),g(x)+f(y)):M^2\to N^2$. More precisely, this condition implies the injectivity of  the map
   $
   E\times E\to N\times N, \ \ (\xi,\eta)\mapsto \xi+\sigma(\eta)
   $. 
   The failure of this condition ought to imply the impossibility  of a factorization ``monomorphism $\circ$ epimorphism" but one needs to be careful since the above definition of $E$ uses the diagonal.

 To test our guess, we reconsider the simplest Example \ref{example not mono}. 
 \begin{example} \label{notmonoepi}
 With the notations of Example \ref{example not mono}, the morphism $\phi:M \stackbin[p_2]{p_1}{\rightrightarrows}N$  does not admit a factorization ``monomorphism $\circ$ epimorphism" in the category $\bm2$. 	
   \end{example}
\proof Assume that $\phi= \psi\circ \eta$ in $\bm2$. By \eqref{iso8} one has  $Z(\eta)\subset Z(\psi\circ \eta)=Z(\phi)$ and $Z(\phi)=\Delta$ since the sequence $0 \stackbin[0]{0}{\rightrightarrows}M\stackbin[p_2]{p_1}{\rightrightarrows}N$ is \strgly exact. Thus $Z(\eta)=\Delta$ for any choice of $\eta$, and if moreover $\eta$ is an epimorphism then Proposition \ref{epimorphism} shows, together with Proposition \ref{iso9}, that it is an isomorphism.   The $\B$-module $M$ does not admit a non-trivial decomposition as a product, and thus the only possible choices for $\eta$ are $M \stackbin[0]{\rho}{\rightrightarrows}M$ or $M \stackbin[\rho]{0}{\rightrightarrows}M$, where $\rho$ is an isomorphism. This reduces the possible factorizations ``mono $\circ$ epi" to $\phi$ itself and since this map fails to be  a monomorphism, there is no factorization ``$\phi=$ mono $\circ$ epi" in the category $\bm2$. \endproof

This example points out a serious issue which is in fact independent of the definition of exactness and is formulated simply in terms of the category $\bm2$. This defect is resolved by the extension of the category $\bm2$ performed below in \S \ref{sect yoneda}.
\begin{rem}\label{hom functors} For $M$ a $\B$-module, let $M^*:=\Hom_\B(M,\B)$, and for $f\in \Hom_\B(M,N)$ let $f^*:N^*\to M^*$ be given by $f^*(\phi):=\phi\circ f$. This defines a contravariant endofunctor of $\bmod$ and also of $\bm2$ using $(f,g)^*:=(f^*,g^*)$ and the compatibility of the composition law \eqref{paircompbm2} with $f \mapsto f^*$. One shows 
\newline
  $(i)$~If the sequence 
  $N^*\stackbin[g^*]{f^*}{\rightrightarrows}M^*\stackbin[0]{0}{\rightrightarrows}0$ is  \strgly exact so is the sequence $0 \stackbin[0]{0}{\rightrightarrows}M\stackbin[g]{f}{\rightrightarrows}N$.\newline
  $(ii)$~If the sequence 
 $M\stackbin[g]{f}{\rightrightarrows}N\stackbin[0]{0}{\rightrightarrows}0$ is \strgly exact so is the dual sequence $0\stackbin[0]{0}{\rightrightarrows}N^*\stackbin[g^*]{f^*}{\rightrightarrows}M^*$.\newline
	Example \ref{example not mono} shows that the converse of both statements fails.
\end{rem}

\subsection{Quotients in $\bm2$}\label{sectquotbm2}

The issue of the existence of quotients in the category $\bm2$ arises naturally because one would like to  define the cohomology of the sequence 
\begin{equation}	
\label{exactseqter}
M \stackbin[\alpha_2]{\alpha_1}{\rightrightarrows}N\stackbin[\beta_2]{\beta_1}{\rightrightarrows}P
\end{equation}
 at $N$ as the quotient $H_N:=Z(\beta_1,\beta_2)/(\Delta+B(\alpha_1,\alpha_2))$. The submodule $\Delta+B(\alpha_1,\alpha_2)\subset Z(\beta_1,\beta_2)$ is not hereditary in general and the  quotient is ill defined. Thus  we proceed as follows by using functors and the Yoneda embedding. 	
\begin{prop}\label{repfunctor} Let $M$ be a $\B$-module and $N\subset M$ a submodule. \newline
$(i)$~The following equality defines a covariant additive functor $F=M/N:\bm2\longrightarrow {\rm Bmod}$
\begin{equation}\label{covfunc}
	M/N(X):=\{(f,g)\in \Hom_{\bm2}(M,X)\mid f(x)=g(x)\qqq x\in N\}.
\end{equation}
$(ii)$~The functor $F=M/N$ is a subfunctor of the representable functor $y_M:=\Hom_\bm2(M,-)$.\newline
$(iii)$~The map $(f,g)\mapsto (g,f)$  induces an involution $\sigma$ on $M/N(X)$ and $\sigma$ is the identity for any $X$ if and only if $N=M$.
\end{prop}
\proof $(i)$~The pairs $(f,g)$ fulfilling $f(x)=g(x), \forall x\in N$ form a $\B$-submodule of the $\B$-module $\Hom_{\bm2}(M,X)$. Let $X\stackbin[\beta]{\alpha}{\rightrightarrows} Y$ be a morphism in $\Hom_{\bm2}(X,Y)$. The composition $(\alpha,\beta)\circ (f,g)$ as in \eqref{paircompbm2} defines a map $F(\alpha,\beta):M/N(X)\to M/N(Y)$ since one has
$$
(\alpha\circ f+\beta\circ g)(x)=(\alpha\circ g+\beta\circ f)(x)\qqq x\in N.
$$
The map $F(\alpha,\beta)$ is additive by construction, \ie it defines a morphism of $\B$-modules. Moreover it depends additively upon the pair $(\alpha,\beta)$. \newline
$(ii)$~By construction the functor $F=M/N$ is a subfunctor of the Yoneda functor $y_M$ which associates to any $\B$-module $X$ the $\B$-module $\Hom_\bm2(M,X)$.\newline
$(iii)$~Lemma \ref{imlem} states that $N\neq M$ iff there exists a $\B$-module $E$ and two distinct morphisms $f,g\in \Hom_\B(M,E)$
which agree on $N$. 
\endproof
Proposition \ref{repfunctor} gives a meaning to the quotient $M/N$ as a covariant additive functor $F=M/N:\bm2\longrightarrow {\rm Bmod}$.  By the Yoneda Lemma  the opposite of the category $\bm2$ embeds fully and faithfully  in the category of covariant additive functors $\bm2\longrightarrow {\rm Bmod}$. The inclusion of $F=M/N$ as a subfunctor of the representable functor $y_M:\bm2\longrightarrow {\rm Bmod}$ corresponds to the ``quotient map" $M\to M/N$. The following example shows that the functor $M/N:\bm2\longrightarrow {\rm Bmod}$ is not representable in general.
\begin{example}\label{example not representable}
 We let $M=\{0,1,2\}$ with the operation of max, and take $N=\{0,2\}\subset M$. The functor $F$ associates, to a $\B$-module $X$, two morphisms $u,v\in \Hom_\B(M,X)$ which agree on $N$. With $a=u(1)$, $b=v(1)$, $c=u(2)=v(2)$ this corresponds to the subset of $X^3$ formed of triples $(a,b,c)$ such that $a\leq c$ and $b\leq c$. Let then $Q$ be the $\B$-module 
$$
Q:=\{0,\alpha,\beta,\alpha\vee \beta,\gamma\}, \ \ x\vee \gamma=\gamma\qqq x\in Q
$$
One has a natural identification $F(X)\simeq \Hom_\B(Q,X)$ by sending $\alpha\mapsto a$, $\beta\mapsto b$, $\gamma\mapsto c$. Assume that the functor $M/N:\bm2\longrightarrow {\rm Bmod}$ is represented by a $\B$-module $Z$. Then  $F(X)\simeq \Hom_\bm2(Z,X)=\Hom_\B(Z,X)^2$ and one gets a contradiction for $X=\B$ since $F(\B)$ has five elements.\end{example}

The above example suggests that  given a sub-$\B$-module $N\subset M$, one can define an analogue of the above module $Q$ as follows. 
\begin{prop}\label{cokerneldef} Let $N\subset M$ be a submodule. On $E=M\times M$ define the following relation 
\begin{equation}\label{represent functor}
(x,y)=z\sim z'=(x',y')\iff f(x)+g(y)=f(x')+g(y')\qqq X, \ (f,g)\in F(X)	
\end{equation}
where $F(X)$ is defined in \eqref{covfunc}. \newline 
$(i)$~The quotient $Q=E/\!\!\sim$ is the cokernel pair $\Cpr(\iota)$  of the inclusion  $\iota:N\to M$.\newline
$(ii)$~The canonical maps $\gamma_j:M\to \Cpr(\iota)$ of the sequence \begin{equation}\label{defncokerused}
	N\stackrel{\iota}{\to} M \stackbin[\gamma_2]{\gamma_1}{\rightrightarrows}\Cpr(\iota)\end{equation}
	are injective and the intersection of the $\gamma_j(M)$ is $\gamma_1(N)=\gamma_2(N)$.\newline
$(iii)$~One has   canonical isomorphisms of endofunctors  of $\bmod$
$$r:F\circ \kappa(X)\to \Hom_\B(Q,X), \  \  \gamma:\Hom_\B(Q,X)\to F\circ \kappa(X),$$
where $\kappa:\bmod\longrightarrow \bm2$ is defined in \eqref{functoriota}.
\end{prop}
\proof By construction \eqref{represent functor} defines an equivalence relation that is compatible with the addition, \ie $z_j\sim z'_j (j=1,2)\implies z_1+z_2\sim z'_1+z'_2$. Thus the quotient $Q=E/\!\!\sim$ is a well defined $\B$-module. \newline
$(i)$~By Definition \ref{propimagedefn}, $\Cpr(\iota)$ is the coequalizer $\coequ (\iota_{(2)})$ of the two maps $s_j\circ \iota$ where the $s_j:M\to M\oplus M$ are the canonical inclusions from $M$ to $M\oplus M$. A morphism $\rho\in \Hom_\B(M\oplus M,X)$ is given by a pair of morphisms $f,g\in \Hom_\B(M,X)$ satisfying the following rule
$$
\rho \circ s_1\circ \iota=\rho \circ s_2\circ \iota\iff (f,g)\in F(X).
$$
Thus the quotient $Q=E/\!\!\sim$ coincides with the construction of the coequalizer of Lemma \ref{doublelem0}. \newline
$(ii)$~Let $x\in N$, then $f(x)=g(x)$ for $(f,g)\in F(X)$ and one gets $(x,0)\sim (0,x)$. Note also that the map $\gamma_1:M\to Q$,
$\gamma_1(x)=(x,0)$ is injective since the pair $(\id_M,\id_M)$ belongs to $F(M)$. The same applies to $\gamma_2:M\to Q$,
$\gamma_2(x)=(0,x)$. By construction any element of $Q$ is a sum of two elements of $\gamma_j(M)$. Let $z\in \gamma_1(M)\cap \gamma_2(M)$, $z=\gamma_1(x)=\gamma_2(y)$. Since the pair $(\id_M,\id_M)$ belongs to $F(M)$ one gets $x=y$. By Lemma \ref{imlem} one gets $x\in N$ which shows that the intersection of the $\gamma_j(M)$ is $\gamma_1(N)=\gamma_2(N)$.\newline  
$(iii)$~We  compare $\Hom_\B(Q,X)$ and $F\circ \kappa(X)=F(X)$. First let $(f,g)\in F(X)$, then the map $E=M\times M\to X$, $(x,y)\mapsto f(x)+g(y)$ is compatible with the equivalence relation $\sim$ and thus induces a map $r(f,g):Q\to X$. One has $r(f,g)\circ \gamma_1=f$ and $r(f,g)\circ \gamma_2=g$. Conversely, let $\rho\in \Hom_\B(Q,X)$, then let $f=\rho\circ \gamma_1$, $g=\rho\circ \gamma_2$. The equivalence $(x,0)\sim (0,x)$ for $x\in N$ shows that $f(x)=g(x)$ for $x\in N$, \ie $(f,g)\in F(X)$. Denote $\gamma(\rho):=(\rho\circ \gamma_1,\rho\circ \gamma_2)$. Consider the maps $r:F(X)\to \Hom_\B(Q,X)$ and $\gamma:\Hom_\B(Q,X)\to F(X)$. For $\rho\in \Hom_\B(Q,X)$ one has $\rho=r(\gamma(\rho))$ since any $z\in Q$ is of the form $z=(x,y)=\gamma_1(x)+\gamma_2(y)$ while $\rho(z)=f(x)+g(y)$, where $f=\rho\circ \gamma_1$, $g=\rho\circ \gamma_2$. Let similarly $(f,g)\in F(X)$
then $r(f,g)\circ \gamma_1=f$ and $r(f,g)\circ \gamma_2=g$, so that $(f,g)=\gamma(r(f,g))$. Thus the maps $r$ and $\gamma$ are inverse of each other and give a natural identification 
$F(X)\simeq \Hom_\B(Q,X)$.
\endproof

\begin{rem} Consider the sequence \eqref{exactseqter}. As in Definition \ref{paircomp1},  the equalizer of the $\beta_j$, 
 $\iota:\equ(\beta_1,\beta_2)\to N$  is a submodule of $N$
 and the 	coequalizer of the $\alpha_j$  is a quotient of $N$, \ie $\gamma:N \to \coequ(\alpha_1,\alpha_2)$. One can thus consider the composition $\gamma\circ \iota:\equ(\beta_1,\beta_2)\to \coequ(\alpha_1,\alpha_2)$ which is a morphism of $\B$-modules, and  define the {\em  weak cohomology} $\hweak_N$ at $N$ of the sequence \eqref{exactseqter} as:
\begin{equation}	
\label{defnweakcohom}
\hweak_N:=\Range(\gamma\circ \iota) , \ \ \gamma\circ \iota:\equ(\beta_1,\beta_2)\to \coequ(\alpha_1,\alpha_2).
\end{equation}
If the sequence \eqref{exactseqter} is  \strgly exact at $N$ then $\hweak_N=0$. 
Indeed, if $\hweak_N\neq 0$
	 there exists a $\B$-module $E$ and $\phi\in \Hom_\B(N,E)$ such that $\phi\circ \alpha_1=\phi\circ \alpha_2$, while 
	the restriction of $\phi$ to $\equ (\beta_1,\beta_2)$ is non-zero, \ie there exists $t\in N$, $\beta_1(t)=\beta_2(t)$ and $\phi(t)\neq 0$. One then has $(t,0)\in Z(\beta_1,\beta_2)$ while $(t,0)\notin B(\alpha_1,\alpha_2)+\Delta$. Indeed otherwise let $x,y\in M$, $z\in N$ such that
	$
	(t,0)=(\alpha_1(x)+\alpha_2(y)+z,\alpha_2(x)+\alpha_1(y)+z)
	$. 
	By applying $\phi$ one gets that $(\phi(t),0)$ is diagonal which contradicts $\phi(t)\neq 0$. \newline
The converse does not hold since for $\psi$ the maximal element of $M^*$ the sequence:
$
0 \stackbin[0]{0}{\rightrightarrows}M\stackbin[0]{\psi}{\rightrightarrows}\B\stackbin[0]{0}{\rightrightarrows}0
$
fulfills $\hweak_M=0$ and $\hweak_\B=0$ but 
 by Corollary \ref{isotest}  this sequence is  \strgly exact  only when  $M=\B$.
 \end{rem}

\section{The  Eilenberg-Moore category $\b2$ of the comonad $\perp$}\label{secthomologbmods}

In this Section (\cf~\S\S\ref{sect yoneda},\ref{sect yoneda1}) we take up the problem of the representability of the functor associated to quotients of objects of $\bmod$ and we provide a solution by considering the natural extension of $\bm2$ to the Eilenberg-Moore category $\b2$ of Proposition \ref{Kleisli} which is simply the category $\bmod$ in the topos of ``sets endowed with an involution" (and as such shares with $\bmod$ most of its abstract categorical properties).
\begin{defn}\label{catb2} Let $\b2$ be the category of $\B$-modules endowed with an involution $\sigma$. The morphisms in $\b2$ are the morphisms of $\B$-modules commuting with $\sigma$, \ie equivariant for the action of $\Z/2\Z$. 
\end{defn} 	
 In \S\ref{grandis remark} we provided the conceptual meaning of the construction of  $\bm2$ and $\b2$ from the category $\bmod$.  In \S \ref{sect yoneda} we provide another direct construction of $\b2$ based on the need to represent the functor associated to quotients. In \S \ref{sect yoneda1} we prove the required representability.  In \S \ref{sectmonad} we analyze the  monad $T$ on $\b2$ corresponding to the adjunction  and show its simplifying role. The notion of \strg exact sequence extends naturally from $\bm2$ to $\b2$ and we end the section with a table comparing these two categories.
   
 \subsection{Extending the category $\bm2$}\label{sect yoneda}

For any $\B$-module $M$ one gets the covariant functor $y_M(-):\bm2\longrightarrow \bmod$, $y_M(N)= \Hom_\bm2(M,N)$. For $M=\B$   the functor  $y_\B(-)$ lands in a finer category than $\bmod$ since one can use the group $S:=\Aut_\bm2(\B)$ to act on the right on morphisms. The only non-trivial element of $S$ is $\sigma=(0,\id)$ and it has order two. One has $(f,g)\circ (0,\id)=(g,f)$,  thus this action exchanges the two copies of $N$ in $y_\B(N)=\Hom_\bm2(\B,N)=N\times N$. One thus obtains a functor $y_\B(-):\bm2\longrightarrow \b2$ which embeds 
the category $\bm2$ as a full subcategory of  $\b2$.
 \begin{lem}\label{yb} $(i)$~The functor $y_\B(-):\bm2\longrightarrow \b2$ associates to a $\B$-module $N$ the square $y_\B(N)=\Hom_\bm2(\B,N)=N\times N$ endowed with the involution which exchanges the two copies of $N$ and to a morphism $\phi=(f,g)\in \Hom_\bm2(N,N')$ the following map $y_\B(\phi):y_\B(N)\to y_\B(N')$
\begin{equation}\label{myyb}
 N\times N\ni (x,y)\mapsto (f(x)+g(y),f(y)+g(x))\in N'\times N'.
\end{equation}	
$(ii)$~The functor $y_\B(-):\bm2\longrightarrow \b2$ is fully faithful.  	\newline
$(iii)$~For any morphism $\phi=(f,g)\in \Hom_\bm2(N,N')$ one has 
 \begin{equation}\label{yb2}
B(\phi)={\rm Range}(y_\B(\phi)), \ \ \  Z(\phi)=y_\B(\phi)^{-1}(\Delta)
\end{equation}
where $\Delta$ is the diagonal in $N'\times N'$.
 \end{lem}
  \proof $(i)$~An element of $y_\B(N)=\Hom_\bm2(\B,N)$ is given by a pair of morphisms of $\B$-modules  from $\B$ to $N$ and hence characterized by a pair of elements of $N$. The composition of morphisms in $\bm2$ can be seen  easily  to give \eqref{myyb}.\newline
 $(ii)$~Let $\phi=(f,g), \phi'=(f',g')$ be elements of $\Hom_\bm2(N,N')$  such that $y_\B(\phi)=y_\B(\psi)$. One then gets, taking $y=0$ in \eqref{myyb}  
 $
 (f(x),g(x))=(f'(x),g'(x))\qqq x\in N
 $.
 This proves that $\phi=\phi'$.  We prove that $y_\B(-)$ is full \ie that for any morphism $\rho: y_\B(N)\to y_\B(N')$ commuting with $\sigma$ there exists a morphism $\phi\in \Hom_\bm2(N,N')$ such that $\rho=y_\B(\phi)$. The restriction of $\rho$ to $N\times \{0\}\subset N\times N=y_\B(N)$ is given by a pair of morphisms $f,g\in \Hom_\B(N,N')$. Since $\rho$ commutes with $\sigma$ one has
  $
  \rho((x,y))=\rho((x,0))+\rho(\sigma((y,0)))=(f(x),g(x))+\sigma((f(y),g(y))
  $
which shows by \eqref{myyb} that $\rho=y_\B(\phi)$ for $\phi=(f,g)$.\newline
 $(iii)$~Follows from \eqref{myyb} and Propositions \ref{ker1} and \ref{coker1}.
  \endproof 
  Thus the notion of \strg exact sequence in $\bm2$ can be re-interpreted, through the functor  $y_\B(-):\bm2\longrightarrow \b2$, in terms of the notions of range and of inverse image in the category $\b2$ but involving the additional structure given by the diagonal in the square which replaces the zero element for abelian groups. This implies that rather than asking the composition of two consecutive maps in a complex to be zero, one requires instead  ${\rm Range}(y_\B(\phi\circ \psi))\subset \Delta$ or equivalently  ${\rm Range}(y_\B(\psi))\subset y_\B(\phi)^{-1}(\Delta)$, \ie $B(\psi)\subset Z(\phi)$. Moreover  Proposition \ref{monomorphism} can be now re-interpreted stating that a morphism $\phi\in \Hom_\bm2(N,N')$ is a monomorphism iff $y_\B(\phi)$ is injective, \ie a monomorphism in $\b2$. Similarly, Proposition \ref{epimorphism} now states that a morphism $\phi\in \Hom_\bm2(N,N')$ is an epimorphism in $\bm2$ iff $y_\B(\phi)$ is surjective, \ie an epimorphism in $\b2$.

Next lemma states that the forgetful functor $I:\b2\longrightarrow \bmod$ is left adjoint to the squaring functor  $\yb=y_\B\circ \kappa:\bmod\to\b2$,  $\yb(M)=(M^2,\sigma)$, which is obtained by composition of $y_\B$ with the functor $\kappa:\bmod\longrightarrow \bm2$ of \eqref{functoriota}.
\begin{lem}\label{forgetful}Let $I:\b2\longrightarrow \bmod$ be the forgetful functor. The composition with the first projection gives, for any object $M$ of $\b2$ and any object $N$ of $\bmod$, a canonical isomorphism
\begin{equation}	
\label{canonicaliso}
\pi:\Hom_\b2(M,\yb(N))\to \Hom_\B(I(M),N).
\end{equation}
which is natural in $M$ and $N$.  \end{lem}
  \proof Let $\phi\in \Hom_\b2(M,\yb(N))$, then the composition $p_1\circ \phi$  with the first projection fulfills  
$p_1\circ \phi\in \Hom_\B(I(M),N)$. Let $\psi\in \Hom_\B(I(M),N)$, define $
\tilde \psi \in \Hom_\b2(M,\yb(N))$ by $\tilde \psi(z)=(\psi(z),\psi(\sigma(z)))$ for $z\in M$. It commutes with $\sigma$ by construction and $p_1\circ\tilde \psi=\psi$. Moreover, for $\phi\in \Hom_\b2(M,\yb(N))$ one has for any $z\in M$, 
$
\phi(z)=(p_1\circ \phi(z),p_2\circ \phi(z))=(p_1\circ \phi(z),p_1\circ \phi(\sigma(z)))
$, 
so that $\phi=\widetilde{p_1\circ \phi}$. By construction the isomorphism $\pi$ is natural in $M$. The naturality in $N$ follows from $g\circ p_1=p_1\circ \yb(g)$, $\forall g\in \Hom_\B(N,N')$. \endproof 
To the adjunction $I\dashv \yb$ corresponds the  comonad  $\perp$ for the category $\bmod$ of Proposition \ref{comonad}. Indeed, one has $\perp=I\circ \yb$ which proves 1 (of Proposition \ref{comonad}). The counit of the adjunction is the first projection and this gives 2. Using the unit $\eta$ of the adjunction, one gets the coproduct $\delta$ as $I\circ \eta\circ \yb:I\circ \yb\longrightarrow I\circ \yb\circ I\circ \yb$ and since $\sigma((x,y))=(y,x)$ one obtains 3. The adjunction\footnote{This remark was pointed to us by M. Grandis} $I\dashv \yb$  is comonadic and this corresponds to the construction of $\b2$ in \S \ref{grandis remark}.
  
  \subsection{Representability in $\b2$}\label{sect yoneda1}
  
   Next, for a submodule $N\subset M$ of a $\B$-module $M$ we shall show that the   covariant functor $F=M/N:\bm2\longrightarrow {\rm Bmod}$ (\cf~\eqref{covfunc})
\begin{equation}	
\label{cohomseqter1}
F(X)=\{(f,g)\mid f,g\in \Hom_\B(M,X),\ f\vert_N=g\vert_N\}.
\end{equation}
extends to a representable functor $\Hom_\b2(Q,-):\b2\longrightarrow {\rm Bmod}$. Thus we look for an object $Q$ of $\b2$ such that, for any $\B$-module $X$, one has a natural isomorphism,  
\begin{equation}	
\label{cohomseqter1bis}
F(X) \simeq \Hom_\b2(Q,y_\B(X))
\end{equation}

We first reconsider Example \ref{example not representable} and take the category $\b2$ into account.
\begin{example}\label{example not representable1}
 Let $M=\{0,1,2\}$, $N=\{0,2\}\subset M$ and $Q$  as in Example \ref{example not representable}.   A morphism $\phi=(f,g)\in\Hom_\bm2(X,X')$ maps $(a,b,c)$ to $(f(a)+g(b),f(b)+g(a), f(c)+g(c))$ and this can be re-interpreted in the following form 
$$
\Hom_\B(Q,X)\ni \tau \mapsto f\circ \tau+ g\circ \tau\circ \sigma
$$
where $\sigma\in \Aut_\bmod(Q)$ is the involution which interchanges $\alpha$ and $\beta$ and fixes $\gamma$. This suggests to view $Q$ as an object of $\b2$ and to compare $F(X)$ with $\Hom_\b2(Q,y_\B(X))$. Let $\rho\in\Hom_\b2(Q,y_\B(X))$, then with $\rho(\alpha)=(a,b)\in X\times X$, one has  $\rho(\beta)=\rho(\sigma(\alpha))=\sigma((a,b))=(b,a)\in X\times X$, while $\rho(\gamma)$ is fixed by $\sigma$ and hence diagonal, \ie of the form $\rho(\gamma)=(c,c)$ with $a\vee c=c$, $b\vee c=c$. Thus one obtains: $F(X)\simeq \Hom_\b2(Q,y_\B(X))$, and the composition $\tau \mapsto f\circ \tau+ g\circ \tau\circ \sigma$ corresponds to \eqref{myyb}. This shows the representability of the functor $F$ in the category $\b2$.
\end{example}
The above example suggests that the quotient functor $M/N$ ought to be representable in the category $\b2$.   We take the notations of Proposition \ref{cokerneldef}.
\begin{prop}\label{repb2}
$(i)$~The involution  $(x,y)\mapsto (y,x) $ on $E=M\times M$ induces an involution $\sigma$ on $Q\simeq \Cpr(\iota)$ which turns $\Cpr(\iota)$ into an object of $\b2$.	\newline
 $(ii)$~The functor  $F:\bm2\longrightarrow \bmod$,  associated by \eqref{cohomseqter1} to the quotient $``M/N"$ is represented  in the category $\b2$ 	by the object $\Cpr(\iota: N\to M)=(Q=E/\!\!\sim,\sigma)$. 
  \end{prop}
  \proof $(i)$~The congruence $\sim$ on $E$ is compatible with the involution $\sigma((x,y))=(y,x)$, \ie $z\sim z'\Rightarrow \sigma(z)\sim \sigma(z')$ since $(f,g)\in F(X)\iff (g,f)\in F(X)$. Thus one obtains an involution, still denoted by $\sigma$, on the $\B$-module $Q$ and this turns it into an object of the category $\b2$. \newline
 $(ii)$~Proposition \ref{cokerneldef} shows  that one has $F(X)\simeq \Hom_\B(Q,X)$, with  $Q=I((Q,\sigma))$. Thus by Lemma \ref{forgetful} one deduces a canonical isomorphism
    $
  F(X)  \simeq \Hom_\b2((Q,\sigma),y_\B(X))
  $.
It remains to show that this identification is compatible with the functoriality. To this end, one needs to keep track on the description of both sides after composing with a morphism $y_\B(\phi):y_\B(X)\to y_\B(X')$ where by \eqref{myyb} and with $\phi=(\alpha,\beta)$, the map $y_\B(\phi)$ is given by $(x,y)\mapsto (\alpha(x)+\beta(y),\alpha(y)+\beta(x))$. Let then $(f,g)\in F(X)$, then the corresponding pair $(f',g')\in F(X')$ is given by 
$
f'=\alpha\circ f+\beta\circ g, \ \ g'=\alpha\circ g+\beta \circ f
$. 
Now, to $(f,g)\in F(X)$ corresponds the morphism $r(f,g):Q\to X$ with $r(f,g)((x,y))=f(x)+g(y)$. Moreover the associated morphism $\rho(f,g)$ in $\Hom_\b2((Q,\sigma),y_\B(X))$ is given by 
$
\rho(f,g)((x,y))=(f(x)+g(y),f(y)+g(x))
$. 
Thus one has 
$$
y_\B(\phi)(\rho(f,g)((x,y)))=(\alpha(f(x)+g(y))+\beta(f(y)+g(x)),\alpha(f(y)+g(x))+\beta(f(x)+g(y))).
$$
One checks that this is the same as 
$$
\rho(f',g')((x,y))=((\alpha\circ f+\beta\circ g)(x)+(\alpha\circ g+\beta \circ f)(y),
(\alpha\circ f+\beta\circ g)(y)+(\alpha\circ g+\beta \circ f)(x))
$$
This proves that the two functors are the same. \endproof

\subsection{The monad $T=\yb\circ I$ on the category $\b2$ and \strg exactness} \label{sectmonad}


We consider the monad associated to the adjunction $I\dashv \yb$ displayed in Lemma \ref{forgetful}. We determine the counit and the unit of this adjunction. The counit is a natural transformation $I\circ \yb\longrightarrow\id_\bmod$ which associates to an object $M$ of $\bmod$ a morphism $\epsilon_M\in \Hom_\B(I\circ \yb(M),M)$. By Lemma \ref{forgetful}, $\epsilon_M$ corresponds to the identity $\id\in \Hom_\b2(\yb(M),\yb(M))$. In fact the proof of that lemma shows that $\epsilon_M$ is the first projection $p_1:M\times M\to M$. 
\begin{prop}\label{monad} The monad associated to the adjunction $I\dashv \yb$ is described by 
\begin{enumerate}
	\item The endofunctor $T:\b2\longrightarrow \b2$,  $T(M,\sigma)=(M^2,\sigma_M)$, $\sigma_M(x,y)=(y,x)$, and for $f\in \Hom_\b2(M,N)$, $T(f)=(f,f)$.
	\item The unit $\eta: 1_\b2\to T$,  $\eta_{(M,\sigma)}=((M,\sigma)\to (M^2,\sigma_M),~  a\mapsto (a,\sigma(a)))$.
	\item The product $\mu:T^2\to T$, $\mu_{(M,\sigma)}= (((M^2)^2,\sigma_{M^2})\to (M^2,\sigma_M),~((x,y),(x',y'))\mapsto (x,x'))$. 
\end{enumerate}	
\end{prop}
\proof From the general construction of a monad associated to an adjunction, the endofunctor is of the form $G\circ F$, where $F\dashv G$. The adjunction of Lemma \ref{forgetful} is $I \dashv \yb$ and thus $T=\yb\circ I$ which shows 1. The unit $\eta$ of this adjunction is a natural transformation $\id_\b2\longrightarrow \yb \circ I$. It associates to an object $(M,\sigma)$ of $\b2$ a morphism $\eta_M\in \Hom_\b2(M,\yb \circ I(M))$. This morphism  corresponds to the identity $\id\in \Hom_\B(I(M),I(M))$ and the proof of Lemma \ref{forgetful} shows that it is given, for every object $(M,\sigma)$ of $\b2$, by the morphism
$
M \ni a\mapsto (a,\sigma(a))\in T(M)=(M^2,\sigma_M)
$.
The product $\mu$ is a natural transformation  $T\circ T\longrightarrow T$ and for an adjunction $F\dashv G$ is given by $G\epsilon F$ where $\epsilon$ is the counit of the adjunction. The counit of the adjunction $I \dashv \yb$ is given by the first projection $p_1:M\times M\to M$. Thus $\mu=\yb p_1 I$ as a natural transformation $T\circ T(M)=\yb\circ I\circ \yb\circ I(M)\to  T(M)=\yb\circ I(M)$. To obtain it one applies the functor $\yb$  to the morphism $\epsilon_{I(M)}\in \Hom_\B(I\circ \yb(I(M)),I(M))$. Given a morphism $f\in \Hom_\B(N,N')$ the morphism $\yb(f)\in \Hom_\b2(\yb(N),\yb(N'))$ acts diagonally as $(z,z')\mapsto (f(z),f(z'))$. This gives the description 3. of the product. \endproof 
To clarify the fact that the cokernel of a morphism of $\bmod$ should be viewed as an object of $\b2$, we continue with the investigation of the monad $T=\yb \circ I$ defined in Proposition \ref{monad}. 
There is a natural notion of an injective object in a category $\cC$ endowed with a monad $T$, it is given by the following
\begin{defn}\label{T-injective def} An object $I$ of $\cC$ is $T$-injective if the unit map $\eta_I:I\to T(I)$ has a retraction, \ie there is a map $f:T(I)\to I$ such that $f\circ \eta_I=\id_I$.	
\end{defn}
In our setup we derive
\begin{lem}\label{T-injective} Any object of $\b2$ is $T$-injective for the monad $T=\yb \circ I$.	
\end{lem}
\proof Let $(M,\sigma)$ be an object of $\b2$. Let $r:T(M)\to (M,\sigma)$ be defined by
$
r((x,y))=x+\sigma(y),~ \forall x,y\in M
$. 
 One has $r(\sigma_M(x,y))=r((y,x))=y+\sigma(x)=\sigma(x+\sigma(y))=\sigma(r((x,y)))$. Thus $r$ is a morphism in $\b2$. Moreover one also has
 $
 r(\eta_M(a))=r((a,\sigma(a)))=a+a=a,~ \forall a\in M
 $. This shows that 
$r$ is a retraction of $\eta_M$.\endproof 
Lemma \ref{T-injective} will be used below in the proof of Lemma \ref{T-injectiveused}.\vspace{.05in} 

In the following part, we shall recast  our previous main constructions and results in terms of the category $\b2$ and its monad structure. For example, Lemma \ref{yb} shows that the notion of \strg exactness  in the category $\bm2$ as provided in Definition \ref{exsequ}, corresponds to the following definition in $\b2$. 
\begin{defn}\label{exactnessb2} A sequence $L\stackrel{f}{\to} M\stackrel{g}{\to}N$ in $\b2$ is \strgly exact at $M$ if 
\begin{equation}\label{exactnessb21}
	{\rm Range}(f)+M^\sigma=g^{-1} (N^\sigma). 
\end{equation}
 \end{defn}
 Notice that this definition implies the weaker condition ${\rm Range}(g\circ f)\subset N^\sigma$. The heuristic behind is that the fixed points $M^\sigma$ of an object $M$ of $\b2$  play the role of the zero element.  The kernel and cokernel of a morphism in $\b2$ are defined as follows 
 
 \begin{defn}\label{ker&coker} For $h\in \Hom_\b2(L,M)$  one sets
 $	\Ker(h):=h^{-1} (M^\sigma) $ 
endowed with the induced involution.  One also sets $\coker(h):=M/\!\!\sim$ ~where for $a,a'\in M$
 \begin{equation}\label{cokernelb2}
 a\sim a'\iff f(a)=f(a')~\ ~ \forall X, \forall f\in \Hom_\b2(M,X)~{\rm s.t.}~ {\rm Range}(h)\subset \Ker(f).
\end{equation}
\end{defn}
The above formula defines an equivalence relation which is compatible with the addition and the  involution and thus the quotient $\coker(h):=M/\!\!\sim$ is an object of $\b2$. We shall show below in Proposition \ref{grandis} that the above notions of kernel and cokernel are the natural ones in the context of homological categories of \cite{gra}. What remains to be seen is whether this notion of cokernel is compatible with the notion given earlier on, in Definition \ref{propimagedefn} \ie as $\Cpr h = \coequ (s_1\circ h,s_2\circ h)$. The issue is to control the morphisms $f\in \Hom_\b2(y_\B(M),X)$ coming from objects $X$ of $\b2$ which are not in the range of $y_\B$. This is indeed possible thanks to Lemma \ref{T-injective} and one obtains 
 \begin{lem}\label{T-injectiveused} Let $h\in \Hom_\b2(L,M)$, then the equivalence relation $\sim$ of \eqref{cokernelb2} defining its cokernel is the same as the following, for $a,a'\in M$
 \begin{equation}\label{cokernelb2bis}
 a\sim a'\iff f(a)=f(a')~\ ~ \forall X, \forall f\in \Hom_\b2(M,T(X))~{\rm s.t.}~ {\rm Range}(h)\subset \Ker(f).
\end{equation}
 \end{lem}
 \proof Let $X$ be an object of $\b2$. We show that if $g(a)=g(a')\ \forall g\in \Hom_\b2(M,T(X))$~ s.t. ${\rm Range}(h)\subset \Ker(g)$, one has $f(a)=f(a')\ \forall f\in \Hom_\b2(M,X)~{\rm s.t.}~ {\rm Range}(h)\subset \Ker(f)$. Let $\eta_X$  and $r$ be as in Lemma \ref{T-injective}. Let $f\in \Hom_\b2(M,X)~{\rm s.t.}~ {\rm Range}(h)\subset \Ker(f)$ and $g=\eta_X\circ f$. One has $g\in \Hom_\b2(M,T(X))~{\rm s.t.}~ {\rm Range}(h)\subset \Ker(g)$, since $\Ker(f)\subset \Ker(\eta_X\circ f)$. Thus one has by hypothesis $g(a)=g(a')$. Since $r\circ g=r\circ\eta_X\circ f=f$, one concludes that $f(a)=f(a')$. \endproof  
 \begin{prop}\label{cokernelequiv} Let $f:M\to N$ be a morphism in $\bmod$.\newline
 $(i)$~$\Kpr f=(\equ f^{(2)},\iota_1,\iota_2)$  as in \eqref{defnker}, endowed with its canonical involution is  isomorphic as an object of $\b2$ to $\ker(\yb f)$ as in  Definition \ref{ker&coker}.\newline
 $(ii)$~$
\Cpr f = \coequ (s_1\circ f,s_2\circ f)$ as in Definition \ref{propimagedefn}, endowed with its canonical involution, is isomorphic as an object of $\b2$ to $\coker(\yb f)$ as in  \eqref{cokernelb2}. 
 \end{prop}
 \proof $(i)$~ By definition $\equ f^{(2)}$ is the equalizer of $(f\circ p_1,f\circ p_2)$: 
$
\equ f^{(2)} \to M^2 \stackbin[f\circ p_2]{f\circ p_1}{\rightrightarrows} N
$, and coincides as a subobject of $M^2$ with $(\yb f)^{-1}(\Delta)$ where $\Delta= (N^2)^\sigma$ is the diagonal.\newline
 $(ii)$~By Lemma \ref{T-injectiveused}, and using $T=\yb\circ I$,    $\coker(\yb f)$ of \eqref{cokernelb2} is the quotient of $\yb N=N^2$ by the equivalence relation, where $X$ varies among objects of $\bmod$, and  $a,a'\in \yb N$
 \begin{equation}\label{cokernelb2ter}
 a\sim a'\iff \phi(a)=\phi(a')~\ ~ \forall X, \forall \phi\in \Hom_\b2(\yb (N),\yb (X))~{\rm s.t.}~ {\rm Range}(\yb f)\subset \Ker(\phi).
\end{equation}
 By Lemma \ref{yb} one has a natural isomorphism 
 $\Hom_\b2(\yb (N),\yb (X)) \simeq \Hom_\bm2(N,X)$  defined by the map which associates to a morphism  $\phi\in\Hom_\b2(\yb (N),\yb (X))$ given by the  matrix $\left(
\begin{array}{cc}
 u & v  \\
 v & u 
\end{array}
\right)$ the morphism in $\bm2$ given by $N\stackbin[v]{u}{\rightrightarrows} X$. The condition 
${\rm Range}(\yb f)\subset \Ker(\phi)$ is equivalent to $u\circ f=v\circ f$, since it means that $\phi\circ \yb f$ is null (\ie fixed by the involution). This condition means exactly that the pair $(u,v)$ defines a morphism $N^2\to X$ which coequalizes the $s_j\circ f$. Thus the equivalence relation
\eqref{cokernelb2ter} in $\yb N=N^2$ takes the form
$$
(x,y)\sim (x',y')\iff u(x)+v(y)=u(x')+v(y')\qqq u,v\in \Hom_\B(N,X)~{\rm s.t.}~ u\circ f=v\circ f.
$$
This is exactly the same as the equivalence relation which defines the coequalizer of the $s_j\circ f$.
\endproof 
When $\sigma_M=\id$, the \strg exactness in $\b2$ of the sequence $L\stackrel{f}{\to} M\stackrel{g}{\to}N$ is automatic. This shows that the \strg exactness of $0 \to L\stackrel{f}{\to}M$  in  $\b2$  does not imply that $f$ is a monomorphism, and similarly that the \strg exactness of $M\stackrel{f}{\to}N\to 0$ in  $\b2$ does not imply that $f$ is an epimorphism.  The application of the monad $T$ improves these statements considerably.
 \begin{prop}\label{use of T} Let $f\in \Hom_\b2(L,M)$, then
 
1. $f$ is a monomorphism  $\iff$ $f$ is injective $\iff$ $0 \to TL\stackrel{Tf}{\to}TM$ is \strgly exact.

2. $f$ is an epimorphism  $\iff$ $f$ is surjective $\iff$ $TL\stackrel{Tf}{\to}TM\to 0$ is \strgly exact.

\end{prop}
\proof 1. $f$ is a monomorphism if the underlying map of $\B$-modules is injective. Conversely, for $a\in L$ there exists a unique morphism $\xi_a:(\B\times \B,\sigma_\B)\to (L,\sigma)$ such that $\xi_a((1,0))=a$. Thus if for some $a\neq b\in L$ one has $f(a)=f(b)$, $f$ cannot be a monomorphism. Hence  $f$ is a monomorphism in $\b2$ $\Leftrightarrow$ $I(f)$ is a monomorphism in $\bmod$ $\Leftrightarrow$ $f$ is injective. Moreover, the \strg exactness of $0 \to TL\stackrel{Tf}{\to}TM$ is equivalent to the injectivity of $f$.  Indeed one has $T(f)=y_\B((I(f),0))$ since $T=\yb\circ I=y_\B\circ \kappa\circ I$ and the equivalence follows from $(ii)$ of Proposition \ref{injtest} since the faithful functor $y_\B$ preserves \strg exactness. \newline
2. $f$ is an epimorphism if the underlying map of $\B$-modules is surjective. Conversely, if the underlying map of $\B$-modules is not surjective there exists by Lemma \ref{imlem} a $\B$-module $X$ and a pair $h\neq k$ of morphisms of $\B$-modules, $h,k\in \Hom_\B(M,X)$, such that $h\circ f=k\circ f$.  Then the corresponding morphisms $\tilde h(x)=(h(x),h(\sigma(x))$, $\tilde k(x)=(k(x),k(\sigma(x))$ fulfill 
$$
\tilde h,\, \tilde k\in \Hom_\b2(M,(X\times X,\sigma_X)), \ \ \tilde h\circ f=\tilde k\circ f, \ \ \tilde h\neq \tilde k
$$
where the equality $h(\sigma(f(x)))=k(\sigma(f(x)))$ follows from $\sigma(f(x))=f(\sigma(x))$. This shows that $f$ is an epimorphism in $\b2$ $\Leftrightarrow$ $I(f)$ is an epimorphism in $\bmod$ $\Leftrightarrow$ $f$ is surjective. Finally, the \strg exactness of $TL\stackrel{Tf}{\to}TM\to 0$ is equivalent to the surjectivity of $f$ by $(ii)$ of Proposition \ref{surjtest}. \endproof
In the next table we compare the interpretation of several definitions in the categories $\bm2$ and $\b2$.

 \begin{table}[H]
\begin{center}
\label{tab:1}       
%
%
\begin{tabular}{| c | c |}
\hline
&\\
\begin{Large}{\bf Category $\bm2$}\end{Large} &\begin{Large}\bf{ Category $\b2$}\end{Large} \\ &\\
\hline &\\
$\phi=(u,v)\in \Hom_\bm2(M,N)$ &  $h\in \Hom_\b2(M,N)$\\  &\\
\hline &\\
$Z(\phi)=Z(u,v)$ & $\Ker(h)=h^{-1} (N^\sigma)$ \\
&\\
\hline &\\
$B(\phi)=\{(u(x)+v(y),u(y)+v(x))\}$& ${\rm Range}(h)=h(M)$\\ 
&\\
\hline &\\
Diagonal $\Delta\subset M\times M$ & Fixed points $M^\sigma$\\  & \\ 
\hline &\\
\strg exactness at $M$ of $L \stackbin[\alpha_2]{\alpha_1}{\rightrightarrows}M\stackbin[\beta_2]{\beta_1}{\rightrightarrows}N$& \strg exactness at $M$ of $L\stackrel{f}{\to} M\stackrel{g}{\to}N$\\ &\\$ B(\alpha_1,\alpha_2)+\Delta=Z(\beta_1,\beta_2)$ & ${\rm Range}(f)+M^\sigma=\Ker(g)$\\ &\\
\hline &\\
$\beta_1\circ \alpha_1+\beta_2\circ \alpha_2=\beta_2\circ \alpha_1+\beta_1\circ \alpha_2$ & ${\rm Range}(f)\subset \Ker(g)$\\ &\\
\hline
\end{tabular}
\end{center}
\end{table}

\section{The category $\b2$ and the null morphisms}\label{sectnullmorphisms}
By construction,  $\b2$ is the category $\bmod$ in the topos of sets endowed with an involution. The category $\b2$ contains, as a subcategory stable under retract, the category $\bmod$ of those objects of $\b2$, called ``null objects", whose involution is the identity. 
The morphisms which factor through a null object are called the null morphisms. They form   an ideal 
\cite{Ehresmann, Lavendhomme} in the category $\b2$. In \S\ref{secthomol} we show
that one obtains in this way a homological category in the sense of \cite{gra1}. In \S\ref{sectkerb2} we prove two results which control the least normal subobject containing a given subobject: \cf~Proposition \ref{existphi1} and Proposition \ref{extreme and range}. Finally, in \S\ref{sectcokerb2} (Proposition \ref{compute cokernel}) we provide  an explicit description of the cokernel of morphisms in $\b2$.

\subsection{$\b2$ as a semiexact homological category}\label{secthomol} 
  The general notion of semiexact category has been developed by Grandis in \cite{gra,gra1}. We view $\b2$ as the category of $\bs$-modules where $\bs$ is the semiring generated over $\B$ by $s$, $s^2=1$ (\!\cite{Gaubert}, p 71 and \cite{Akian}, Definition 2.18). One has $\bs=\{0,1,s,p\}$ where $p=1+s$, $p^2=p$. We use the notion of null morphisms associated to the ideal $N=\{0,p\}\subset \bs$. We shall  show that the category $\b2$ is semiexact. 
 \begin{prop}\label{grandis} The pair given by the category $\b2$ and the null morphisms: $\cN\subset \Hom_\b2(L,M)$
 $$
 f\in\cN \iff f(x)=\sigma(f(x))\qqq x\in L
 $$
 forms a semiexact category in the sense of \cite{gra,gra1}. The corresponding notions of kernel and cokernel are the same as in Definition \ref{ker&coker}.	
 \end{prop}
\proof First we show that $\cN$ is a closed ideal in $\b2$. For $f\in \cN$ one has $g\circ f\in \cN$, $\forall g\in \Hom_\b2(M,N)$ since $g(M^\sigma)\subset N^\sigma$. One has also $f\circ h\in \cN$, $\forall h\in\Hom_\b2(E,L)$. The closedness means that any null morphism factors through a null identity. This is the case here since any $f\in\cN\subset \Hom_\b2(L,M)$ factors through $M^\sigma$ which is a null object (the identity morphism belongs to $\cN$). In fact any $f\in\cN\subset \Hom_\b2(L,M)$ also factors through $L^\sigma$ using the projection $p:L\to L^\sigma$, $p(x):=x+\sigma(x)$.  Next we prove that every morphism has a kernel and cokernel with respect to $\cN$. In general the kernel $\sker(f):\ker(f)\to L$ is characterized by 
$$
f\circ \sker(f)\in \cN, \  \  f\circ g\in \cN\iff \exists!\, h~{\rm s.t.}~ g=\sker(f)\circ h.
$$
  For $\b2$ the condition $f\circ g\in \cN$ means that the range of $g$ is contained in $f^{-1}(M^\sigma)$ and this gives the required unique factorization $g=\sker(f)\circ h$, where $\sker(f):f^{-1}(M^\sigma)\to L$ is the inclusion. Thus one gets the agreement with  Definition \ref{ker&coker}. \newline
  The cokernel $\scoker(f):M\to \coker(f)$ is characterized in turn by 
$$
\scoker(f)\circ f\in \cN, \  \  g\circ f\in \cN\iff \exists!\, h~{\rm s.t.}~ g=h\circ \scoker(f).
$$
 We have defined in \eqref{cokernelb2} the  cokernel as the quotient $\coker(f):=M/\!\!\sim$ where for $b,b'\in M$
 \begin{equation}\label{cokernelb2bis1}
 b\sim b'\iff g(b)=g(b')\ \  \forall X, \forall g\in \Hom_\b2(M,X)~{\rm s.t.}~ g\circ f\in \cN.
\end{equation}
 Let $\scoker(f):M\to \coker(f)$ be the quotient map. One has $\scoker(f)\circ f\in \cN$, since for $b=f(a)$ one has $b\sim \sigma(b)$ as $g(b)=g(\sigma(b))$ for any $g\mid g\circ f\in \cN.$ Indeed 
 $$
 g(b)=(g\circ f)(a)=\sigma((g\circ f)(a))=g(\sigma(f(a))=g(\sigma(b)).
 $$
  Moreover, by construction, any $g$ s.t. $g\circ f\in \cN$ factors uniquely as $g=h\circ \scoker(f)$.
  We have thus shown that every morphism has a kernel and cokernel with respect to $\cN$ and that they agree with Definition \ref{ker&coker}.\endproof
  \begin{defn}\label{exactnessb2bis}  $(i)$~Let $f\in \Hom_\b2(L,M)$, then the normal image $\imm(f)\subset M$ is the kernel of the  cokernel $\scoker(f)$.\newline
  $(ii)$~We say that a sequence  of $\b2$: $L\stackrel{f}{\to} M\stackrel{g}{\to}N$ is {\em  exact} at $M$ if $\imm(f)=\ker(g)$.
 \end{defn}
  Notice that this definition of the normal image as the kernel of the cokernel corresponds to the definition of the image in an abelian category. By definition the cokernel of $f$ is the quotient of $M$ by the relation \eqref{cokernelb2}. Thus the kernel of the cokernel is given by the elements of $M$ whose image in the cokernel is fixed under $\sigma$. Thus one has 
  \begin{equation}\label{imagedefn}
 	b\in \imm(f)\iff g(b)=\sigma(g(b))\qqq X, \ g\in \Hom_\b2(M,X)~{\rm s.t.}~ {\rm Range}(f)\subset \ker(g).
 \end{equation}
  \begin{prop}\label{tentative1} Let $f\in\Hom_\b2(L,M)$.\newline
  $(i)$~The sequence $0\to L\stackrel{f}{\to}M$ is exact at $L$ if and only if it is \strgly exact at $L$.\newline
$(ii)$~${\rm Range}(f)+M^\sigma\subset \imm(f)$.\newline
 $(iii)$~$\xi+\sigma(\xi)\in {\rm Range}(f)+\sigma(\xi)$, $\forall \xi\in \imm(f)$.\newline
 $(iv)$~\Strg exactness (in the sense of Definition \ref{exactnessb2}) implies exactness.\newline
 $(v)$~Proposition \ref{use of T} continues to hold if one replaces \strg exactness by exactness.\newline
  $(vi)$~The following sequence is exact in $\b2$ 
 \begin{equation}\label{5exactnice}
 	0 \to \Ker(f)\to L\stackrel{f}{\to}M\to \coker(f)\to 0.
 \end{equation}
\end{prop}
 \proof $(i)$~\Strg exactness of $0 \to L\stackrel{f}{\to}M$ means that  $L^\sigma=\ker(f)$,
 while exactness of $0 \to L\stackrel{f}{\to}M$ means that $\imm(0)=\ker(f)$. Thus it is enough to show that for any object $F$ of $\b2$ one has  $\imm(0)=F^\sigma$. One has $ F^\sigma\subset \imm(0)$. Moreover the kernel of the identity map $\id:F\to F$ is $F^\sigma $ and thus $\imm(0)\subset F^\sigma$. \newline
 $(ii)$~This follows from  \eqref{imagedefn}.\newline
 $(iii)$~Let $E={\rm Range}(f)$ and $\zeta=\xi+\sigma(\xi)$. Assume $\zeta=\xi+\sigma(\xi)\notin \xi+E $. We show that $\xi\notin \imm(f)$. The result will then follow since $\imm(f)$ is $\sigma$-invariant. Let
$
P=\{\alpha \in M\mid \exists \eta\in \xi+E,\ \alpha\leq \eta\leq \zeta\}.
$
  By construction $P$ is hereditary and is a sub-$\B$-module since for $\alpha\leq \eta\leq \zeta$, $\alpha'\leq \eta'\leq \zeta$ one has 
  $
  \alpha+\alpha'\leq \eta+\eta'\leq \zeta, \  \  \eta+\eta'\in \xi+E
  $.
  Let $\omega\in \Hom_\B(M,\B)$ be such that $\omega(\alpha)=0\Leftrightarrow \alpha \in P$. One has $P\cap E=[0,\zeta]\cap E$, since for $\alpha\in E$, $\alpha\leq \zeta$,   $\eta=\xi +\alpha$ one has $\eta\in \xi+E$, $\alpha\leq \eta\leq \zeta$. Thus  one gets $\omega(\sigma(\alpha))=\omega(\alpha)$ for any $\alpha \in E$. Since $\xi \in P$  (with $\eta=\xi$) one has $\omega(\xi)=0$. One has by hypothesis that $\zeta=\xi+\sigma(\xi)\notin \xi+E $. Thus $\omega(\zeta)=1$ and it follows that $\omega(\sigma(\xi))=1$.   
  Let then $h\in \Hom_\b2(M,\yb\B)$ be given by 
  $
  h(\alpha)=(\omega(\alpha),\omega(\sigma(\alpha))~ \forall\alpha \in M
  $.
  Since $\omega(\sigma(\alpha))=\omega(\alpha)$ for any $\alpha \in E$, one has $h(E)\subset (\yb\B)^\sigma$ while $h(\xi)=(0,1)\notin (\yb\B)^\sigma$.\newline
 $(iv)$~\Strg exactness means ${\rm Range}(f)+M^\sigma=\ker(g)$ and by $(ii)$ ${\rm Range}(f)+M^\sigma\subset \imm(f)$ with $\imm(f)\subset \ker(g)$, thus it follows that $\imm(f)=\ker(g)$.\newline
 $(v)$~By $(i)$  statement 1. of Proposition \ref{use of T} continues to hold if one replaces \strg exactness by exactness. To see that this also holds for statement 2. it is enough to show that exactness of $TL\stackrel{Tf}{\to}TM\to 0$ implies that $f$ is surjective. But if $\imm(Tf)=TM$ one has $\xi+\sigma(\xi)\in {\rm Range}(Tf)+\sigma(\xi)$ for any $\xi\in TM$, using $(iii)$, and taking $\xi=(x,0)$ for $x\in M$  gives $x\in {\rm Range}(f)$ and the required surjectivity of $f$. \newline
   $(vi)$~We first show exactness at $\ker(f)$. The normal image $\imm(0)=\ker(f)^\sigma$. Since the inclusion $\iota:\ker(f)\to L$ is injective and compatible with $\sigma$, one has 
 $
 z\in \ker(\iota)=\iota^{-1}(L^\sigma)\Leftrightarrow z\in \ker(f)^\sigma
 $.
 Since $\imm(0)=\ker(f)^\sigma$ this gives exactness at $\ker(f)$. We now show exactness at $L$. By $(ii)$ ${\rm Range}(\iota)\subset \imm(\iota)$. By construction ${\rm Range}(\iota)=\ker(f)$. The other inclusion $\imm(\iota)\subset \ker(f)$ follows from ${\rm Range}(f\circ \iota)\subset M^\sigma$. Exactness at $M$ follows from the definition $\imm(f)=\ker(\coker(f))$. Finally concerning exactness at $\coker(f)$, note that the kernel of the last map (to $0$) is
 $\coker(f)$ while the quotient map $\eta:M\to\coker(f)$ is surjective by construction and hence $\coker(f)\subset {\rm Range}(\eta)\subset \imm(\eta)$  (using $(ii)$) so that 
 the equality $\coker(f)=\imm(\eta)$ is proved. 
\endproof 
  In \cite{gra} a sequence $L\stackrel{f}{\to}M\stackrel{g}{\to}N$ of morphisms in a semiexact category is called of order two if $g\circ f\in \cN$. The sequence is said to be exact if $\Ker(\scoker(f))=\Ker(g)$. This definition agrees with Definition \ref{exactnessb2bis}. The notion of normal image introduced in \cite{gra} is the same as $\imm(f)=\Ker(\scoker(f))$. The notion of normal coimage again as in \cite{gra} is $\coimm(g)=\Coker(\sker(g))$ and the condition $\coker(f)=\coimm(g)$ is equivalent to exactness because one has in general the equalities $\scoker(\sker(\scoker(f)))=\scoker(f)$ and $\sker(\scoker(\sker(f)))=\sker(f)$. The notion of an exact morphism is defined in \cite{gra}  based on a diagram such as
  \begin{equation}\label{canonical factorization}
\xymatrix@C=5pt@R=25pt{
\ker(f) \ar[rr]^{\ \sker(f)}  &&
L \ar[d]_{\scoker(\sker(f))}\ar[rr]^{f} && M\ar[rr]^{\!\!\!\!\scoker(f)}&&\, \, \coker(f)
\\
&&\coimm(f) \ar[rr]^{\ \ \ \tilde f}&& \imm(f)\ar[u]_{\sker(\scoker(f))}}
\end{equation}
  where for $f$ to be exact one requires that the map $\tilde f$ is an isomorphism. 
  
  When $M=0$, one has $\sker(f)=\id_L$ and the cokernel $\scoker(\sker(f))$ is the quotient of $L$ by the equivalence relation \eqref{cokernelb2bis1} for $f=\id_L$. One checks that this equivalence relation is given by $b\sim b'\iff p(b)=p(b')$ where $p(x)=x+\sigma(x)$ gives the projection  $p: L\to L^\sigma$. Thus one has 
   \begin{equation}\label{coimagezero}
 \coimm(L\stackrel{0}{\to}M)=L^\sigma, \  \ \scoker(\sker(0))=L\stackrel{p}{\to}L^\sigma
\end{equation}
and this shows that the zero map $L\stackrel{0}{\to}M$ is never exact unless $L=0$. \vspace{.05in}

The following example shows that $\b2$ is not generalized exact in the sense of \cite{gra} \S 1.3.6 (\ie a semi-exact category in which every morphism is exact). 

 \begin{example}\label{excok}  Let $M=\B^2$. It has 4 elements: $M=\{0,\ell,m,\ell\vee m=n\}$. Let $\iota:N=\{0,m,n\}\to M=\{0,\ell,m,n\}$ be the inclusion.  We consider  the morphism $f=\yb\iota:\yb N=N^2\to \yb M=M^2$.  The normal image  $\imm f=\ker(\scoker(f))$ is given by the elements of $\yb M=M^2$ which are fixed points of the involution in   the quotient $\coker(f)$ of $\yb M$ by the equivalence relation 
 \begin{equation}\label{equivcoke}
(x,y)\sim (x',y')\iff h(x)+g(y)=h(x')+g(y')  \end{equation}
 $$
 \forall X,~\forall (h,g)\in\Hom_\B(M,X)~{\rm s.t.}~ h\vert_N=g\vert_N.
 $$
 One has $(\ell,m)\in M^2\setminus (N^2+\Delta)$ and  $(\ell,m)\in \imm f$, \ie  $(\ell,m)\sim (m,\ell)$ that is 
 $$
h,g \in\Hom_\B(M,X),\ \ h\vert_N=g\vert_N \Rightarrow h(\ell)+g(m)= h(m)+g(\ell).
 $$	
 Indeed,  since $\ell\vee m=n\in N$ one has $h(\ell)+h(m)= g(m)+g(\ell)$, but since $h(m)=g(m)$ this gives $h(\ell)+g(m)= h(m)+g(\ell)$.
 The quotient $\coker(f)$ of $M^2$ by the equivalence relation \eqref{equivcoke} is the object of $\b2$ obtained by adjoining to $N$ (endowed with $\sigma=\id$) three elements $\alpha,\sigma(\alpha), \alpha+\sigma(\alpha)$ with addition given by $\alpha+z=n$, $\forall z\neq 0$, $z\in N$. The quotient map $\eta:M^2\to M^2/\!\!\sim$ is such that 
 $$
 \eta((u,v))=u+v\qqq u,v\in N, \ \ \eta((\ell,0))=\alpha.
 $$
  \end{example}
 
In the above example, since $f$ is injective, its kernel is a null object and the cokernel of its kernel is the identity map (using \eqref{cokernelb2} and the identity map $\id\in \Hom_\b2(\yb N,\yb N)$ to show that the equivalence relation is trivial).  Thus $\tilde f$ is the inclusion of $N^2$ in $\imm(T\iota)$ and is not an isomorphism.  The above example also shows that $f=\yb\iota$ is not a normal monomorphism in the sense of \cite{gra}. In fact kernels and cokernels establish an anti-isomorphism between the ordered sets of normal subobjects (kernels of morphisms in $\b2$) and normal quotients as one uses in the context of abelian categories. Normal subobjects and normal quotients form anti-isomorphic lattices. Thus normal subobjects give the relevant notion of subobject. It follows from \cite{gra1} 1.3.3. that in a semiexact category, every morphism has a normal factorization through its normal coimage and its normal image as in \eqref{canonical factorization}. In \cite{gra1} one denotes a normal monomorphism by the symbol
  $L\stackrel{f}{\rightarrowtail} M$ while for a normal epimorphism one uses
  $L\stackrel{f}{\xtwoheadrightarrow{}} M$. The next Lemma holds in any semiexact category and we include its proof for completeness. 
  \begin{lem}\label{exact examples} Let $L\stackrel{f}{\to}N$ be a morphism with $f=k\circ h$, where   $h:L \stackrel{}{\xtwoheadrightarrow{}}  M$ is a normal epimorphism and $k:M \stackrel{}{\rightarrowtail} N$ is a  normal monomorphism. Then $f$ is exact. \end{lem}
 \proof Since $k$ is injective  $k^{-1}(N^\sigma)=M^\sigma$ and the kernel of $f$ is equal to the kernel of $h$. By hypothesis $h$ is the cokernel of its kernel and thus $h=\scoker(\sker(f))$ and $M=\coimm(f)$. Since $h$ is surjective  $h(L^\sigma)=M^\sigma$ and the condition $\phi\circ f$ null is equivalent to $\phi\circ k$ null. Thus the cokernel of $f$ is the same as the cokernel of $k$ and the factorization \eqref{canonical factorization} takes the form
   \begin{equation}\label{canonical factorization bis}
\xymatrix@C=5pt@R=25pt{
\ker(h) =\ker(f)\ar[rr]^{}  &&
L \ar[d]_{h=\scoker(\sker(f))}\ar[rr]^{f} && N\ar[rr]^{}&&\, \, \coker(f)=\coker(k)
\\
&&M \ar[rr]^{\id_M}&& M\ar[u]_{k=\sker(\scoker(f))}}
\end{equation}
 which shows that $f$ is exact.\endproof

 Exact sequences $L\stackrel{f}{\to}M\stackrel{g}{\to}N$ are defined by $\Ker(\scoker(f))=\Ker(g)$. In \cite{gra1} 1.3.5. the following notion of short exact sequence\footnote{Such short exact sequences are the analogues of the kernel-cokernel pairs in additive categories: \cite{Buhler}} is introduced. 
 
 \begin{defn}\label{defn shortexact}
 	 A short \exxx sequence is given by a pair of maps $$L\stackrel{m}{\rightarrowtail}M\stackrel{p}{\xtwoheadrightarrow{}}N,~~ m= \sker(p)\ \& \  p= \scoker(m).$$
 \end{defn} 
  This means that 
  $L$ is a normal subobject of $M$ and $N$ is a normal quotient of $M$. As pointed out in \cite{gra1} \S 1.5.2, one has
  \begin{prop}\label{shortex} let $f:L\to M$ and $g:M\to N$ be morphisms in $\b2$. The following conditions are equivalent 
 \begin{enumerate}
  \item The sequence
  $L\stackrel{f}{\rightarrowtail}M\stackrel{g}{\xtwoheadrightarrow{}}N$  is  short \exxx\!.
 \item The sequence $0\to L\stackrel{f}{\to}M\stackrel{g}{\to}N\to 0$ is exact and the  morphisms $f$ and $g$ are exact.
\end{enumerate}  	
  \end{prop}

      Direct and inverse images of normal subobjects can be organized by a transfer functor  with values in the category of lattices and Galois connections. Exact functors are  introduced in \cite{gra1} 1.7. \vspace{.05in}
  
  Next, we test in our context the axiom $ex2$ of  \cite{gra1} 1.3.6. \ie the stability under composition of the normal monos and normal epis. 
  \begin{lem}\label{normalsubinj} Let $L \subset M$ be a normal subobject of the object $M$ of $\b2$. There exists an injective object $F$ of the category $\b2$ and a morphism $\phi:M\to F$ such that $L=\ker(\phi)$.
\end{lem}
 \proof By hypothesis there exists $f\in \Hom_\b2(M,X)$ such that $L=\ker(f)$.  Let $\eta_X$  and $r$ be as in Lemma \ref{T-injective}. One has $L=\ker(\eta_X\circ f)$ since $\eta_X^{-1}(\Delta)=X^\sigma$ by construction. Let $E$ be an injective object of $\bmod$ and $\iota:I(X)\to E$ an embedding as in Proposition \ref{inj3}. Then $g=\yb(\iota):T X\to \yb E=F$ is a morphism in $\b2$ such that $g^{-1}(F^\sigma)=T(X)^\sigma$. Thus one has $L=\ker(g\circ \eta_X\circ f)=\ker \phi$, with $\phi= g\circ \eta_C\circ f$. By Lemma \ref{injective b2},  $\yb E=F$ is an injective object of $\b2$ since the underlying $\B$-module is $E^2$.  \endproof 
  \begin{lem}\label{normalmonos} The normal monomorphisms in the category $\b2$ are stable under composition.
\end{lem}
 \proof Let $L \subset M$ be a normal subobject of an object $M$ of $\b2$ and similarly let 
 $M \subset N$ be a normal subobject of an object $N$ of $\b2$. By Lemma \ref{normalsubinj} there exists an injective object $X$ of the category $\b2$ and a morphism $\phi:M\to X$ such that $L=\ker(\phi)$. Since $X$ is injective we can extend $\phi$ to a morphism $\psi\in \Hom_\b2(N,X)$. Let $f\in \Hom_\b2(N,Y)$ such that $M=\ker(f)$. We let $\rho=(\psi,f)\in \Hom_\b2(N,X\times Y)$. The inverse image by $\rho$ of $(X\times Y)^\sigma=X^\sigma\times Y^\sigma$ is contained in $M=\ker(f)$ and coincides with the inverse image of $X^\sigma$ by the restriction of $\psi$ to $M$, \ie with $L=\ker(\phi)$. Thus $L=\ker(\rho)$. \endproof
 
 \begin{lem}\label{stablenormepis} In the category $\b2$ the normal epis are stable under composition.	
\end{lem}
\proof Let  $\alpha:L\to M$ and $\beta:M\to N$ be normal epis. One needs to show that 
\begin{equation}\label{toshow}
\beta\circ \alpha(a)=\beta\circ \alpha(b)\iff f(a)=f(b)\qqq f:L\to Z,\ \ker(f)\supset \ker(\beta\circ\alpha).	
\end{equation}
We first show, using the fact that $\alpha$ is a normal epi, that one has the equivalence $$\ker(f)\supset \ker(\beta\circ\alpha)\iff \exists \psi:M\to Z\mid f=\psi\circ\alpha, \ \ker(\psi)\supset \ker(\beta).$$
Indeed, the implication $\Leftarrow$ is immediate. Conversely, if $\ker(f)\supset \ker(\beta\circ\alpha)$
one has $\ker(f)\supset \ker(\alpha)$ and since $\alpha$ is normal one can factorize $f=\psi\circ\alpha$. Moreover, for any element $u$ of $\ker\beta$ one can find using the surjectivity of $\alpha$ an element $v\in L$ such that $\alpha(v)=u$. One has $v\in \ker(\beta\circ\alpha)$ since $\alpha(v)\in \ker(\beta)$ and thus since $\ker(f)\supset \ker(\beta\circ\alpha)$ one has $v\in \ker(f)$ and $f(v)\in Z^\sigma$ which shows that $\psi(u)\in Z^\sigma$ and $\ker(\psi)\supset \ker(\beta)$. Thus one gets the required equivalence. Next, since $\beta$ is a normal epi one has 
$$
\beta\circ \alpha(a)=\beta\circ \alpha(b)\iff \psi(\alpha(a))=\psi(\alpha(b))\qqq\psi:M\to Z\mid \ker(\psi)\supset \ker(\beta).
$$
 Thus one obtains \eqref{toshow} and that $\beta\circ \alpha$ is a normal epi. \endproof 
  The axiom $ex3$ of  \cite{gra1} 1.3.6, (sub-quotient axiom, or homology axiom) which if satisfied, defines a  homological category is stated as follows:\newline
  $ex3$: Given a normal mono $M\stackrel{m}{\rightarrowtail} N$ and a normal epi $N\stackrel{q}{\xtwoheadrightarrow{}} Q$ with $m\geq \ker(q)$,  the morphism $q\circ m$ is exact.\vspace{.05in}
  
  In \cite{gra1} the basic example of a homological category is the category $\Se_2$ of pairs of sets $(X,X_0)$ with $X_0\subset X$ and maps $f:X\to Y$ such that $f(X_0)\subset Y_0$. The  null maps are those with $f(X)\subset Y_0$. This category is shown to be semiexact and homological. A morphism of $\Se_2$ is exact iff $f$ is injective and $Y_0\subset f(X)$. This example suggests that for $\b2$ a necessary condition for $f:L\to M$ to be exact should be that $M^\sigma\subset f(L)$. Indeed
  \begin{lem}\label{exactneccond} In the category $\b2$ if $f:L\to M$ is exact then  $M^\sigma\subset f(L)$.
  \end{lem}
\proof In the factorization \eqref{canonical factorization}, the left vertical arrow is surjective and the right vertical arrow is a kernel so that $\imm(f)=\ker(\scoker(f))$ automatically contains $M^\sigma$. Thus if the horizontal arrow is an isomorphism one gets $M^\sigma\subset f(L)$.\endproof
Now, given $m$ and $q$ as in $ex3$ we have  $\ker(q)=q^{-1}(Q^\sigma)$ and since  $N\stackrel{q}{\xtwoheadrightarrow{}} Q$ is a normal epi, it is surjective and thus the map $q$ restricted to $\ker(q)$ surjects onto $Q^\sigma$. The condition $m\geq \ker(q)$ means that the normal subobject $M\stackrel{m}{\rightarrowtail} N$ contains $\ker(q)$ and it follows that the map $q\circ m$  surjects onto $Q^\sigma$, \ie it fulfills the necessary condition of  Lemma \ref{exactneccond}. In fact, since $M\stackrel{m}{\rightarrowtail} N$ is a normal subobject let $f:N\to Z$ be a morphism of $\b2$ such that $M=\ker(f)$. One has $\ker(q)\subset \ker(f)$ and since  $N\stackrel{q}{\xtwoheadrightarrow{}} Q$ is a normal epi it follows that $q=\scoker(\ker(q))$ which according to \eqref{cokernelb2} means  $$
 q(b)= q(b')\iff h(b)=h(b')\ \forall X, \forall h\in \Hom_\b2(N,X)~{\rm s.t.}~ \ker(q)\subset \Ker(h).$$
 This shows, taking $h=f$, that $f(b)$ only depends upon $q(b)$ and that there exists a morphism $\phi\in \Hom_\b2(Q,Z)$ such that $f=\phi\circ q$. Thus one has $\ker(f)=q^{-1}(\ker(\phi))$ and the normal subobject $S=\ker(\phi)\subset Q$ is the natural candidate to make the following diagram commutative
  \begin{equation}\label{guess for S}
\xymatrix@C=25pt@R=35pt{
 M \ \ar@{->>}[d]_{q'}\ar@{>->}[rr]^{m} && N \ar@{->>}[d]^{q}
\\
S \ \ar@{>->}[rr]^{k}&& Q}
\end{equation}
Since $S=\ker(\phi)$ its inclusion in $Q$ is a normal mono $S\stackrel{k}{\rightarrowtail} Q$. The map $q\circ m$ has range in $S$ and one needs to show that the induced map $q':M\to S$ is a normal epi. Since $q$ is a normal epi it is surjective and thus $q':q^{-1}(\ker(\phi))\to \ker(\phi)$ is also surjective. To show that $q'$ is a normal epi one needs to prove that it is the cokernel of its kernel. Its kernel is $\ker(q)\subset \ker(f)$. The equivalence relation defining the cokernel of $\ker(q)\subset \ker(f)$ is, for $b,b'\in M =\ker(f)$
$$
 b\sim_1 b'\iff g(b)=g(b')\ \forall X, \forall g\in \Hom_\b2(\ker(f),X)~{\rm s.t.}~ \ker(q)\subset \Ker(g).$$
 We need to show that this equivalence relation  is the same as that defined by $q(b)=q(b')$ which in turns is 
 $$
 b\sim_2 b'\iff h(b)=h(b')\ \forall X, \forall h\in \Hom_\b2(N,X)~{\rm s.t.}~ \ker(q)\subset \Ker(h).$$
Now, any such $h$ gives a $g$ by restriction to $\ker(f)$ and thus for for $b,b'\in M =\ker(f)$, one has $b\sim_1 b'\Rightarrow b\sim_2 b'$. Conversely we need to show that if $g(b)\neq g(b')$ for some $g\in \Hom_\b2(\ker(f),X)$ s.t. $\ker(q)\subset \Ker(g)$, then the same holds for some $h$.  Let $\eta_X$  and $r$ be as in Lemma \ref{T-injective}. One has $\Ker(g)=\ker(\eta_X\circ g)$ since $\eta_X^{-1}(\Delta)=X^\sigma$ by construction. Let $E$ be an injective object of $\bmod$ and $\iota:I(X)\to E$ an embedding as in Proposition \ref{inj3}. Then $u=\yb(\iota):T(X)\to \yb(E)=F$ is a morphism of $\b2$ such that $u^{-1}(F^\sigma)=T(X)^\sigma$. Thus replacing $g$ by $g'=u\circ \eta_X\circ g$ one still has $g'(b)\neq g'(b')$, $\ker(q)\subset \Ker(g')$ while now the object $\yb(E)=F$ is injective in the category $\b2$. Thus one can extend $g'\in \Hom_\b2(\ker(f),F)$ to $h\in \Hom_\b2(N,F)$ and one has $\ker(q)\subset \Ker(h)$ since $\ker(q)\subset \Ker(g')$. Moreover since $b,b'\in M =\ker(f)$ one has $h(b)\neq h(b')$ since $g'(b)\neq g'(b')$. Thus we have shown that the equivalence relation defining the cokernel of $\ker(q)\subset \ker(f)$ is the same as $q(b)=q(b')$. This shows that the restriction $q'$ of $q$ to $\ker(f)$ is a normal epimorphism and hence that the vertical arrow $q'$ of the diagram \eqref{guess for S} is a normal epimorphism. Thus $q\circ m$ factorizes as $k\circ q'$ and is exact by Lemma \ref{exact examples}. 
We can therefore state our main result
\begin{thm}\label{homological} The category $\b2$ is homological.	
\end{thm}
\proof By Proposition \ref{grandis} the category $\b2$ is semiexact. It follows from Lemma \ref{normalmonos} and Lemma \ref{stablenormepis} that $\b2$ satisfies $ex2$. We show that  it also satisfies $ex3$. Given a normal mono $M\stackrel{m}{\rightarrowtail} N$ and a normal epi $N\stackrel{q}{\xtwoheadrightarrow{}} Q$ with $m\geq \ker(q)$, the above construction determines a factorization \eqref{guess for S}, and one has $q\circ m=k\circ q'$ where $k$ is a normal mono and $q'$ a normal epi. Thus Lemma \ref{exact examples} applies and shows that $q\circ m=k\circ q'$ is exact. Thus the category $\b2$ is homological. \endproof

   Next, we investigate how morphisms in $\b2$ act on normal subobjects. 
   
\begin{figure}[H]
\begin{center}
\includegraphics[scale=0.4]{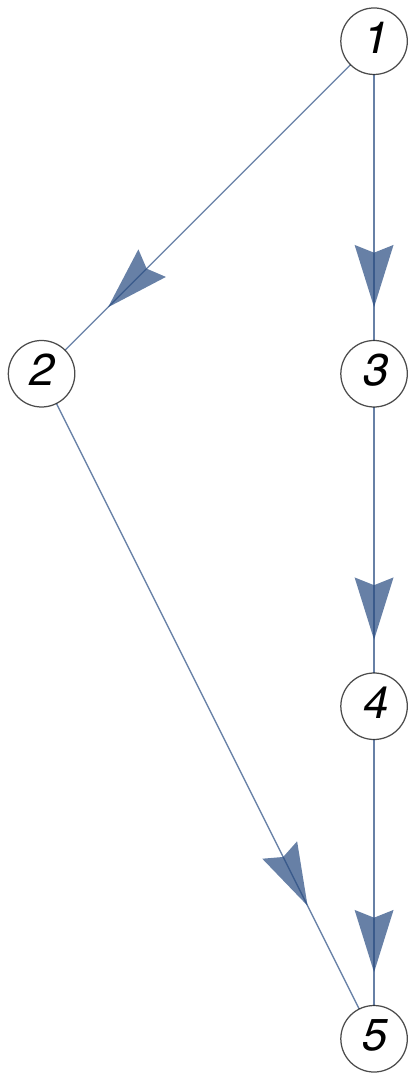}
\end{center}
\caption{The lattice of normal submodules of $\yb N=N^2$}\label{normsubmodn}
\end{figure}

   For any object $E$ of $\b2$ we let $\nsb(E)$ be the lattice of normal subobjects of $E$. For $N,M\in \nsb(E)$ one has $N\cap M\in \nsb(E)$ since the intersection of two kernels $\ker(f)\cap\ker(g)$ is the kernel of the map $(f,g)$ to the product of the codomains. Thus the lattice operation $\wedge$ on $\nsb(E)$ coincides with the intersection. The operation $\vee$ of the lattice is more delicate. By construction $E_1\vee E_2$ contains $E_1+E_2$ but next example shows the existence of two normal submodules $E_j\in \nsb(E)$ whose sum $E_1+E_2$ is not normal.
   \begin{example}\label{excok1}   
 Let $M=\{0,m,\ell,n\}$ be as in Example \ref{excok} and take $E=\yb M$. One easily checks  that the following submodules are normal submodules $E_j\in \nsb(E)$
$$
E_1=\{0,m\}\times \{0,m\}+\Delta, \  \ E_2=\{0,n\}\times \{0,n\}+\Delta.
$$
 Example \ref{excok} shows that the smallest element of $\nsb(E)$ which contains $E_1+E_2$ is $\imm(T\iota)=M^2\setminus \{(\ell,0),(0,\ell)\}$ and this normal submodule $E_1\vee E_2$ is strictly larger than $E_1+E_2$ since it contains $(\ell,m)\notin E_1+E_2$. 
\end{example} 
We recall that the modular condition for a lattice states that: $(E\vee F)\wedge G=E\vee(F\wedge G)$, for $E\subset G$. Figure \ref{normsubmodn} provides the graph of the lattice of normal submodules of $\yb N=N^2$ for $N=\{0,m,n\}$.. 
 The figure shows clearly that this lattice is not modular since taking $E=e(3)\subset G=e(4)$ and $F=e(2)$ one gets $(E\vee F)\wedge G=G$ while $E\vee(F\wedge G)=E$. \vspace{.05in}

 Following \cite{gra1} 1.5.6., a morphism in a semiexact category determines direct and inverse images of normal subobjects. In our setup this gives
 \begin{prop}\label{dirinvim} Let $f\in \Hom_\b2(E,F)$.\newline
 $(i)$~For $m\in \nsb(E)$, the direct image $f_*(m)$  is equal to $\imm(f\circ m)$.\newline
 $(ii)$~For $M\in \nsb(F)$, the inverse image $f^*(M)$  is equal to  $f^{-1}(M)$. 	
 \end{prop}
\proof $(i)$~This is by definition \cf~  (1.69) of \cite{gra1} 1.5.6.\newline
$(ii)$~It is enough, using equation (1.70) of \cite{gra1} 1.5.6., to show that for $M=\ker(\scoker(n))$ one has $\ker(\scoker(n)\circ f )=f^{-1}(M)$: this is immediate.\endproof 

\subsection{Kernels and normal subobjects in $\b2$}\label{sectkerb2}
In this subsection we show that for $E\subset F$ a subobject in $\b2$, the normal image $\imm(E)$ of the inclusion is determined in a local manner. More precisely, let $\xi\in F$ and assume that $\xi+\sigma(\xi)\in E$ (this is possible since adding $F^\sigma$ to $E$ does not change $\imm(E)$).  We look for a morphism $L\in \Hom_\b2(F,\yb\B)$ such that $L(E)\subset (\yb\B)^\sigma$ and $L(\xi)\notin (\yb\B)^\sigma$. The composition $\phi=p_j\circ L$ with the first or second projection is an element of $\Hom_\B(F,\B)$ such that 
 \begin{equation}\label{conditions on phi}
 \phi(\sigma(x))=\phi(x)\qqq x\in E, \ \  \phi(\xi)=0, \  \ \phi(\sigma(\xi))=1.
 \end{equation}
 Conversely, the existence of such an element of $\Hom_\B(F,\B)$ suffices to reconstruct $L=(\phi,\phi\circ\sigma)$. Note now that the restriction $\rho$ of $\phi$ to $E^\sigma$ uniquely determines $\phi$ on $E+\B\xi\subset F$ by means of the equality
 $$
 \phi(x+\epsilon \xi)=\rho(x+\sigma(x))\qqq x\in E,  \ \epsilon \in \B.
 $$
 The relative position of $\xi$ with respect to $E$ determines a triple $(\cR,\cS,\alpha)$ where 
 \begin{enumerate}
 \item $\cR\subset E\times E$ is a submodule (in $\bmod$) of $E\times E$ described by an equivalence relation.
 \item	$\cS\subset E$ is an additive subset.
 \item $\alpha=\sigma(\alpha)\in \cS$.
 \end{enumerate}
 Indeed, one lets $\cR=\{(x,y)\mid x+\xi=y+\xi\}$. This is a submodule (in $\bmod$) of $E\times E$. $\cS=\{z\mid z+\xi=z\}$  is an additive subset of $E$, and $\alpha=\xi+\sigma(\xi)$ belongs to $\cS$, since $\alpha+\xi=\alpha$. The existence of $\phi$ verifying \eqref{conditions on phi} is determined by the triple  $(\cR,\cS,\alpha)$ as follows, where $p:E\to E^\sigma$ is the projection $p(x)=x+\sigma(x)$.
  \begin{lem}\label{existphi} The existence of $\phi$ verifying \eqref{conditions on phi} is
  equivalent to the existence of $\rho\in \Hom_\B(E^\sigma,\B)$ such that 
  $$
  \rho(p(x))=\rho(p(y))\qqq (x,y)\in \cR, \ \ \rho(\alpha)=1.
  $$  
    \end{lem}
\proof The conditions are necessary as follows from 
$$
\rho(p(x))=\phi(x)=\phi(x+\xi)=\phi(y+\xi)=\phi(y)=\rho(p(y))
$$
and from
$$
s+\xi=s\Rightarrow p(s)+\alpha=p(s), \ \ \rho(\alpha)=1 \Rightarrow \rho(p(s))=1.
$$
Conversely, define $\phi$ on $E+\B\xi\subset F$ by the equality
 $$
 \phi(x+\epsilon \xi)=\rho(p(x))\qqq x\in E,  \ \epsilon \in \B.
 $$
If $x+\epsilon\xi=x'+\epsilon'\xi$ then $x+\xi=x'+\xi$ since one can add $\xi$ to both terms. Thus $(x,x')\in \cR$ and by hypothesis $\rho(p(x))=\rho(p(x'))$. This shows that $\phi$ is well defined on $E+\B\xi\subset F$ and $\phi\in \Hom_\B(E+\B\xi,\B)$. Moreover $\phi(\xi)=0$.  By Proposition \ref{inj1} $(iii)$ one can extend $\phi$ to an element $\tilde \phi\in \Hom_\B(F,\B)$. Since $\phi(\alpha)=1$ by hypothesis one has $$
\tilde\phi(\sigma(\xi))=\tilde\phi(\xi+\sigma(\xi))=\phi(\alpha)=1.
$$
 Thus $\tilde \phi$ fulfills \eqref{conditions on phi}.\endproof 
 Note that any pair $(a,b)$ in the projection $p(\cR)$ of the equivalence relation $\cR$ fulfills $\alpha+a=\alpha+b$ since $a=x+\sigma(x)$, $b=y+\sigma(y)$ with $x+\xi=y+\xi$. 
 Thus $p(\cR)\subset \cR_\alpha$. To $p(\cR)$ corresponds a canonical element $\psi\in \Hom_\B(E^\sigma,\B)$, namely the largest element compatible with the relation. This is obtained as $$\psi=\vee\{\rho\in \Hom_\B(E^\sigma,\B)\mid \rho(p(x))=\rho(p(y))\qqq (x,y)\in \cR\}.$$ Thus the existence of $\rho$ as in Lemma \ref{existphi} is equivalent to the condition $\psi(\alpha)=1$.

\begin{prop}\label{existphi1} Let $E\subset F$ be a subobject in $\b2$ containing $F^\sigma$. Then for $\xi \in F$ one has  $\xi \in \imm(E)$ (the normal image of the inclusion) if and only if there exists a finite sequence $a_0,a'_0,a_1,a'_1,\ldots, a_n,a'_n$ of elements of $E$ such that 
 $$
 \xi=a_0+\xi,\ p(a_0)=p(a'_0),\ a'_0+\xi=a_1+\xi,\ p(a_1)=p(a'_1), \ldots, p(a_n)=\xi+\sigma(\xi).
 $$
  \end{prop}
\proof Let us assume that there exists a finite sequence $a_0,a'_0,a_1,a'_1,\ldots, a_n,a'_n$ of elements of $E$ as in the statement. Let $f\in \Hom_\b2(F,C)$ be such that $E\subset \ker(f)$ and let us show that $\xi\in \ker(f)$. One has, using the equality $f(a)=f(p(a))$ for all $a\in E$ which implies $f(a_j)=f(a'_j)$ for all $j$
$$
f(\xi)=f(a_0+\xi)=f(a'_0)+f(\xi)=f(a'_0+\xi)=f(a_1+\xi)=\ldots =f(a_n+\xi)=f(p(a_n)+\xi)
$$
 and since $p(a_n)+\xi=\xi+\sigma(\xi)\in F^\sigma$ one has $f(\xi)\in C^\sigma$, \ie $\xi\in \ker(f)$. Thus the existence of the sequence ensures that $\xi \in \imm(E)$. We now assume that no such sequence exists and we define a submodule $Z\subset E^\sigma$ as the equivalence class $\cU(0)$ of $0$ for the equivalence relation $\cU$ on $E^\sigma$ generated by the relation $p(\cR)$ where $\cR=\{(x,y)\mid x+\xi=y+\xi\}\subset E\times E$ was defined above. The non-existence of the sequence means that $\alpha\notin \cU(0)$ where $\alpha=\xi+\sigma(\xi)$. Let then $\rho\in \Hom_\B(E^\sigma,\B)$ be defined by 
 $
 u\in \rho^{-1}(0)\Leftrightarrow \exists v\in \cU(0)\mid u\leq v
 $. 
 This condition defines an hereditary submodule, thus $\rho$ is well defined. Since $\cU(0)$ is saturated for the relation $p(\cR)$, $\rho$ fulfills the first condition of Lemma \ref{existphi} and it remains to show that $\rho(\alpha)=1$. Note that the subset  $V\subset E^\sigma$ determined as $V=\{v\in E^\sigma\mid v+\alpha=\alpha\}$ contains $0$ and is such that if $(v,v')\in p(\cR)$ then $v\in V\iff v'\in V$. Indeed for $v=p(x)$, $v'=p(x')$ and $x+\xi=x'+\xi$ one has $p(x)+\alpha=p(x')+\alpha$. It follows that $\cU(0)\subset V$. Now if $\alpha\in \rho^{-1}(0)$ one has $\exists v\in \cU(0)\mid \alpha\leq v$ and since $v\in V$ one has $v\leq \alpha$ and thus $v=\alpha$ which is a contradiction since $\alpha\notin \cU(0)$. Thus one has $\rho(\alpha)=1$ and 
 by Lemma \ref{existphi} there exists $L\in \Hom_\b2(F,\yb\B)$ such that $E\subset \ker(L)$ but $\xi\notin\ker(L)$. \endproof 
Proposition \ref{existphi1}  provides a natural filtration of $\imm(E)$ as: $\imm(E)=\cup_n \imm^{(n)}(E)$. Thus, we set 
 $$
 \imm^{(1)}(E):=\{\xi\in F\mid \exists a\in E, \ a+\xi=\xi,\ p(a)=\xi+\sigma(\xi)\}.
 $$
 One has $E\subset \imm^{(1)}(E)$ since for $\xi\in E$ one can take $a=\xi$. Also $\imm^{(1)}(E)$ is a subobject since it is invariant under $\sigma$ and with $a,a'$ associated to $\xi,\xi'$ the sum $a+a'$ is associated to $\xi+\xi'$. The next level is $\imm^{(2)}(E)$ which is defined by
 $$
 \{\xi\in F\mid \exists a_0,a'_0,a_1\in E, \  \xi=a_0+\xi,\ p(a_0)=p(a'_0),\ a'_0+\xi=a_1+\xi,\ p(a_1)=\xi+\sigma(\xi)\}.
 $$
Taking $a_0=a'_0=0$ one sees that $\imm^{(1)}(E)\subset\imm^{(2)}(E)$. \vspace{.05in}

Next, we show with an example that this inclusion is strict in general. \vspace{.05in}

One lets $F=\yb N$, $N=\{0,m,n\}$, $m<n$, and considers the inclusion $E\subset F$ with
 $$
E=\{ (0, 0), (0, m), (m, 0), (m, m), (m, n), (n, m), (n, n)\}\subset F.
$$
For $\xi=(0,n)$ the relation $p(\cR)$ on $F^\sigma=N$ is in fact given by the same symmetric set 
$$
p(\cR)=\{ (0, 0), (0, m), (m, 0), (m, m), (m, n), (n, m), (n, n)\}\subset N\times N
$$
and one sees that it is not an equivalence relation since it contains $(0, m)$ and $(m, n)$ but not $(0, n)$. One has $\xi+\sigma(\xi)=n$ when viewed as an element of $F^\sigma=N$ and the composition $p(\cR)\circ p(\cR)$ is required to get $n$ in the equivalence class of $0$. This shows that $\xi \in \imm^{(2)}(E)$ but $\xi \notin \imm^{(1)}(E)$.
\begin{prop}\label{extreme and range}$(i)$~Let $\phi\in \Hom_\b2(L,M)$ and $\xi\in M$, $\xi \neq \sigma(\xi)$ be an indecomposable element, \ie such that the subset $\{0,\xi\}\subset M$ is a hereditary sub $\B$-module. Then $\xi\in \phi(L)\Leftrightarrow \xi\in \imm(\phi)$. \newline
 $(ii)$~Assume that $M$ is generated by its indecomposable elements, then one has
 \begin{equation}\label{nice image}
 	\imm(\phi)=M\iff \phi(L)+M^\sigma=M.
 \end{equation}	
 $(iii)$~Let $\phi\in \Hom_\b2(L,M)$ and $\xi\in M$,  be  minimal among non-null elements, \ie such that $\xi$ is the only non-null element in  $[0,\xi]\subset M$. Then $\xi\in \phi(L)\iff \xi\in \imm(\phi)$. 
 \end{prop}
\proof $(i)$~This follows from Proposition \ref{existphi1} but we provide here a direct proof. One has to show that if $\xi\notin \phi(L)$  there exists a morphism $\psi\in \Hom_\b2(M,X)$ such that $\phi(L)\subset \ker(\psi)$ while $\xi\notin \ker(\psi)$. Let $f\in \Hom_\B(M,\B)$ be defined by  $f(b)=0\iff b\in \{0,\xi\}$. One has $f\in \Hom_\B(M,\B)$ since the subset $\{0,\xi\}\subset M$ is an hereditary submodule by hypothesis. Take $X=\yb\B$ and define $\psi$ as 
$$
\psi(b)=(f(b),f(\sigma(b)))\qqq b\in M.
$$
Since $\xi \neq \sigma(\xi)$ one has $f(\sigma(\xi))=1$ and thus $\xi\notin \ker(\psi)$. Since $\xi \notin \phi(L)$ one has $f(u)=1$, $\forall u\in \phi(L)$, $u\neq 0$ and since $\phi(L)$ is globally invariant under $\sigma$ one gets that $\phi(L)\subset \ker(\psi)$.\newline
 $(ii)$~We show the implication $\Rightarrow$. By $(i)$ if $\xi\in M$ is indecomposable, one has either $\xi\in M^\sigma$ or $\xi\in \phi(L)$.  Since $M$ is generated by its indecomposable elements one gets $\phi(L)+M^\sigma=M$.\newline
 $(iii)$~Let $f\in \Hom_\B(M,\B)$ be defined by  $f(b)=0\iff b\in [0,\xi]$, one has $f\in \Hom_\B(M,\B)$. Define $\psi$ as above  by 
$
\psi(b)=(f(b),f(\sigma(b)))\qqq b\in M
$. Then as above $\xi\notin \ker(\psi)$. Moreover  if $\xi\notin \phi(L)$ one has for $a\in L$ that $\phi(a)\neq \xi$ and $\phi(\sigma(a))\neq \xi$, thus  either $f(\phi(a))=1=f(\phi(\sigma(a)))$, or $\phi(a)<\xi$ (or $\phi(\sigma(a))<\xi$) in which case $\phi(a)$ is null and hence $f(\phi(a))=0=f(\phi(\sigma(a)))$. This shows that $\phi(L)\subset \ker(\psi)$, and hence $\xi\notin \imm(\phi)$. 
  \endproof 
  
  \subsection{Cokernels  in $\b2$}\label{sectcokerb2}

Proposition \ref{existphi1} gives a good control on the normal image and we now discuss the cokernel. The first guess for the cokernel of a normal submodule $E\subset F$ is the disjoint union of $F^\sigma$ with the complement of $E$ in $F$, while the map $\tilde p:F\to (F^\sigma\cup E^c)$ is defined as 
$$
\tilde p(x)=p(x)=x+\sigma(x)\qqq x\in E, \  \ \tilde p(x)=x\qqq x\in E^c.
$$
 In \S\ref{sectcokerdiag} we shall provide several examples where 
  this rule defines the cokernel. For this to hold one needs to show that the operation in $F^\sigma\cup E^c$ given by
  \begin{equation}\label{operation single0}
(\xi,\eta)\mapsto \{\tilde p(u+v)\mid \tilde p (u)=\xi,\ \tilde p(v)=\eta\}
\end{equation}
is single valued. This is clear when $\xi,\eta\in E^c$, since then the fibers have one element. It also holds for $\xi,\eta\in F^\sigma$ since then one has $u,v\in E$ and by linearity of $p$ one has $u+v\in E$, $p(u+v)=p(u)+p(v)=\xi+\eta$. One can thus assume that $\xi \in E^c$ and $\eta\in E$.  One then has to show that 
\begin{equation}\label{single0} 
	\#\{\tilde p(u+\xi)\mid u\in E, \ p(u)=\eta\}=1.
\end{equation}
  When all of the $u+\xi$ which appear are in $E$, the uniqueness follows from the linearity of $p$. Thus the interesting case to study is when $\xi \in E^c$, $u+\xi \in E^c$. \vspace{.05in}
  
  The following example provides a case where uniqueness fails.
  
  \begin{example}  \label{example coker}	
  Take $F=\yb M$ with $M=\{0,m,\ell,n\}$ and consider the submodule
  $$
  E_1=\{(0, 0), (0, m), (m, 0), (m, m), (\ell, \ell), (\ell, n), (n, \ell), (n, n)\}.
  $$
  Take $\xi=(0,\ell)$ and $\eta=(m,m)$. The fiber $\{u\in E_1, \ p(u)=\eta\}$ consists of the three elements $(0, m), (m, 0), (m, m)$. The elements $u+\xi$ with $u$ in the fiber are all in $E^c_1$ and are the three distinct elements $(0,n),(m,n),(m,\ell)$. It follows that these three elements have the same image in the cokernel.
   \end{example}
      
      \begin{figure}[H]
\begin{center}
\includegraphics[scale=1]{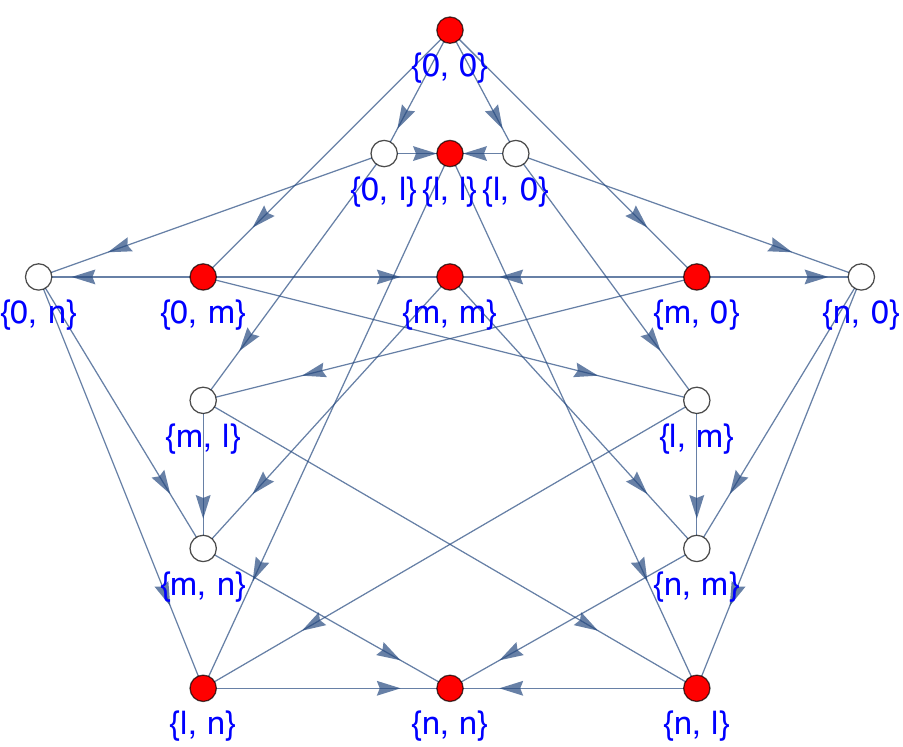}
\end{center}
\caption{The submodule $E_1\subset F=\yb M$ \label{cokerconstruct} }
\end{figure}

    The next proposition shows that in general, to obtain the cokernel, it suffices to divide $E^c$ by the equivalence relation generated by \eqref{simple add kernel} which is the analog in characteristic $1$ of the operation of quotient by the submodule.
      \begin{prop}\label{compute cokernel} Let $E\subset F$ be a subobject in $\b2$ containing $F^\sigma$. Then the cokernel of the inclusion $E\subset F$ is the quotient of $F^\sigma\cup E^c$ by the smallest equivalence relation such that 
   \begin{equation}\label{simple add kernel}
   \xi \in E^c,\ u,v\in E,\ p(u)=p(v)\Rightarrow \xi +u\sim \xi +v.
  \end{equation} 
   The cokernel map is the quotient map on $E^c$ and the projection $p$ on $E$. 
   \end{prop}
  \proof  By construction the cokernel of the inclusion $E\subset F$ is the quotient of $F$ by the equivalence relation
  $$
 \alpha \sim_\scoker \beta\iff f(\alpha)=f(\beta)\qqq f~{\rm s.t.}~  E\subset \ker(f).
  $$ 
  Let $f\in \Hom_\b2(F,X)$ be such that $E\subset \ker(f)$. Then for $u\in E$ one has $f(u)\in X^\sigma$ and hence $f(u)=f(p(u))$. This shows that 
  $
  u\sim_\scoker p(u), ~\forall u\in E
  $ and hence that the cokernel is the quotient of $F^\sigma\cup E^c$ by the restriction to $F^\sigma\cup E^c$ of the equivalence relation $\sim_\scoker$. This equivalence relation fulfills \eqref{simple add kernel} and we need to show that it coincides with the equivalence relation $\sim $ generated by \eqref{simple add kernel}. The latter  equivalence relation is given by $\xi\sim \eta\iff  \exists (\xi_j,u_j,v_j)$, 
  $\xi_j\in F$, $u_j,v_j\in E$
  $$
 \xi_1+u_1=\xi,\ \xi_n+v_n=\eta,
\ \ p(u_j)=p(v_j), \  \xi_j+v_j=\xi_{j+1}+u_{j+1}\ \forall j, 1\leq j \leq n.
  $$
  This equivalence relation is compatible with the addition, \ie if $\xi\sim \eta$ and $\xi'\sim \eta'$ one has $\xi+\xi'\sim \eta+\eta'$ as one gets by adding term by term the sequences $(\xi_j,u_j,v_j)$, $(\xi'_j,u'_j,v'_j)$. Moreover $\xi\sim \eta\Rightarrow p(\xi)=p(\eta)$. Now as above define a map $\tilde p:F\to Z:=(F^\sigma\cup E^c)/\!\!\sim$, by 
$$
\tilde p(x)=p(x)\qqq x\in E, \  \ \tilde p(x)=x/\sim\qqq x\in E^c.
$$
We have shown that $x\sim y\Rightarrow x \sim_\scoker y$ and to show the converse it is enough to prove that the addition in $F$ descends to $Z$, while $\tilde p:F\to Z$ is additive. This follows from the compatibility 
$$
\xi\sim \eta \ \& \ \xi'\sim \eta'\Rightarrow \xi+\xi'\sim \eta+\eta'.
$$
 The commutation with $\sigma$ is automatic. \endproof 
 
 In the statement of Proposition \ref{compute cokernel} we did not assume that $E$ is normal as a submodule of $F$. It is important to relate Proposition \ref{existphi1} with  Proposition \ref{compute cokernel}. We check that the existence of a finite sequence $a_0,a'_0,a_1,a'_1,\ldots, a_n,a'_n$ of elements of $E$ such that 
 $$
 \xi=a_0+\xi,\ p(a_0)=p(a'_0),\ a'_0+\xi=a_1+\xi,\ p(a_1)=p(a'_1), \ldots, p(a_n)=\xi+\sigma(\xi)
 $$
implies that $\tilde p(\xi)=\sigma(\tilde p(\xi))$ with the above notations. Indeed, to show that $\xi\sim p(\xi)$ we can take all $\xi_j=\xi$, $u_1=a_0$, $v_1=a'_0$, $u_j=a_{j-1}$, $v_j=a'_{j-1}$. The equalities $p(a_j)=p(a'_j)$ and  $a'_{j-1}+\xi=a_j+\xi$ mean $p(u_j)=p(v_j), \  \xi_j+v_j=\xi_{j+1}+u_{j+1}$ which gives the equivalence of $\xi$ with the last term $a_n+\xi\sim p(a_n)+\xi=\xi+\sigma(\xi)$. This shows that $\xi\in \imm$, as required. 

As a corollary  of Proposition \ref{compute cokernel} we obtain
\begin{prop}\label{fixed cokernel} Let $E\subset F$ be a subobject in $\b2$ containing $F^\sigma$. Let $\coker(\iota)$ be the cokernel of the inclusion $\iota:E\to F$.  Then there are canonical isomorphisms of $\B$-modules 
   \begin{equation}\label{cokernel fixed}
   F^\sigma\stackrel{\scoker(\iota)}{\simeq}(\coker(\iota))^\sigma \stackrel{\scoker(\id_F)}{\simeq}F^\sigma
  \end{equation} 
  whose composition is the identity on $F^\sigma$.
   \end{prop}
   \proof $\scoker(\iota):F\to \coker(\iota)$ restricts to a map of $\B$-modules $\gamma:F^\sigma\to (\coker(\iota))^\sigma$. This map is surjective because $\scoker(\iota):F\to \coker(\iota)$ is surjective by construction, and the projection $p$ on the $\sigma$-fixed points commutes with morphisms in $\b2$. To show that $\gamma$ is injective it is enough to check that the composition $\scoker(\id_F)\circ \gamma$ is the identity on $F^\sigma$. One has by Proposition \ref{compute cokernel},  $\coker(\id_F)=F^\sigma$ and $\scoker(\id_F)$ is the projection $p:F\to F^\sigma$, thus one gets the required result. It follows that the two maps of \eqref{cokernel fixed} give an isomorphism and its inverse. \endproof
   
     \section{Analogy with operators}\label{section analogy operators}
   The category $\b2$ is homological but not modular. As exlained in \cite{gra1} 2.3.5, the modularity property is used as a replacement of the standard argument in diagram chasing for additive categories, \eg $f(a)=f(b)\Rightarrow a-b\in \ker(f)$. In \S\ref{section duality b2} we develop the analogy of the category $\b2$ with the category of operators in Hilbert spaces  and in \S\ref{section kernel and injective} we introduce a substitute for the modularity property by showing in Theorem \ref{kernel and injective} that for a large class of objects $E$ of $\b2$ the nullity of the kernel of a morphism with domain $E$ is equivalent to the injectivity of the restriction of the morphism to each fiber of the projection on the null elements. Finally, in \S\ref{sectinjproj} we discuss injective and projective objects of $\b2$ and analyze the simplest example of a non-trivial short \exxx sequence of finite objects. 
    \subsection{Duality in $\b2$}\label{section duality b2}

In this subsection we define a duality in the category $\b2$. One first notes that there is an internal Hom functor $\Homi_\b2(L,M)$ obtained by involving the symmetry $\sigma(f):=\sigma\circ f=f\circ \sigma$ on the $\B$-module $\Hom_\b2(L,M)$. To dualize we use the object $\yb\B$ of $\b2$ and we define the dual $E^*:=\Homi_\b2(E,\yb\B)$. By Lemma \ref{forgetful} one has a canonical isomorphism $I(E^*)\simeq I(E)^*$ of the dual $E^*$ in the above sense with the dual of $I(E)$ in the sense of \S\ref{galois connect} \ie $I(E)^*=\Hom_\B(I(E),\B)$. It is given explicitly by the formula
\begin{equation}\label{iso of duals}
	E^*\ni \phi\mapsto \psi =p_1\circ \phi\in I(E)^*, \ \ p_1((x,y)):=x\qqq (x,y)\in \yb\B.
\end{equation}
Conversely, given $\psi\in I(E)^*$ the associated $\phi\in E^*$ is obtained by symmetrization \ie as $\phi(x)=(\psi(x),\psi(\sigma(x)))$.
The additional structure given by the involution corresponds to the involution $\psi \mapsto \psi\circ \sigma$, for  $\psi\in \Hom_\B(I(E),\B)$. In particular, the results of \S\ref{galois connect} apply and allow one to reconstruct $E$ from $E^*$ by biduality  as follows. One defines the normal dual as 
$$
E^*_{\rm norm}:=\{\phi\in \Hom_\b2(E,\yb\B)\mid \phi(\vee x_\alpha)=\vee(\phi(x_\alpha)\}.
$$
Note that $\id(I(E))$ involves hereditary submodules of $E$ which are not in general $\sigma$-invariant.  The map $J\mapsto \sigma(J)$ endows $\id(I(E))$ with an involution which coincides with the involution of the dual $E^*$ under the identifications $\id(I(E))\simeq I(E)^*\simeq E^*$. 
The isomorphism  $\tilde \epsilon: \id(I(E))\simeq (I(E)^*)^*_{\rm norm}$ of Proposition \ref{bidual} is compatible with the involutions and one obtains
 \begin{prop}\label{bidualbmods} The object $E$ of $\b2$ is the subobject of $(E^*)^*_{\rm norm}$ given by the compact elements of this latter complete algebraic lattice.
 \end{prop}
 \proof The evaluation map $\phi\mapsto \phi(x)\in \yb\B$ determines an embedding  $\rho:E\to (E^*)^*_{\rm norm}$. This embedding is compatible with the embedding $\epsilon:I(E)\subset (I(E)^*)^*_{\rm norm}$ of Proposition \ref{bidual} using the isomorphisms of \eqref{iso of duals}  
 $$
  I((E^*)^*_{\rm norm})\simeq (I(E^*))^*_{\rm norm}\simeq  (I(E)^*)^*_{\rm norm}
 $$
 and by applying the equality $(p_1\circ \phi)(x)=p_1(\phi(x))$.
 \endproof 
 
 One defines the notion of the orthogonal of a subobject  by using the natural pairing between $E$ and its dual $E^*:=\Homi_\b2(E,\yb\B)$,  \ie
$$
\langle x,y\rangle_\sigma:=y(x)\in \yb\B\qqq x\in E,y\in E^*, \ \ F^\perp:=\{y\mid \langle x,y\rangle_\sigma \ {\rm null}\qqq x\in F\}.
$$ 
We recall that the term ``null" means fixed under $\sigma$, so that $\langle x,y\rangle_\sigma \ {\rm null}\Leftrightarrow \langle x,y\rangle_\sigma \in (\yb\B)^\sigma$.
  \begin{lem}\label{normal condition} Let $F\subset E$ be a subobject of an object $E$ of $\b2$. Then the least normal subobject containing $F$ is $(F^\perp)^\perp$.
    \end{lem}
\proof Note first that $(F^\perp)^\perp$ is normal since it is a kernel by construction. Moreover it contains $F$. Let $F\subset \ker(f)$ with $f\in \Hom_\b2(E,X)$. Using the injection $\eta_X:X\to T(X)$ one can replace $X$ by $T(X)$ and using an embedding $I(X)\to \B^Y$ (Proposition \ref{inj3}) one can replace $X$ by $\yb(\B^Y)=(\yb\B)^Y$. Then one sees that $f(z)$ is null iff  $\forall x\in Y$ the component $f_x(z)\in \yb\B$ is null. Since $F\subset \ker(f)$ one has $F\subset \ker(f_x)$ and thus $f_x\in F^\perp$ $\forall x$. Hence, for $z \in (F^\perp)^\perp$ the component $f_x(z)\in \yb\B$ is null $\forall x\in Y$ and one gets  $(F^\perp)^\perp\subset \ker(f)$. \endproof 
Proposition \ref{bidualbmods} shows that an object of $\b2$ is uniquely determined by its dual. The following result (to be compared with Proposition \ref{compute cokernel}) thus gives an efficient way to determine the cokernel of any morphism $\phi\in \Hom_\b2(E,F)$.
\begin{prop}\label{bidualcoker} Let $\phi\in \Hom_\b2(E,F)$. The dual of the cokernel of $\phi$ is canonically isomorphic to the kernel of $\phi^*\in \Hom_\b2(F^*,E^*)$:~ $\coker(\phi)^*=\Ker(\phi^*)$.
 \end{prop}
 \proof By construction  $\coker(\phi)$ is the quotient of $F$ by the equivalence relation
 $$
 \alpha \sim_\scoker \beta\iff f(\alpha)=f(\beta)\qqq f~{\rm s.t.}~ \ \phi(E)\subset \ker(f).
  $$ 
  Since the map $\scoker: F\to \coker(\phi)$ is surjective, the map $\scoker^*:  \coker(\phi)^*\to F^*$ is injective and it remains to show that its range is the kernel of $\phi^*$. The map $\theta=\scoker\circ \phi$ is null, \ie $\sigma\circ\theta=\theta \circ \sigma=\theta$, thus the same holds for $\theta^*=\phi^*\circ \scoker^*$.  This shows that   $\scoker^*(\coker(\phi)^*)\subset \ker(\phi^*)$. Conversely, let $f\in \ker(\phi^*)\subset F^*$. One has $f\in \Hom_\b2(F,\yb\B)$ and since $f\in \ker(\phi^*)$ then $f\circ \phi$ is null and $\phi(E)\subset \ker(f)$. Thus $\alpha \sim_\scoker \beta\Rightarrow  f(\alpha)=f(\beta)$ and $f$ induces a morphism $g: \coker(\phi)\to  \yb\B$ such that $g\circ \scoker =f$. This means  $\scoker^*(g)=f$ and hence  $f\in \scoker^*(\coker(\phi)^*)$. Thus we get  
  $\scoker^*(\coker(\phi)^*)= \ker(\phi^*)$.\endproof 
 \begin{rem}\label{normal closure rem}{\rm $(i)$~Proposition \ref{bidualcoker} implies that the cokernel of the inclusion $j:E\to F$ of a subobject is the same as the cokernel of the identity $\id_F:F\to F$ iff the kernel of $j^*:F^*\to E^*$ is null. One has $\ker(j^*)=\{\phi\in F^*\mid \phi\circ j \ {\rm null}\}=E^\perp$. Thus $\ker(j^*)=(F^*)^\sigma \iff E^\perp=(F^*)^\sigma$ and the above statement follows directly from Lemma \ref{normal condition}.\newline
 $(ii)$~When  the  object $E$ of $\b2$ is finite, Proposition 	\ref{bidualbmods} simplifies and gives a canonical isomorphism of biduality $E\simeq (E^*)^*$.\newline
 $(iii)$~The subtlety arising from the existence of non-normal subobjects in the category $\b2$ is  analogous to the existence of non-closed subspaces of a Hilbert space in the category of Hilbert spaces and  linear operators. In the latter category the range of a morphism $T:\cH_1\to \cH_2$ is in general not closed and it is natural to define the cokernel as
 $$
 \coker(T):=\cH_2/\!\!\sim,\  \  \alpha\sim \beta\iff f(\alpha)=f(\beta)\qqq f:\cH_2\to \cH_3~{\rm s.t.}~ f\circ T=0.
 $$
 Since $f$ is continuous, $f=0$ on the {\em closure} of the range of $T$ and thus  $\coker(T)=\cH_2/(\overline{{\rm Range}(T)})$. The equality $\coker(T)^*=\Ker(T^*)$ holds for any morphism $T$. 
 \newline
 $(iv)$~ In the category of operators in Hilbert spaces one introduces (\!\!\cite{bbk}) the notion of {\em strict morphism} $f$. These are morphisms with closed range, and they are characterized by the existence of a quasi-inverse $g$ such that 
  $f=fgf$ and $g=gfg$. The morphism $g$ is obtained as the composition of the orthogonal projection on the closed range of $f$ with the inverse of $f_{\vert \ker(f)^\perp}$. Thus $fg$ is the orthogonal projection on the closed range of $f$  and $gf$ is the  orthogonal projection on the support of $f$, \ie  $\ker(f)^\perp$, hence  $f=fgf$ and $g=gfg$. However in general it is not true that the composition of strict morphisms is strict.  Indeed, a selfadjoint idempotent is strict, but in general the product of two selfadjoint idempotents is not strict since the operator valued angle can have arbitrary spectrum.
 }\end{rem}
 
   In the Hilbert space context any operator admits a canonical decomposition as the product of its restriction to its support, which is the orthogonal of its kernel, with the orthogonal projection on its support. We give an analogous statement for morphisms in $\b2$.
    \begin{lem}\label{ortho proj} Let $f\in \Hom_\b2(E,F)$, then $\ker(f)={\rm Range}(f^*)^\perp$.
      \end{lem}
   \proof  One has $\ker(f)=\{x\mid f(x)=\sigma(f(x))\}$, moreover since the duality between $F$ and $F^*$ is separating one has $f(x)=\sigma(f(x))\iff \phi(f(x))=\phi(\sigma(f(x)))$ $\forall \phi\in F^*$. Thus 
   $$
   x\in \ker(f)\iff \phi\circ f(x)\ {\rm null}\qqq \phi\in F^*\iff x\in {\rm Range}(f^*)^\perp.
   $$
      \endproof    
   Next, we apply the canonical decomposition of Proposition \ref{support} to obtain the
   following analogue in the category $\b2$ of  the above decomposition of operators in Hilbert space    
   \begin{prop}\label{support bmods} Let $f\in \Hom_\b2(E,F)$ be a morphism of finite objects in $\b2$. Then\newline
   $(i)$~The involution $\sigma$ on $E$ restricts to the support $S={\rm Support}(I(f))$ of $I(f)\in \Hom_\B(I(E),I(F))$ and the morphisms $q_S$ and $f\vert_S$ are morphisms in $\b2$.\newline
   $(ii)$~The canonical factorization $f=f\vert_S\circ q_S$ holds in $\b2$.\newline
   $(iii)$~The kernel of $f$ is the orthogonal of its support:
   $
   \ker(f)=\hat S^\perp
   $, 
   where $\hat S$ is a submodule of $E^*$ using the canonical isomorphism $E^{\rm op}\simeq E^*$.   	
  \end{prop}
   \proof $(i)$~The definition of the support as $S=\{z\in E\mid f(y)\leq f(z)\Rightarrow y\leq z\}$ shows that $z\in S\iff \sigma(z)\in S$. Similarly the definition of $q_S$ shows that it commutes with $\sigma$ and since by hypothesis $f$ commutes with $\sigma$ one gets $(i)$.\newline
   $(ii)$~Follows from Proposition \ref{support}.   
   \newline
   $(iii)$~Since  $f\vert_S$ is injective, one has $f\vert_S(u)=\sigma(f\vert_S(u))\iff u=\sigma(u)$ and thus the inverse image of null elements under $f\vert_S(u)$ is $S^\sigma$. Thus the kernel of $f=f\vert_S\circ q_S$ is the same as the kernel of $q_S$.
   By Lemma \ref{ortho proj} one has 
   $
   q_S(\xi)=q_S(\sigma(\xi))\Leftrightarrow \langle \xi,\phi_s\rangle_\sigma=\langle\sigma(\xi),\phi_s\rangle_\sigma~ \forall s\in S
   $. 
   The equality on the right hand side shows that $\xi\in \hat S^\perp
   $ 
   where $\hat S$ is viewed as a submodule of $E^*$ using the canonical isomorphism $E^{\rm op}\simeq E^*$. The equality on the left hand side is the definition of $\ker(q_S)$ thus one gets the required equality $ \ker(f)=\hat S^\perp$. \endproof
   
      \subsection{Nullity of kernel and injectivity}\label{section kernel and injective}

 As a general tool to prove injectivity of a morphism $f\in \Hom_\b2(E,F)$ one has
\begin{thm}\label{kernel and injective strong} Let $E$ be an object of $\b2$ and $x,y\in E$ such that  $x+\sigma(y)$ is not a null element (\ie is not fixed by $\sigma$). Then for any morphism $f\in \Hom_\b2(E,F)$ with null kernel, one has $f(x)\neq f(y)$.
\end{thm}
\proof If $f(x)=f(y)$ one has $f(x)+\sigma(f(y))$ null and hence $x+\sigma(y)\in \ker(f)$, thus we get
\begin{equation}\label{kernel significance}
	f(x)=f(y)\Rightarrow x+\sigma(y)\in \ker(f).
\end{equation}
When $\ker(f)$ is null we obtain the required statement. \endproof 
\begin{rem}\label{relation equiv} $(i)$~The relation ``$x \cR y\iff x+\sigma(y)\in E^\sigma$", between elements $x,y$ of an object $E$ of $\b2$ is symmetric and reflexive but in general not transitive. The simplest example of a non transitive relation is for $E=\yb\B$. Here one sees that
$
(1,0)\cR(1,1), (0,1)\cR (1,1)
$ while $(1,0)\cR(0,1)$ does not hold.\newline
$(ii)$~In fact, we show that for any object $E$ of $\b2$, and $x,y\in E$, if $p(x)=p(y)$ then $x$ and $y$ belong to the same class of the equivalence relation $\sim_\cR$ generated by the relation $\cR$. Indeed one has 
$$
x+\sigma(x)=y+\sigma(y)\Rightarrow x+\sigma(x+y)=(x+y)+\sigma(x).
$$
This shows that $x+\sigma(x+y)\in E^\sigma$, hence $x\cR(x+y)$.  Thus $x\sim_\cR x+y$. Similarly $y\sim_\cR x+y$ so that $x\sim_\cR y$ by transitivity of $\sim_\cR$.	
\end{rem}
   Next theorem shows that, by implementing a technical hypothesis on an object $E$ of $\b2$, the kernel of a morphism $f:E\to F$ in $\b2$ plays a role similar to the kernel of a linear map inasmuch as its nullity is equivalent to injectivity. Since on the null objects the nullity of the kernel is automatic, the statement of injectivity is relative to the projection $p:E\to E^\sigma$, $p(x)=x+\sigma(x)$. 
   \begin{defn}\label{sigma inject} A morphism $f\in \Hom_\b2(E,F)$ is said to be $\sigma$-{\em injective} if the restriction of $f$ to each fiber of $p$ is injective.	
\end{defn}
Notice that a morphism $f\in \Hom_\b2(E,F)$ is injective if and only if it is $\sigma$-injective and its restriction to the null elements is injective.
   \begin{thm}\label{kernel and injective} Let $E$ be a finite object of $\b2$ whose dual $E^*$ is generated by its minimal non-zero elements. Then for $f\in \Hom_\b2(E,F)$ the following statement holds \vspace{.05in}
   
   \centerline{$f$ is $\sigma$-injective $\iff$ $\ker(f)$ is null.}   	
   \end{thm}
   \proof First assume that the restriction of $f$ to the fibers of $p$ is injective, then since $x$ and $p(x)$ are on the same fiber of $p$, the equality $f(x)=\sigma(f(x))$ implies $f(x)=f(p(x))$ and hence $x=p(x)$ which shows that the kernel of $f$ is null. Conversely, assume that $\ker(f)$ is null.  By Lemma  \ref{support bmods} $(ii)$ one has $\ker(f)={\rm Range}(f^*)^\perp$ and hence by Lemma \ref{normal condition} the least normal subobject of $E^*$ containing ${\rm Range}(f^*)$ is $({\rm Range}(f^*)^\perp)^\perp=(\ker(f)^\perp=E^*$ since $\ker(f)$ is null. Note that we use here the finiteness of $E$ to apply Lemma \ref{normal condition} to $E^*$, using the identification of $E$ with the dual of $E^*$. Thus the normal image of the inclusion ${\rm Range}(f^*)\subset E^*$ is $E^*$ and  Proposition \ref{extreme and range} $(ii)$ shows that ${\rm Range}(f^*)+(E^*)^\sigma=E^*$. Let now $\xi, \eta\in E$ be such that $p(\xi)=p(\eta)$ and $f(\xi)=f(\eta)$. We prove that $\xi=\eta$. It is enough to show that one has $L(\xi)=L(\eta)$ for any $L\in E^*$ since the duality is separating (both in $\bmod$ and in $\b2$). Since ${\rm Range}(f^*)+(E^*)^\sigma=E^*$ it is enough to show that $L(\xi)=L(\eta)$ when $L\in {\rm Range}(f^*)$ and when $L=L\circ \sigma$. In the latter case the equality follows from $p(\xi)=p(\eta)$. For $L\in {\rm Range}(f^*)$ the equality follows from $f(\xi)=f(\eta)$. \endproof 
   \begin{rem}\label{trivial kernel} It is not true in general that for $f\in \Hom_\b2(E,F)$ and $\ker(f)$ null, one has $f^{-1}(\{0\})=\{0\}$. However the nullity of $\ker(f)$ implies  $f^{-1}(\{0\})\subset E^\sigma$, and that the fiber of $p$ above any  $a\in f^{-1}(\{0\})$ is reduced to $a$. Thus the restriction of $f$ to these fibers is injective.   	
   \end{rem}

   The condition that the  dual $E^*$ is generated by its minimal non-zero elements can be weakened still ensuring the conclusion of Theorem \ref{kernel and injective}
   \begin{cor}\label{trivial kernel1}  Let $E$ be a finite object of $\b2$ such  that for any fiber $F$ of $p$ the dual of the $\B$-module  	$\{0\}\cup F$ is generated by its minimal non-zero elements. Then for  $f\in \Hom_\b2(E,F)$ the following statement holds \vspace{.05in}
   
   \centerline{$f$ is $\sigma$-injective $\iff$ $\ker(f)$ is null.}   	
   \end{cor}
   \proof This follows from Theorem \ref{kernel and injective} applied to the restriction of $f$ to the object $\{0\}\cup F$ of $\b2$. The kernel of this restriction is null if $\ker(f)$ is null. \endproof 
   
    The hypothesis of Corollary \ref{trivial kernel1} holds for the  object $E=\coker(F(\alpha''))$ described in Figure \ref{graphcokernel} while  the  dual $E^*$ is not generated by its minimal non-zero elements (see Remark \ref{trivial kernel2}).\vspace{.05in}

Next, we show using Theorem \ref{kernel and injective}, that the condition that the dual $M^*$ of a $\B$-module $M$ is generated by its minimal non-zero elements suffices to derive the implication 

\centerline{$0 \stackbin[0]{0}{\rightrightarrows}M\stackbin[g]{f}{\rightrightarrows}N$  exact at $M$ $\Rightarrow$ $(f,g)$ is a monomorphism in $\bm2$.}
This statement was discussed at length in \S\ref{sectmeaningst}  where several counterexamples have been given none of which though fulfills the above condition. Moreover, the second statement of next Corollary \ref{exactrad} exhibits a large class of objects $X$ in $\b2$ which fulfill the analogue of the fundamental property holding for morphisms in an abelian category \ie
\begin{equation}\label{nullmono}
	f\in \Hom_\b2(X,Y) \ \text{monomorphism}\ \iff \ker(f)\ \text{null}.
\end{equation}
\begin{cor}\label{exactrad} Let $M$ be a finite object of $\bmod$ whose dual $M^*$ is generated by its minimal non-zero elements. \newline
$(i)$~The sequence  $0 \stackbin[0]{0}{\rightrightarrows}M\stackbin[g]{f}{\rightrightarrows}N$ is \strgly exact at $M$   if and only if $(f,g)$ is a monomorphism in $\bm2$.	\newline
$(ii)$~$\phi\in \Hom_\b2(\yb M,X)$ is a monomorphism if and only if its kernel is null.
\end{cor}
\proof   $(i)$~By Proposition \ref{monomorphism}, we need to show   that, assuming \strg exactness, the map $\iota:M^2\to N^2$, $\iota(x,y)=(f(x)+g(y), g(x)+f(y))$ of \eqref{137} is injective. Lemma \ref{injN} $(ii)$ shows that its restriction to the null elements (\ie to the diagonal) is injective. Moreover by hypothesis the kernel of $\iota$ (viewed as a morphism in $\b2$) is null and the dual of $\yb M=M^2$ is generated by its minimal elements. Thus Theorem \ref{kernel and injective} shows that $\iota$ is $\sigma$-injective and hence injective since its restriction to the null elements  is injective.  \newline
$(ii)$~The map $\phi$ is a monomorphism iff it is injective. If it is injective its kernel is null since $\phi\circ \sigma=\sigma\circ \phi$. Conversely, assume that the kernel of $\phi$ is null. Let $\eta_X:X\to TX$ be the unit, $\eta_X(b)=(b,\sigma(b))$. The kernel of $\eta_X\circ \phi$ is null.  By Lemma \ref{yb}, the morphism $\eta_X\circ \phi\in \Hom_\b2(\yb M,TX)$ derives from a morphism $(f,g)$ in $\bm2$. Moreover the normal image of the morphism $0:0\to \yb M$ in $\b2$ is the diagonal $\Delta\subset \yb M$ and thus exactness in $\b2$ at $\yb M$ is equivalent to \strg exactness. Then  the required injectivity follows from $(i)$. \endproof

\begin{figure}[H]
\begin{center}
\includegraphics[scale=.8]{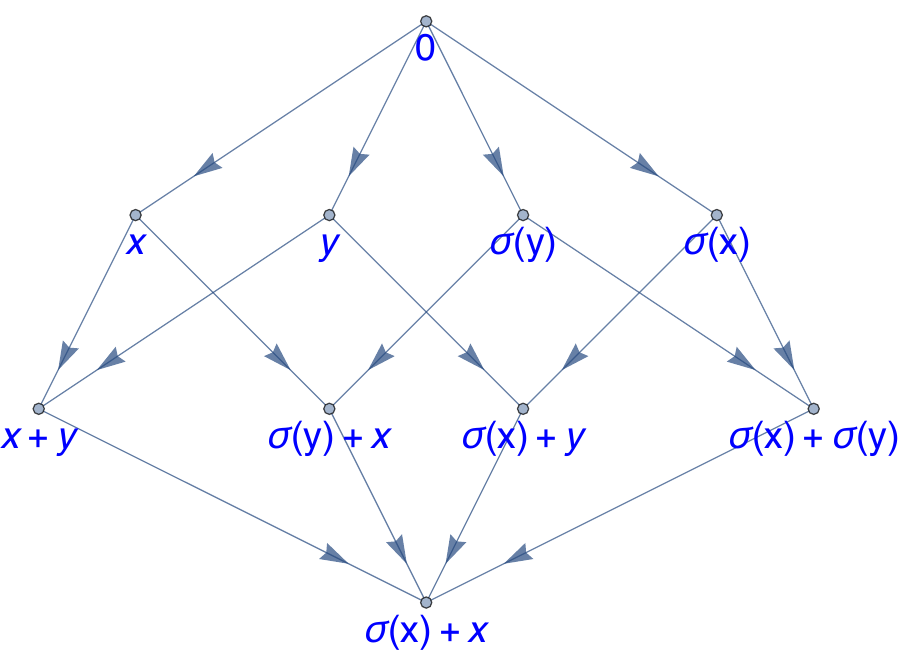}
\end{center}
\caption{The structure of the object $R$ of $\b2$. \label{submodrel} }
\end{figure}

In order to analyse the relations between two elements $x,y\in E$ such that $p(x)=p(y)$, we define a  finite object $R$ of $\b2$ whose dual $R^*$ is generated by its minimal elements, and a morphism $g:R\to E$ whose range  contains $x$ and $y$. Then we determine the normal subobjects of $R$.
The list of elements of $R$ is $\{0,a,c,d,b,a+c,a+d,b+c,b+d,a+b\}$ and the generators $a,b,c,d$ fulfill the relation $a+b=c+d$. The involution $\sigma$ interchanges $a$ with $b$ and $c$ with $d$. The only non-zero fixed element of $\sigma$ is $a+b=c+d$. Figure \ref{submodrel} shows the graph of the order relation. The notation in Figure \ref{submodrel} is to label the vertices by their image under the map $g:R\to E$ which is defined as follows
$$
g(a)=x,\ g(b)=\sigma(x), \ g(c)=y, \ g(d)=\sigma(y). 
$$
The relation $p(x)=p(y)$ corresponds to the relation $a+b=c+d$ and shows that $g$ is well defined. The first statement of the next proposition then shows that $R$ also fulfills  \eqref{nullmono}.

\begin{figure}[H]
\begin{center}
\includegraphics[scale=.8]{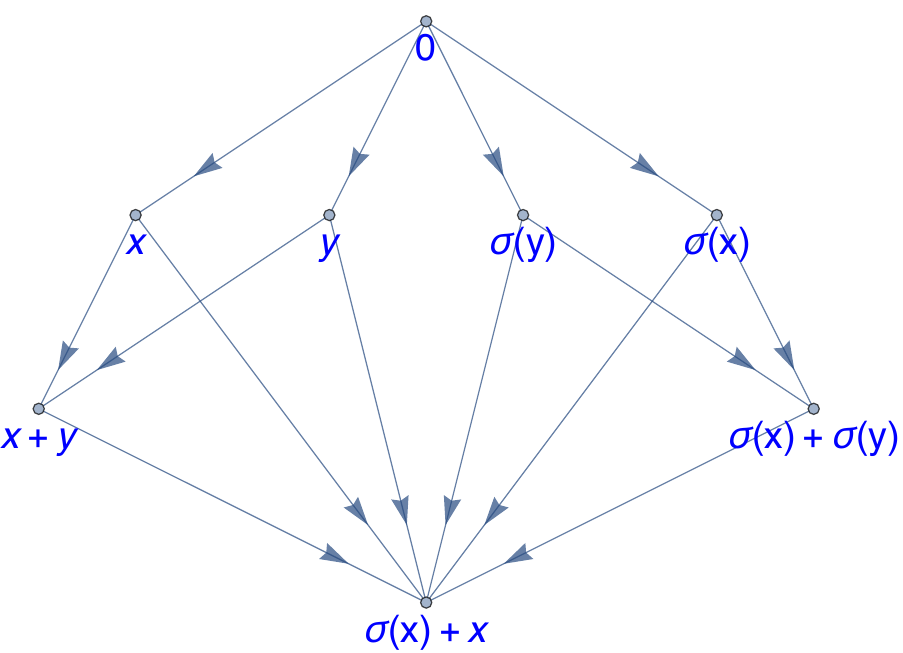}
\end{center}
\caption{The object $R'$ and the  effect of the further relation $x+\sigma(y)=\sigma(x)+y$. \label{submodrel2} }
\end{figure}

\begin{prop}\label{kernel and injective strong bis} $(i)$~$\phi\in \Hom_\b2(R,X)$ is a monomorphism if and only if its kernel is null.\newline
$(ii)$~Let $E$ be an object of $\b2$ and $x,y\in E$ such that $p(x)=p(y)$ and $x+y$ and $x+\sigma(y)$ are not null elements. Then for any morphism $f\in \Hom_\b2(E,F)$ with null kernel, the following ten elements are distinct
$$
0, f(x), \sigma(f(x)), f(y), \sigma(f(y)), f(x+y),\sigma(f(x+y)), f(x+\sigma(y)), f(y+ \sigma(x)),f(x+\sigma(x)).
$$	
\end{prop}
\proof $(i)$~The dual $R^*$ is generated by its minimal elements  so Theorem \ref{kernel and injective} applies and shows that the restriction of $\phi$ to the non-zero elements of $R$ is injective if $\ker(\phi)$ is null. Moreover under this hypothesis  one has $\phi(x)\neq 0$ $\forall x\neq 0$ since otherwise $\ker(\phi)$ would contain all of $R$ (since $x+\sigma(x)$ is the largest element of $R$ for any $x\in R$, $x\neq 0$). Thus $\phi$ is injective. The converse is clear.\newline
$(ii)$~We show that the condition that both $x+y$ and $x+\sigma(y)$ are not null elements, implies that the kernel of $g:R\to E$ is null. It is enough to prove that any normal submodule $N$ of $R$ which is not null contains $a+c$ or $b+d$. One has $R^\sigma=\{0,a+b\}$, thus Proposition \ref{existphi1} shows that for  $u\in N$, $u\neq 0$, any $\xi\in R$ such that $u+\xi=\xi$ satisfies $\xi\in N$. Indeed, one takes $a_0=u$ and one has $a_0+\xi=\xi$ and $p(a_0)=\xi+\sigma(\xi)$. Next, let $N$ be a submodule  of $R$ which is $\sigma$-invariant and fulfills:  
$u\in N, u\neq 0 \Rightarrow v\in N\qqq v\geq u$.
Then if $N\neq \{0\}$ one has $a+b\in N$ and if $N\neq \{0, a+b\}$ it contains $a+c$ or $b+d$ as shown by Figure \ref{submodrel}. we have shown that  the kernel of $g:R\to E$ is null, thus since $\ker(f)$ is null by hypothesis, the kernel of $f\circ g$ is also null. It thus follows from $(i)$ that the map $f\circ g:R\to F$ is injective.\endproof 
Let now $E$ be an object of $\b2$ and $x,y\in E$ such that $p(x)=p(y)$ and $x+\sigma(y)=\sigma(x)+y$. To understand the meaning of this relation we consider the finite object $R'$ of $\b2$ obtained by imposing on $R$ the further relation $a+d=b+c$. Then it follows that these two (equal) elements are both equal to $a+b$ and Figure \ref{submodrel2} shows the graph of the order on $R'$. This graph shows in particular that  $R'^*$ is not generated by its minimal elements. In fact the two minimal elements of $R'^*$ correspond to the segments $[0,a+c]$ and $[0,b+d]$.

Next, we investigate the meaning of the  technical hypothesis  that the dual $E^*$ of an object $E$ of $\b2$ is generated by its minimal elements. Next Example shows a class of cases where it fails.
\begin{example} \label{not sigma inject}
	Let $X$ be an arbitrary $\B$-module. We define an object $E=X\vee X$ of $\b2$ as follows. The underlying $\B$-module is  $X\cup_0 X'\cup \{t\}$, where $X'$ is a second copy of $X$ and the zero elements are identified, while the additional element $t$ is the largest element in $E$. The addition restricts to the addition in each copy of $X$ and is otherwise defined by $x+x'=t$ whenever $x\in X$, $x'\in X'$ are non-zero. The element $t$  is fixed by the involution $\sigma$ which interchanges the two copies of $X$. The elements $0,t$ are the only null elements. One has $\B\vee \B=\yb\B$ and moreover the maximal element $\mu$ of $X^*=\Hom_\B(X,\B)$, \ie the element such that $\mu^{-1}(\{0\})=\{0\}$, defines a morphism  $\mu\vee \mu:X\vee X\to \B\vee \B=\yb\B$ in $\b2$. By construction the kernel of $\mu\vee \mu$ is null but as soon as $X$ has more than one element $\mu\vee \mu$ fails to be $\sigma$-injective.
\end{example}
 To understand what happens in Example \ref{not sigma inject}, we assume for simplicity that $X$ is finite. The support $S$ of $f=\mu\vee \mu$ is by definition $S=\{z\in E\mid f(y)\leq f(z)\Rightarrow y\leq z\}$ and this selects the subset $S=\{0,u,u',t\}\subset E=X\vee X$ where $u:=\sum_X j$ is the largest element of $X$ (which exists since $X$ is assumed to be finite). When viewed as a subset of $E^*\simeq E^{\rm op}$, the sub-object $\hat S$ contains the null elements and we seek to understand why its normal closure is $E^*$. This follows from Proposition \ref{existphi1} since for any $\xi \in X^{\rm op}\subset E^{\rm op}$, non null, one has with $a=u^{\rm op}$: $\xi+a=\xi$ (since $a$ is minimal) and $a+\sigma(a)=\xi+\sigma(\xi)$. Thus the condition  of Proposition \ref{existphi1} holds with $n=0$ and $a_0=a$. When $X$ is not finite, $u:=\sum_X j$ no longer exists in general as an element of $X$ but it still makes sense as an ideal $J=\{x\in E\mid x\leq u\}$ so that the associated element $\phi_u\in E^*$ is meaningful. 
 \begin{rem}\label{strong exact not iso} Let $\mu\vee \mu:X\vee X\to \B\vee \B=\yb\B$ be as in Example \ref{not sigma inject}. Then the sequence 
 \begin{equation}\label{strong exact not iso1}
 	 0\to X\vee X\stackrel{\mu\vee \mu}{\to} \B\vee \B\to 0
 \end{equation}
 	is \strgly exact since the kernel of $\mu\vee \mu$ is null and its range is $\B\vee \B$. This result is in sharp contrast with Proposition \ref{iso9} which forbids for any object of $\bm2$ not isomorphic to $\B\vee \B$ to be on the left hand side of \eqref{strong exact not iso1}.
 \end{rem}
 
We can now refine Theorem \ref{kernel and injective}.
\begin{thm}\label{kernel and injective bis} Let $E$ be a finite object of $\b2$,  $f\in \Hom_\b2(E,F)$ a morphism with null kernel and $x,y\in E$ such that $p(x)=p(y)$ and $f(x)=f(y)$.  Then  $x\sim_{\rad(E)} y$.  
\end{thm}
 \proof As in the proof of Theorem \ref{kernel and injective} one gets that the normal image of the inclusion ${\rm Range}(f^*)\subset E^*$ is $E^*$. Let then $\xi\in E^*$ be a minimal element, then we show that $\xi(x)=\xi(y)$. If $\sigma(\xi)=\xi$ this follows from $p(x)=p(y)$. If $\sigma(\xi)\neq\xi$, Proposition \ref{extreme and range} states that $\xi\in {\rm Range}(f^*)$, \ie $\xi=L\circ  f$ for some $L\in F^*$, but then $\xi(x)=L(f(x))=L(f(y))=\xi(y)$. \endproof    
 In general (\cf~Example \ref{rad not functorial})  a morphism $f\in \Hom_\B(E,F)$ does not induce a morphism on the quotients by the radical congruence. One may wonder if the hypothesis of Theorem \ref{kernel and injective bis} may define a finer notion compatible with morphisms. More precisely, we let $h\in \Hom_\b2(E,Z)$ be a morphism with null kernel and $x,y\in E$ such that $p(x)=p(y)$ and $h(x)=h(y)$.  Then given a morphism $f\in \Hom_\b2(E,F)$ we complete the square below 
\begin{equation}\label{attempt functorial}
\xymatrix@C=25pt@R=25pt{
 E \ar[d]_{h}\ar[rr]^{f} && F\ar[d]^{\tilde h}
\\
 Z \ar[rr]^{\tilde f}&& F\oplus_{E}Z}
\end{equation}
by using the coequalizer $F\oplus_{E}Z$ of the two morphisms $E\to F\oplus Z$ obtained from  $f$ and $h$. By construction $\tilde h(f(x))=\tilde h(f(y))$  (since the diagram is commutative) and $p(f(x))=p(f(y))$. Thus if one could prove that $\ker(\tilde h)$ is null one could apply Theorem \ref{kernel and injective bis}. 
 On the other hand, the next example, related to  Example \ref{rad not functorial} shows that $\ker(\tilde h)$ is not null in general.
 
\begin{figure}[H]
\begin{center}
\includegraphics[scale=0.8]{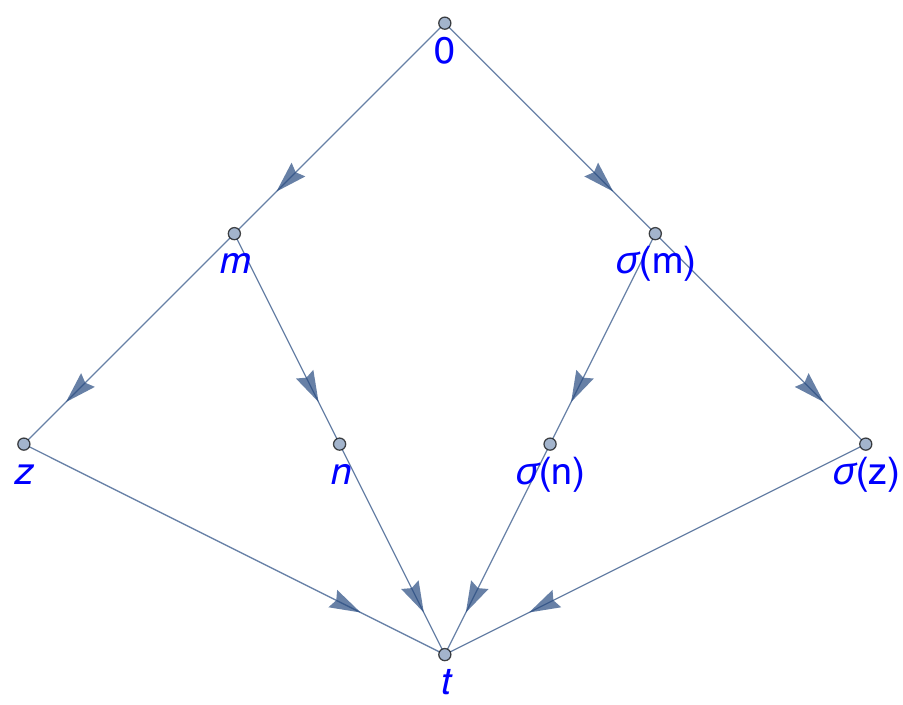}
\end{center}
\caption{The object $F$ and the subobject $E=N \vee N$. \label{submodnm} }
\end{figure}

 \begin{example}\label{rad not functorial bis} Let $N=\{0,m,n\}$ with $0<m<n$ and $E=N \vee N$, $Z=\B\vee\B$, $h=\mu\vee\mu$ (with the notations of Example \ref{not sigma inject}). We let $f:E\to F$ be the inclusion of $E$ in the object $F$ of $\b2$ obtained by adjoining a pair of elements $z, \sigma(z)$ as in Figure \ref{submodnm}. In $F$ the radical congruence is trivial so $f(n)$ and $f(m)$ are not equivalent. One can see directly that $\ker(\tilde h)$ is not null as follows. In $F\oplus_{E}Z$ one has, using $\tilde h(n)=\tilde h(m)$:  
 $
\tilde h(t)=\tilde h(z\vee n)=\tilde h(z\vee m) =\tilde h(z)
 $ 
 	and thus $z\in \ker(\tilde h)$.
 \end{example}
Example \ref{rad not functorial bis} suggests to refine the notion of kernel of a morphism in $\b2$ by taking into account the kernel of the push-outs. We first show how to enlarge a finite object of $\b2$ so that the extended object satisfies the hypothesis of Theorem \ref{kernel and injective}.
\begin{lem}\label{embed in good}
	Let $E$ be a finite object of $\b2$. Then there exists a finite object $F$ of $\b2$ whose dual $F^*$ is generated by its minimal elements and an embedding $f:E\to F$ of $E$ as a subobject of $F$.
\end{lem}
\proof Using the unit of the monad as in \S\ref{grandis remark} we can assume that the $\B$-module $E$ is of the form $A^2$ endowed with the involution $\sigma_A$, $\sigma_A(x,y):=(y,x)$. By Proposition \ref{inj3} there exists a finite $n\in \N$ and an embedding $A\subset \B^n$. This determines an embedding $f:E\to F$ of $E$ as a subobject of $F=\B^n\times \B^n$ where the dual $F^*$ is generated by its minimal elements since $F^*\simeq \B^n\times \B^n$. \endproof

Next, assume that in the push-out diagram \eqref{attempt functorial} the map $h$ is injective. Then we show that $\tilde h$ is injective. 

Let $u,v\in F$ with $u\neq v$. Let $L\in F^*$ such that $L(u)\neq L(v)$. Then $L\circ f\in E^*$ and since $ h$ is injective there exists $L'\in Z^*$ such that $L'\circ h=L\circ f$. Thus the morphism $(L,L')\in \Hom_\b2(F\oplus Z,\B)$ coequalizes the two morphisms from $E$ to $F\oplus Z$ and hence induces a morphism $L'':F\oplus_{E}Z\to \B$. One has $L''(\tilde h(u))=L(u)\neq L(v)=L''(\tilde h(v))$ so that $\tilde h(u)\neq \tilde h(v)$. Assume now that the map $h$ is only $\sigma$-injective. Then one can replace $Z$ by $Z':=Z \times E^\sigma$
and $h$ by $h'=(h,p)$. One has $\ker(h')=\ker(h)$. Moreover $h$ is  $\sigma$-injective iff $h'$ is injective. This leads us to consider the following push-out diagram
\begin{equation}\label{attempt functorial2}
\xymatrix@C=25pt@R=25pt{
 E \ar[d]_{h'=h\times p}\ar[rr]^{f} && F\ar[d]^{\tilde h'}
\\
 Z\times E^\sigma \ar[rr]^{\tilde f}&& F\oplus_{E}(Z\times E^\sigma)}
\end{equation}
\begin{lem}\label{embed in good bis}
	Let $E,F$ and $f:E\to F$ be as in Lemma \ref{embed in good} and let  $h\in \Hom_\b2(E,Z)$. Let $\tilde h'$ be as in the push-out diagram
\ref{attempt functorial2}. Then the following equivalence holds\vspace{.05in}

   \centerline{$h$  is $\sigma$-injective $\iff$ $\ker(\tilde h')$ is null.} 
\end{lem}
\proof Assume that $h$  is $\sigma$-injective, then $h'$ is injective and thus by the above reasoning $\tilde h'$ is injective. Thus $\ker(\tilde h')$ is null. Conversely, assume that $\ker(\tilde h')$ is null. Then by Theorem \ref{kernel and injective} the morphism $\tilde h'$ is $\sigma$-injective. This implies that $h$ itself is $\sigma$-injective. Indeed, for $x,y\in E$ with $p(x)=p(y)$ one has $p(f(x))=p(f(y))$ and if $h(x)=h(y)$ one gets using the commutativity of the diagram \eqref{attempt functorial2}, that 
$\tilde h'(f(x))=\tilde h'(f(y))$ so that $f(x)=f(y)$ and $x=y$.\endproof
In order to understand the push-out diagram \eqref{attempt functorial2} we consider the simplest case \ie $Z=0$. Then the diagram simplifies as
\begin{equation}\label{attempt functorial3}
\xymatrix@C=25pt@R=25pt{
 E \ar[d]_{ p}\ar[rr]^{f} && F\ar[d]^{\tilde p}
\\
  E^\sigma \ar[rr]^{\tilde f}&& F\oplus_{E} E^\sigma}
\end{equation} 
\begin{lem}\label{coker pushout}
	With the above notations one has  $F\oplus_{E} E^\sigma=\coker(f)$ and $\tilde p=\scoker(f)$.
\end{lem}
\proof The map $u\mapsto (u,0)$ from $F$ to $F\oplus_{E} E^\sigma$ is surjective since for $v\in E^\sigma$ one has $(0,v)=(0,p(v))\sim (f(v),0)$. Thus $F\oplus_{E} E^\sigma$ is a quotient of $F$ and for any object $C$ of $\b2$ the morphisms $\phi\in \Hom_\b2(F\oplus_{E} E^\sigma,C)$ are those morphisms 
$\psi\in \Hom_\b2(F,C)$ such that $\psi\circ f$ factors through $p$. This means that $\psi\circ f$ is null and hence as a quotient of $F$, $F\oplus_{E} E^\sigma$ coincides with $\coker(f)$, while the map $u\mapsto (u,0)$ is the same as $\scoker(f)$.\endproof 

\begin{rem}\label{pushout converse} The proof of Lemma \ref{coker pushout} shows that if in the  push-out diagram \eqref{attempt functorial} $\ker(\tilde h)$ is null then $h$  is $\sigma$-injective. However the converse is not true in general. For instance, when $Z=0$ the  push-out diagram \eqref{attempt functorial}  gives $F\oplus_{E}0=0$ as soon as $f(E)$ contains the largest element of $F$. If moreover $E$ is null while $F$ is not null, one would obtain that $h=0$ is $\sigma$-injective while $\ker(\tilde h)$ being not null.	
\end{rem}
Thus one can weaken the conditions on the morphism $f:E\to F$ so that Lemma \ref{embed in good bis} still holds. Notice first that one does not need to assume that $f$ is an embedding in fact it is enough to assume that $f$ is $\sigma$-injective. Moreover the condition that the dual $F^*$ is generated by its minimal elements is not necessary, rather one simply requires that any morphism with domain $F$ having null kernel is $\sigma$-injective. Indeed, Corollary \ref{trivial kernel1} states the existence  of objects $F$ which fulfill this last condition while $F^*$ is not generated by its minimal elements. Finally notice that using the diagram \eqref{attempt functorial2} one has $f(\ker(h))\subset \ker(\tilde h')$ and moreover when $f$ is $\sigma$-injective a non-null element of $\ker(h)$ is mapped by $f$ to a non-null element of $\ker(\tilde h')$.

\subsection{Injective and Projective objects in  $\b2$}\label{sectinjproj}

Proposition \ref{inj3} states that the product $\B^X$ of any number of copies of $\B$ is an injective object in the category of $\B$-modules and that any $\B$-module  is isomorphic to a submodule of a product $\B^X$. Next lemma relates the study of injective/projective objects in $\b2$ to the corresponding one in $\bmod$ by applying the (forgetful) functor $I: \b2\longrightarrow \bmod$.	
\begin{lem}\label{injective b2}
	$(i)$~An object $E$ of $\b2$ is injective if and only if the underlying $\B$-module $I(E)$ is injective in $\bmod$.\newline
	$(ii)$~An object $E$ of $\b2$ is projective if and only if the underlying $\B$-module $I(E)$ is projective in $\bmod$.\newline
	$(iii)$~For the {\em finite} objects in $\b2$ the properties of being injective/projective are equivalent and they mean that an object is a retract of a finite product $(\yb\B)^n$. 
\end{lem}
\proof $(i)$~Assume that $I(E)$ is injective in $\bmod$, consider an inclusion $L\subset M$ in $\b2$ and let $f:L\to E$ be a morphism in $\b2$. Since $I(E)$ is injective, let $g:M\to E$ be a morphism in $\bmod$ extending $f$. Then $h=g+\sigma\circ g\circ \sigma$ agrees with $f$ on $L$ and commutes with $\sigma$ so that it defines an extension of $f$ to a morphism $M\to E$ in $\b2$. Thus $E$ is injective in $\b2$. Conversely, let 
$E$ be an object  of $\b2$, then by Proposition \ref{inj3} there exists an embedding $\iota:I(E)\subset \B^X$, and the map $x\mapsto u(x)=(\iota(x),\iota(\sigma(x)))\in \yb\B^X$ gives an embedding in $\b2$ of $E$ as a subobject of $\yb\B^X$. If the object $E$ of $\b2$ is injective there exists a retraction of $u$, \ie a morphism  of $\b2$, $v:\yb\B^X\to E$ such that $v\circ u=\id$. It follows that the same equality holds in $\bmod$ and since $I(\yb\B^X)$ is injective in $\bmod$, so is $I(E)$.\newline
	$(ii)$~Assume that $I(E)$ is projective in $\bmod$, and let $M\to N$ be an epimorphism in $\b2$ and $f:E\to N$ a morphism in $\b2$. Since $I(E)$ is projective, let $g:E\to M$ be a morphism in $\bmod$ lifting $f$. Then as for $(i)$ $h=g+\sigma\circ g\circ \sigma$ determines a lift of $f$ in $\b2$. For the converse one uses Proposition \ref{inj3} and the existence of a surjective morphism $k:\B^{(X)}\to I(E)$ which extends to a surjective morphism of $\b2$ from $\yb\B^{(X)}$ to $E$. \newline
	$(iii)$~By Proposition \ref{inj3} for finite object of $\bmod$ the properties of being injective/projective are equivalent and they mean that an object is a retract of a finite product $\B^n$. Thus by 	
	$(i)$ and $(ii)$ one obtains the required statement. \endproof 
	\subsubsection{Normal subobjects of injective objects in $\b2$}	
	We now investigate the condition on an object $A$ of $\b2$ that $A$ embeds as a {\em normal} subobject of an injective object of $\b2$. 
	\begin{lem}\label{injective fixed2} Let $E$ be an object of $\b2$ such that $E^\sigma=\B$. Then if $E$ is projective one has $E=\B$ or $E=\yb \B$.		
	\end{lem}
\proof Let $X$ be the complement of $E^\sigma=\B$ in $E$, and $Y=X/\sigma$ the orbit space of the action of $\sigma$ on $X$ and $\iota:Y\to X$ an arbitrary section of the canonical surjection $X\to Y$. Let $\phi: (\yb \B)^{(Y)}\to E$ be the morphism in $\b2$ which is given on the copy of $\yb \B$ corresponding to $y\in Y$ by 
$$
\phi((1,0)_y):=\iota(y), \ \ \phi((0,1)_y):=\sigma(\iota(y))
$$
The morphism $\phi: (\yb \B)^{(Y)}\to E$ is surjective by construction. 	Let $p:(\yb \B)^{(Y)}\to \B^{(Y)}$ be the projection on the fixed points of the involution, we identify $\B^{(Y)}$ with the boolean $\B$-module of finite subsets of $Y$. Let $\psi: E\to (\yb \B)^{(Y)}$ be a morphism in $\b2$ such that $\phi\circ \psi=\id_E$. Since $E^\sigma=\B$ contains only one non-zero element $\tau$ one has $p(\psi(\xi))=\psi(\tau)$ for any non-zero $\xi\in E$. Let $S\subset Y$ be the finite subset of $Y$ corresponding to $\psi(\tau)$. Let us show that the existence of $\psi$ implies the following property of $E$:
\begin{equation}\label{easy sum}
	\xi+\eta=\tau \qqq \xi,\eta\in E\setminus \{0\}, \xi\neq \eta
\end{equation}
To prove this we can assume that $\xi,\eta\in E\setminus \{0,\tau\}$. For $\xi \in E\setminus \{0,\tau\}$ the components $\psi(\xi)_y$ are $0$ for $y\notin S$ and are either $(1,0)$, in which case we write $y\in S_\xi$, or $(0,1)$ for $y\in S$ since they are non-zero and if one component is $(1,1)$ one would get $\phi(\psi(\xi)))=\tau$. Thus to each $\xi \in E\setminus \{0,\tau\}$ corresponds a partition of $S$ as a disjoint union of $S_\xi$ and its complement. For two distinct elements $\xi,\eta\in E\setminus \{0,\tau\}$ these partitions are different and thus there exists $y\in S$ such that the components $\psi(\xi)_y$ and $\psi(\eta)_y$ are different and this gives $\phi\circ \psi(\xi+\eta)=\tau$ and thus \eqref{easy sum}. We have shown that if $E$ is projective it fulfills \eqref{easy sum}, and hence that the underlying $\B$-module $I(E)$ is obtained from the set $X$ by adjoining $0$ and the element $\tau$ while the addition is idempotent and uniquely specified by \eqref{easy sum}. Let then $\rho:\B^{(X)}\to I(E)$ the surjective morphism in $\bmod$  which associates  to any finite subset $F\subset X$ the sum $\sum_F\xi$. If $E$ is projective so is $I(E)$ by Lemma \ref{injective b2}, and thus there exists a section $\delta\in \Hom_\B(I(E), \B^{(X)})$, $\rho\circ \delta=\id_E$. For any $x\in X$ the element $\delta(x)$ is given by the subset $F=\{x\}$ and since the range of $\delta$ must be a submodule of  $\B^{(X)}$ one gets that $X$ has at most two elements. This shows that one is in one of the two cases $E=\B$ or $E=\yb \B$.	\endproof 
\begin{prop}\label{injective fixed2bis} Let $E$ be a finite object of $\b2$ such that $E^\sigma=\B$. Then if $E$ is isomorphic to a normal subobject of a finite injective object of $\b2$  one has $E=\B$ or $E=\yb \B$.		
\end{prop}
\proof Consider an embedding $E\subset I$ of $E$ as a normal subobject of a finite injective object of $\b2$. The equality $E^\sigma=I^\sigma$ shows that $I^\sigma=\B$. Since $I$ is finite and injective it is also projective by  Lemma \ref{injective b2}. Thus by Lemma \ref{injective fixed2} one has $I=\B$ or $I=\yb \B$ and hence $E=\B$ or $E=\yb \B$.	\endproof 
	As a corollary of Proposition \ref{injective fixed2bis} we get that among the examples \ref{not sigma inject} only the trivial ones give a normal subobject of a finite injective object.
	\subsubsection{A kernel-cokernel sequence in $\b2$}\label{ad-ex}
\begin{figure}[H]
\begin{center}
\includegraphics[scale=0.4]{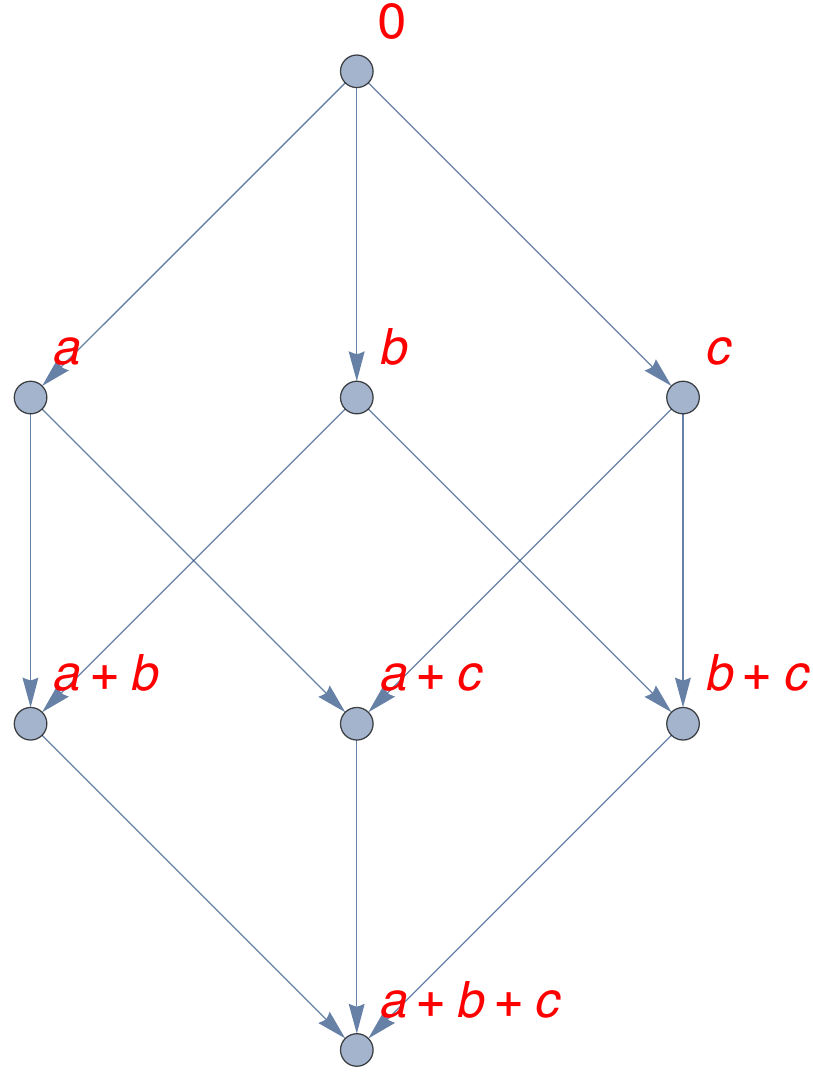}
\end{center}
\caption{Boolean object of $\bmod$ freely generated by $a,b,c$. \label{boole3} }
\end{figure}

\begin{figure}[H]
\begin{center}
\includegraphics[scale=0.5]{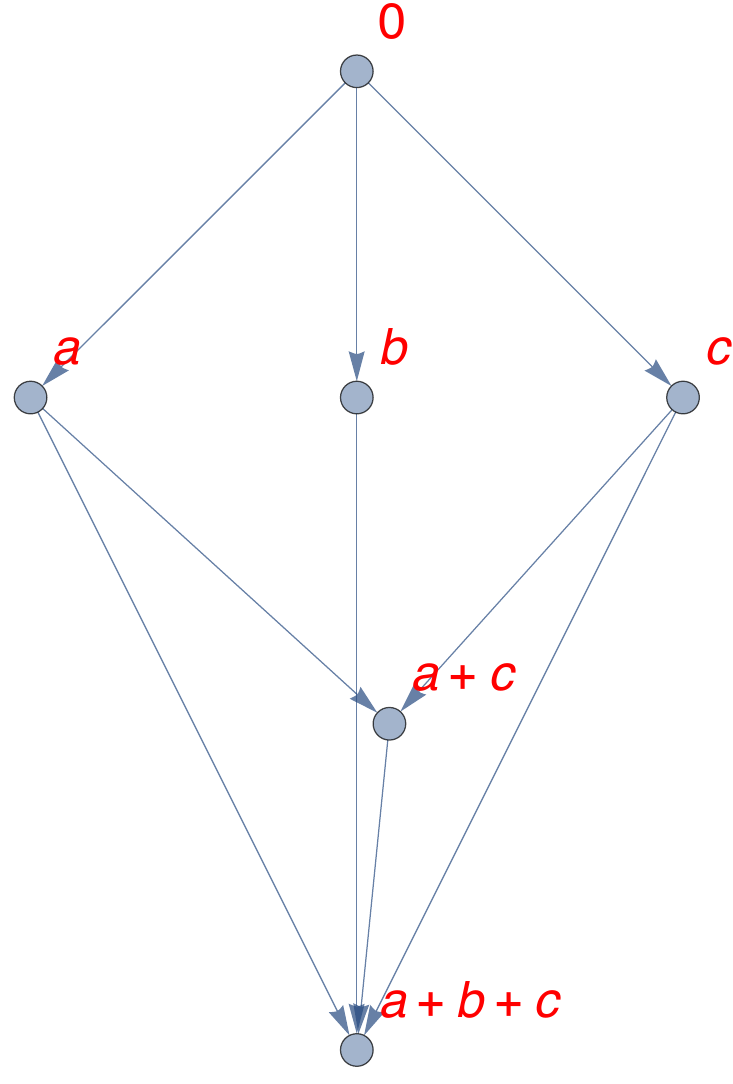}
\end{center}
\caption{The object $S$ of $\bmod$. \label{boole4} }
\end{figure}

	We investigate a specific  example of an object of $\b2$ which is not projective. 	
	Let $S$ be the object of $\bmod$ with three generators $a,b,c$ such that $a+b=b+c$. It follows that $a+b=b+c=a+b+c$. Thus in the Boolean object $\B^3$ of $\bmod$ freely generated by $a,b,c$ (Figure \ref{boole3}) one identifies $a+b=b+c=a+b+c$.
The object $\B^3$ of $\bmod$  becomes an object of $\b2$ when  endowed with the involution $\sigma(a)=c$, $\sigma(b)=b$. Let $S$ (see Figure \ref{boole4}) be the quotient of $\B^3$ by the relation $a+b=b+c$ and endowed with the induced involution. Let $\phi:\B^3\to S$ the quotient map. Consider the functor $H:=\Homi_\b2(S,-)$ viewed as a covariant endofunctor on $\b2$, using the natural internal Hom.

\begin{prop}\label{shortexactsimple} $(i)$~The kernel  of $\phi:\B^3\to S$ is the subobject $J\subset \B^3$ which is the complement of $a,c$ in $\B^3$.\newline
$(ii)$~The sequence $s:~J\stackrel{}{\rightarrowtail{}} \B^3 \stackrel{\phi}{\xtwoheadrightarrow{}} S$ is a short \exxx sequence and is non split in $\b2$.\newline
$(iii)$~$\id_S,\sigma_S\in H(S)=\Homi_\b2(S,S)$  are the only elements  which do not belong to the range of $H(\phi):H(\B^3)\to H(S)$. 
\newline
$(iv)$~The cokernel of the morphism $H(\phi):H(\B^3)\to H(S)$ is given by the projection $p$ on the range of $H(\phi)$ and by the identity on the complement of the range.
   \end{prop}
   \begin{figure}[H]
\begin{center}
\includegraphics[scale=0.5]{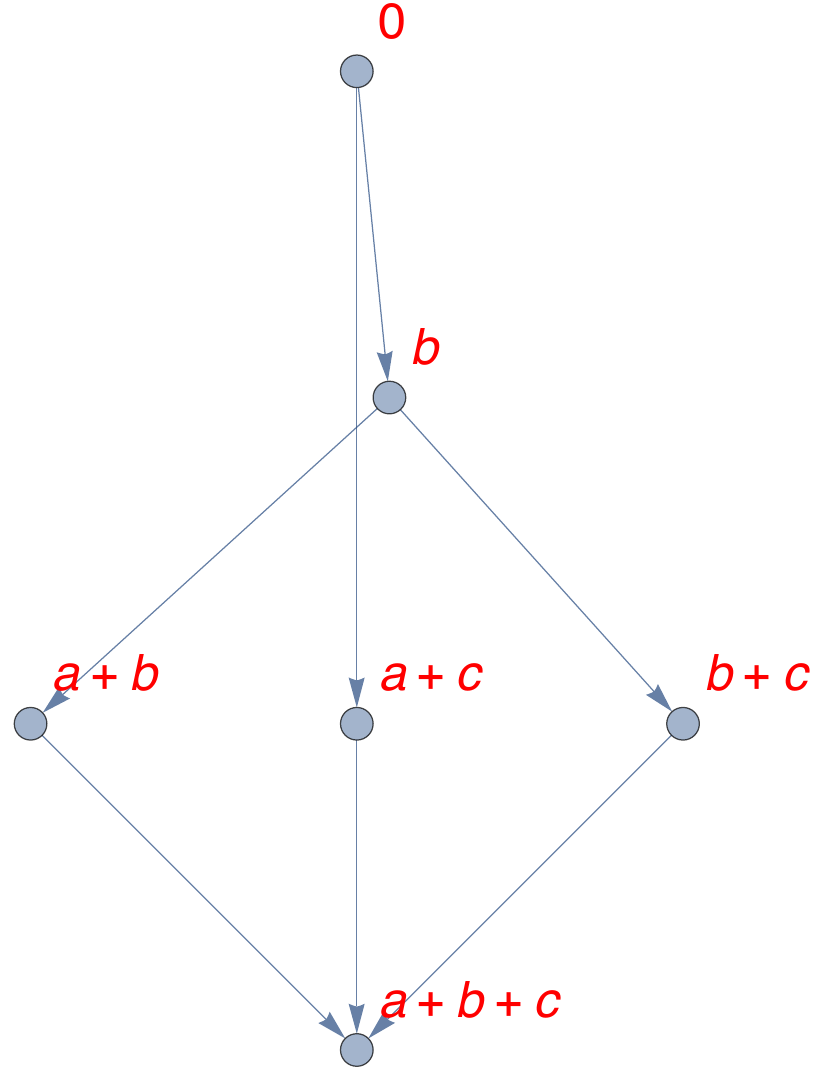}
\end{center}
\caption{Kernel $J$ of the map $\phi$ from $\B^3$ to $S$. \label{boole5} }
\end{figure}
   \proof $(i)$~The map $\phi:\B^3\to S$ is defined as 
   $$
   0\to 0,a\to a,b\to b,c\to c,a+b\to a+b+c,a+c\to a+c,b+c\to a+b+c,a+b+c\to a+b+c.
   $$
   and the submodule $S^\sigma\subset S$ of null elements is $S^\sigma=\{0,b,a+c,a+b+c\}$. It follows that the only elements of $\B^3$ which do not belong to the kernel of $\phi$ are $a,c$.\newline
$(ii)$~We show that the cokernel of the inclusion $J\subset \B^3$ is given by the map $\phi$. By Proposition \ref{compute cokernel}, this cokernel is the  quotient of $(\B^3)^\sigma\cup J^c$ by the smallest equivalence relation fulfilling the rule 
  $$
   \xi \in J^c,\ u,v\in J,\ p(u)=p(v)\Rightarrow \xi +u\sim \xi +v.
  $$
$\scoker(\phi)$ agrees with $p$ on $J$ and with the quotient map on $J^c$. For $\xi \in J^c=\{a,c\}$ one has $\xi +u\in J$ for any non-zero element $u\in J$. This shows that the equivalence relation is trivial, thus $\coker(\phi)=(\B^3)^\sigma\cup J^c$. Since the  kernel  of $\phi:\B^3\to S$ is the subobject $J\subset \B^3$, the map $\phi$ factors through $\coker(\phi)$ and hence agrees with it. The proof of $(iii)$ below shows that $s$ is not split. \newline
$(iii)$~We first show that the element $\id_S\in H(S)=\Homi_\b2(S,S)$ does not belong to the range of the morphism $H(\phi):H(\B^3)\to H(S)$. For any object $M$ of $\b2$, an element  $f \in H(M)$ is given by a pair $(\alpha,\beta)$ of elements of $M$ such that $\beta\in M^\sigma$ and $\alpha + \beta\in M^\sigma$. One lets $\alpha=f(a)$ and $\beta=f(b)$. For the element $\id_S\in H(S)=\Homi_\b2(S,S)$ one has $\alpha=a$ and $\beta=b$. These two elements lift uniquely to elements of $\B^3$ but the lifts $\alpha'=a, \beta'=b$ no longer satisfy  $\alpha' + \beta'\in (\B^3)^\sigma$. In the same way one sees that the element $\sigma_S\in H(S)=\Homi_\b2(S,S)$ does not belong to the range of the morphism $H(\phi):H(\B^3)\to H(S)$. There are $6$ endomorphisms of $S$ which are not null, they correspond to the values of $(\alpha,\beta)$ reported here below
$$
\left(
\begin{array}{cc}
 a & b \\
 a & a+c \\
 a & a+b+c \\
 c & b \\
 c & a+c \\
 c & a+b+c \\
\end{array}
\right)
$$
 Then, $\id_S$ which corresponds to $(a,b)$, and $\sigma_S$ which correponds to $(c,b)$ are both not liftable. But the other $4$ endomorphisms are liftable since the sum $\alpha+\beta$ is null already in $\B^3$. 
\newline
$(iv)$~Let $\xi=\id_S\in H(S)=\Homi_\b2(S,S)$. We determine the interval $[0,\xi]\subset  H(S)$. One has $f\in [0,\xi]$ iff $\alpha\leq a$ and $\beta \leq b$. For $\alpha= a$ and $\beta =0$ this does not fulfill $\alpha + \beta\in S^\sigma$. Thus the only non-trivial element $f\in [0,\xi]$ corresponds to
$\alpha= 0$ and $\beta =b$. This element is null, and thus Proposition \ref{extreme and range} $(iii)$ applies to show that $\xi\notin \imm(H(\phi))$.  This shows that the cokernel of the morphism $H(\phi):H(\B^3)\to H(S)$ is non null.
The result then follows from Proposition \ref{compute cokernel}.\endproof 
Next, we consider the endomorphisms of the short exact  sequence $s:~J\stackrel{}{\rightarrowtail{}} \B^3 \stackrel{\phi}{\xtwoheadrightarrow{}} S$, \ie the endomorphisms 
$f\in \End_\b2(\B^3)$ such that $f(J)\subset J$. They define a correspondence between $\End_\b2(J)$ and  $\End_\b2(S)$ displayed in the following diagram
\begin{equation}\label{quotient maph}
	\End_\b2(J)\stackrel{{ res}}{\leftarrow} \End(s)\stackrel{{ quot}}{\to} \End_\b2(S).
\end{equation}
The left arrow is given by restriction of $f$ to $J\subset \B^3$ and the right one is defined by the induced morphism on the cokernel.
\begin{prop}\label{shortexactsimple1} $(i)$~$f\in \End_\b2(\B^3)$ is uniquely specified by the pair $(f(a),f(b))\in \B^3\times (\B^3)^\sigma$.\newline
$(ii)$~Among the $32$ endomorphisms $f\in \End_\b2(\B^3)$ only two do not fulfill $f(J)\subset J$, they correspond to the pairs $(a,0)$ and $(c,0)$.\newline
$(iii)$~For any $h\in \End_\b2(S)$ such that $h\notin \imm(H(\phi))$ the restriction  of $h$ to $S^\sigma$ is an automorphism.\newline
 $(iv)$~Let $v\in\End_\b2(J)$  admit more than one extension to $\B^3$. Then for any of these extensions $w$ the restriction to $S^\sigma$ of the induced morphisms  $w''\in \End_\b2(S)$ fails to be surjective and  $w''\in \imm(H(\phi))$.
   \end{prop}
\proof $(i)$~The commutation of $f$ with the involution $\sigma$ of $\B^3$ entails that $f(c)=\sigma(f(a))$ and that $f(b)\in (\B^3)^\sigma$. Conversely, any pair $(f(a),f(b))\in \B^3\times (\B^3)^\sigma$ uniquely defines an $f\in \End_\b2(\B^3)$.\newline
$(ii)$~Let  $f\in \End_\b2(\B^3)$. One has $f(b)\in (\B^3)^\sigma\subset J$. Thus if $f(a)=0$ then $f(\B^3)\subset (\B^3)^\sigma\subset J$. Also, if $f(a)\in J$ then
$f(\B^3)\subset J$. Thus, by Proposition \ref{shortexactsimple} $(i)$, we just need to consider the cases $f(a)=a$ and $f(a)=c$.  Assume $f(a)=a$. Then if $f(b)\neq 0$ one gets $f(J)\subset J$ since $f(c)=c$ and for any subset $Y\subset 
\{a,b,c\}$ not reduced to $a$ or $c$ one has $\sum_Y f(y)\in J$. If $f(b)= 0$ one has $f(a+b)=a$ and this contradicts $f(J)\subset J$.\newline
$(iii)$~By Proposition \ref{shortexactsimple} $(iii)$, the only such $h$ are $\id_S$ and $\sigma_S$.\newline
 $(iv)$~The restriction map of \eqref{quotient maph} is surjective since $\B^3$ is injective. By $(ii)$ there are $30$ endomorphisms of the short \exxx sequence $s$ and when we take their restriction to $J$ one obtains the following $4$ non-trivial fibers with three elements each (where we list the corresponding pairs $(f(a),f(b))\in \B^3\times (\B^3)^\sigma$): 
 $$
 \left(
\begin{array}{cc}
 a & a+c \\
 c & a+c \\
 a+c & a+c \\
\end{array}
\right), \  \ \left(
\begin{array}{cc}
 a+b & a+c \\
 b+c & a+c \\
 a+b+c & a+c \\
\end{array}
\right), \ \ \left(
\begin{array}{cc}
 a & a+b+c \\
 c & a+b+c \\
 a+c & a+b+c \\
\end{array}
\right),\ \ \left(
\begin{array}{cc}
 a+b & a+b+c \\
 b+c & a+b+c \\
 a+b+c & a+b+c \\
\end{array}
\right)
$$ 
In all these cases the restriction of $f$ to $(\B^3)^\sigma$ fails to be surjective. This restriction is the same as the restriction to $S^\sigma$ of the induced morphisms $w''$ on $S$ by Proposition \ref{fixed cokernel}. Thus by $(iii)$ one gets $w''\in \imm(H(\phi))$.
There are $22$ elements in $\End_\b2(J)$ and since the  restriction map \eqref{quotient maph} is surjective one concludes that all other fibers are reduced to a single element. \endproof

\section{Homological algebra in a homological category}\label{homological grandis}
 
 For our long term purpose of obtaining 
  a computable sheaf cohomology of $\B$-modules over a topos, the case of interest is the right derived functors $\ext^n$ of the functor $F:=\Hom(L,-)$ for a fixed object $L$. The  category in which this process should take place, in the case of the one point topos, is the category $\b2$. In this section we apply the construction of chapter 4 of \cite{gra1} of satellite functors as Kan extensions. After recalling in \S\ref{sect chain complexes} the framework of homological algebra from \cite{gra1}, we give a straightforward adaptation of the construction to left exact functors and right satellite in \S \ref{right satellite}. In \S \ref{sect condition a}  we explain the key condition (a) of \cite{gra1} which allows one to compute the satellite functor $SF$. The latter is defined, on an object $X$ as a colimit indexed by a  comma  category $J=\sequ(\cE)\downarrow_{P'} X$ built from  short \exxx sequences. This indexing category does not in general admit a final object but under good circumstances it admits a weakly final object coming from a semi-resolution $i$ of $X$. Condition (a) of \cite{gra1}, Theorem 4.2.2, is the requirement that such weakly final object of the indexing category $J$ becomes final after applying the functor $\alpha\longrightarrow \coker(F(\alpha''))$ which defines the colimit. In other words applying this functor should erase the ambiguity created by the weak finality of the semi-resolution. We explain 
  the link between condition (a) and the radical  in the auxiliary section \S \ref{sect a and radical}.
Our main result in this development is explained in \S \ref{sect reduction to endo}, it is Theorem \ref{reduction toendo}  which, under the hypothesis that the middle term of the semi-resolution is both injective and projective, reduces the proof of condition (a) to endomorphisms of the weakly final  short \exxx sequence. This result is applied in Theorem \ref{reduction toendo cor1} to show that the satellite functor of the hom functor is not null for a specific finite object of $\b2$.

\subsection{Short \exxx sequences of chain complexes}\label{sect chain complexes}
We first recall that in the case of an abelian category the functor $\Hom(A,-)$ (for fixed $A$) is left exact and hence admits right derived functors $R^nT=\ext^n(A,-)$ which are computed, by using an injective resolution $I^-$, as the cohomology of the complex  
 $$
 \cdots\to\Hom(A,I^j)\to \Hom(A,I^{j+1})\to \Hom(A,I^{j+2}) \to\cdots
 $$
We start this section by comparing two notions of ``left exactness" in this setup. The first notion we take up is the purely categorical notion
\begin{defn}\label{categorical left exact}
A functor between finitely complete categories is called left exact (or flat) if it preserves finite limits.	
\end{defn}
 Let $F\dashv G$ be a pair of functors 
 with $F$
 left adjoint to $G$. Then $F$ 
preserves all colimits and $G$
 preserves all limits. For the category $\b2$ and a fixed object $A$ we consider the internal $\Hom$ functor $\Homi(A,-)$ as a covariant functor $G:\b2\longrightarrow \b2$, then  one has the adjunction 
 \begin{equation}\label{basic adjunction}
 \Hom_\b2(X,\Homi(A,Y))\simeq \Hom_\b2(X\otimes A,	Y)
 \end{equation}
where $\otimes$ denotes the tensor product over $\B$ endowed with the involution, \ie the tensor product of the two involutions. This suffices to show that $G$ is left  exact in the categorical sense of Definition \ref{categorical left exact}. 
 We shall now review some results taken from \cite{gra1} and explain their expected role in the development of a computable sheaf cohomology of $\B$-modules over a topos. Chapter III, 3.3 of \cite{gra1} develops the homology sequence associated to a short \exxx sequence of chain complexes. In a homological category $\cE$  one defines unbounded chain complexes as sequences
 $$
 \cdots \to A_{n+1}\stackrel{\partial_{n+1}}{\longrightarrow}A_n \stackrel{\partial_{n}}{\longrightarrow} A_{n-1}\to\cdots
 $$
 indexed by $n\in \Z$ and satisfying the condition that $\partial_n\circ \partial_{n+1}$ is always null. A morphism of complexes is a sequence of morphisms $f_n:A_n\to B_n$ such that $\partial_n\circ f_n=f_{n-1}\circ \partial_n$. One defines the null morphisms of complexes as those for which all the components $f_n$ are null. This yields the category $\chain(\cE)$ of chain complexes over $\cE$ and the subcategory $\chainp(\cE)$ of positive chain complexes (a positive chain complex is completed in negative indices by null objects). These categories are also homological.

One of the key results of \cite{gra1} is the following statement (homology sequence) \cf~ 3.3.5.
\begin{thm}\label{homology sequence} Let $\cE$ be a homological category and consider a short \exxx sequence of chain complexes
$$
U\stackrel{m}{\rightarrowtail{}} V \stackrel{p}{\xtwoheadrightarrow{}} W, \ \ m=\sker(p), \ p=\scoker(m).
$$
	$(a)$~There is a homology sequence of order two, natural for morphisms of short \exxx sequences	\begin{equation}\label{long exact sequence}
	\cdots\longrightarrow  H_n(V)\stackrel{H_n(p)}{\longrightarrow}H_n(W)\stackrel{\partial_n}{\longrightarrow}H_{n-1}(U)\stackrel{H_{n-1}(m)}{\longrightarrow}H_{n-1}(V)\longrightarrow\cdots
	\end{equation}
where  $\partial_n$ is induced by the differential $\partial_n^V$
of the complex $V$.\newline
$(b)$~If the differential $\partial_n^V$ of the central complex $V$ is an exact morphism, so is the differential $\partial_n$ of the homology sequence; moreover, the sequence itself is exact in the domain of $\partial_n$ (\ie $H_n(W)$) and in its codomain (i.e. $H_{n-1}(U)$).\newline
$(c)$~If the following conditions hold for every $n\geq 0$, the homology sequence is exact
\begin{equation}\label{long exact sequence1}
	(B_nV\vee U_n)\wedge Z_nV=B_nV\vee(U_n\wedge Z_nV)
\end{equation}
\begin{equation}\label{long exact sequence2}
\partial^*\partial_*(U_n)=	U_n\vee Z_nV, \ \ \partial_*\partial^*(U_{n-1})=U_{n-1}\wedge B_{n-1}V.
\end{equation}
\end{thm}

The central exactness stated in $(b)$  
 is the key to prove the universality of chain homology for non-exact categories (cf. \cite{gra1}  Section 4.5).
The exact couples are treated in 3.5 of \opcit

\subsection{The right satellite of a left exact functor}\label{right satellite}
In this subsection we transpose the treatment of the left satellite of a right exact functor explained in \cite{gra1}  Section 4.1, to the construction of the right satellite of a left exact functor such as the internal $\Hom$ functor $\Homi(A,-)$. The right satellite is constructed as a left Kan extension in the sense of \cite{MacLane} p. 240. One considers a homological category $\cE$ and the category $\sequ(\cE)$ of short \exxx sequences of $\cE$ \ie of sequences of the form
$$
A'\stackrel{a'}{\rightarrowtail}A\stackrel{a''}{\xtwoheadrightarrow{}}A''~{\rm s.t.}~ a'= \sker(a'')\ \& \  a''= \scoker(a').
$$
 The morphisms between two such sequences of $\cE$ are the same as the morphisms 
 for the corresponding 2-truncated chain complexes (\ie one has morphisms $A'\to B'$ etc. commuting with the $a',b'$ etc). Then, one considers a left exact functor $F:\cE\longrightarrow \cB$ where $\cB$ is also a homological category and cocomplete. One seeks to construct a long sequence of derived functors $F^n$ where $F^0=F$ and where each short \exxx sequence of $\cE$ gives rise to a long sequence of the form
 $$
 \cdots \to F^nA'\stackrel{F^n a'}{\to}F^nA\stackrel{F^na''}{\to}F^nA''\stackrel{d^n}{\to} F^{n+1}A'\stackrel{F^{n+1} a'}{\to}F^{n+1}A\stackrel{F^{n+1}a''}{\to}F^{n+1}A''\to\cdots
 $$
 In order to construct the first right satellite $SF=F^1$ one associates to a short \exxx sequence of $\cE$ the object of $\cB$ given by $\coker(F(a''))$. In fact, one should define
 $SF$ in such a way that all these objects map to $SF(A')$. Thus  $SF(A')$ is defined as a colimit of the $\coker(F(a''))$'s. One fixes an object $X$ of $\cE$ and considers the comma  category $\sequ(\cE)\downarrow_{P'} X$, where  $P':\sequ(\cE)\longrightarrow \cE$ is the functor of projection on $A'$. Thus, an object of this comma category is of the form
 \begin{equation}\label{Comma category}
X\stackrel{x}{\leftarrow} A'\stackrel{a'}{\rightarrowtail}A\stackrel{a''}{\xtwoheadrightarrow{}}A''.
\end{equation}  
We introduce the following notion to handle the preservation of null objects by satellite functors.
\begin{defn}\label{N-retraction}  An  $\cN$-retraction $\rho$ of an homological category  $\cE$ is provided by an endofunctor $\rho:\cE\longrightarrow \cE$ which projects on the subcategory of null objects, and by a natural transformation $p$ from the identity functor to $\rho$.
\end{defn}
 In fact the retraction exists in general for any homological category  $\cE$ and it is defined by the functor 
$$
\rho:\cE\longrightarrow \cE,\ \ \rho(X):=\coker(\id_X), \ \ p_X:=\scoker(\id_X):X\to \coker(\id_X).
$$
For  the category $\b2$ the endofunctor $\rho$ associates to an object $N$ of $\b2$  the fixed points $N^\sigma$ under $\sigma$. The natural transformation $p$ is defined by $p(x)=x+\sigma(x)$.
 \begin{lem}\label{first satellite} Let $\sequ_{\rm small}$ be a small subcategory of the category $\sequ(\cE)$ of short \exxx sequences in the homological category $\cE$. Let $F:\cE\longrightarrow \cB$ be a covariant functor where $\cB$ is a cocomplete  homological category. The following colimit is meaningful for any object $X$ of $\cE$ and defines a covariant functor
 \begin{equation}\label{Kan extension}
	SF:\cE\longrightarrow \cB, \ \ SF(X):=\varinjlim_I \coker(F(a'')), \ \ I=\sequ_{\rm small}(\cE)\downarrow_{P'} X.
\end{equation} 
Moreover, assuming  that $F$ sends null objects to null objects, one obtains for any object of $\sequ_{\rm small}$,  an order two sequence 
 \begin{equation}\label{satellite sequ}
	FA'\stackrel{F a'}{\to}FA\stackrel{Fa''}{\to}FA''\stackrel{d}{\to} SF A'\stackrel{SF a'}{\to}SF A\stackrel{SFa''}{\to}SFA''.\end{equation} 
  \end{lem} 
 \proof The smallness of $\sequ_{\rm small}$ ensures that the comma category $\sequ_{\rm small}(\cE)\downarrow_{P'} X$ is small. Since $\cB$ is cocomplete, the colimit $SF(X):=\varinjlim_I \coker(F(a''))$ is meaningful. Let $f:X\to Y$ be a morphism in $\cE$. The following assignment  defines a functor $I\longrightarrow J=\sequ_{\rm small}(\cE)\downarrow_{P'} Y$
 $$
 \left(X\stackrel{x}{\leftarrow} A'\stackrel{a'}{\rightarrowtail}A\stackrel{a''}{\xtwoheadrightarrow{}}A''\right)\mapsto \left(Y\stackrel{f\circ x}{\leftarrow} A'\stackrel{a'}{\rightarrowtail}A\stackrel{a''}{\xtwoheadrightarrow{}}A''\right).
 $$
 This  yields a natural morphism $SF(f)\in\Hom_\cB(SFX,SFY)$,  where $SF(X):=\varinjlim_I \coker(F(a''))$  and $SF(Y):=\varinjlim_J \coker(F(a''))$. Thus $SF:\cE\longrightarrow \cB$ is a covariant functor.   To obtain the map $d$ of \eqref{satellite sequ} one uses the object $i$ of $\sequ_{\rm small}(\cE)\downarrow_{P'} A'$ given by $i:~A'\stackrel{\id}{\leftarrow} A'\stackrel{a'}{\rightarrowtail}A\stackrel{a''}{\xtwoheadrightarrow{}}A''$
 and the natural morphism $\phi_i: \coker(F(a''))\to SF(A')$ given by the index $i$ and the construction of the colimit. After composition with $\scoker(F(a'')):A''\to\coker(F(a''))$, one obtains the morphism  $FA''\stackrel{d}{\to} SF A'$, $d=\phi_i\circ \scoker(F(a''))$. One has to show that $d\circ F(a'')$ is null and  this follows from 
 $\scoker(F(a'')\circ F(a'')$ null. Next, we consider $SFa'\circ d$. By the functoriality of $SF$ one gets the map $\phi_j\circ  \scoker(F(a''))$, where the object $j$ of $\sequ_{\rm small}(\cE)\downarrow_{P'} A$ is given by composition with $i$ as
 $$
 i=\left(A'\stackrel{\id}{\leftarrow} A'\stackrel{a'}{\rightarrowtail}A\stackrel{a''}{\xtwoheadrightarrow{}}A''\right)\mapsto j=\left(A\stackrel{a'}{\leftarrow} A'\stackrel{a'}{\rightarrowtail}A\stackrel{a''}{\xtwoheadrightarrow{}}A''\right).
 $$
 The following commutative diagram defines a morphism in the comma category $\sequ_{\rm small}(\cE)\downarrow_{P'} A$ from the object $j$ to the object $j'$ represented by the lower horizontal line
\begin{equation}\label{colimit compat}
\xymatrix@C=20pt@R=35pt{
&  A' \ar[ld]_{a'}\ar[d]_{a'}\ar@{>->}[rr]^{a'} && A\ar[d]_{\id} \ar@{->>}[rr]^{a''}&& A''\ar[d]
\\
A & \ar[l]_{\id}A \ \ar@{>->}[rr]^{\id}&& A \ar@{->>}[rr]^{\scoker(\id)}&& \coker(\id) }
\end{equation}
 By construction  $\coker(\id)$ is a null object and it follows that $SFa'\circ d$ is null, since in the colimit one has  $\phi_j\circ  \scoker(F(a''))\sim \phi_{j'}\circ  \scoker(\id)$ which factorizes through the null object $\coker(\id)$. It remains to show that the composition $SFa''\circ SFa'$ is null. This will follow if we show that the image by $SF$ of a null object $N$ is a null object. In order to prove this statement we use an  $\cN$-retraction $\rho$ of $\cE$ 
\begin{equation}\label{using retraction}
\xymatrix@C=20pt@R=35pt{
&  A' \ar[ld]_{x}\ar[d]^{p_{A'}}\ar@{>->}[rr]^{a'} && A\ar[d]_{p_A} \ar@{->>}[rr]^{a''}&& A''\ar[d]^{p_{A''}}
\\
N &\ar[l]_{\rho(x)} \rho(A') \ \ar@{>->}[rr]^{\rho(a')}&& \rho(A) \ar@{->>}[rr]^{\rho(a'')}&& \rho(A'') }
\end{equation}
The  diagram \eqref{using retraction} is commutative when $N$ is a null object and this  shows that $SF(N)$ is a colimit of null objects and hence a null object. \endproof 
The notion of normally injective object in a homological category is introduced in \cite{gra1} 4.2.1. An object $I$ is said to be normally injective if for any normal monomorphism $m:A\to B$, every morphism $f:A\to I$ extends to a morphism $g:B\to I$ such that $f=g\circ m$. 
\begin{lem}\label{inject vanishing}Let $I$ be a normally injective object of $\cE$ then, with the notations of Lemma \ref{first satellite},  $SF(I)$ is a null object of $\cB$.	
\end{lem}
\proof Let $i$ be an object of the comma category $\sequ_{\rm small}(\cE)\downarrow_{P'} I$
as in the upper horizontal line of the diagram
\begin{equation}\label{colimit compat1}
\xymatrix@C=20pt@R=35pt{
I&  A' \ar[l]_{x}\ar[d]^{a'}\ar@{>->}[rr]^{a'} && A\ar[d]_{\id} \ar@{->>}[rr]^{a''}&& A''\ar[d]
\\
 & \ar[ul]^{\tilde x}A \ \ar@{>->}[rr]^{\id}&& A \ar@{->>}[rr]^{\scoker(\id)}&& \coker(\id) }
\end{equation}
Since $I$ is normally injective the morphism $x:A'\to I$ extends to a morphism $\tilde x:A\to I$ such that the diagram \eqref{colimit compat1} is commutative. By construction of the colimit it follows that $i$ gives a null object in $SF(I)$. Indeed, the cokernel $\coker(\id)$ is a null object and hence both $F(\coker(\id))$
and $\coker(F(\scoker(\id)))$ are null. \endproof 

  \subsection{Condition (a)}\label{sect condition a}

Consider two morphisms of short \exxx sequences in $\b2$ of the form
\begin{equation}\label{weakly final abelian}
\xymatrix@C=20pt@R=35pt{
  C'\ar[d]^{v}\ar@{>->}[rr]^{c'} && C \dar[d]_{w_1}^{\ \, w_2} \ar@{->>}[rr]^{c''}&& C''\dar[d]_{w''_1}^{\ \, w''_2}
\\
I' \ \ar@{>->}[rr]^{\alpha'}&& I \ar@{->>}[rr]^{\alpha''}&& I''}
\end{equation}
so that both morphisms coincide on $v:C'\to I'$. Thus one starts with two extensions $w_j:C\to I$ of the morphism $v:C'\to I'$ using the fact that $I$ is assumed to be injective.  Given a covariant functor 
$F:\b2\longrightarrow \b2$ one wishes to show the equality 
\begin{equation}\label{grandis condition}
F(w''_1)=F(w''_2): \coker(F(c''))\to \coker(F(\alpha'')).	
\end{equation}
In order to understand the inherent difficulty to obtain this condition (labelled condition (a) in \cite{gra1}, Theorem 4.2.2) in relation with the functoriality, we first recall the proof that such condition holds for the functor $F:=\Hom(Q,-)$ in abelian categories. In this classical setup  
$$
\coker(F(c''))=\Hom(Q,C'')/(c''\circ \Hom(Q,C)), \ 
\coker(F(\alpha''))=\Hom(Q,I'')/(\alpha''\circ \Hom(Q,I)).
$$
In the abelian context, one can take $w=w_1-w_2$ and this induces $w''=w''_1-w''_2$ in $\Hom(C'',I'')$. Since $w_1\circ c'=w_2\circ c'$, one has $w\circ c'=0$ and this determines (diagram chasing) a lift $\tilde w\in \Hom(C'',I)$ such that $\alpha''\circ \tilde w=w''$. It follows that
$F(w''): \coker(F(c''))\to \coker(F(\alpha''))$ is zero since for any $\phi\in \Hom(Q,C'')$
the composition $w''\circ \phi$ lifts to $\tilde w\circ \phi\in \Hom(Q,I)$ and thus vanishes in the quotient $\coker(F(\alpha''))=\Hom(Q,I'')/(\alpha''\circ \Hom(Q,I))$.\vspace{.05in}

In $\b2$ the cokernel of a morphism $\phi:L\to M$ is never reduced to $\{0\}$ and maps always surjectively  onto the cokernel of the identity map $\id_M$ and $\scoker(\id_M)$ is the projection $p:M\to M^\sigma$. In fact, Proposition \ref{fixed cokernel} states a canonical isomorphism $(\coker(\phi))^\sigma\simeq M^\sigma$. Thus when we consider the functor $F:=\Hom(Q,-)$ and test the functoriality of $F(w'')$ as in \eqref{grandis condition}, we can first replace the involved cokernels $\coker(F(c''))$ and $\coker(F(\alpha''))$ by their null submodules and test the functoriality there.
\begin{lem}\label{functoriality on fixed}
$(i)$~Let $N$, $M$ be two objects of $\b2$. Then one has a canonical isomorphism 
$$
r:\Homi_\b2(N,M)^\sigma\simeq \Hom_\B(N^\sigma,M^\sigma).
$$
$(ii)$~Let $Q$ be an object of $\b2$ and $F:=\Homi_\b2(Q,-)$.	Let $w_j:C\to I$ be two extensions  of the morphism $v:C'\to I'$ for short \exxx sequences as in \eqref{weakly final abelian}. Then the restrictions to the null objects $\coker(F(c''))^\sigma$ of the two morphisms $F(w''_1), F(w''_2)$ are equal.
\end{lem}
\proof $(i)$~For any $\phi \in \Hom_\b2(N,M)$ the restriction of $\phi$ to $N^\sigma$ gives an element  $r(\phi)\in \Hom_\B(N^\sigma,M^\sigma)$. When $\phi\in \Homi_\b2(N,M)^\sigma$ one has $\phi=\phi\circ \sigma=\sigma\circ \phi$ and thus $\phi=r(\phi)\circ p$.
Moreover for any $\psi\in \Hom_\B(N^\sigma,M^\sigma)$ one has $\psi\circ p\in \Homi_\b2(N,M)^\sigma$ which shows that $r$ is an isomorphism.\newline
$(ii)$~Proposition \ref{fixed cokernel} gives canonical isomorphisms 
$$
F(C'')^\sigma\to\coker(F(c''))^\sigma, \ \ F(I'')^\sigma\to\coker(F(\alpha''))
$$
and using $(i)$ one obtains, by composition, canonical isomorphisms 
$$
\Hom_\B(Q^\sigma,C''^\sigma)\to\coker(F(c''))^\sigma, \ \ \Hom_\B(Q^\sigma,I''^\sigma)\to\coker(F(\alpha'')).
$$
Under these isomorphisms the restrictions to the null objects $\coker(F(c''))^\sigma$ of the two morphisms $F(w''_1), F(w''_2)$ are given by composition on the left with the restrictions $w''_j:C''^\sigma\to I''^\sigma$ which give maps 
$$
 \Hom_\B(Q^\sigma,C''^\sigma)\ni \phi\mapsto w''_j\circ \phi \in \Hom_\B(Q^\sigma,I''^\sigma).
$$
Thus in order to prove $(ii)$ it is enough to show that  the restrictions $w''_j:C''^\sigma\to I''^\sigma$ are equal. Since the sequence $c$ is short exact, the  term $C\stackrel{c''}{\to} C''$ is the cokernel of the first term and Proposition \ref{fixed cokernel} applies to give a canonical isomorphism $C''^\sigma\simeq C^\sigma$. The same result applies to the short \exxx sequence $i$ and gives a canonical isomorphism $I''^\sigma\simeq I^\sigma$. In the two short \exxx sequences $c$ and $i$ the left term is the kernel of the last one, and thus contains the null elements. This gives natural inclusions $C^\sigma\subset C'$ and $I^\sigma\subset I'$. Moreover by construction the maps $w_j$ both induce the same map $v:C'\to I'$ and thus the same map on the submodules 
$C^\sigma\subset C'$. This shows that  the restrictions $w''_j:C''^\sigma\to I''^\sigma$ are equal.\endproof
It follows from Lemma \ref{functoriality on fixed} that if we let $p:F(I'')\to F(I'')^\sigma$ be the projection, one has  $p\circ F(w''_1)=p\circ F(w''_2)$, since for any morphism $\psi:L\to M$ in $\b2$ the composition $p\circ \psi=\psi\circ p$ is determined by the restriction to $L^\sigma$. This equality suffices to get $ F(w''_1)(u)= F(w''_2)(u)$ when the $ F(w''_j)(u)$ belong to the image of $F(\alpha'')$, since by Proposition \ref{compute cokernel} the cokernel coincides with the projection $p$ on that image.
Thus the interesting case to consider is when for a $u\in F(C'')=\Homi_\b2(Q,C'')$ the compositions $w''_j\circ u\in \Homi_\b2(Q,I'')$ do not lift to $\Homi_\b2(Q,I)$. \vspace{.05in}

Let us take the notations of Section~\ref{ad-ex}. 
 \begin{thm}\label{sequence s} The short \exxx sequence $s:J\stackrel{}{\rightarrowtail{}} \B^3 \stackrel{\phi}{\xtwoheadrightarrow{}} S$ satisfies condition $(a)$ of \cite{gra1} (\ie \eqref{grandis condition}) with respect to the functor $H:=\Homi_\b2(S,-)$  and all endomorphisms of $J=\ker(\phi)$.
 \end{thm}
\proof By Lemma \ref{functoriality on fixed} it is enough to show that if $v\in\End_\b2(J)$  admits more than one extension to $\B^3$ then for any extension $w$ of $v$ to $\B^3$ the action of $w''$ by left multiplication on $\coker(H(\phi))$ is null. By $(iv)$ of Proposition \ref{shortexactsimple1} the restriction to null elements of $w''$ fails to be surjective and thus the same holds for any $w''\circ u$, $\forall u\in \End_\b2(S)$ which shows that the range of left multiplication by $w''$ is contained in $\imm(H(\phi))$ and is hence null in  $\coker(H(\phi))$. \endproof
The next step is to extend Theorem \ref{sequence s} to general morphisms from a short \exxx sequence $c$  as in \eqref{weakly final abelian}. We first determine  the freedom in extending the morphism $v:C'\to J$ to a morphism $w:C\to \B^3$. 
 \begin{lem}\label{sequence s1} Let $C'\stackrel{c'}{\rightarrowtail{}} C \stackrel{c''}{\xtwoheadrightarrow{}} C''$ be a short \exxx sequence and $v:C'\to J$ be a morphism. The extensions of $v$ to a morphism $w:C\to \B^3$ are uniquely determined by the element $L_w\in C^*$ defined by 
 $$
 L_w=\epsilon_a \circ w, \qquad \epsilon_a(z)=0 \iff z\leq  b+c.
 $$ 	
 \end{lem}
\proof The three elements $\epsilon_a, \epsilon_b, \epsilon_c$ generate the dual of $ \B^3$ thus $w$ is uniquely determined by the composition of these elements with $w$. But $\epsilon_c=\epsilon_a \circ \sigma$ so that $\epsilon_a \circ w$ determines $\epsilon_c \circ w=L_w\circ \sigma$. Let us show that $\epsilon_b \circ w$ is uniquely determined by $v$. One has $\epsilon_b \circ\sigma=\epsilon_b$ and thus
$\epsilon_b \circ w=\epsilon_b \circ w\circ p$ where $p(x)=x+\sigma(x)$ is the projection on $C^\sigma$. One has $C^\sigma=C'^\sigma$ and  the restriction of  $w$ to $C^\sigma$ is uniquely determined by $v$. \endproof 
Lemma \ref{sequence s1} provides a second interpretation of the proof Proposition \ref{shortexactsimple1} $(iv)$. Indeed, one considers the short \exxx sequence $C'\stackrel{c'}{\rightarrowtail{}} C \stackrel{c''}{\xtwoheadrightarrow{}} C''$ identical to $s$. The issue of extending elements of $J^*=\Hom_\B(J,\B)$ to elements of $(\B^3)^*$ is related to the map $\phi:\B^3\to S$, using the natural isomorphism $J^*\simeq S$ visible in Figures \ref{boole4} and \ref{boole5}. It follows that the only element of $J^*=\Hom_\B(J,\B)$ which does not extend uniquely to $\B^3$ is the maximal element $\tau$ which takes the value $1\in \B$ on any non-zero element of $J$.  Thus  from Lemma \ref{sequence s1} one derives that the only  $v\in\End_\b2(J)$ which admit more than one extension to $\B^3$ fulfills  $\epsilon_a \circ v=\tau$. This entails that the restriction of $v$ to $J^\sigma$ cannot be surjective.

\subsection{Condition (a) and the radical}\label{sect a and radical}
 In this subsection we investigate  properties of general morphisms of short \exxx sequences.
We keep the notations of Lemma \ref{functoriality on fixed}.
Let $Q$ be an object of $\b2$ and $F:=\Homi_\b2(Q,-)$.	Let $w_j:C\to I$ be two extensions  of the morphism $v:C'\to I'$ for short \exxx sequences as in \eqref{weakly final abelian}.  It follows from Lemma \ref{functoriality on fixed} that if we let $p:F(I'')\to F(I'')^\sigma$ be the projection, one has  $p\circ F(w''_1)=p\circ F(w''_2)$. We now want to apply Theorem \ref{kernel and injective} to conclude that 
$F(w''_1)=F(w''_2): \coker(F(c''))\to \coker(F(\alpha''))$, where 
$
\coker(F(c''))=\Hom(Q,C'')/(c''\circ \Hom(Q,C))$,  
$\coker(F(\alpha''))=\Hom(Q,I'')/(\alpha''\circ \Hom(Q,I))$. In order to prove this we construct an object $E$ of $\b2$ together with two morphisms 
$$
\Homi_\b2(C',I')\stackrel{\rho}{\longleftarrow} E\stackrel{\theta}{\longrightarrow}
\Homi_\b2(\coker(F(c'')),\coker(F(\alpha''))).
$$
Thus an element of $E$ is an equivalence class of elements $w\in \Homi_\b2(C,I)$ compatible with the short \exxx sequences, \ie such that $w(\ker(c''))\subset \ker(\alpha'')$. We let $\rho(w)$ be the restriction $\rho(w)\in \Homi_\b2(C',I')$ and we let $\theta(w)$ be the action $F(w''): \coker(F(c''))\to \coker(F(\alpha''))$. We set
$$
w_1\sim w_2\iff \rho(w_1)=\rho(w_2)\ \& \ \theta(w_1)=\theta(w_2).
$$
This means that the quotient $E$ is the range of the map 
\begin{equation}\label{map rho theta}
	(\rho,\theta):\Homi_\b2^{\rm comp}(C,I)\to \Homi_\b2(C',I')\times \Homi_\b2(\coker(F(c'')),\coker(F(\alpha'')))
\end{equation}
 where the upper index in $\Homi_\b2^{\rm comp}(C,I)$ refers to the compatibility with the short \exxx sequences. One has the following general fact
 \begin{lem}\label{kernel of rho} With the above notations, the kernel of the map $\rho:E\to 
 	\Homi_\b2(C',I')$ is null.
 \end{lem}
 \proof We need to show that if $w\in\Homi_\b2^{\rm comp}(C,I)$ fulfills $\rho(w)$ null, then $\theta(w)$ is also null. Assume $\rho(w)$ null, then $w\circ c'=\alpha'\circ \rho(w)$ is null. Thus since $c$ is a short \exxx sequence there exists a unique induced morphism $u\in \Homi_\b2(C'',I)$ such that $w=u\circ c''$. One has $w''=\alpha''\circ u$. Let us consider the diagram, with $F:=\Homi_\b2(Q,-)$,
  \begin{equation}\label{applying F again}
\xymatrix@C=20pt@R=35pt{
 F(C)\ar[d]_{F(w)} \ar[rr]^{F(c'')}&& F(C'')\ar[lld]_{u\circ - }\ar[d]^{w''\circ - }\ar[rr]^{\!\!\!\scoker(F(c''))}&& \ \coker(F(c''))\ar[d]^{\theta(w)}
\\
 F(I) \ar[rr]^{F(\alpha'')}&& F(I'') \ar[rr]^{\!\!\!\!\scoker(F(\alpha''))}&&\ \coker(F(\alpha''))}\end{equation}
 The middle vertical arrow is given by left composition with $w''$, \ie by $t\mapsto w''\circ t$ for $t\in \Homi_\b2(Q,C'')$. Since $w''=\alpha''\circ u$ one has 
 $
 w''\circ t=\alpha''\circ u\circ t
 $
which gives a factorization of the middle vertical arrow as $F(\alpha'')\circ F(u)$. This shows that the range of the middle vertical arrow is contained in the range of $F(\alpha'')$ and hence in the kernel of $\scoker(F(\alpha''))$. Using the commutativity of the diagram \eqref{applying F again} this proves that $\theta(w)$ is null. \endproof 
We now involve Definition \ref{sigma inject} and obtain
 \begin{lem}\label{sigma inject hyp}
 $(i)$~With the above notations, the restriction of the map $\rho:E\to 
 	\Homi_\b2(C',I')$ to null elements is injective.\newline
 	$(ii)$~If the map $\rho:E\to 
 	\Homi_\b2(C',I')$ is $\sigma$-injective then for any two extensions $w_j$ of the morphism $v:C'\to I'$, one has $F(w''_1)=F(w''_2): \coker(F(c''))\to \coker(F(\alpha''))$.
 \end{lem}
\proof $(i)$~By construction $E$ is a sub-object of the product which is the right hand side of \eqref{map rho theta}, thus $E^\sigma$ is a submodule of  $\Homi_\b2(C',I')^\sigma\times \Homi_\b2(\coker(F(c'')),\coker(F(\alpha'')))^\sigma$ which by Lemma \ref{functoriality on fixed} is $\Hom_\B(C'^\sigma,I'^\sigma)\times \Hom_\B(\coker(F(c''))^\sigma,\coker(F(\alpha''))^\sigma)$. Lemma \ref{functoriality on fixed} gives moreover natural isomorphisms $C'^\sigma=C''^\sigma$,  
$\coker(F(c''))^\sigma\simeq \Hom_\B(Q^\sigma,C''^\sigma)$, $I'^\sigma=I''^\sigma$, and $\coker(F(\alpha''))^\sigma \simeq \Hom_\B(Q^\sigma,I''^\sigma)$. Thus $E^\sigma$ is a submodule of 
$$
\Hom_\B(C''^\sigma,I''^\sigma)\times \Hom_\B(\Hom_\B(Q^\sigma,C''^\sigma),\Hom_\B(Q^\sigma,I''^\sigma)).
$$
In fact Lemma \ref{functoriality on fixed} shows that $E^\sigma=\{(v,L(v))\mid v\in \Hom_\B(C''^\sigma,I''^\sigma)\} $, where $L(v)$ denotes the left composition with $v$, \ie the map $L(v)(u):= v\circ u\in \Hom_\B(Q^\sigma,I''^\sigma)$, $\forall u \in \Hom_\B(Q^\sigma,C''^\sigma)$. The restriction of $\rho:E\to 
 	\Homi_\b2(C',I')$ to $E^\sigma$ is, under the identification $\Homi_\b2(C',I')^\sigma\simeq \Hom_\B(C''^\sigma,I''^\sigma)$, the projection $(v,L(v))\mapsto v$ and is thus injective. 
 	\newline
 	$(ii)$~Since the map $\rho:E\to 
 	\Homi_\b2(C',I')$ is assumed to be $\sigma$-injective it follows from $(i)$ that it is injective.
 	We now show that this implies  $F(w''_1)=F(w''_2): \coker(F(c''))\to \coker(F(\alpha''))$. The extensions $w_j$ of $v$ define elements $w_j\in \Homi_\b2^{\rm comp}(C,I)$ and, using \eqref{map rho theta}, we let $\xi_j=(\rho,\theta)(w_j)\in E$. Since the map $\rho:E\to 
 	\Homi_\b2(C',I')$ is injective, and since $\rho(\xi_1)=v=\rho(\xi_2)$ one has $\xi_1=\xi_2$ and hence $\theta(\xi_1)=\theta(\xi_2)$ which implies  $F(w''_1)=F(w''_2): \coker(F(c''))\to \coker(F(\alpha''))$.\endproof 
 	\begin{cor}\label{functoriality on fixed applied} $(i)$~With the above notations,  for any two extensions $w_j$ of the morphism $v:C'\to I'$, the morphisms $F(w''_j): \coker(F(c''))\to \coker(F(\alpha''))$ fulfill $p(F(w''_1))=p(F(w''_2))$.\newline
 	$(ii)$~If the short \exxx sequence  $c:~C'\stackrel{c'}{\rightarrowtail{}} C \stackrel{c''}{\xtwoheadrightarrow{}} C''$ is split then the morphisms $F(w''_j): \coker(F(c''))\to \coker(F(\alpha''))$ are null and equal. 		
 	\end{cor}
 	\proof $(i)$~The extensions $w_j$ of $v$ define elements $w_j\in \Homi_\b2^{\rm comp}(C,I)$ and, as above, using \eqref{map rho theta}, we let $\xi_j=(\rho,\theta)(w_j)\in E$. By Lemma \ref{sigma inject hyp} the restriction of the map $\rho:E\to 
 	\Homi_\b2(C',I')$ to null elements is injective. Thus, since $\rho(\xi_1)=v=\rho(\xi_2)$, one has $\rho(p(\xi_1))=p(v)=\rho(p(\xi_2))$ and
 	$p(\xi_1)=p(\xi_2)$. Hence $\theta(p(\xi_1))=\theta(p(\xi_2))$ which implies  $p(F(w''_1))=p(F(w''_2))$.\newline
 	$(ii)$~It is enough to show that the object $\coker(F(c''))$ is null since this implies that the morphisms $F(w''_j)$ are null and hence equal by $(i)$. By hypothesis the short \exxx sequence $c$ is split, thus there exists a section $s:C''\to C$ such that $c''\circ s=\id$. It follows that $F(c'')\circ F(s)=\id$ and  the map  $F(c'')$ is surjective. Thus its cokernel $\coker(F(c''))$ is null. \endproof 

 Lemma \ref{kernel of rho} shows that the kernel of the morphism $\rho:E\to 
 	\Homi_\b2(C',I')$ is null.	Thus if one could apply  Theorem \ref{kernel and injective} one would conclude that $\rho$ is $\sigma$-injective, and hence by Lemma \ref{sigma inject hyp} that the condition $(a)$ of \cite{gra1} holds. Instead, we apply Theorem \ref{kernel and injective strong} and obtain the following
 	\begin{prop}\label{sigma inject general} With the above notations,  for any two extensions $w_j$ of  $v:C'\to I'$, the morphisms $F(w''_j): \coker(F(c''))\to \coker(F(\alpha''))$ fulfill the relation
 	\begin{equation}\label{pseudo equal}
 		F(w''_1)(u)+\sigma (F(w''_2)(u))\ \text{null}\qqq u\in \coker(F(c'')).
 	\end{equation}
 \end{prop}
	\proof We keep the same notations as in the proof of Lemma \ref{sigma inject hyp}. Thus we let $\xi_j=(\rho,\theta)(w_j)\in E$. We apply Theorem \ref{kernel and injective strong} to $\rho:E\to 
 	\Homi_\b2(C',I')$. Lemma \ref{kernel of rho} shows that the kernel of $\rho$ is null, thus since $\rho(\xi_1)=\rho(\xi_2)$, the sum $\xi_1+\sigma(\xi_2)$ is null. It follows that $\theta(\xi_1)+\sigma(\theta(\xi_2))$ is null and hence $F(w''_1)+\sigma (F(w''_2))$ is null which gives \eqref{pseudo equal}. \endproof 
 	 Figure \ref{graphcokernel} shows that for $\coker(F(\alpha''))$ there are two non-trivial fibers of the projection $p$ to the sixteen null elements. They are the fibers of $(\alpha+\gamma,\beta+\delta)$ and $(\beta+\delta,\alpha+\gamma)$. The fiber of $(\alpha+\gamma,\beta+\delta)$ is given by
 	 $$
 	p^{-1}((\alpha+\gamma,\beta+\delta))= \{(\alpha,\beta),(\gamma,\delta),(\alpha+\gamma,\beta),(\alpha+\gamma,\delta),(\alpha,\beta+\delta),(\gamma,\beta+\delta),(\alpha+\gamma,\beta+\delta)\}.
 	 $$
 	  Inside these fibers one can find non-trivial maximal submonoids $M$ such that
 	 $z_1+\sigma(z_2)\ \text{null}, \forall z_j\in M$. For instance, one can take for the fiber of $(\alpha+\gamma,\beta+\delta)$
 	 $$
 	 M_1=\{(\alpha,\beta),(\alpha+\gamma,\beta),(\alpha,\beta+\delta),(\alpha+\gamma,\beta+\delta)\}
 	 $$ 
 	 or 
 	 $$
 	 M_2=\{(\gamma,\delta),(\alpha+\gamma,\delta),(\gamma,\beta+\delta),(\alpha+\gamma,\beta+\delta)\}.
 	 $$ 
 	 This shows that the information one gets from Proposition \ref{sigma inject general}
	is not sufficient to get condition $(a)$.

	By applying Theorem \ref{kernel and injective bis} one obtains the following information on the injectivity of the  morphism  $\rho:E\to 
 	\Homi_\b2(C',I')$. 
	\begin{lem}\label{inject bis} With the above notations, for any two extensions $w_j$ of the morphism $v:C'\to I'$ one has
	$
	(v,F(w''_1))\sim_{\rad(E)}(v,F(w''_2))
	$, where $F(w''_j): \coker(F(c''))\to \coker(F(\alpha''))$.
 \end{lem}
	\proof By construction $E$ is a sub-object of the product which is the right hand side of \eqref{map rho theta}. Let $\xi_j=(v,F(w''_j))\in E$, one has $\rho(\xi_1)=\rho(\xi_2)$. The proof of Lemma \ref{sigma inject hyp} shows that the restriction of $\rho$ to the null elements $E^\sigma$ is injective. Thus one has $p(\xi_1)=p(\xi_2)$. By Lemma \ref{kernel of rho}, the kernel of the map $\rho:E\to 
 	\Homi_\b2(C',I')$ is null. Thus Theorem \ref{kernel and injective bis} applies and one gets the required equivalence $
	(v,F(w''_1))\sim_{\rad(E)}(v,F(w''_2))
	$.\endproof 

\subsection{Reduction of condition (a) to endomorphisms}\label{sect reduction to endo}

We consider a short \exxx sequence $\iota:I'\stackrel{i'}{\rightarrowtail{}} I \stackrel{i''}{\xtwoheadrightarrow{}} I''$ where we now make the following ``Frobenius\footnote{By reference to the notion of Frobenius algebra, \ie of algebras for which projective modules are the same as injective modules. Note that in $\b2$ this property holds for finite objects.} hypothesis" that the middle term $I$ is an object of $\b2$ which is {\em injective and projective}. 
We consider the functor $F:=\Homi_\b2(I'',-)$ and we use the fact that the object $I$ is projective as follows 
 \begin{lem}\label{sequence s2}
 With the above notations, let 	$C'\stackrel{c'}{\rightarrowtail{}} C \stackrel{c''}{\xtwoheadrightarrow{}} C''$ be a short \exxx sequence  and let $x\in \coker(F(c''))$. Then there exists a morphism of short \exxx sequences $\xi:\iota\to c$ such that 
 \begin{equation}\label{gc trick}
F(\xi'')(\scoker(F(i'')(\id_{I''}))=x,\ \  F(\xi''): \coker(F(i''))\to \coker(F(c'')).	
\end{equation}
 \end{lem}
\proof Since $\scoker(F(c''))$ is surjective let $z\in F(C'')$ be such that $\scoker(F(c''))(z)=x$. One has  $z\in F(C'')=\Homi_\b2(I'',C'')$. Since $I$ is projective and the morphism $c''$ is surjective one can lift the morphism $z\circ i'':I\to C''$ to a morphism $\xi:I\to C$ such that $c''\circ \xi= z\circ i''$. For $u \in \ker(i'')$ one has that $c''\circ \xi(u)$ is null and thus $\xi(u)\in \ker(c'')$ so that $\xi$ restricts to a morphism $\xi':I'\to C'$ and with $\xi'':=z$ one obtains 
a morphism of short \exxx sequences in $\b2$ of the form
\begin{equation}\label{use of proj}
\xymatrix@C=20pt@R=35pt{
  I'\ar[d]^{\xi'}\ar@{>->}[rr]^{i'} && I \ar[d]^{\xi} \ar@{->>}[rr]^{i''}&& I''\ar[d]^{\xi''=z}
\\
C' \ \ar@{>->}[rr]^{c'}&& C \ar@{->>}[rr]^{c''}&& C''}
\end{equation}
This gives a commutative square using the functor $F=\Homi_\b2(I'',-)$
\begin{equation}\label{applying F to sequ}
\xymatrix@C=20pt@R=35pt{
 F(I)\ar[d]_{F(\xi)} \ar[rr]^{\!\!\!\!\!\!\!\!\!\!\!\! F(i'')}&& \Homi_\b2(I'',I'')\ar[d]^{\xi''\circ - }\ar[rr]^{\;\;\scoker(F(i''))}&& \coker(F(i''))\ar@{-->}[d]
\\
 F(C) \ar[rr]^{\!\!\!\!\!\!\!\!\!\!\!\!  F(c'')}&& \Homi_\b2(I'',C') \ar[rr]^{\;\;\scoker(F(c''))}&& \coker(F(c''))}\end{equation}
 The commutativity of the diagram given by the left square in \eqref{applying F to sequ} insures that the dotted vertical arrow $F(\xi''): \coker(F(i''))\to \coker(F(c''))$ exists and that the diagram given by the right square is commutative. Moreover the image of $\scoker(F(i'')(\id_{I''})$ by $F(\xi'')$  is then the same as the  image by $\scoker(F(c''))$ of $\xi''\circ \id_{I''}=z$ and is thus equal to $x$ which gives \eqref{gc trick}.\endproof 
 \begin{thm}\label{reduction toendo} Let $\iota:I'\stackrel{i'}{\rightarrowtail{}} I \stackrel{i''}{\xtwoheadrightarrow{}} I''$ be a short \exxx sequence such that the middle term $I$ is an object of $\b2$ which is {\em injective and projective}. 
Then the functor $F:=\Homi_\b2(I'',-)$ satisfies condition $(a)$  with respect to an arbitrary morphism of short \exxx sequences $c\to \iota$, provided this holds for all  endomorphisms of $I'$.
 \end{thm}
 \proof Let $c:C'\stackrel{c'}{\rightarrowtail{}} C \stackrel{c''}{\xtwoheadrightarrow{}} C''$ be a short \exxx sequence, $v:C'\to I'$ be a morphism, and $w_j$ be extensions of $v$ to morphisms of short \exxx sequences from $c$ to $\iota$. We need to show the equality 
 \begin{equation}\label{grandis condition1}
F(w''_1)=F(w''_2): \coker(F(c''))\to \coker(F(i'')).	
\end{equation} 
Let $x\in \coker(F(c''))$. Then by Lemma \ref{sequence s2},  there exists a morphism of short \exxx sequences $\xi:\iota\to c$ such that \eqref{gc trick} holds, \ie 
$F(\xi'')(\scoker(F(i'')(\id_{I''}))=x$. We now consider the composition $w_j\circ \xi$  of the morphism $\xi:\iota\to c$ with either of the extensions  $w_j$  of $v$ to morphisms of short \exxx sequences from $c$ to $\iota$. These provide two extensions of the endomorphism $v\circ \xi': I'\to I'$. 
 \begin{equation}\label{combinazione}
\xymatrix@C=20pt@R=35pt{
  I'\ar[d]^{\xi'}\ar@{>->}[rr]^{i'} && I \ar[d]^{\xi} \ar@{->>}[rr]^{i''}&& I''\ar[d]^{\xi''=z}
\\
C'\ar[d]^{v}\ar@{>->}[rr]^{c'} && C \dar[d]_{w_1}^{\ \, w_2} \ar@{->>}[rr]^{c''}&& C''\dar[d]_{w''_1}^{\ \, w''_2}
\\
I' \ \ar@{>->}[rr]^{i'}&& I \ar@{->>}[rr]^{i''}&& I''}
\end{equation}
 Since by hypothesis the short \exxx sequence $\iota$ satisfies condition $(a)$  with respect to  all  endomorphisms of $I'$, we get the equality $F(w''_1\circ \xi'')=F(w''_2\circ \xi''): \coker(F(i''))\to \coker(F(i''))$. We apply this equality to the element $\scoker(F(i'')(\id_{I''})\in \coker(F(i''))$. By \eqref{gc trick} we derive $F(\xi'')(\scoker(F(i'')(\id_{I''}))=x$ and hence we get $F(w''_1)(x)=F(w''_2)(x)$ so that \eqref{grandis condition1}
 holds. \endproof 
 \begin{thm}\label{reduction toendo cor1} Let $S$ be as in Theorem \ref{sequence s}. The satellite functor $SH$ of  $H:=\Homi_\b2(S,-)$ is non-null and $SH(J)$  is the cokernel of the morphism $H(\phi):H(\B^3)\to H(S)$ of Proposition \ref{shortexactsimple}, $(iv)$.\end{thm}
\proof  Let $\sequ_{\rm small}$ be any small subcategory of the category $\sequ(\b2)$ of short \exxx sequences in the homological category $\b2$, with $\sequ_{\rm small}$  large enough to contain all short \exxx sequences of finite objects. By construction  $H:\b2\longrightarrow \b2$ is a covariant functor and $\b2$ is  a cocomplete  homological category. By \eqref{Kan extension} the following colimit makes sense for any object $X$ of $\b2$, and defines the  covariant functor $SH$
 \begin{equation}\label{Kan extension1}
	SH:\b2\longrightarrow \b2, \ \ SH(X):=\varinjlim_I \coker(H(a'')), \ \ I=\sequ_{\rm small}(\b2)\downarrow_{P'} X
\end{equation} 
We proceed as in the proof of Theorem 4.2.2 of \cite{gra1} and show that the object $$j=\left(J\stackrel{\id}{\leftarrow} J\stackrel{\subset}{\rightarrowtail}\B^3\stackrel{\phi}{\xtwoheadrightarrow{}}S\right)$$
is a  final object for the functor $H$ in the comma  category $I=\sequ_{\rm small}(\b2)\downarrow_{P'} J$. Indeed, since  $\B^3$ is injective, for any object $J\stackrel{v}{\leftarrow} C'\stackrel{c'}{\rightarrowtail}C\stackrel{c''}{\xtwoheadrightarrow{}}C''$ of $I$ one gets 
  a morphism in $I$
\begin{equation}\label{colimit compatib}
\xymatrix@C=20pt@R=35pt{
&  C' \ar[ld]_{v}\ar[d]_{v}\ar@{>->}[rr]^{c'} && C\ar[d]_{w} \ar@{->>}[rr]^{c''}&& C''\ar[d]_{w''}
\\
J & \ar[l]_{\id}J \ \ar@{>->}[rr]^{\subset}&& \B^3 \ar@{->>}[rr]^{\phi}&& S }
\end{equation}
where the choice of $w$ is not unique in general but where, by Theorem \ref{reduction toendo} combined with Theorem \ref{sequence s}, the induced action $H(w''): \coker(H(c''))\to \coker(H(\phi))$ is unique. It follows that the colimit \eqref{Kan extension1} is simply the evaluation on the final object, \ie is given by $\coker(H(\phi))$ and Proposition \ref{shortexactsimple}, $(iv)$ shows that it is non-null. \endproof 

 \section{The cokernel of the diagonal}\label{sectcokerdiag}
 
 In Proposition \ref{cokernelequiv} $(ii)$, we have shown  that the functor $F$ associated to the  cokernel of a morphism of $\B$-modules $f:L\to M$ becomes representable \ie of the form $\Homi_\b2(Q,-)$ in the category $\b2$. 
 In this section we prove that the satellite functor Ext of the functor $F$ is non null for the simplest natural example where the morphism $f:L\to M$ is the diagonal $\Delta: \B\to \B\times \B$. 
  
  The cokernel of the diagonal enters in the construction of the \v Cech version of sheaf cohomology. Given a constant sheaf $\cF$ with fiber an abelian group $H$ on a topological space $X$, and an open cover $\cU$ by open sets $U_j$ one has an exact sequence of sheaves of abelian groups of the form
\begin{equation}\label{abelian exact}
	\cF\stackrel{a'}{\rightarrowtail} \prod \cF\vert_{U_j}\stackrel{a''}{\xtwoheadrightarrow{}}\cF''
\end{equation}
where the map $a'$ is obtained from the restriction maps $\cF(U)\to \cF(U\cap U_j)$ and  the last arrow $a''$ is the cokernel of $a'$. This dictates the construction of the sheaf $\cF''$ as the cokernel of $a'$. At a point $s$ belonging to two open sets the stalk of $\prod \cF\vert_{U_j}$  is $H\times H$ and the map $a'$ is the diagonal $h\mapsto (h,h)$. Thus the stalk of $\cF''$ at $s$ is the cokernel of the diagonal and hence gets identified with  $H$. More precisely, in this abelian case the short \exxx sequence 
\begin{equation}\label{abelian exact bis}
	H\stackrel{\Delta}{\rightarrowtail} H\times H\stackrel{\scoker(\Delta)}{\xtwoheadrightarrow{}}H, \  \ \Delta(h):=(h,h)
\end{equation}
splits and  the identification of the cokernel with $H$ can be done using the difference, \ie the map $(h_1,h_2)\mapsto h_2-h_1$. This becomes instrumental in defining the coboundary in the \v Cech complex associated to the open cover.
In order to treat the more general case where the abelian group $H$ is replaced by a $\B$-module $L$ we  investigate the cokernel of the diagonal map $\Delta:L\to L\times L$, $h\mapsto (h,h)$  which appears in \eqref{abelian exact bis}. We consider the simplest case \ie $L=\B$. \vspace{.05in}

The steps which allow us to compute the satellite functor $SF$ of the cokernel of the diagonal are the following:
\begin{enumerate}
\item In \S \ref{sectrepfunctor} we determine the object $Q=\coker(\yb\Delta)$ of $\b2$ representing the cokernel of the diagonal $\Delta:\B\to \B\times \B$.
\item In Lemma \ref{condition to lift} we show  that, with $F=\Homi_\b2(Q,-)$, the cokernel $\coker(F(\scoker(\yb\Delta)))$ is non-null.
\item  In \S \ref{new sect duality} we determine  $\coker(F(\scoker(\yb\Delta)))$ using duality.
\item  In \S \ref{sect comma cat} we analyze the elements of the comma category $\cI=\sequ(\b2)\downarrow_{P'} K$ associated to endomorphisms of the kernel of $\scoker(\yb\Delta)$.
\item The correspondence underlying the multiple extensions of endomorphisms of the kernel is discussed in \S \ref{subsubcorrespondence} as a preparation for the proof of property (a) for endomorphisms.
\end{enumerate}
The main result which computes the satellite functor $SF$ of the functor $F:=\Homi_\b2(Q,-)$ as the cokernel $\coker(F(\alpha''))$ is obtained in \S \ref{subsecfunctact}, Theorem \ref{reduction toendo cor2}.

\subsection{The cokernel pair of $\Delta:\B\to \B\times \B$ as a representable functor} \label{sectrepfunctor} 
The cokernel pair of the diagonal $\Delta:\B\to \B\times \B$ is first comprehended as the covariant functor \eqref{covfunc} of Lemma \ref{repfunctor}
$$
F(X)=	\B^2/\B(X):=\{(f,g)\in \Hom_{\bm2}(\B^2,X)\mid f(x)=g(x)\qqq x\in \Delta\}.
$$
The pair $(f,g)\in \Hom_{\bm2}(\B^2,X)$ is characterized by the $4$ elements of $X$
$$
a=f((1,0)), \ b=f((0,1)), \ c=g((1,0)), \ d=g((0,1))
$$
and the only constraint is that $a+b=c+d$. Given a morphism $(h,k)\in \Hom_\bm2(X,X')$, one obtains the morphism $(f',g')\in \Hom_{\bm2}(\B^2,X')$ corresponding to 
$$
a'=h(a)+k(c),\ b'=h(b)+k(d),\ c'=h(c)+k(a), \ d'=h(d)+k(b).
$$
By Proposition \ref{repb2} the functor $F$ is represented in $\b2$ and the object $\coker(\yb\Delta)$ representing $F$ is the quotient $Q$ of $\B^2\times \B^2$ by the equivalence relation \eqref{represent functor}. We view the elements of $\B^2\times \B^2$ as sums over the subsets of the set $\{\alpha ,\beta ,\gamma ,\delta \}$. One finds that the first ten elements of the list 
$$
0,\alpha ,\beta ,\gamma ,\delta ,\alpha +\beta ,\alpha +\gamma ,\alpha +\delta ,\beta +\gamma ,\beta +\delta,$$ $$ \gamma +\delta ,\alpha +\beta +\gamma ,\alpha +\beta +\delta ,\alpha +\gamma +\delta ,\beta +\gamma +\delta ,\alpha +\beta +\gamma +\delta 
$$
represent all the elements in the quotient, while all the others project on $\alpha +\beta\sim \gamma +\delta$.
\begin{figure}[H]
\begin{center}
\includegraphics[scale=1]{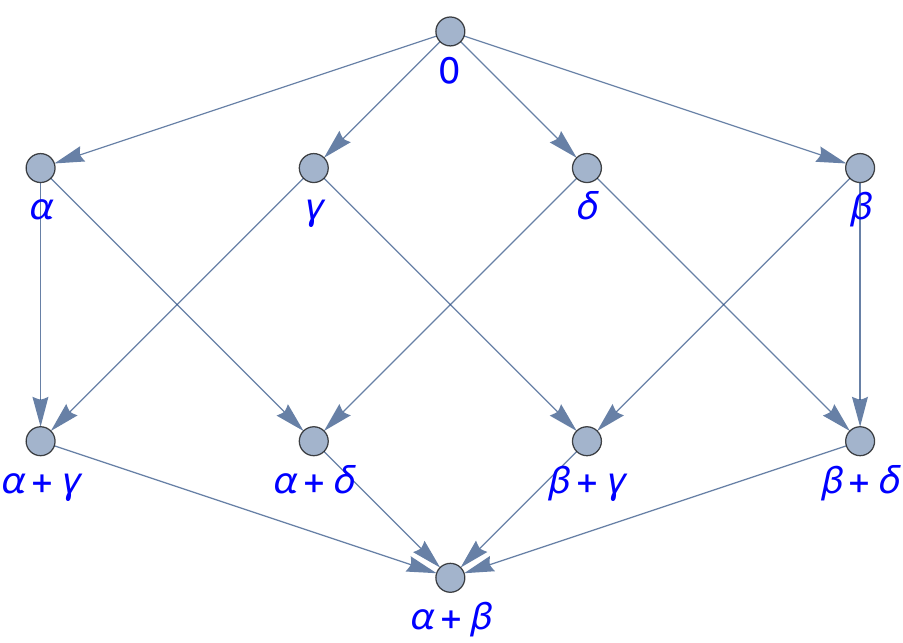}
\end{center}
\caption{The structure of the $\B$-module $Q$ \label{newnormsubmodn} }
\end{figure}

The involution $\sigma$ is given by  $\sigma(\alpha)=\gamma$, $\sigma(\beta)=\delta$. Its fixed points in $Q$ are 
$
0,\alpha +\beta,\alpha +\gamma,\beta +\delta
$
and the preimages of the fixed points are the following $10$ elements of $\B^2\times \B^2$
$$
0,\alpha +\beta,\alpha +\gamma,\beta +\delta,\gamma +\delta ,\alpha +\beta +\gamma ,\alpha +\beta +\delta ,\alpha +\gamma +\delta ,\beta +\gamma +\delta ,\alpha +\beta +\gamma +\delta 
$$
which form the kernel $K:=\ker(\scoker(\yb\Delta))$. The intersection of this kernel with the sums not invoking $\gamma$ or $\delta$ is as expected the diagonal $\{0,\alpha +\beta\}\subset \B^2$. By Proposition \ref{cokernelequiv} $(ii)$, one derives the following short \exxx sequence 
\begin{equation}\label{short diag exact}
	K\stackrel{\sker(\scoker(\yb\Delta))}{\rightarrowtail} \B^2\times \B^2\stackrel{\scoker(\yb\Delta)}{\xtwoheadrightarrow{}}Q.
\end{equation}
We want to detect the obstruction to lift back from $Q$ to $\B^2\times \B^2$. Thus we  consider the represented functor $F:=\Homi_\b2(Q,-)$ which we view as a covariant endofunctor $\b2\longrightarrow\b2$ (using the natural internal Hom). We need to understand 
$$
F(\B^2\times \B^2)=\Homi_\b2(Q,\B^2\times \B^2).
$$
Any given $\phi \in \Homi_\b2(Q,\B^2\times \B^2)$ is uniquely determined by $a=\phi(\alpha)$ and $b=\phi(\beta)$ which must fulfill the condition that $a+b$ is fixed under $\sigma$.
There are only $4$ elements of $\B^2\times \B^2$ which are invariant under $\sigma$ and they are: 
$
0, \ \alpha+\gamma, \ \beta+\delta, \ \alpha +\beta +\gamma +\delta 
$.
\begin{lem}\label{cokernel of diag} The object $\Homi_\b2(Q,\B^2\times \B^2)$ of $\b2$ is the square of $Q^*=\Homi_\b2(Q,\yb\B)$.
\end{lem}
\proof As an object of $\b2$, $\B^2\times \B^2$ is the product of the two copies of $\yb\B$ given by 
$$
\yb\B\sim \{0,\alpha ,\gamma,\alpha +\gamma\},\ \sigma(\alpha)=\gamma,\ \ 
\yb\B'\sim \{0,\beta ,\delta,\beta +\delta\},\ \sigma(\beta)=\delta.
$$
\endproof
By Lemma \ref{forgetful} one has a canonical isomorphism in $\bmod$ 
\begin{equation}\label{canonicalisobis}
\pi:\Hom_\b2(Q,\yb\B)\to \Hom_\B(I(Q),\B).
\end{equation}
A morphism $\phi\in \Hom_\B(I(Q),\B)$ is given by $4$ elements $a,b,c,d\in \B$ such that $a\vee b=c\vee d$. Either all are $0$ or $a\vee b=c\vee d=1$, and this condition has $9$ solutions. One thus gets the ten possibilities given by the lines of the matrix

$$
\left(
\begin{array}{cccc}
 0 & 0 & 0 & 0 \\
 0 & 1 & 0 & 1 \\
 0 & 1 & 1 & 0 \\
 0 & 1 & 1 & 1 \\
 1 & 0 & 0 & 1 \\
 1 & 0 & 1 & 0 \\
 1 & 0 & 1 & 1 \\
 1 & 1 & 0 & 1 \\
 1 & 1 & 1 & 0 \\
 1 & 1 & 1 & 1 \\
\end{array}
\right)
$$
\begin{figure}[H]
\begin{center}
\includegraphics[scale=1]{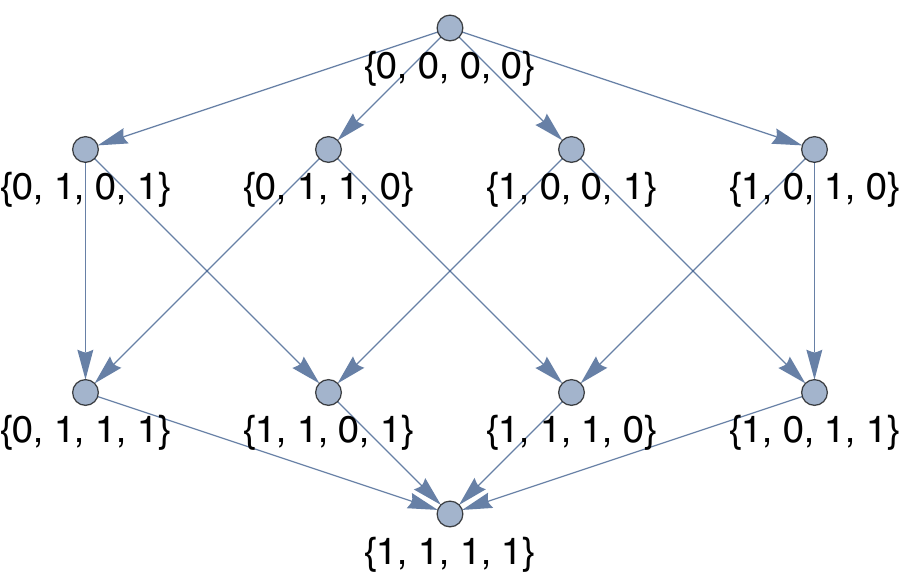}
\end{center}
\caption{The structure of the $\B$-module $Q^*$ \label{bmodtstar} }
\end{figure}
For an object $X$ of $\b2$,  $\phi\in\Hom_\b2(Q,X)$ is specified by the pair $(\phi(\alpha),\phi(\beta))\in X\times X$, with the only constraint that $\phi(\alpha)+\phi(\beta)$ is fixed under $\sigma$. For $X=\yb\B$ we thus specify two elements $\phi(\alpha),\phi(\beta)\in \yb\B$. The map $\pi$ of \eqref{canonicalisobis} is such that
$$
\pi(\phi)=(a,b,c,d)\iff \phi(\alpha)=(a,c),\ \ \phi(\beta)=(b,d).
$$

The involution on $\Homi_\bmod(I(Q),\B)$ coming from the internal Hom structure of $\Homi_\b2(Q,\yb\B)$ by the isomorphism \eqref{canonicalisobis}, is given by composition with $\sigma_Q$ and hence by the action of the permutation $(a,b,c,d)\mapsto (c,d,a,b)$ on the labels. Thus while the $\B$-modules underlying $Q$ and $Q^*$ are isomorphic, the objects $Q$ and $Q^*$ of $\b2$ are distinct since the involutions are not the same (the two elements $(1,0,1,0)$ and $(0,1,0,1)$ are fixed under the involution).\vspace{.05in}

Next, we investigate the map $F(\scoker(\yb\Delta))$ obtained by composition with $\scoker(\yb\Delta)$
$$
F(\scoker(\yb\Delta)): \Homi_\b2(Q,\B^2\times \B^2)\to \Endi_\b2(Q),\  \   \phi\mapsto \scoker(\yb\Delta)\circ \phi=F(\scoker(\yb\Delta))(\phi).
$$
We know that this map cannot be surjective since the identity $\id\in \Endi_\b2(Q)$ cannot be lifted. We view $F(\scoker(\yb\Delta))$ as a morphism in $\b2$ and want to compute its cokernel. An element $\psi\in \Endi_\b2(Q)$ is encoded by the values $(\psi(\alpha),\psi(\beta))\in Q\times Q$ and the only constraint on the pair is that $\psi(\alpha)+\psi(\beta)$ is $\sigma$-invariant \ie belongs to $
\{0,\alpha +\beta,\alpha +\gamma,\beta +\delta\}
$. One obtains 70 solutions which are listed below. The list $L_{\alpha,\gamma}$ of nine solutions with $\psi(\alpha)+\psi(\beta)=\alpha+\gamma$ is
$$
((0,\alpha+\gamma),(\alpha,\gamma),(\alpha,\alpha+\gamma),(\gamma,\alpha),(\gamma,\alpha+\gamma),(\alpha+\gamma,0),(\alpha+\gamma,\alpha),(\alpha+\gamma,\gamma),(\alpha+\gamma,\alpha+\gamma))
$$
The list $L_{\beta,\delta}$ of nine solutions with  $\psi(\alpha)+\psi(\beta)=\beta +\delta$ is
$$
((0,\beta+\delta),(\beta,\delta),(\beta,\beta+\delta),(\delta,\beta),(\delta,\beta+\delta),(\beta+\delta,0),(\beta+\delta,\beta),(\beta+\delta,\delta),(\beta+\delta,\beta+\delta))
$$
\begin{lem}\label{condition to lift0}  $\psi\in \End_\b2(Q)$ can be lifted to $\phi \in  \Hom_\b2(Q,\B^2\times \B^2)$ if and only if the pair $(\psi(\alpha),\psi(\beta))\in Q\times Q$ is a sum of two elements of the  lists $L_{\alpha,\gamma} \cup\{(0,0)\}$ and  $L_{\beta,\delta} \cup\{(0,0)\}$.
\end{lem}
\proof Notice first that each element $\psi$ of the two lists is uniquely liftable to $\phi \in  \Hom_\b2(Q,\B^2\times \B^2)$, where $\phi(\alpha)=\psi(\alpha)$ and $\phi(\beta)=\psi(\beta)$. This follows from the $\sigma$-invariance in $\B^2\times \B^2$ of the elements $\alpha+\gamma$ and $\beta+\delta$. Thus any element of $\End_\b2(Q)$ which is a sum of elements of the  lists $L_{\alpha,\gamma} \cup\{(0,0)\}$ and  $L_{\beta,\delta} \cup\{(0,0)\}$ can be lifted. Conversely, Lemma \ref{cokernel of diag} shows that any element of $\Hom_\b2(Q,\B^2\times \B^2)$ belongs to the sums of the two lists and thus only these sums in 
$\End_\b2(Q)$ are liftable.
\endproof 
The elements of $\End_\b2(Q)$ which do not belong to $L_{\alpha,\gamma} \cup L_{\beta,\delta} \cup\{(0,0)\}$ are those coming 
from the  solutions in $Q\times Q$ of the equation $\psi(\alpha)+\psi(\beta)=\alpha+\beta$. There are $51$ of them
$$
 (0 , \alpha +\beta  ),\ (
 \alpha  , \beta  ),\ (
 \alpha  , \alpha +\beta  ),\ (
 \alpha  , \beta +\gamma  ),\ (
 \alpha  , \beta +\delta  ),\ (
 \beta  , \alpha  ),\ (
 \beta  , \alpha +\beta  ),\ $$ $$ (
 \beta  , \alpha +\gamma  ),\ (
 \beta  , \alpha +\delta  ),\  (
 \gamma  , \delta  ),\   (
 \gamma  , \alpha +\beta  ),\ (
 \gamma  , \alpha +\delta  ),\ (
 \gamma  , \beta +\delta  ),\ (
 \delta  , \gamma  ),\ $$ $$  (
 \delta  , \alpha +\beta  ),\ (
 \delta  , \alpha +\gamma  ),\  (
 \delta  , \beta +\gamma  ),\ (
 \alpha +\beta  , 0 ),\ (
 \alpha +\beta  , \alpha  ),\ (
 \alpha +\beta  , \beta  ),\  (
 \alpha +\beta  , \gamma  ),\ $$ $$  (
 \alpha +\beta  , \delta  ),\ (
 \alpha +\beta  , \alpha +\beta  ),\   (
 \alpha +\beta  , \alpha +\gamma  ),\ (
 \alpha +\beta  , \alpha +\delta  ),\ (
 \alpha +\beta  , \beta +\gamma  ),\ (
 \alpha +\beta  , \beta +\delta  ),\  $$ $$  (
 \alpha +\gamma  , \beta  ),\ (
 \alpha +\gamma  , \delta  ),\  (
 \alpha +\gamma  , \alpha +\beta  ),\ (
 \alpha +\gamma  , \alpha +\delta  ),\   (
 \alpha +\gamma  , \beta +\gamma  ),\ (
 \alpha +\gamma  , \beta +\delta  ),\ $$ $$ (
 \alpha +\delta  , \beta  ),\ (
 \alpha +\delta  , \gamma  ),\ (
 \alpha +\delta  , \alpha +\beta  ),\   (
 \alpha +\delta  , \alpha +\gamma  ),\   (
 \alpha +\delta  , \beta +\gamma  ),\ (
 \alpha +\delta  , \beta +\delta  ),\ $$ $$(
 \beta +\gamma  , \alpha  ),\ (
 \beta +\gamma  , \delta  ),\ (
 \beta +\gamma  , \alpha +\beta  ),\ (
 \beta +\gamma  , \alpha +\gamma  ),\    (
 \beta +\gamma  , \alpha +\delta  ),\ (
 \beta +\gamma  , \beta +\delta  ),\ $$ $$ (
 \beta +\delta  , \alpha  ),\ (
 \beta +\delta  , \gamma  ),\   (
 \beta +\delta  , \alpha +\beta  ),\ (
 \beta +\delta  , \alpha +\gamma  ),\ (
 \beta +\delta  , \alpha +\delta  ),\ (
 \beta +\delta  , \beta +\gamma  ).
 $$
Among them only the following $23$  are  liftable as one sees using Lemma \ref{condition to lift0} 
$$
(0,\alpha +\beta ),(\alpha ,\alpha +\beta ),(\beta ,\alpha +\beta ),(\gamma ,\alpha +\beta ),(\delta ,\alpha +\beta ),(\alpha +\beta ,0),(\alpha +\beta ,\alpha ),$$ $$(\alpha +\beta ,\beta ),(\alpha +\beta ,\gamma ),(\alpha +\beta ,\delta ),(\alpha +\beta ,\alpha +\beta ),(\alpha +\beta ,\alpha +\gamma ),(\alpha +\beta ,\alpha +\delta ),$$ $$(\alpha +\beta ,\beta +\gamma ),(\alpha +\beta ,\beta +\delta ),(\alpha +\gamma ,\alpha +\beta ),(\alpha +\gamma ,\beta +\delta ),(\alpha +\delta ,\alpha +\beta ),$$ $$(\alpha +\delta ,\beta +\gamma ),(\beta +\gamma ,\alpha +\beta ),(\beta +\gamma ,\alpha +\delta ),(\beta +\delta ,\alpha +\beta ),(\beta +\delta ,\alpha +\gamma ).
$$
 Thus the range of $F(\scoker(\yb\Delta)): \Homi_\b2(Q,\B^2\times \B^2)\to \Endi_\b2(Q)$ consists of $1+9+9+23$ elements and for instance fails to contain the identity map. \vspace{.05in}
 
  Next, we investigate $\coker(F(\scoker(\yb\Delta)))$. We apply Proposition \ref{extreme and range} to show that it is not a null object of $\b2$.   
\begin{lem}\label{condition to lift} The class of $\id\in \Endi_\b2(Q)$ is non-null  in the cokernel $\coker(F(\scoker(\yb\Delta)))$.
\end{lem}
\proof We apply Proposition \ref{extreme and range} to the element $\xi =\id_Q\in \End_\b2(Q)$. 
It is given by the pair $(\alpha,\beta)$ and does not lift to an element of $\Hom_\b2(Q,\B^2\times \B^2)$ because $\alpha+\beta$ is not $\sigma$-invariant in $\B^2\times \B^2$.
Thus $\xi\notin F(\scoker(\yb\Delta))(\Hom_\b2(Q,\B^2\times \B^2))$. The element $\sigma(\xi)$ corresponds to the pair $(\gamma,\delta)$ and in fact it represents the endomorphism $\sigma\in \End_\b2(Q)$. It remains to show that $\xi$ is indecomposable. The only non-trivial way to write the pair $(\alpha,\beta)$ as a sum in $Q\times Q$ involves $(\alpha,0)$ or $(0,\beta)$ but none of the two terms comes from an element of $\End_\b2(Q)$ since it does not fulfill the condition of $\sigma$-invariance. 
 \endproof 
 In fact, by applying Proposition \ref{existphi1} one determines all elements of $\Endi_\b2(Q)$ whose class is non-null (\ie not $\sigma$-invariant) in the cokernel $\coker(F(\scoker(\yb\Delta)))$. They form the list
 \begin{equation}\label{complement list}
 C=\{ (\alpha ,\beta ),(\alpha ,\beta +\delta ),(\beta ,\alpha ),(\beta ,\alpha +\gamma ),(\gamma ,\delta ),(\gamma ,\beta +\delta ),$$ $$(\delta ,\gamma ),(\delta ,\alpha +\gamma ),(\alpha +\gamma ,\beta ),(\alpha +\gamma ,\delta ),(\beta +\delta ,\alpha ),(\beta +\delta ,\gamma )\}.
 \end{equation}

 The complement of this list is the submodule $\imm(F(\scoker(\yb\Delta)))\subset \Endi_\b2(Q)$.
 
 Lemma \ref{condition to lift} shows that the short \exxx sequence \eqref{short diag exact}, \ie
$	K\stackrel{\subset}{\rightarrowtail} \B^2\times \B^2\stackrel{\scoker(\yb\Delta)}{\xtwoheadrightarrow{}}Q
$, gives rise after applying the functor  $F=\Homi_\b2(Q,-)$ to a non-null cokernel for the map $F(a'')$ with $a''=\scoker(\yb\Delta)$. This yields an element of the satellite functor $SF(K)$. Since $F$ is a $\Hom$-functor the first satellite functor here is the analogue of $\ext^1(Q,K)$. It remains to show that this element remains non-null in the colimit \eqref{satellite sequ}, with respect to a suitable small class of exact sequences. Here for instance one can take the class $\sequ_{\rm small}$ of exact sequences only involving finite objects of $\b2$. This ensures the smallness of the involved categories. 
\begin{figure}[H]
\begin{center}
\includegraphics[scale=1]{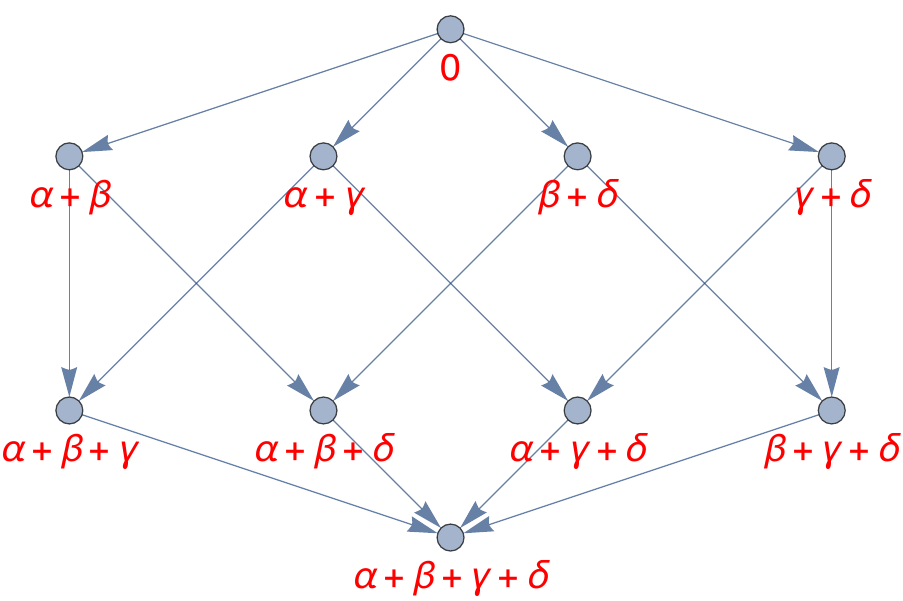}
\end{center}
\caption{The structure of  $K=\ker(\scoker(\yb\Delta))$ \label{kermod0} }
\end{figure}

The structure of $K=\ker(\scoker(\yb\Delta))$ is displayed in Figure \ref{kermod0}. In particular one sees that it is isomorphic to the dual $Q^*$ of $Q$ (see Figure \ref{bmodtstar}) since the involution fixes $\alpha+\gamma$ and $\beta+\delta$ and interchanges $\alpha+\beta$ with $\gamma+\delta$. We are thus dealing with a specific element of $\ext^1(Q,Q^*)$. Inside $K=\ker(\scoker(\yb\Delta))$ one has the range of the diagonal map $\Delta:\B\to \B\times\B$ once lifted to the square. Here it gives the submodule generated by $\alpha+\beta$ and $\gamma+\delta$. It follows from Proposition \ref{extreme and range} that $K={\rm Range}(\Delta)+K^\sigma$. 
 \begin{lem}\label{structure of K} $(i)$~The projection $p:K\to K^\sigma$, $p(x):=x+\sigma(x)$ admits only one non-trivial fiber which is $p^{-1}(\{\tau\})=\{\tau\}\cup (K^\sigma)^c$, where $\tau=\alpha+\beta+\gamma +\delta$. \newline
 $(ii)$~The fiber of the projection $p:\End_\b2(K)\to \End_\b2(K)^\sigma$, $p(\phi):=p\circ \phi=\phi\circ p$ over $\psi\in \End_\bmod(K^\sigma)$ is reduced to $\psi\circ p$ if $\psi(\tau)\neq \tau$.\newline
 $(iii)$~Let $\psi\in \End_\bmod(K^\sigma)$ with $\psi(\tau)=\tau$, then the following map is bijective
 $$
 \{\phi\in \End_\b2(K)\mid p(\phi)=\psi\}\to p^{-1}(\{\tau\}),\ \ \phi\mapsto \phi(\alpha+\beta).
 $$ 	
 \end{lem}
 \proof $(i)$~Using Figure \ref{kermod0} one sees that there is only one non-trivial fiber and it has seven elements. \newline
 $(ii)$~$\alpha+\beta$, $\sigma(\alpha+\beta)$ generate $K$ when taken together with $K^\sigma$. Thus $\phi\in \End_\b2(K)$ is uniquely determined by $\psi=p\circ\phi$ and $\phi(\alpha+\beta)\in p^{-1}(\{\psi(\tau)\})$. By $(i)$ this shows that $\psi$ uniquely determines $\phi$ when $\psi(\tau)\neq \tau$.\newline
    $(iii)$~By the proof of $(ii)$ it follows that $\phi$ is uniquely determined by $\psi$ and by $\phi(\alpha+\beta)\in p^{-1}(\{\tau\})$ using $\psi(\tau)=\tau$. Conversely, any value  $\phi(\alpha+\beta)=\xi\in p^{-1}(\{\tau\})$ defines an extension of $\psi$ as an  endomorphism of $K$ such that $\phi(\gamma+\delta)=\sigma(\xi)$.\endproof
    It follows from Lemma \ref{structure of K} that  $\End_\b2(K)$ has $70$ elements among which $7$ correspond to the seven elements of $\End_\bmod(K^\sigma)$ such that $\psi(\tau)\neq \tau$ and $63=9\times 7$ correspond to the nine elements of $\End_\bmod(K^\sigma)$ such that $\psi(\tau)= \tau$.

 \begin{figure}[H]
\begin{center}
\includegraphics[scale=1]{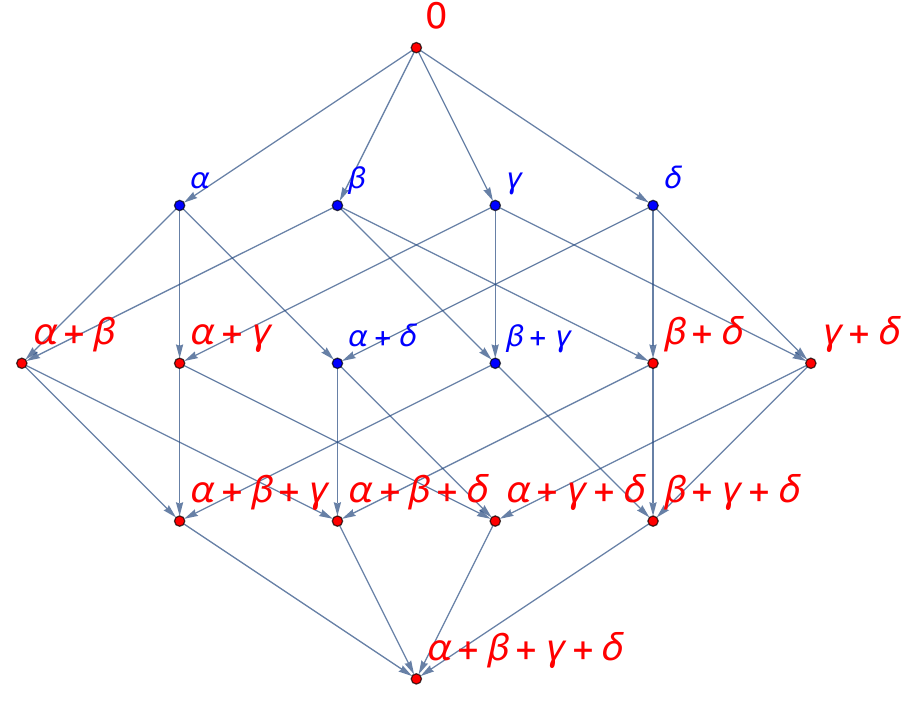}
\end{center}
\caption{The inclusion  $K=\ker(\scoker(\yb\Delta))\subset \B^2\times \B^2$ \label{kermod} }
\end{figure}
\subsection{Duality and the cokernel $\coker(F(\scoker(\yb\Delta)))$}\label{new sect duality}
 In order to show that the satellite functor $SF(K)$ is non-trivial we use a dual theory. We develop it in the explicit example of  $\coker(F(\scoker(\yb\Delta)))$. \vspace{.05in}
 
 We have seen that there are twelve elements of $\Endi_\b2(Q)$ whose image in  $\coker(F(\scoker(\yb\Delta)))$ are non-null, but this does not imply that they are all distinct in the quotient. To determine $\coker(F(\scoker(\yb\Delta)))$ one has to find all linear forms $\phi\in \Hom_\b2(\Endi_\b2(Q),\yb\B)$ whose restriction to the submodule $\imm(F(\scoker(\yb\Delta)))\subset \Endi_\b2(Q)$ is null. One notices that this restriction is entirely specified by a linear form $\ell\in \Hom_\B(E,\B)$, where $E=(\Endi_\b2(Q))^\sigma=\End_\bmod(Q^\sigma)$. Any linear form $\phi$ is uniquely determined by $p_1\circ \phi$,  where $p_1:\yb\B\to \B$ is the first projection. Moreover $p_1\circ \phi$ is determined by $\ell$ and the $4$ values

 $$
 p_1\circ \phi((\alpha ,\beta ))\in \B,\  p_1\circ \phi((\beta ,\alpha ))\in \B,\ p_1\circ \phi((\gamma ,\delta ))\in \B,\ p_1\circ \phi((\delta ,\gamma ))\in \B $$
 since the values of $p_1\circ \phi$ on the remaining $8$ elements of the complement of $\imm(F(\scoker(\yb\Delta)))\subset \Endi_\b2(Q)$ are then uniquely determined using the relations 
 $$
 \left(
\begin{array}{ccc}
 \left(
\begin{array}{cc}
 \alpha  & \beta  \\
 0 & 0 \\
 \alpha  & \beta  \\
\end{array}
\right) & \left(
\begin{array}{cc}
 \alpha  & \beta  \\
 0 & \beta +\delta  \\
 \alpha  & \beta +\delta  \\
\end{array}
\right) & \left(
\begin{array}{cc}
 \alpha  & \beta  \\
 \alpha +\gamma  & 0 \\
 \alpha +\gamma  & \beta  \\
\end{array}
\right) \\
 \left(
\begin{array}{cc}
 \beta  & \alpha  \\
 0 & 0 \\
 \beta  & \alpha  \\
\end{array}
\right) & \left(
\begin{array}{cc}
 \beta  & \alpha  \\
 0 & \alpha +\gamma  \\
 \beta  & \alpha +\gamma  \\
\end{array}
\right) & \left(
\begin{array}{cc}
 \beta  & \alpha  \\
 \beta +\delta  & 0 \\
 \beta +\delta  & \alpha  \\
\end{array}
\right) \\
 \left(
\begin{array}{cc}
 \gamma  & \delta  \\
 0 & 0 \\
 \gamma  & \delta  \\
\end{array}
\right) & \left(
\begin{array}{cc}
 \gamma  & \delta  \\
 0 & \beta +\delta  \\
 \gamma  & \beta +\delta  \\
\end{array}
\right) & \left(
\begin{array}{cc}
 \gamma  & \delta  \\
 \alpha +\gamma  & 0 \\
 \alpha +\gamma  & \delta  \\
\end{array}
\right) \\
 \left(
\begin{array}{cc}
 \delta  & \gamma  \\
 0 & 0 \\
 \delta  & \gamma  \\
\end{array}
\right) & \left(
\begin{array}{cc}
 \delta  & \gamma  \\
 0 & \alpha +\gamma  \\
 \delta  & \alpha +\gamma  \\
\end{array}
\right) & \left(
\begin{array}{cc}
 \delta  & \gamma  \\
 \beta +\delta  & 0 \\
 \beta +\delta  & \gamma  \\
\end{array}
\right) \\
\end{array}
\right)
$$
One can then determine which among the $16\times 16$ possible choices give rise to a linear form and one finds $28$ solutions whose values on the $16+12$ elements of $E\cup C$ are the columns of the large matrix displayed below.
\begin{lem}\label{coker compute} The cokernel of $F(\scoker(\yb\Delta))$ is the map $\tilde p:\Endi_\b2(Q)\to E\cup C$ which is the identity on $C$ and the projection $p(\xi):=\xi+\sigma(\xi)$ on the complement of $C\subset \Endi_\b2(Q)$.
\end{lem}
\proof The map $\tilde p$ is well defined and surjective by construction, and what is non-trivial is the statement that
\begin{equation}\label{operation single}
(\xi,\eta)\mapsto \{\tilde p(u+v)\mid \tilde p (u)=\xi,\ \tilde p(v)=\eta\}
\end{equation}
is single valued. This is clear when $\xi,\eta\in E$ by linearity of $p$ and also for $\xi,\eta \in C$. One can thus assume that $\xi \in C$ and $\eta\in E$.  One then has to show that 
\begin{equation}\label{single} 
	\#\{\tilde p(u+\xi)\mid u\notin C, \ p(u)=\eta\}=1.
\end{equation}
If $u+\xi\notin C$ then $\tilde p(u+\xi)=p(u+\xi)=p(u)+p(\xi)=\eta+p(\xi)$. Thus the case one needs to consider with care is when $u+\xi\in C$. But in our case, the only values of $u$ such that $C\cap (C+u)\neq \emptyset$ are $(\alpha +\gamma,0)$, $(0,\alpha+\gamma)$, $(\beta+\delta,0)$, $(0,\beta+\delta)$ and  all these elements belong to $E$ and the fiber $p^{-1}(\{u\})$ in the complement of $C$ is reduced to the single element $u\in E$. This proves \eqref{single} and shows that the operation \eqref{operation single} is well defined. It is automatically associative and commutative and $\tilde p$ is a morphism by construction. \endproof

The following huge matrix describes the $28$ linear forms $\phi\in \Hom_\b2(\Endi_\b2(Q),\yb\B)$ whose restriction to the submodule $\imm(F(\scoker(\yb\Delta)))\subset \Endi_\b2(Q)$ is null. They are defined by the composite $p_1\circ \phi$ which form the columns of the matrix.  \vspace{0.15in}
 \begin{small}
 $$
\left(
\begin{array}{cc}
 (0,0) & 0,0,0,0,0,0,0,0,0,0,0,0,0,0,0,0,0,0,0,0,0,0,0,0,0,0,0,0 \\
 (0,\alpha +\beta ) & 0,0,0,0,1,1,1,1,1,1,1,1,1,1,1,1,1,1,1,1,1,1,1,1,1,1,1,1 \\
 (0,\alpha +\gamma ) & 0,0,0,0,0,0,0,0,0,0,1,1,1,1,1,1,1,1,1,1,1,1,1,1,1,1,1,1 \\
(0,\beta +\delta ) & 0,0,0,0,1,1,1,1,1,1,0,0,0,0,0,0,1,1,1,1,1,1,1,1,1,1,1,1 \\
 (\alpha +\beta ,0) & 0,1,1,1,0,1,1,1,1,1,0,1,1,1,1,1,0,1,1,1,1,1,1,1,1,1,1,1 \\
 (\alpha +\beta ,\alpha +\beta ) & 0,1,1,1,1,1,1,1,1,1,1,1,1,1,1,1,1,1,1,1,1,1,1,1,1,1,1,1 \\
 (\alpha +\beta ,\alpha +\gamma ) & 0,1,1,1,0,1,1,1,1,1,1,1,1,1,1,1,1,1,1,1,1,1,1,1,1,1,1,1 \\
 (\alpha +\beta ,\beta +\delta ) & 0,1,1,1,1,1,1,1,1,1,0,1,1,1,1,1,1,1,1,1,1,1,1,1,1,1,1,1 \\
 (\alpha +\gamma ,0) & 0,0,1,1,0,0,1,1,1,1,0,0,1,1,1,1,0,0,0,0,1,1,1,1,1,1,1,1 \\
 (\alpha +\gamma ,\alpha +\beta ) & 0,0,1,1,1,1,1,1,1,1,1,1,1,1,1,1,1,1,1,1,1,1,1,1,1,1,1,1 \\
 (\alpha +\gamma ,\alpha +\gamma ) & 0,0,1,1,0,0,1,1,1,1,1,1,1,1,1,1,1,1,1,1,1,1,1,1,1,1,1,1 \\
 (\alpha +\gamma ,\beta +\delta ) & 0,0,1,1,1,1,1,1,1,1,0,0,1,1,1,1,1,1,1,1,1,1,1,1,1,1,1,1 \\
 (\beta +\delta ,0) & 0,1,0,1,0,1,0,1,1,1,0,1,0,1,1,1,0,1,1,1,0,0,0,1,1,1,1,1 \\
 (\beta +\delta ,\alpha +\beta ) & 0,1,0,1,1,1,1,1,1,1,1,1,1,1,1,1,1,1,1,1,1,1,1,1,1,1,1,1 \\
 (\beta +\delta ,\alpha +\gamma ) & 0,1,0,1,0,1,0,1,1,1,1,1,1,1,1,1,1,1,1,1,1,1,1,1,1,1,1,1 \\
 (\beta +\delta ,\beta +\delta ) & 0,1,0,1,1,1,1,1,1,1,0,1,0,1,1,1,1,1,1,1,1,1,1,1,1,1,1,1 \\
 (\alpha ,\beta ) & 0,0,1,1,1,1,1,1,1,1,0,0,1,0,1,1,1,0,1,1,1,1,1,0,1,1,1,1 \\
 (\alpha ,\beta +\delta ) & 0,0,1,1,1,1,1,1,1,1,0,0,1,0,1,1,1,1,1,1,1,1,1,1,1,1,1,1 \\
 (\beta ,\alpha ) & 0,1,0,1,0,1,0,0,1,1,1,1,1,1,1,1,1,1,1,1,0,1,1,1,0,1,1,1 \\
 (\beta ,\alpha +\gamma ) & 0,1,0,1,0,1,0,0,1,1,1,1,1,1,1,1,1,1,1,1,1,1,1,1,1,1,1,1 \\
 (\gamma ,\delta ) & 0,0,1,1,1,1,1,1,1,1,0,0,1,1,0,1,1,1,0,1,1,1,1,1,1,0,1,1 \\
 (\gamma ,\beta +\delta ) & 0,0,1,1,1,1,1,1,1,1,0,0,1,1,0,1,1,1,1,1,1,1,1,1,1,1,1,1 \\
 (\delta ,\gamma ) & 0,1,0,1,0,1,0,1,0,1,1,1,1,1,1,1,1,1,1,1,1,0,1,1,1,1,0,1 \\
 (\delta ,\alpha +\gamma ) & 0,1,0,1,0,1,0,1,0,1,1,1,1,1,1,1,1,1,1,1,1,1,1,1,1,1,1,1 \\
 (\alpha +\gamma ,\beta ) & 0,0,1,1,1,1,1,1,1,1,0,0,1,1,1,1,1,0,1,1,1,1,1,1,1,1,1,1 \\
 (\alpha +\gamma ,\delta ) & 0,0,1,1,1,1,1,1,1,1,0,0,1,1,1,1,1,1,0,1,1,1,1,1,1,1,1,1 \\
 (\beta +\delta ,\alpha ) & 0,1,0,1,0,1,0,1,1,1,1,1,1,1,1,1,1,1,1,1,0,1,1,1,1,1,1,1 \\
 (\beta +\delta ,\gamma ) & 0,1,0,1,0,1,0,1,1,1,1,1,1,1,1,1,1,1,1,1,1,0,1,1,1,1,1,1 \\
\end{array}
\right)$$
\end{small}
\vspace{0.4in}

Thus one is in fact computing a kernel in the dual theory. A more refined study is needed to state the non-vanishing of the classes in the colimit which defines the satellite functor $SF$. The natural framework for discussing the duality is as in Lemma \ref{first satellite}, where we assume now that both categories $\cE$ and $\cB$ are endowed with a duality which given by a contravariant functor of the form $X\longrightarrow X^*=\Hom(X,\beta)$, where $\beta$ is an injective object so that the functor is exact. We then associate to the covariant functor $F:\cE\longrightarrow \cB$  the new covariant functor obtained by conjugation
 \begin{defn}\label{conjugate functor}  The functor $F^*:\cE\longrightarrow \cB$ is defined as  the composition $F^*(X):=(F(X^*))^*$.
 \end{defn} 
 In the case of the functor $F=\Homi_\b2(Q,-)$ the guess for the conjugate functor is $G:=Q\otimes -$, in view of the adjunction
 \begin{equation}\label{hom tensor}
 G(M)^*=\Homi_\b2(Q\otimes M,\yb\B)\simeq \Homi_\b2(Q,\Homi_\b2(M,\yb\B))=F(M^*)
 \end{equation}
One cannot assert that the conjugate of $F$ is $G$ except for finite objects since for these  the duality is involutive  as shown in \S \ref{section duality b2}. Thus we should expect that the conjugate functor of $F=\Homi_\b2(Q,-)$ admits the analogue of the Tor functor as a left satellite. The short \exxx sequence \eqref{short diag exact}
	$$K\leftarrow \ker(\scoker(\yb\Delta))\stackrel{\subset}{\rightarrowtail} \B^2\times \B^2\stackrel{\scoker(\yb\Delta)}{\xtwoheadrightarrow{}}Q$$ together with the isomorphism  $K\simeq Q^*$ provides an element of the comma category $\cI=\sequ_{\rm small}(\b2)\downarrow_{P'} K$ which enters in the construction of $SF(K)$. However, we can view the short \exxx sequence \eqref{short diag exact} as self-dual and use it together with the identity map $\id:Q\to Q$ as
	$$
	Q^*\stackrel{\subset}{\rightarrowtail} \B^2\times \B^2\stackrel{\scoker(\yb\Delta)}{\xtwoheadrightarrow{}}Q\stackrel{\id}{\longleftarrow} Q
	$$
	as an element of the comma category $\cJ=Q\downarrow_{P''} \sequ_{\rm small}(\b2)$ 
	which enters in the construction of $SG(Q)$. We note that the short \exxx sequence 
	\begin{equation}\label{exact alpha}
		\alpha: Q^*\simeq K\stackrel{\alpha'}{\rightarrowtail} \B^2\times \B^2\stackrel{\alpha''}{\xtwoheadrightarrow{}}Q, \ \ \alpha'=\subset, \ \ \alpha''=\scoker(\yb\Delta)
	\end{equation}
	is both a left semiresolution of $Q$ and a right semiresolution of $Q^*$ in the sense of 
	\cite{gra1}, 4.2.1, (4.27).

 \begin{figure}[H]
\begin{center}
\includegraphics[scale=2]{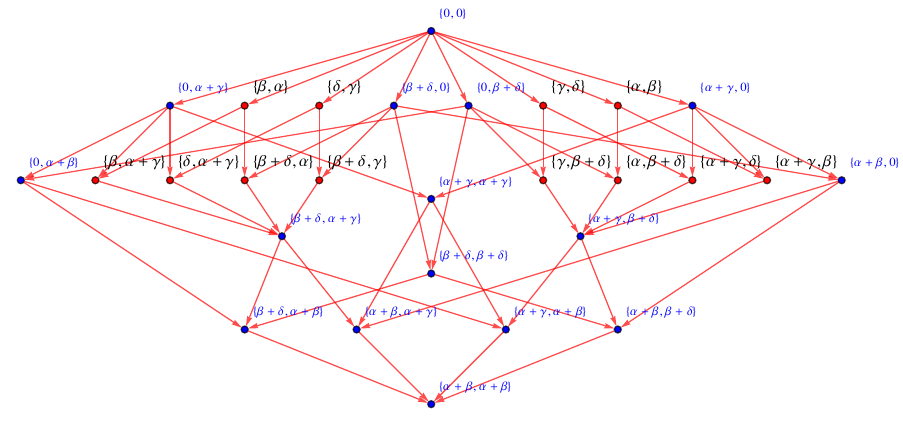}
\end{center}
\caption{Graph of the cokernel $\coker(F(\alpha''))$. \label{graphcokernel} }
\end{figure}
 
  Indeed, the object $\B^2\times \B^2$ of $\b2$ is injective and projective. If one could apply Theorem 4.2.2 of \cite{gra1}, \S 4.2.1, one would conclude that both satellite functors $SF$ and $SG$ can be computed using these semiresolutions as
$$
SF(K)=\coker(F(\alpha'')),\ \ SG(Q)=\ker(G(\alpha')).
$$ 
This raises the problem of proving that condition (a) of Theorem 4.2.2 of \cite{gra1} is valid. This issue already arises with the  short \exxx sequence $\alpha$ as in \eqref{exact alpha} for which the problem of checking condition (a) arises for both functors $F$ and $G$. We first consider $F$ and the 
comma category $\cI=\sequ_{\rm small}(\b2)\downarrow_{P'} K$ 
	which enters in the construction of $SF(K)$. The first question is whether one can determine an action of the group $\Aut_\b2(K)$ on the object $\coker(F(\alpha''))$. This action should reflect the functoriality of the satellite functor $SF$. The group $\Aut_\b2(K)$ is given by the permutations of the $4$ generators $\alpha +\beta ,\alpha +\gamma ,\beta +\delta ,\gamma +\delta$ which respect the partition in the two subsets of $\sigma$-fixed points, \ie $\alpha +\gamma ,\beta +\delta$ and its complement. This gives the following permutations
	$$
	\left(
\begin{array}{cccc}
 \alpha +\beta  & \alpha +\gamma  & \beta +\delta  & \gamma +\delta  \\
 \alpha +\beta  & \beta +\delta  & \alpha +\gamma  & \gamma +\delta  \\
 \gamma +\delta  & \alpha +\gamma  & \beta +\delta  & \alpha +\beta  \\
 \gamma +\delta  & \beta +\delta  & \alpha +\gamma  & \alpha +\beta  \\
\end{array}
\right)
$$
These permutations extend to the automorphisms of $\B^2\times \B^2$ associated to the permutations
 \begin{equation}\label{liftable auto}
 	\left(
\begin{array}{cccc}
 \alpha  & \beta  & \gamma  & \delta  \\
 \beta  & \alpha  & \delta  & \gamma  \\
 \gamma  & \delta  & \alpha  & \beta  \\
 \delta  & \gamma  & \beta  & \alpha  \\
\end{array}
\right)
 \end{equation}
which form a subgroup of order $4$ of the group $\Aut_\b2(\B^2\times \B^2)$ which is of order $8$ and  consists of all permutations of $(\alpha,\beta,\gamma,\delta)$ which commute with $\sigma$. There are $8$ such permutations 
\begin{equation}\label{full auto}
\left(
\begin{array}{cccc}
 \alpha  & \beta  & \gamma  & \delta  \\
 \alpha  & \delta  & \gamma  & \beta  \\
 \beta  & \alpha  & \delta  & \gamma  \\
 \beta  & \gamma  & \delta  & \alpha  \\
 \gamma  & \beta  & \alpha  & \delta  \\
 \gamma  & \delta  & \alpha  & \beta  \\
 \delta  & \alpha  & \beta  & \gamma  \\
 \delta  & \gamma  & \beta  & \alpha  \\
\end{array}
\right)
 \end{equation}
An element of the group $\Aut_\b2(Q)$ is determined by its action on the generators (see Figure \ref{normsubmodn}) and the associated permutation of $(\alpha,\beta,\gamma,\delta)$ must commute with $\sigma$ and also respect the equation $\alpha+\beta=\gamma+\delta$. This reduces the group \eqref{full auto} to its subgroup \eqref{liftable auto}.
We thus get:
\begin{lem}\label{covariance} $(i)$~The functors $P'$ and $P''$ establish  isomorphisms $P':\Aut(\alpha)\to \Aut_\b2(K)$ and $P'':\Aut(\alpha)\to \Aut_\b2(Q)$ where $\alpha$ is the exact sequence \eqref{exact alpha}.\newline
$(ii)$~This determines a canonical functorial action $SF(u)\in \Aut(\coker(F(\alpha''))$ for $u\in \Aut_\b2(K)$.\newline
$(iii)$~Similarly one gets a canonical functorial action $SG(u)\in \Aut(\ker(G(\alpha'))$ for $u\in \Aut_\b2(Q)$.\newline
\end{lem}
\proof $(i)$~The subgroup \eqref{liftable auto} of $\Aut_\b2(\B^2\times \B^2)$ preserves globally the subobject of $\B^2\times \B^2$, $K=\ker(\scoker(\yb\Delta))\simeq Q^*$, and thus acts as automorphisms of the short \exxx sequence $\alpha$. The above computation shows that the functors $P'$ and $P''$ establish  isomorphisms. \newline
$(ii)$~To compute $SF(u)\in \Aut(\coker(F(\alpha''))$ one applies the permutations of  \eqref{liftable auto} to the twelve elements of $C$ of \eqref{complement list}. One checks that one obtains a subgroup of order $4$ of $\Aut(\coker(F(\alpha''))$ whose non-trivial elements act without fixed points on $C$. \newline
$(iii)$~The proof is the same as for $(ii)$.\endproof
The short \exxx sequence $\alpha$ is self-dual, \ie $\alpha^*=\alpha$, where the duality is defined by applying the functor $\Hom_\b2(-,\yb\B)$. Moreover, the satellite functors $SF$ and $SG$ ought to be related by  the duality of \eqref{hom tensor}. \vspace{.05in}

\subsection{The comma category $\cI=\sequ(\b2)\downarrow_{P'} K$}\label{sect comma cat}

Next, we investigate the more delicate aspect of functoriality involving all endomorphisms  $u\in \End_\b2(K)$. The general idea is to use the injectivity of $\B^2\times \B^2$ in the short \exxx sequence $\alpha$ to show that $(\alpha,\id_{K})$ is a weakly final object in the comma category $\cI=\sequ_{\rm small}(\b2)\downarrow_{P'} K$ 
	which enters in the construction of $SF(K)$. Thus let $c$ be an object of the comma category
	$$
	c=\{K\stackrel{v}{\leftarrow} C'\stackrel{c'}{\rightarrowtail}C \stackrel{c''}{\xtwoheadrightarrow{}}C''\}.
	$$ 
Since $\B^2\times \B^2$ is injective, one can extend the morphism $\alpha'\circ v:C'\to \B^2\times \B^2$ to a morphism $w:C\to \B^2\times \B^2$ such that $w\circ c'=\alpha'\circ v$. 
\begin{equation}\label{weakly final}
\xymatrix@C=20pt@R=35pt{
K\ar[d]_\id &  C' \ar[l]_{v}\ar[d]^{v}\ar@{>->}[rr]^{c'} && C\ar[d]_{w} \ar@{->>}[rr]^{c''}&& C''\ar@{-->}[d]
\\
K & \ar[l]_{\id_{K}}K \ \ar@{>->}[rr]^{\alpha'}&& \B^2\times \B^2 \ar@{->>}[rr]^{\alpha''}&& Q}
\end{equation}
In order to obtain a morphism  $c \to (\alpha,\id_{T^*})$ in the comma category it remains to fill the dotted vertical arrow of the diagram \eqref{weakly final}. In fact this is automatic using the proof of Theorem 4.2.2 of \cite{gra1}.
\begin{lem}\label{lem weakly final} The object $(\alpha,\id_{K})$ is  weakly final  in the comma category $\cI=\sequ(\b2)\downarrow_{P'} K$.
\end{lem}
\proof One just needs to fill the dotted vertical arrow of the diagram \eqref{weakly final}, but since the map $c''$ is surjective it is enough to show that $\alpha''\circ w:C\to Q$ is compatible with the equivalence relation associated to $c''$. Since $c$ is a short \exxx sequence this compatibility is equivalent to showing that $\alpha''\circ w\circ c'$ is null, and since $w\circ c'=\alpha'\circ v$ this follows from the nullity of $\alpha''\circ \alpha'$. \endproof 
We now apply this result to test the functoriality of $SF$ with respect to the endomorphisms 
$\End_\b2(K)$. Lemma \ref{lem weakly final} shows that any endomorphism $v\in \End_\b2(K)$ lifts to an endomorphism $w\in \End(\alpha)$ of the short \exxx sequence $\alpha$. For any of the lifts $w\in P'^{-1}(\{v\})$ one obtains an associated element of $\End(\coker(F(\alpha'')))$. More precisely the endomorphism $w \in \End(\alpha)$ induces an endomorphism $w''\in \End_\b2(Q)$ making the following diagram commutative
\begin{equation}\label{w double prime}
\xymatrix@C=20pt@R=35pt{
 \B^2\times \B^2\ar[d]_{w} \ar@{->>}[rr]^{\alpha''}&& Q\ar[d]^{w''}
\\
 \B^2\times \B^2 \ar@{->>}[rr]^{\alpha''}&& Q }
\end{equation} 
This gives a commutative square using the functor $F=\Homi_\b2(Q,-)$
\begin{equation}\label{applying F}
\xymatrix@C=20pt@R=35pt{
 F(\B^2\times \B^2)\ar[d]_{F(w)} \ar[rr]^{F(\alpha'')}&& \Homi_\b2(Q,Q)\ar[d]^{w''\circ - }\ar[rr]^{\scoker(F(\alpha''))}&& \coker(F(\alpha''))\ar@{-->}[d]
\\
 F(\B^2\times \B^2) \ar[rr]^{F(\alpha'')}&& \Homi_\b2(Q,Q) \ar[rr]^{\scoker(F(\alpha''))}&& \coker(F(\alpha''))}\end{equation} 
 and one needs to compute the dotted vertical arrow whose existence follows as in Lemma \ref{lem weakly final}. Note that this vertical arrow is uniquely determined by the functorial action $F(w'')$ of $w''$ on $F(Q)=\Homi_\b2(Q,Q)$, and thus it only depends upon $w''$ rather than on $w$. Thus to compute the functoriality one just needs the correspondence from $v\in \End_\b2(K)$ to $w''\in \End_\b2(Q)$. An explicit computation shows that there are $19$ elements of $\End_\b2(K)$ which give ambiguous values of $w''$ and the $w''$ is unique for the remaining $51$ elements of $\End_\b2(K)$. We provide below the list of the ambiguous elements of $\End_\b2(K)$ by specifying how they act on the three elements $\alpha+\beta$, $\alpha+\gamma$ and $\beta+\delta$. On the right hand side we give the associated choices for $w''$ by their action on the two elements $\alpha,\beta$ of $Q$. 
 $$
 \{\alpha +\gamma ,\alpha +\gamma ,\alpha +\gamma \}\to \left(
\begin{array}{cc}
 \alpha  & \gamma  \\
 \alpha  & \alpha +\gamma  \\
 \gamma  & \alpha  \\
 \gamma  & \alpha +\gamma  \\
 \alpha +\gamma  & \alpha  \\
 \alpha +\gamma  & \gamma  \\
 \alpha +\gamma  & \alpha +\gamma  \\
\end{array}
\right)
$$
$$
\{\beta +\delta ,\beta +\delta ,\beta +\delta \}\to \left(
\begin{array}{cc}
 \beta  & \delta  \\
 \beta  & \beta +\delta  \\
 \delta  & \beta  \\
 \delta  & \beta +\delta  \\
 \beta +\delta  & \beta  \\
 \beta +\delta  & \delta  \\
 \beta +\delta  & \beta +\delta  \\
\end{array}
\right)
$$
$$
\{\alpha +\beta +\gamma ,\alpha +\gamma ,\alpha +\beta +\gamma +\delta \}\to \left(
\begin{array}{cc}
 \alpha  & \alpha +\beta  \\
 \alpha  & \beta +\gamma  \\
 \gamma  & \alpha +\beta  \\
 \alpha +\gamma  & \alpha +\beta  \\
 \alpha +\gamma  & \beta +\gamma  \\
\end{array}
\right)
$$
$$
\{\alpha +\beta +\gamma ,\alpha +\beta +\gamma +\delta ,\alpha +\gamma \}\to \left(
\begin{array}{cc}
 \alpha +\beta  & \alpha  \\
 \alpha +\beta  & \gamma  \\
 \alpha +\beta  & \alpha +\gamma  \\
 \beta +\gamma  & \alpha  \\
 \beta +\gamma  & \alpha +\gamma  \\
\end{array}
\right)
$$
$$
\{\alpha +\beta +\gamma ,\alpha +\beta +\gamma +\delta ,\alpha +\beta +\gamma +\delta \}\to \left(
\begin{array}{cc}
 \alpha +\beta  & \alpha +\beta  \\
 \alpha +\beta  & \beta +\gamma  \\
 \beta +\gamma  & \alpha +\beta  \\
\end{array}
\right)
$$
$$
\{\alpha +\beta +\delta ,\beta +\delta ,\alpha +\beta +\gamma +\delta \}\to \left(
\begin{array}{cc}
 \beta  & \alpha +\beta  \\
 \beta  & \alpha +\delta  \\
 \delta  & \alpha +\beta  \\
 \beta +\delta  & \alpha +\beta  \\
 \beta +\delta  & \alpha +\delta  \\
\end{array}
\right)
$$
$$
\{\alpha +\beta +\delta ,\alpha +\beta +\gamma +\delta ,\beta +\delta \}\to \left(
\begin{array}{cc}
 \alpha +\beta  & \beta  \\
 \alpha +\beta  & \delta  \\
 \alpha +\beta  & \beta +\delta  \\
 \alpha +\delta  & \beta  \\
 \alpha +\delta  & \beta +\delta  \\
\end{array}
\right)
$$
$$
\{\alpha +\beta +\delta ,\alpha +\beta +\gamma +\delta ,\alpha +\beta +\gamma +\delta \}\to \left(
\begin{array}{cc}
 \alpha +\beta  & \alpha +\beta  \\
 \alpha +\beta  & \alpha +\delta  \\
 \alpha +\delta  & \alpha +\beta  \\
\end{array}
\right)
$$
$$
\{\alpha +\gamma +\delta ,\alpha +\gamma ,\alpha +\beta +\gamma +\delta \}\to \left(
\begin{array}{cc}
 \alpha  & \alpha +\beta  \\
 \gamma  & \alpha +\beta  \\
 \gamma  & \alpha +\delta  \\
 \alpha +\gamma  & \alpha +\beta  \\
 \alpha +\gamma  & \alpha +\delta  \\
\end{array}
\right)
$$
$$
\{\alpha +\gamma +\delta ,\alpha +\beta +\gamma +\delta ,\alpha +\gamma \}\to \left(
\begin{array}{cc}
 \alpha +\beta  & \alpha  \\
 \alpha +\beta  & \gamma  \\
 \alpha +\beta  & \alpha +\gamma  \\
 \alpha +\delta  & \gamma  \\
 \alpha +\delta  & \alpha +\gamma  \\
\end{array}
\right)
$$
$$
\{\alpha +\gamma +\delta ,\alpha +\beta +\gamma +\delta ,\alpha +\beta +\gamma +\delta \}\to \left(
\begin{array}{cc}
 \alpha +\beta  & \alpha +\beta  \\
 \alpha +\beta  & \alpha +\delta  \\
 \alpha +\delta  & \alpha +\beta  \\
\end{array}
\right)
$$
$$
\{\beta +\gamma +\delta ,\beta +\delta ,\alpha +\beta +\gamma +\delta \}\to \left(
\begin{array}{cc}
 \beta  & \alpha +\beta  \\
 \delta  & \alpha +\beta  \\
 \delta  & \beta +\gamma  \\
 \beta +\delta  & \alpha +\beta  \\
 \beta +\delta  & \beta +\gamma  \\
\end{array}
\right)
$$
$$
\{\beta +\gamma +\delta ,\alpha +\beta +\gamma +\delta ,\beta +\delta \}\to \left(
\begin{array}{cc}
 \alpha +\beta  & \beta  \\
 \alpha +\beta  & \delta  \\
 \alpha +\beta  & \beta +\delta  \\
 \beta +\gamma  & \delta  \\
 \beta +\gamma  & \beta +\delta  \\
\end{array}
\right)
$$
$$
\{\beta +\gamma +\delta ,\alpha +\beta +\gamma +\delta ,\alpha +\beta +\gamma +\delta \}\to \left(
\begin{array}{cc}
 \alpha +\beta  & \alpha +\beta  \\
 \alpha +\beta  & \beta +\gamma  \\
 \beta +\gamma  & \alpha +\beta  \\
\end{array}
\right)
$$
$$
\{\alpha +\beta +\gamma +\delta ,\alpha +\gamma ,\alpha +\beta +\gamma +\delta \}\to \left(
\begin{array}{cc}
 \alpha  & \alpha +\beta  \\
 \gamma  & \alpha +\beta  \\
 \alpha +\gamma  & \alpha +\beta  \\
\end{array}
\right)
$$
$$
\{\alpha +\beta +\gamma +\delta ,\beta +\delta ,\alpha +\beta +\gamma +\delta \}\to \left(
\begin{array}{cc}
 \beta  & \alpha +\beta  \\
 \delta  & \alpha +\beta  \\
 \beta +\delta  & \alpha +\beta  \\
\end{array}
\right)
$$
$$
\{\alpha +\beta +\gamma +\delta ,\alpha +\beta +\gamma +\delta ,\alpha +\gamma \}\to \left(
\begin{array}{cc}
 \alpha +\beta  & \alpha  \\
 \alpha +\beta  & \gamma  \\
 \alpha +\beta  & \alpha +\gamma  \\
\end{array}
\right)
$$
$$
\{\alpha +\beta +\gamma +\delta ,\alpha +\beta +\gamma +\delta ,\beta +\delta \}\to \left(
\begin{array}{cc}
 \alpha +\beta  & \beta  \\
 \alpha +\beta  & \delta  \\
 \alpha +\beta  & \beta +\delta  \\
\end{array}
\right)
$$
Besides the above $18$ ambiguous elements there is a unique element of $\End_\b2(K)$ which extends in $49$ different ways as an endomorphism of the short \exxx sequence $\alpha$. But it gives rise to only $7$ choices of $w''$ as follows
$$
\{\alpha +\beta +\gamma +\delta ,\alpha +\beta +\gamma +\delta ,\alpha +\beta +\gamma +\delta \}\to \left(
\begin{array}{cc}
 \alpha +\beta  & \alpha +\beta  \\
 \alpha +\beta  & \alpha +\delta  \\
 \alpha +\beta  & \beta +\gamma  \\
 \alpha +\delta  & \alpha +\beta  \\
 \alpha +\delta  & \beta +\gamma  \\
 \beta +\gamma  & \alpha +\beta  \\
 \beta +\gamma  & \alpha +\delta  \\
\end{array}
\right)
$$
We now describe in our case the explicit computation of the action of $w''\in \End_\b2(Q)$ by naturality on the cokernel $\coker(F(\alpha''))$ as in \eqref{w double prime}.

One gets seven fibers of seven elements
\bigskip
$$
\left(
\begin{array}{cc}
 \alpha +\beta  & \alpha +\beta  \\
 \alpha +\beta  & \alpha +\delta  \\
 \alpha +\beta  & \beta +\gamma  \\
 \alpha +\delta  & \alpha +\beta  \\
 \alpha +\delta  & \beta +\gamma  \\
 \beta +\gamma  & \alpha +\beta  \\
 \beta +\gamma  & \alpha +\delta  \\
\end{array}\right),
\left(
\begin{array}{cc}
 \alpha +\beta  & \alpha  \\
 \alpha +\beta  & \gamma  \\
 \alpha +\beta  & \alpha +\gamma  \\
 \alpha +\delta  & \gamma  \\
 \alpha +\delta  & \alpha +\gamma  \\
 \beta +\gamma  & \alpha  \\
 \beta +\gamma  & \alpha +\gamma  \\
\end{array}
\right),\left(
\begin{array}{cc}
 \alpha +\beta  & \beta  \\
 \alpha +\beta  & \delta  \\
 \alpha +\beta  & \beta +\delta  \\
 \alpha +\delta  & \beta  \\
 \alpha +\delta  & \beta +\delta  \\
 \beta +\gamma  & \delta  \\
 \beta +\gamma  & \beta +\delta  \\
\end{array}
\right),\left(
\begin{array}{cc}
 \alpha  & \alpha +\beta  \\
 \alpha  & \beta +\gamma  \\
 \gamma  & \alpha +\beta  \\
 \gamma  & \alpha +\delta  \\
 \alpha +\gamma  & \alpha +\beta  \\
 \alpha +\gamma  & \alpha +\delta  \\
 \alpha +\gamma  & \beta +\gamma  \\
\end{array}
\right),$$
$$\left(
\begin{array}{cc}
 \beta  & \alpha +\beta  \\
 \beta  & \alpha +\delta  \\
 \delta  & \alpha +\beta  \\
 \delta  & \beta +\gamma  \\
 \beta +\delta  & \alpha +\beta  \\
 \beta +\delta  & \alpha +\delta  \\
 \beta +\delta  & \beta +\gamma  \\
\end{array}
\right),\left(
\begin{array}{cc}
 \alpha  & \gamma  \\
 \alpha  & \alpha +\gamma  \\
 \gamma  & \alpha  \\
 \gamma  & \alpha +\gamma  \\
 \alpha +\gamma  & \alpha  \\
 \alpha +\gamma  & \gamma  \\
 \alpha +\gamma  & \alpha +\gamma  \\
\end{array}
\right),\left(
\begin{array}{cc}
 \beta  & \delta  \\
 \beta  & \beta +\delta  \\
 \delta  & \beta  \\
 \delta  & \beta +\delta  \\
 \beta +\delta  & \beta  \\
 \beta +\delta  & \delta  \\
 \beta +\delta  & \beta +\delta  \\
\end{array}
\right)
$$

and the homomorphism is injective on the remaining $21$ elements of $\End_\b2(Q)$.\vspace{.05in}

In the next subsection we shall  describe  how things work for the short \exxx sequence $\alpha$ in order to understand, using Lemma \ref{functoriality on fixed}, why condition $(a)$ holds there. In fact it is important to distinguish two levels.
\subsection{The correspondence $v\mapsto w''$ at the level of $\alpha$}\label{subsubcorrespondence}
In this Section we describe the correspondence between  $\End_\b2(K)$ and $\End_\b2(Q)$ coming from the ambiguity in the extension of an endomorphism of $K$ to an endomorphism of the short \exxx sequence $\alpha$.
We ignore the functor $F$ and work directly at the level of the short \exxx sequence $\alpha$ using the two morphisms 
\begin{equation}\label{quotient map}
	\End_\b2(K)\stackrel{{ res}}{\leftarrow} \End(\alpha)\stackrel{{ quot}}{\to} \End_\b2(Q)
\end{equation}
where $K\simeq Q^*$, $\End_\b2(K)$ is determined by Lemma \ref{structure of K}, $\End(\alpha)$ is the subobject of the endomorphisms $\End_\b2(\B^2\times \B^2)$  which map $K$ to $K$, the left arrow is the restriction to $K$ and the right arrow is the quotient action. Let $X=\{\alpha,\beta,\gamma,\delta\}$ be the set of minimal elements of $\B^2\times \B^2$.   An element $z\in\B^2\times \B^2$ is given uniquely by the subset $Z\subset X$ such that $z=\sum_Z \epsilon$. We denote by $\vert z\vert$ the cardinality of $Z$. An element $\phi\in\End_\b2(\B^2\times \B^2)$ is uniquely determined by the two subsets $A,B\subset X$ such that $\phi(\alpha)=\sum_A \epsilon$ and $\phi(\beta)=\sum_B \epsilon$. The restriction of $\phi$ to $K$ (if it exists) determines $A\cup B$, $A\cup \sigma(A)$ and $B\cup \sigma(B)$ using $\phi(\alpha+\beta)$, $\phi(\alpha+\gamma)$, $\phi(\beta+\delta)$. If $A\cup \sigma(A)$ or $B\cup \sigma(B)$ is empty the union $A\cup B$ determines both $A$ and $B$ uniquely.
  To deal with the lack of uniqueness of the extension of $\phi\in\End_\b2(K)$ as an endomorphism in $\b2$ of $\B^2\times \B^2$ we can thus assume that   $\phi(\tau_j)\neq 0$ $\forall j$, where we let $\tau_1=\alpha+\gamma\in K^\sigma$, $\tau_2=\beta+\delta\in K^\sigma$. We encode the elements $\phi\in\End_\b2(K)$ by the triple $(\phi(\alpha+\beta),\phi(\tau_1),\phi(\tau_2))$ which fulfills the conditions of Lemma \ref{structure of K} (so that in particular $\phi(\tau_j)\in K^\sigma=\{0,\tau_1,\tau_2,\tau\}$).
  \begin{lem}\label{lifting to alpha}
$(i)$~Among the seven elements of $\End_\b2(K)$ such that $\phi(\tau)\neq \tau$ only two, namely $(\tau_j,\tau_j,\tau_j)$ for $j=1,2$, admit more than one extension  to $\B^2\times \B^2$.\newline
 $(ii)$~Let $\phi\in\End_\b2(K)$ with $\phi(\tau)= \tau$. Then $\phi$ admits more than one extension to $\B^2\times \B^2$ if and only if either $\{\vert\phi(\tau_1)\vert,\vert\phi(\tau_2)\vert\}=\{2,4\}$ and $\phi(\alpha+\beta)>\phi(\tau_j)\mid \vert\phi(\tau_j)\vert=2$, or if $\phi(\tau_j)=\tau$, $\forall j$, and $\vert \phi(\alpha+\beta)\vert\geq 3$.
\end{lem}
\proof 
$(i)$~Under the hypothesis $(i)$, the $\sigma$-invariant set $A\cup B \cup \sigma(A)\cup \sigma(B)$ is $\neq X$ and is thus  a single $\sigma$-orbit. If $A\cup B$ has one element this uniquely determines $A$ and $B$ thus the only cases left are when $A\cup B=A\cup \sigma(A)=B\cup \sigma(B)$ is a single orbit. There are two such cases given by 
$\phi(\alpha+\beta)=\phi(\alpha+\gamma)=\phi(\beta+\delta)=\alpha+\gamma$ and $\phi(\alpha+\beta)=\phi(\alpha+\gamma)=\phi(\beta+\delta)=\beta+\delta$. In each case the extensions correspond to the seven choices of two non-empty subsets whose union is a given set with two elements. \newline
$(ii)$~Under the hypothesis $(ii)$, the $\sigma$-invariant set $A\cup B \cup \sigma(A)\cup \sigma(B)$ is $X$. Let us show that if $A\cup B$ has two elements the extension is unique. Indeed $A\cup B$ intersects each $\sigma$-orbit in a single point and thus one has $A=(A\cup B)\cap (A\cup  \sigma(A))$ and $B=(A\cup B)\cap (B\cup  \sigma(B))$. Similarly if $A\cup  \sigma(A)$ and $B\cup  \sigma(B)$ are both $\neq X$ they must be disjoint orbits, thus $A\cap B=\emptyset$ and  one has $A=(A\cup B)\cap (A\cup  \sigma(A))$ and $B=(A\cup B)\cap (B\cup  \sigma(B))$. We have thus shown that if  $\phi$ admits more than one extension one has $\vert \phi(\alpha+\beta)\vert\geq 3$ and either $\{\vert\phi(\tau_1)\vert,\vert\phi(\tau_2)\vert\}=\{2,4\}$ or $\phi(\tau_j)=\tau$, $\forall j$. It remains to exclude the case where $A\cup  \sigma(A)$ is a single $\sigma$-orbit and is not contained in $A\cup B$. In that case $A=(A\cup  \sigma(A))\cap(A\cup B)$ is uniquely determined and has one element. Also $B\cap \sigma(A)=\emptyset$ since otherwise one would have $A\cup  \sigma(A)\subset A\cup B$.  Moreover since $\#(A\cup B)\geq 3$ and $B\cup \sigma(B)=X$ one gets $B=\sigma(A)^c$. This shows the uniqueness of the extension in that case. \endproof

We thus obtain the  list of the $19$ elements of $\End_\b2(K)$ with multiple extensions
$$
\left(
\begin{array}{ccc}
 \alpha +\gamma  & \alpha +\gamma  & \alpha +\gamma  \\
 \beta +\delta  & \beta +\delta  & \beta +\delta  \\
 \alpha +\beta +\gamma  & \alpha +\gamma  & \alpha +\beta +\gamma +\delta  \\
 \alpha +\beta +\gamma  & \alpha +\beta +\gamma +\delta  & \alpha +\gamma  \\
 \alpha +\beta +\gamma  & \alpha +\beta +\gamma +\delta  & \alpha +\beta +\gamma +\delta  \\
 \alpha +\beta +\delta  & \beta +\delta  & \alpha +\beta +\gamma +\delta  \\
 \alpha +\beta +\delta  & \alpha +\beta +\gamma +\delta  & \beta +\delta  \\
 \alpha +\beta +\delta  & \alpha +\beta +\gamma +\delta  & \alpha +\beta +\gamma +\delta  \\
 \alpha +\gamma +\delta  & \alpha +\gamma  & \alpha +\beta +\gamma +\delta  \\
 \alpha +\gamma +\delta  & \alpha +\beta +\gamma +\delta  & \alpha +\gamma  \\
 \alpha +\gamma +\delta  & \alpha +\beta +\gamma +\delta  & \alpha +\beta +\gamma +\delta  \\
 \beta +\gamma +\delta  & \beta +\delta  & \alpha +\beta +\gamma +\delta  \\
 \beta +\gamma +\delta  & \alpha +\beta +\gamma +\delta  & \beta +\delta  \\
 \beta +\gamma +\delta  & \alpha +\beta +\gamma +\delta  & \alpha +\beta +\gamma +\delta  \\
 \alpha +\beta +\gamma +\delta  & \alpha +\gamma  & \alpha +\beta +\gamma +\delta  \\
 \alpha +\beta +\gamma +\delta  & \beta +\delta  & \alpha +\beta +\gamma +\delta  \\
 \alpha +\beta +\gamma +\delta  & \alpha +\beta +\gamma +\delta  & \alpha +\gamma  \\
 \alpha +\beta +\gamma +\delta  & \alpha +\beta +\gamma +\delta  & \beta +\delta  \\
 \alpha +\beta +\gamma +\delta  & \alpha +\beta +\gamma +\delta  & \alpha +\beta +\gamma +\delta \\
\end{array}
\right)
$$
In each of these cases one finds $7$ extensions of $\phi$ to $\B^2\times \B^2$, except for the last one which admits $49$ extensions. To obtain the $w''$ we recall that the cokernel $\scoker:\B^2\times \B^2\to Q$ of the inclusion $K\subset \B^2\times \B^2$ is obtained (using Proposition \ref{compute cokernel}) as the identity on the first ten elements of the list 
$$
0,\alpha ,\beta ,\gamma ,\delta ,\alpha +\beta ,\alpha +\gamma ,\alpha +\delta ,\beta +\gamma ,\beta +\delta,$$ $$ \gamma +\delta ,\alpha +\beta +\gamma ,\alpha +\beta +\delta ,\alpha +\gamma +\delta ,\beta +\gamma +\delta ,\alpha +\beta +\gamma +\delta 
$$
 while all the others project on $\alpha +\beta\sim \gamma +\delta$. Thus the quotient map
 in \eqref{quotient map} is obtained by applying the map $\scoker$ to the values $\phi(\alpha),\phi(\beta)$ for endomorphisms of $\B^2\times \B^2$ which preserve $K$ globally. To see what happens we review the various cases of Lemma \ref{lifting to alpha}. For the case $(i)$, one finds that the operation of passing from $w$ to $w''$ is injective and one gets (for $\tau_1$)  the corresponding $7$ distinct elements of $\End_\b2(Q)$ 
 $$
 \left(
\begin{array}{cc}
 \alpha  & \gamma  \\
 \alpha  & \alpha +\gamma  \\
 \gamma  & \alpha  \\
 \gamma  & \alpha +\gamma  \\
 \alpha +\gamma  & \alpha  \\
 \alpha +\gamma  & \gamma  \\
 \alpha +\gamma  & \alpha +\gamma  \\
\end{array}
\right)
$$
In the case $(ii)$, the first type is $\phi=(\alpha +\beta +\gamma ,\alpha +\gamma ,\alpha +\beta +\gamma +\delta )$, which admits the seven extensions 
$$
\left(
\begin{array}{cc}
 \alpha  & \beta +\gamma  \\
 \alpha  & \alpha +\beta +\gamma  \\
 \gamma  & \alpha +\beta  \\
 \gamma  & \alpha +\beta +\gamma  \\
 \alpha +\gamma  & \alpha +\beta  \\
 \alpha +\gamma  & \beta +\gamma  \\
 \alpha +\gamma  & \alpha +\beta +\gamma  \\
\end{array}
\right)
$$
but when one takes the associated $w''$ the elements $(\gamma  , \alpha +\beta) $ and $(\gamma  , \alpha +\beta+\gamma) $ for instance give the same result so one finally obtains only five values as follows
$$
 \left(
\begin{array}{cc}
 \alpha  & \alpha +\beta  \\
 \alpha  & \beta +\gamma  \\
 \gamma  & \alpha +\beta  \\
 \alpha +\gamma  & \alpha +\beta  \\
 \alpha +\gamma  & \beta +\gamma  \\
\end{array}
\right)
$$
 The second type is $\phi=(\alpha +\beta +\gamma ,\tau ,\tau)$ which admits the seven extensions 
 $$
 \left(
\begin{array}{cc}
 \alpha +\beta  & \beta +\gamma  \\
 \alpha +\beta  & \alpha +\beta +\gamma  \\
 \beta +\gamma  & \alpha +\beta  \\
 \beta +\gamma  & \alpha +\beta +\gamma  \\
 \alpha +\beta +\gamma  & \alpha +\beta  \\
 \alpha +\beta +\gamma  & \beta +\gamma  \\
 \alpha +\beta +\gamma  & \alpha +\beta +\gamma  \\
\end{array}
\right)
$$
but when one takes the associated $w''$ one obtains only three endomorphisms of $Q$
 $$
 \left(
\begin{array}{cc}
 \alpha +\beta  & \alpha +\beta  \\
 \alpha +\beta  & \beta +\gamma  \\
 \beta +\gamma  & \alpha +\beta  \\
\end{array}
\right)
$$
The third type is $\phi=(\tau,\alpha +\gamma ,\tau)$ which admits the seven extensions 
 $$
 \left(
\begin{array}{cc}
 \alpha  & \beta +\gamma +\delta  \\
 \alpha  & \alpha +\beta +\gamma +\delta  \\
 \gamma  & \alpha +\beta +\delta  \\
 \gamma  & \alpha +\beta +\gamma +\delta  \\
 \alpha +\gamma  & \alpha +\beta +\delta  \\
 \alpha +\gamma  & \beta +\gamma +\delta  \\
 \alpha +\gamma  & \alpha +\beta +\gamma +\delta  \\
\end{array}
\right) 
$$
and which induce the following three $w''$
$$
\left(
\begin{array}{cc}
 \alpha  & \alpha +\beta  \\
 \gamma  & \alpha +\beta  \\
 \alpha +\gamma  & \alpha +\beta  \\
\end{array}
\right)
$$
Finally, there remains the element $\phi=(\tau,\tau,\tau)$ which admits $49$ extensions which correspond to all choices of subsets $A,B\subset X$ such that $A\cup B=A\cup\sigma(A)=B\cup\sigma(B)=X$. One checks that they induce the following seven $w''$
$$
\left(
\begin{array}{cc}
 \alpha +\beta  & \alpha +\beta  \\
 \alpha +\beta  & \alpha +\delta  \\
 \alpha +\beta  & \beta +\gamma  \\
 \alpha +\delta  & \alpha +\beta  \\
 \alpha +\delta  & \beta +\gamma  \\
 \beta +\gamma  & \alpha +\beta  \\
 \beta +\gamma  & \alpha +\delta  \\
\end{array}
\right)
$$
We have now completed the description of the correspondence \eqref{quotient map} and shown that it is univalent on $51$ elements of $\End_\b2(K)$ and $1\to 7$, $1\to 5$ and $1\to 3$ on the remaining $19$ elements of $\End_\b2(K)$.  

\subsection{The functor $F$ and the action on the cokernel}\label{subsecfunctact}
\begin{figure}[H]
\begin{center}
\includegraphics[scale=1.4]{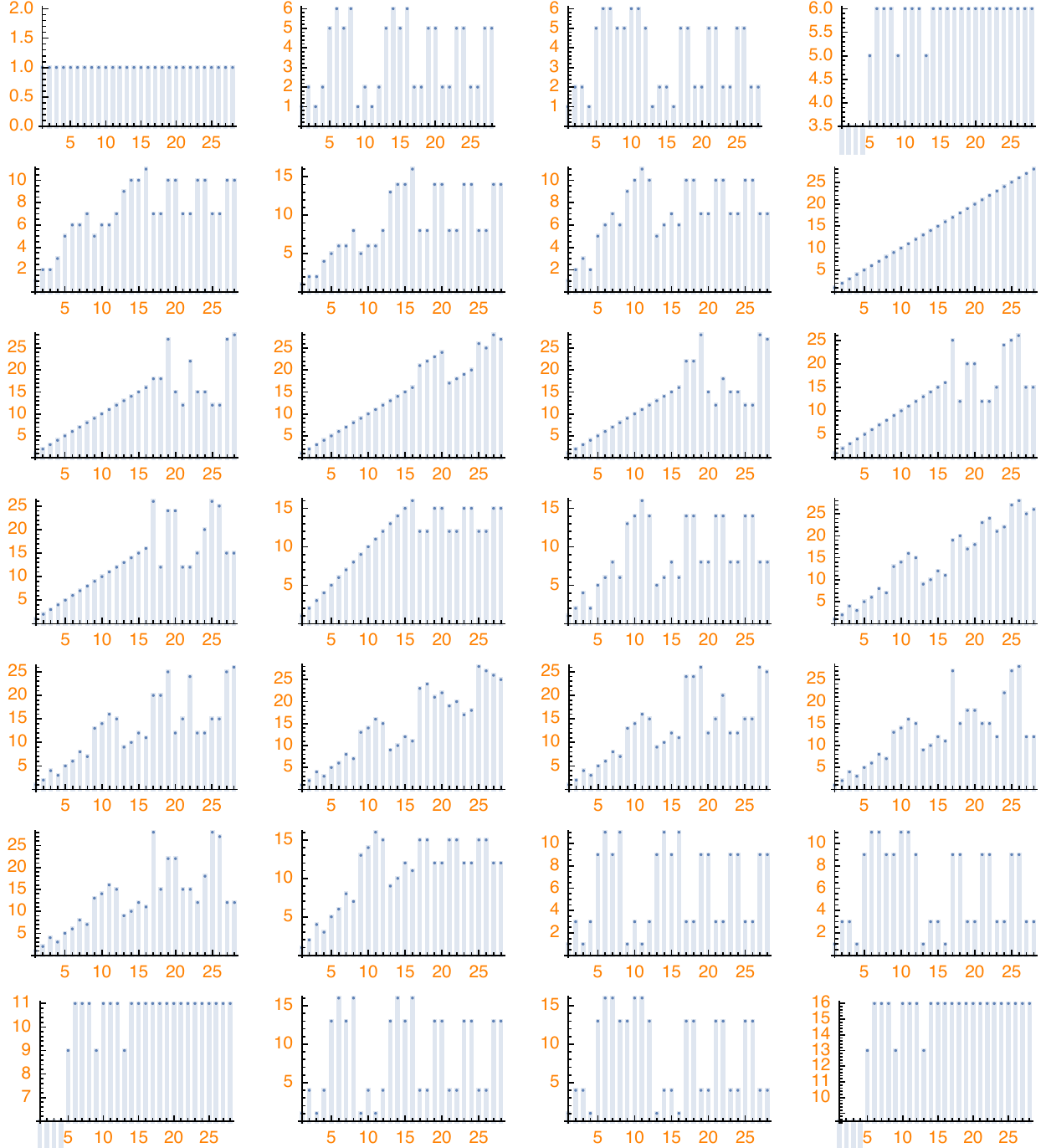}
\end{center}
\caption{Graphs of the 28 transformations of the cokernel induced by endomorphisms of the short \exxx sequence $\alpha$. The elements of the cokernel  labeled in $\{1,2,\ldots,16\}$ are the null elements. Each transformation maps null elements to null elements. There are 16 null transformations, \ie those whose range is formed of null elements.\label{graphoftwentyeigth} }
\end{figure}

The  discussion of \ref{subsubcorrespondence} is independent of the functor $F$ and the key issue is to show that the action of the various $w''$ in the same multivalued piece of the correspondence \eqref{quotient map} all act in the same way on the cokernel $\coker(F(\alpha''))$ of figure \ref{graphcokernel}. One needs to find a conceptual reason why the multivalued pieces all act by null transformations \ie by the sixteen null transformations among the 28 described in Figure \ref{graphoftwentyeigth}. The required independence will then follow from Lemma \ref{functoriality on fixed}
since a morphism $\phi$ in $\b2$ whose range is contained in null elements is determined uniquely as $\psi \circ p$  where $\psi$ is the restriction of $\phi$ to null elements.
In order to simplify the verification that the multivalued pieces all act by null transformations we first investigate the effect of a natural transformation $\mu:F\to F'$ on the computation of the action of $w''$ on the cokernels. With the notations of \eqref{weakly final abelian}, we start with a morphism of short \exxx sequences in $\b2$
\begin{equation}\label{weakly final abelian1}
\xymatrix@C=20pt@R=35pt{
  C'\ar[d]^{v}\ar@{>->}[rr]^{c'} && C \ar[d]_{w} \ar@{->>}[rr]^{c''}&& C''\ar[d]_{w''}
\\
I' \ \ar@{>->}[rr]^{\alpha'}&& I \ar@{->>}[rr]^{\alpha''}&& I''}
\end{equation}
and compare the induced actions on cokernels 
\begin{equation}\label{right action}
F(w''): \coker(F(c''))\to \coker(F(\alpha'')),\ \ F'(w''): \coker(F'(c''))\to \coker(F'(\alpha'')).	
\end{equation}
 \begin{lem}\label{naturality}
$(i)$~Let $w:C\to I$ be a morphism of short \exxx sequences as in \eqref{weakly final abelian1}, and $\mu:F\to F'$ be a natural transformation. The  morphisms induced by 
$F(w''), F'(w''), \mu_{C''}, \mu_{I''}$ form a commutative square.\newline
 $(ii)$~Let $F=F':=\Homi(Q,-)$ and $\mu:F\to F'$ be given by right composition with $\rho\in \End_\b2(Q)$. Then the induced morphism $\tilde \rho:\coker(F(\alpha''))\to \coker(F(\alpha''))$ commutes with the morphism $F(w''): \coker(F(\alpha''))\to \coker(F(\alpha''))$
 for any endomorphism $w$ of the short \exxx sequence $\alpha$.
\end{lem}
\proof $(i)$~One has, using the natural transformation $\mu$, the commutative diagram
\begin{equation*}
\xymatrix{%
&F(C)\ar@{->}[rr]^{F(c'')}\ar@{->}[dl]_{\mu_C}\ar@{->}[dd]_(.65){F(w)}|!{[d];[d]}\hole&& F(C'') %
\ar@{->}[rr]\ar@{->}[dl]_{\mu_{C''}} \ar@{->}[dd]_(.75){F(w'')}|!{[d];[d]}\hole &&
\coker(F(c''))\ar@{-->}[dl]\ar@{->}[dd]
\\
F'(C)\ar@{->}[rr]^(.65){F'(c'')}\ar@{->}[dd]_{F'(w)}&&F'(C'')\ar@{->}[rr]\ar@{->}[dd]_(.75){F'(w'')}&&\coker(F'(c''))
\ar@{->}[dd]\\
&F(I)\ar@{->}[rr]^(.35){F(\alpha'')}|!{[r];[r]}\hole\ar@{->}[dl]_{\mu_I} %
&&F(I'')\ar@{->}[rr]|!{[r];[r]}\hole
\ar@{->}[dl]_{\mu_{I''}}
&&\coker(F(\alpha''))\ar@{-->}[dl]_{}\\
F'(I)\ar@{->}^{F'(\alpha'')}[rr]&&F'(I'')
\ar@{->}[rr]&&\coker(F'(\alpha''))\\
}
\end{equation*}
where the dashed arrows are induced by the natural transformations $\mu_{C''}:F(C'')\to F'(C'')$ and $\mu_{I''}:F(I'')\to F'(I'')$.\newline
$(ii)$~Follows from $(i)$ since right composition with $\rho$ defines a morphism of functors.\endproof 

 \begin{lem}\label{property of complement}
$(i)$~For any $\phi\in \End_\b2(Q)$ such that $\phi\notin \imm(F(\alpha''))$ the restriction $\psi\in \End_\bmod(Q^\sigma)$ of $\phi$ to $Q^\sigma$ is an automorphism.\newline
 $(ii)$~Let $v\in\End_\b2(K)$  admit more than one extension to $\B^2\times \B^2$. Then for any of these extensions $w$ the induced morphism $w''\in \End_\b2(Q)$ fulfills $w''\in \imm(F(\alpha''))$.
 \end{lem}
\proof $(i)$~The twelve elements $\phi\in \End_\b2(Q)$ which do not belong to $\imm(F(\alpha''))$ are given by the list \eqref{complement list}, \ie
\begin{equation}\label{complement list1}
 (\alpha ,\beta ),(\alpha ,\beta +\delta ),(\beta ,\alpha ),(\beta ,\alpha +\gamma ),(\gamma ,\delta ),(\gamma ,\beta +\delta ),$$ $$(\delta ,\gamma ),(\delta ,\alpha +\gamma ),(\alpha +\gamma ,\beta ),(\alpha +\gamma ,\delta ),(\beta +\delta ,\alpha ),(\beta +\delta ,\gamma ).	
 \end{equation}
 The corresponding restrictions $\psi\in \End_\bmod(Q^\sigma)$ are given by applying the projection $p(u)=u+\sigma(u)$ to each term of the list and one gets, with $\tau_1=\alpha+\gamma$, $\tau_2=\beta+\delta$ the list
 $ (\tau_1 ,\tau_2 ),(\tau_1 ,\tau_2 )$, $(\tau_2 ,\tau_1 ),(\tau_2 ,\tau_1 ),(\tau_1 ,\tau_2 ),(\tau_1 ,\tau_2 ),(\tau_2 ,\tau_1 ),(\tau_2 ,\tau_1 ),(\tau_1 ,\tau_2 ),(\tau_1 ,\tau_2 ),(\tau_2 ,\tau_1 ),(\tau_2 ,\tau_1 )
 $ which gives in all cases an automorphism of $Q^\sigma$.\newline
 $(ii)$~Let $v\in\End_\b2(K)$  admit more than one extension to $\B^2\times \B^2$. Then by Lemma \ref{lifting to alpha} the restriction of $v$ to $K^\sigma$ fails to be surjective. For any extension $w$ of $v$ to $\B^2\times \B^2$ the restriction to the null elements is the same and this also holds for the induced morphisms $w''$ on $Q$ as follows from Proposition \ref{fixed cokernel}. Thus $(i)$ shows that the induced morphism $w''\in \End_\b2(Q)$ fulfills $w''\in \imm(F(\alpha''))$.\endproof 
 
 \begin{thm}\label{sequence alpha} The short \exxx sequence $\alpha:K\stackrel{\sker(\scoker(\yb\Delta))}{\rightarrowtail} \B^2\times \B^2\stackrel{\scoker(\yb\Delta)}{\xtwoheadrightarrow{}}Q$ satisfies condition $(a)$  with respect to the functor $F$ and all endomorphisms of $K\simeq Q^*$.  	
 \end{thm}
\proof By Lemma \ref{functoriality on fixed} it is enough to show that if $v\in\End_\b2(K)$  admits more than one extension to $\B^2\times \B^2$ then for any extension $w$ of $v$ to $\B^2\times \B^2$ the action of $w''$ by left multiplication on $\coker(F(\alpha''))$ is null. By $(ii)$ of Lemma \ref{property of complement} the restriction to null elements of $w''$ fails to be surjective and thus the same holds for any $w''\circ u$ $\forall u\in \End_\b2(Q)$ which shows that the range of left multiplication by $w''$ is contained in $\imm(F(\alpha''))$ and is hence null in  $\coker(F(\alpha''))$. \endproof

\begin{rem}\label{rem sequence alpha}{\rm  One may also check directly that none of the $w''$ in the long list of the multiple values for the correspondence \eqref{quotient map} is one of the twelve elements of the complement of $\imm(F(\alpha''))$. One then applies Lemma \ref{naturality} to conclude using the commutation of right multiplications by elements of $\End_\b2(Q)$ to reduce the verification that the left multiplication by multiple values is always null to its verification on $\id\in \End_\b2(Q)$.	
}\end{rem}
\begin{rem}\label{trivial kernel2}{\rm Note  that $\coker(F(\alpha''))$ (see Figure \ref{graphcokernel}) fulfills the condition of Corollary \ref{trivial kernel1}. Indeed, consider a non-trivial fiber of $p$, such as $p^{-1}((\alpha+\gamma,\beta+\delta))$ and let $N=\{0\}\cup p^{-1}((\alpha+\gamma,\beta+\delta))$. Then $N$ has  $4$ maximal ideals, which correspond to the tuples involving each three different letters, and the remaining ideals are obtained by pairwise intersections of the maximal ones. The fact that the dual of $\coker(F(\alpha''))$ is not generated by its minimal elements is obtained by considering the 4 maximal elements (see Figure \ref{graphcokernel}) of this $\B$-module and observing that all of them are fixed under $\sigma$. This shows that the same holds for any element obtained as the top element of the intersection of the corresponding hereditary submodules. Equivalently the submodule of  the dual of $\coker(F(\alpha''))$ which is generated by the minimal elements is null.} 
		\end{rem}
		Let $Q,K\simeq Q^*$ be as in Theorem \ref{sequence alpha}.
\begin{thm}\label{reduction toendo cor2} The satellite functor $SF$ of the functor $F:=\Homi_\b2(Q,-)$ is non-null and $SF(K)$  is the cokernel $\coker(F(\alpha''))$.\end{thm}
\proof   The proof is identical to the proof of Theorem \ref{reduction toendo cor1}, with the role of 
Theorem \ref{sequence s} replaced here by Theorem \ref{sequence alpha}. In the same way $\sequ_{\rm small}$ is any small subcategory of the category $\sequ(\b2)$ of short \exxx sequences in the homological category $\b2$, with $\sequ_{\rm small}$  large enough to contain all short \exxx sequences of finite objects.\endproof

\end{document}